\documentclass[12inch,a4paper]{article}
\usepackage[top=3cm, bottom=2.5cm, left=2.4cm, right=2.7cm]{geometry}
\usepackage{graphicx}
\usepackage{epstopdf}
\usepackage{booktabs}
\usepackage{blkarray}
\usepackage[british]{babel}
\usepackage{amssymb}
\usepackage{setspace}
\usepackage{amsfonts}
\usepackage{amsmath}
\usepackage{float}
\usepackage{color}
\usepackage{url}

\usepackage{epstopdf}

\newtheorem{proposition}{Proposition}

\newcommand{\Z}{ \mathbb{Z} }
\usepackage{tcolorbox}
\usepackage{epstopdf}
\usepackage{subcaption}
\providecommand{\keywords}[1]{\textit{Key words:} #1}
\usepackage[labelfont=bf]{caption}

\usepackage{titlesec}
\titleformat*{\subsection}{\normalfont \large}
\usepackage{abstract}
\usepackage{slashbox}
\begin{document}

% ********************************

% ********************************

\title{A review and comparative study on functional time series techniques}

\author{Javier \'Alvarez-Li\'ebana$^1$}
\maketitle
\begin{flushleft}
$^1$ Department of Statistics and O. R., University of Granada, Spain. 

\textit{E-mail: javialvaliebana@ugr.es}
\end{flushleft}

\doublespacing

% ********************************

% ********************************

\renewcommand{\absnamepos}{flushleft}
\setlength{\absleftindent}{0pt}
\setlength{\absrightindent}{0pt}
\renewcommand{\abstractname}{Summary}
\begin{abstract}
This paper reviews the main estimation and prediction results derived in the context of functional time series, when Hilbert and Banach spaces are considered, specially, in the context of autoregressive processes of order one (ARH(1) and ARB(1) processes, for $H$ and $B$ being a Hilbert and Banach space, respectively). Particularly, we pay attention to the estimation and prediction results, and statistical tests, derived in both parametric  and non-parametric frameworks. A comparative study between different ARH(1) prediction approaches is developed in the simulation study undertaken.  
\end{abstract}


\section{Introduction}
\label{sec:1}
Since the beginning, time series analysis has played a key role in the analysis of temporal correlated data. Due to the huge computing advances, data began to be gathered with an increasingly temporal resolution level. Statistical analysis of temporal correlated high-dimensional data has then  become a very active research area. This paper presents an overview of the main estimation and prediction approaches, formulated in the context of functional time series.

\vspace{0.4cm}
Autoregressive Hilbertian processes of order one (ARH(1) processes) were firstly introduced by  Bosq (1991), where a real separable Hilbert  space $H$ is considered. The functional estimation problem was addressed by a moment-based estimation  of the linear bounded autocorrelation operator, involved in the conditional expectation, providing the least-squares functional predictor (ARH(1) predictor). Projection into the theoretical and empirical eigenvectors of the  autocovariance operator is considered in the computation of a consistent moment-based estimator of the autocorrelation operator. The model introduced in Bosq (1991) has been successfully used by Cavallini et al. (1994), on the forecasting of electricity consumption in Bologna (Italy). Central limit theorems, formulated in Merlev\`ede (1996) and Merlev\`ede et al. (1997), were applied to derive the asymptotic properties of ARH(1)  parameter estimators and predictors (see, e.g., Bosq, 1999a, 1999b). Close graph Theorem allowed Mas (1999) to derive a truncated componentwise estimator of the adjoint of the autocorrelation operator. A central limit theorem for the formulated estimator was also derived. Enhancements to the model firstly established in Bosq (1991), under the Hilbert-Schmidt assumption of the autocorrelation operator, were presented in the comprehensive monograph by Bosq (2000). Specifically, the asymptotic properties of the formulated truncated componentwise parameter estimator of the autocorrelation operator, and of their associated plug-in predictors, are  derived, from the asymptotic behaviour of the empirical eigenvalues and eigenvectors of the autocovariance operator. Improvements of the above-referred results were also provided in  Mas (2000), where an extra regularity condition on the inverse of the autocovariance operator was imposed, to obtain the so-called resolvent class estimators. Efficiency of the componentwise estimator of the autocorrelation operator proposed in Bosq (2000) was studied in Guillas (2001). Mas \& Menneateau (2003a) proved some extra asymptotic results for the empirical functional second-order moments. Mas \& Menneateau (2003b) applied the perturbation theory, in the derivation of the asymptotic properties of the proposed estimators.   Asymptotic behaviour of the ARH(1) estimators was obtained  in Mas (2004, 2007) under  weaker   assumptions (in particular, under the assumption of compactness of the autocorrelation operator).  Menneateau (2005) formulated some laws of the iterated logarithm for an ARH(1) process. Weak-consistency results, in the norm of Hilbert-Schmidt operators, have been recently proposed in \'Alvarez-Li\'ebana et al. (2017), when the covariance and autocorrelation operators  admit a diagonal spectral decomposition, in terms of a common eigenvectors system, under the Hilbert-Schmidt assumption of the autocorrelation operator. Alternative ARH(1) estimation techniques are presented in Besse \& Cardot (1996), where a spline-smoothed-penalized functional principal component analysis (spline-smoothed-penalized FPCA), with rank constraint, was performed. A spline-smoothed-penalized FPCA estimator of the autocorrelation operator of an ARH(p) process was obtained in Cardot (1998), proving  its consistency. A robust spline-smoothed-penalized FPCA estimator of the autocorrelation operator, for a class of ARH(1) processes, was discussed in Besse et al. (2000). 

Statistical tests for the lack of dependence, in the context of linear processes in function spaces, were  derived in Kokoszka et al. (2008). Change point analysis was applied to test the stability and stationarity of an ARH(1) process, in Horv\'ath et al. (2010) and Horv\'ath et al. (2014), respectively. A practical survey about how the ARH(1) framework can be applied to  forecasting electricity consumption was presented in Andersson \& Lillest\o l (2010). The asymptotic normality of the empirical autocovariance operator, and its associated eigenelements, was studied in Kokoszka \& Reimherr (2013a), in the context of non-linear and weakly dependent functional time series. Consistency results of componentwise parameter estimators addressed in H\"ormann \& Kidzi\'nski (2015), in the context of a general functional linear regression, can be particularized to the ARH(1) framework, under weaker assumptions than those ones assumed in Bosq (2000). The case of the autocorrelation operator depending on an unknown real-valued  parameter has been also considered (see Kara-Terki \& Mourid, 2016). This scenario can be applied to the prediction of an Ornstein-Uhlenbeck process in the ARH(1) framework, as developed in \'Alvarez-Li\'ebana et al. (2016). 

\vspace{0.4cm}
An extension of the classical ARH(1) to ARH(p) processes, with $p$ greater than one, was established in Bosq (2000). The crucial choice of the lag order $p$ was discussed in  Kokoszka \& Reimherr (2013b). An extended class of ARH(1) processes, by including exogenous variables in their dependence structure, was formulated  in Guillas (2000). An ARH(p) model, with $p$ greater than one,  to modelling the exogenous random variables, was subsequently proposed by Damon \& Guillas (2002). Indeed, Damon \& Guillas (2005) derived a Markovian representation by its reformulation as a $H^p$-valued ARHX(1) process. Marion \& Pumo (2004) considered the first derivatives of an ARH(1) process, as the exogenous variables to be included in the model, by introducing  the so-called ARHD(1) process. The ARHD(1) process was characterized as a particular ARH(1) model in Mas \& Pumo (2007). Conditional autoregressive Hilbertian processes of order one (CARH(1) processes, also known as doubly stochastic Hilbertian process of order one), were introduced in Guillas (2002), aimed to include exogenous information in a non-additive way (see also Cugliary, 2013, where an underlying  multivariate process  was considered). Mourid (2004) considers  the randomness of the autocorrelation operator  by conditioning to  each element of the sample space. ARH(p) processes with random coefficients (RARH(p) processes) were then defined. Some asymptotic results for the already mentioned RARH(p) processes are provided in Allam \& Mourid (2014). Weakly dependent processes were analysed in H\"ormann \& Kokoszka (2010). Hilbertian periodically correlated autoregressive processes of order one (PCARH(1) processes) have been defined by Soltani \& Hashemi (2011), and later extended to the Banach-valued context by Parvardeh et al. (2017). Spatial extension of the classical ARH(1) models was  firstly proposed in Ruiz-Medina (2011). Their moment-based estimation was detailed in Ruiz-Medina  (2012).  An extension, to the context of processes in function spaces, of the well-known real-valued ARCH model, has been derived in H\"ormann et al. (2013). Recently, Ruiz-Medina \& \'Alvarez-Li\'ebana (2017a) derived the asymptotic efficiency and equivalence of both, classical and Beta-prior-based Bayesian diagonal  componentwise ARH(1) parameter estimators and predictors,  when the autocorrelation operator is not assumed to be a compact operator, in the Gaussian case. Ruiz-Medina \& \'Alvarez-Li\'ebana (2017b) provide sufficient conditions for the strong-consistency, in the trace norm, of the autocorrelation operator of an ARH(1) process, when it does not admit a diagonal spectral decomposition in terms of the eigenvectors of the autocovariance operator. A two-level hierarchical Gaussian model is applied in  Kowal et al. (2017) on the forecasting of ARH(p) processes.

\vspace{0.4cm}
Concerning alternative bases, Grenander's theory of sieves (see Grenander, 1981) was adapted by Bensmain \& Mourid (2001), for the estimation of the autocorrelation operator of an ARH(1) process, from a Fourier-basis-based decomposition in a finite dimensional subspace. An integral form of the autocorrelation operator was assumed in Laukaitis \& Ra\u{c}kauskas (2002), who consider smoothing techniques based on B-spline and Fourier bases. Antoniadis \& Sapatinas (2003) suggested three linear wavelet-based predictors, two of them are built from the componentwise and resolvent estimators of the autocorrelation operator, already established in Bosq (2000) and Mas (2000), respectively. Hyndman \& Ullah (2007) formulated a prediction approach of the mortality and fertility rates based on a real-valued ARIMA forecasting of the FPC scores, for each non-parametric smoothed sample curve. Focusing on the predictor, the idea of replacing the FPC with the directions more relevant to forecasting, by searching a reduced-rank approximation, was firstly exhibited in Kargin \& Onatski (2008) (see also Didericksen et al., 2012, where approaches by Bosq, 2000, and Kargin \& Onatski, 2008, were compared). As an extension of the work by Hyndman \& Ullah (2007), a weighted version of the FPLSR and FPCA approaches was established in Hyndman \& Shang (2009), with a decreasing weighting over time, as often, e.g., in demography. For the purpose of taking into account the information coming from the dynamic dependence, which is usually ignored in the FPCA literature, a dynamic functional principal components approach was simultaneously introduced by Panaretos \& Tavakoli (2013) and H\"orman et al. (2015). The main difference between them lies in the explicit construction of the dynamic scores performed in H\"orman et al. (2015), while a functional process of finite rank was built in Panaretos \& Tavakoli (2013). As discussed in Aue et al. (2015), a multivariate VARMA approach can be used to improve the curve-by-curve approaches above-referred by Hyndman \& Ullah (2007) and Hyndman \& Shang (2009).

\vspace{0.4cm}
Moving-average Hilbertian processes (MAH processes) and ARMAH processes can be defined as a particular case of  Hilbertian general linear processes (LPH). From the previous above-referred works by Bosq (1991) and Mourid (1993), sufficient conditions for the invertibility of LPH were provided in Merlev\`ede (1995, 1996). A Markovian representation of a stationary and invertible LPH, as well as a consistent plug-in predictor, was derived in Merlev\`ede (1997). A few new asymptotic results were derived by Bosq (2000) from the previous theoretical properties. The conditional central limit theorem (see Dedecker \& Merlev\`ede, 2002) was extended to functional stochastic processes in Dedecker \& Merlev\`ede (2003), allowing its application to LPH. The weak convergence for the empirical autocovariance and cross-covariance operators of LPH was proved in Mas (2002). Further results, that those one formulated in Bosq (2000) for LPH, were obtained by Bosq (2007) and Bosq \& Blanke (2007), where the study of a consistent predictor of MAH processes was addressed. Componentwise estimation of a MAH(1) process was studied in Turbillon et al. (2008), under properly assumptions. Tools proposed in Hyndman \& Shang (2008) can be used in the outlier detection of an observed ARMAH(p,q) process. Wang (2008) proposed a real-valued non-linear ARMA model in which functional MA coefficients were considered (see also Chen et al., 2016, where MA coefficients were given by smooth functions). Tensorial products of ARMAH processes have been analysed in Bosq (2010), when innovations are assumed to be a martingale differences sequence.

\vspace{0.4cm}
An extensive literature has been also developed concerning the Banach-valued time series framework. Strong-consistency results on the estimation of a Banach-valued autoregressive process of order one (ARB(1) process) were presented in Pumo (1992, 1998), when $B = \mathcal{C}([0,1]))$ (so-called $AR\mathcal{C}(1)$ process). Its natural extension to $AR\mathcal{C}(p)$ models, with $p$ greater than one, and the characterization of the Ornstein-Uhlenbeck process as an $AR\mathcal{C}(1)$ process, was addressed in Mourid (1993, 1996), respectively. The notion of non-causality was developed in Guillas (2000) for Banach-valued stochastic processes, while the suitable periodicity of the ARB representation was determined by Benyelles \& Mourid (2001). Kuelb's Lemma (see Kuelbs, 1976) played a key role in this Banach-valued framework, since it provides a dense and continuous embedding $B \hookrightarrow H$, where $H$ is the completation of $B$ under a weaker topology (see Labbas \& Mourid, 2002, where  a componentwise estimation of the autocorrelation operator of an ARB(1) process is achieved in a general real separable Banach space). Non-plug-in $AR\mathcal{C}(1)$ prediction is achieved in  Mokhtari \& Mourid (2003), applying the theory of Reproducing Kernel Hilbert Spaces (RKHSs). Dehling \& Sharipov (2005) addressed the asymptotic properties of both empirical mean and empirical autocovariance operator, for an ARB process with weakly dependent innovation process. A wide monograph concerning the estimation of ARB processes, by using the already mentioned sieve estimators, was developed by Rachedi (2005). Recently, for $\mathcal{D} = \mathcal{D}([0,1])$ space (the space of right-continuous functions with left limits), El Hajj (2011, 2013) have derived the asymptotic properties of the parameter estimators and predictors of $AR\mathcal{D}(1)$  processes. The intensity of jumps for these $\mathcal{D}$-valued processes was later analysed in Blanke \& Bosq (2014). Stationarity of ARMAB processes was studied in Spangenberg (2013), under suitable conditions. Central limit theorems for PCARB(1) processes have been recently derived in Parvardeh et al. (2017).

\vspace{0.4cm}
A great amount of authors have been developed alternative non-parametric  prediction techniques in both functional time series and functional  regression frameworks, where the main goal is to  forecast the predictable part of the paths by applying non-parametric methods. Besse et al. (2000) formulated  a functional non-parametric kernel-based predictor of an ARH(1) process. Their work can be seen as a functional extension of the methodology adopted in Poggi (1994), where a non-parametric kernel-based forecasting of the electricity consumption was performed in a multivariate framework. Non-parametric methods were proposed in Cuevas et al. (2002), in the estimation of the underlying linear operator of a functional linear regression, where both explanatory and response variables are valued in a function space. A two-steps prediction approach, based on a non-parametric kernel-based prediction of the scaling coefficients, with respect to a wavelet basis decomposition of a stationary stochastic process with values in a function space, was firstly exhibited in Antoniadis et al. (2006). This method, also so-called kernel wavelet functional (KWF) method, was used, and slightly modified, in Antoniadis et al. (2012), on the forecasting of French electricity consumption, when the hypothesis of stationarity is not held. Improvements in the KWF approach were suggested in Cugliari (2011), where continuous wavelet transforms are also considered.
A functional version of principal component regression (FPCR), and partial least-squares regression (FPLSR), was formulated by Reiss \& Ogden (2007). In the case where the  response is a Hilbert-valued variable, and the explanatory variable takes its values in a general function space, equipped with a semi-metric, Ferraty et al. (2012) obtained a non-parametric kernel estimator of the underlying regression operator. 

\vspace{0.4cm}
The outline of this paper is as follows. References detailed in Sections 2-7 are divided by thematic areas in a chro\-nicle. Section 2 is devoted to the study of the different ARH(1) componentwise estimation frameworks, as well as related results, based on the projections into the theoretical and empirical eigenvectors of the autocovariance operator. Section 3 deals improvements of the classical ARH(1) framework. Parametric forecasting of functional time series, based on the projection into alternative bases, such as Fourier, B-spline or wavelet bases, will be presented in Section 4. Section 5 studies MAH processes, as a particular case of LPH. The Banach-valued context is detailed  in Section 6. Non-parametric techniques, in the context of functional time series and functional linear regression, are described  in Section 7. Section 8 will introduce the main elements involved in the particular context of the ARH(1) diagonal componentwise estimation. In Section 9, a  comparative study is undertaken for  illustration of  the performance of some of the most used  estimation and prediction ARH(1) methodologies. Specifically, the  approaches presented in    \'Alvarez-Li\'ebana et al. (2017), Antoniadis \& Sapatinas (2003), Besse et al. (2000), Bosq (2000) and Guillas (2001) are compared. Proof details and useful theoretical results are provided in the Supplementary Material (see Sections S.1-S.3), where the numerical results obtained in the simulation and comparative studies are outlined as well (see Sections S.4-S.5).

% *************************

% *************************

\section{ARH(1) componentwise estimator, based on the eigenvectors of the autocovariance operator}
\label{sec:2}

ARH(1) processes introduced by Bosq (1991) seek to extend the classical AR(1) model to functional data. In that work, a continuous-time stochastic process $\xi = \left\lbrace \xi_t, \ t \geq 0 \right\rbrace$ is turned into a set of $H$-valued random variables $X = \left\lbrace X_n(t) = \xi_{n\delta + t}, \ n \in \mathbb{Z} \right\rbrace$, with $0 \leq t \leq \delta$. In the sequel, let us consider zero-mean stationary processes, being $H$ a real separable Hilbert. The ARH(1) process was defined by
\begin{equation}
X_n = \rho \left( X_{n-1}\right) + \varepsilon_n, \quad X_n,~\varepsilon_n \in H, \quad n \in \mathbb{Z}, \quad \rho \in \mathcal{L}(H), \label{eq_1}
\end{equation}
\noindent where $\mathcal{L}(H)$ is the space of bounded linear operators on $H$. If $\varepsilon = \left\lbrace \varepsilon_n, \ n \in \mathbb{Z} \right\rbrace$ is assumed to be a $H$-valued strong white noise (and uncorrelated with $X_0$), $\displaystyle \sum_{j=0}^{\infty} \left\| \rho^{j} \right\|_{\mathcal{L}(H)}^{2} < \infty$ is required to the stationarity condition. From the central limit theorems formulated in Merlev\`ede (1996) and Merlev\`ede et al. (1997), the following asymptotic results for the autocovariance operator $C = {\rm E} \left\lbrace X_n \otimes X_n \right\rbrace$, for $n \in \mathbb{Z}$, under ${\rm E}\left\lbrace \left\| X_0 \right \|_{H}^{4} \right\rbrace < \infty$ (so-called \textbf{Assumption A3}) and $\left\| X_0 \right \|_{H} < \infty$, were obtained in Bosq (1999a, 1999b):
\begin{equation}
Z_i = X_i \otimes X_i - C = R \left(Z_{i-1} \right)+E_i, \quad R(z) = \rho z \rho^{*} \in \mathcal{S}(H), \quad \left\| C_n - C \right\|_{\mathcal{S}(H)} = \mathcal{O} \left(\left(\frac{\ln(n)}{n}\right)^{1/2} \right)~a.s., \label{eq_2} 
\end{equation}

\noindent being $C_n = \frac{1}{n} \displaystyle \sum_{i=0}^{n-1} X_i \otimes X_i$, for each $n \geq 2$, $\left\lbrace E_n, \ n \in \mathbb{Z} \right\rbrace$ a martingale difference sequence and $\mathcal{S} (H)$ the class of Hilbert-Schmidt operators on $H$. Since $C$ is compact, from the close graph Theorem, $\rho^{*} = C^{-1} D$ is bounded of the domain of $C^{-1}$, which is a dense subspace in $H$. Then, the autocorrelation operator can be built as $\rho = Ext\left( D C^{-1} \right) = \left(D C^{-1} \right)^{**}$. Mas (1999) provided the asymptotic normality of the formulated  estimator of $\rho^{*}$, under ${\rm E} \left\lbrace \left\| C^{-1} \varepsilon_0 \right\|_{H}^{2} \right\rbrace < \infty$, by projecting into $H_{k_n} = sp\left( \phi_1, \ldots, \phi_{k_n} \right)$, when the eigenvectors $\left\lbrace \phi_j, \ j \geq 1 \right\rbrace$ of $C$ are assumed to be known. Assumption $C_{1}  > \ldots > C_{j} > \ldots > 0$ (\textbf{Assumption A1}) was also imposed, where $\left\lbrace C_j, \ j \geq 1 \right\rbrace$  denote the eigenvalues of $C$.

The asymptotic results formulated in Merlev\`ede (1996) and Merlev\`ede et al. (1997) are crucial in the derivation of some extra asymptotic properties for $C$ and $D= {\rm E} \left\lbrace X_n \otimes X_{n+1} \right\rbrace$ by Bosq (2000). In particular, under \textbf{Assumption A3}, if $D_n = \frac{1}{n-1} \displaystyle\sum_{i=0}^{n-2} X_i \otimes X_{i+1}$ denotes the empirical estimator of the cross-covariance operator $D$, the almost sure convergence to zero of $\frac{n^{1/4}}{\left( \ln(n) \right)^{\beta}} \displaystyle \sup_{j \geq 1} \left\| C_n - C \right\|_{\mathcal{S}(H)}$ and  $\frac{n^{1/4}}{\left( \ln(n) \right)^{\beta}} \displaystyle \sup_{j \geq 1}  \left\| D_n - D \right\|_{\mathcal{S}(H)}$, for any $\beta > 1/2$, was proved. Bosq (2000) also checked that stationarity comes down actually to the existence of an integer $j_0 \geq 1$, with $\left\| \rho^{j_0} \right\|_{\mathcal{L}(H)} < 1$. Under \textbf{Assumptions A1} and \textbf{A3}, as well as the Hilbert-Schmidt assumption of $\rho$, when a spectral decomposition of $C_n$ is achieved in terms of $\left\lbrace C_{n,j}, \ j \geq 1 \right\rbrace$ and $\left\lbrace \phi_{n,j}, \ j \geq 1 \right\rbrace$, the strong-consistency of the following non-diagonal estimator was derived in the above-referred work:
\begin{equation}
\widetilde{\rho}_n = \left( \widetilde{\pi}^{k_n} D_n C_{n}^{-1} \widetilde{\pi}^{k_n} \right), \quad \widetilde{\rho}_n (x) = \displaystyle \sum_{l=1}^{k_n} \widetilde{\rho}_{n,l} (x) \phi_{n,l}, \quad \widetilde{\rho}_{n,l} (x) = \sum_{j=1}^{k_n} C_{n,j}^{-1} \left(\frac{1}{n-1} \displaystyle \sum_{i=0}^{n-2} \widetilde{X}_{i,j} \widetilde{X}_{i+1,l} \right) \langle x, \phi_{n,j} \rangle_H, \label{eq23}
\end{equation}

\noindent assuming that $\left\lbrace \phi_{j}, \ j \geq 1  \right\rbrace$ are unknown, with $\widetilde{X}_{i,j} = \langle X_i, \phi_{n,j} \rangle_H$, for any $j \geq 1$ and $i \in \mathbb{Z}$, being $\widetilde{\pi}^{k_n}$ the orthogonal projector into $\widetilde{H}_{k_n} = sp \left( \phi_{n,1}, \ldots, \phi_{n,k_n} \right)$, for a suitable truncation parameter $k_n$. In the estimation approach formulated in equation (\ref{eq23}), the non-diagonal autocorrelation operator and covariance operator of the error term are defined as follows:
\begin{eqnarray}
\rho \left( X \right) \left( t \right) &=& \displaystyle \int_{a}^{b} \psi \left( t,s \right) X(s) ds, \quad \psi \left( t,s \right) = \displaystyle \sum_{j=1}^{\infty} \displaystyle \sum_{h=1}^{\infty} \rho_{j,h} \phi_j (t) \phi_h (s) \simeq \displaystyle \sum_{j=1}^{M} \displaystyle \sum_{h=1}^{M} \rho_{j,h} \phi_j (t) \phi_h (s), \label{eq_nueva_5} \\
C_{\varepsilon} &=& \displaystyle \sum_{j=1}^{\infty} \displaystyle \sum_{h=1}^{\infty} \sigma_{j,h}^{2} \phi_j \otimes \phi_h  \simeq \displaystyle \sum_{j=1}^{M} \displaystyle \sum_{h=1}^{M} \sigma_{j,h}^{2} \phi_j \otimes \phi_h, \label{eq_nueva_6}
\end{eqnarray}

\noindent for $M$ large enough.

Besides the componentwise estimator of $\rho^{\ast}$, Mas (2000) proposed to approximate $C$ by a linear operator smoothed by a family of  functions $\left\lbrace b_{n,p}(x) = \frac{x^p}{\left(x + b_n \right)^{p+1}}, \ p \geq  0,~n \in \mathbb{N} \right\rbrace$, which converge to $1/x$ point-wise, being $\left\lbrace b_n, \ n \in \mathbb{N} \right\rbrace$ a strictly positive sequence decreasing to zero. The formulated estimators (so-called resolvent class estimators) of $\rho^{\ast}$ were given by Mas (2000) as follows:
\begin{equation}
\widehat{\rho}_{n,p}^{*}= b_{n,p} \left(C_n \right) D_{n}^{*}, \quad b_{n,p} \left(C_n \right) = \left(C_n + b_n I\right)^{-\left(p +1 \right)} C_{n}^{p}, \quad p \geq 0, \quad n \in \mathbb{N},
\end{equation}

\noindent in a manner that $b_{n,p} \left(C_n \right)$ is a compact operator, for each $p \geq 1$ and $n \in \mathbb{N}$, with deterministic norm equal to $b_{n}^{-1}$. Under the non-diagonal scenario in equations (\ref{eq_nueva_5})-(\ref{eq_nueva_6}), a similar philosophy  was adopted by Guillas (2001), in the derivation of the efficiency of the componentwise estimator of $\rho$ formulated in Bosq (2000), in ways that $C_n$ was regularized by a sequence $u=\left\lbrace u_n, \ n \geq 1 \right\rbrace$, with $0 < u_n \leq \beta C_{k_n}$, for $0 < \beta < 1$. Hence,  let us defined $C_{n,j,u}^{-1} =  1 / \displaystyle \max\left(C_{n,j}, u_n \right)$, for any $j \geq 1,~n \geq 2$. An efficient estimator, when $\left\lbrace \phi_j, \ j \geq 1 \right\rbrace$ are unknown, and under \textbf{Assumptions A1} and \textbf{A3}, was stated in Guillas (2001) by
\begin{eqnarray}
 \widetilde{\rho}_{n,u} (x) &=& \displaystyle \sum_{l=1}^{k_n} \widetilde{\rho}_{n,l,u} (x) \phi_{n,l}, \quad \widetilde{\rho}_{n,l,u} (x) = \sum_{j=1}^{k_n} C_{n,j,u}^{-1} \left(\frac{1}{n-1} \displaystyle \sum_{i=0}^{n-2} \widetilde{X}_{i,j} \widetilde{X}_{i+1,l} \right) \langle x, \phi_{n,j} \rangle_H, \label{eq25}
\end{eqnarray}

\noindent for a well-suited truncation parameter, providing the mean-square convergence. Remark that in equation (\ref{eq25}), Hilbert-Schmidt condition over $\rho$ is not needed. We may also cite Mas (2004), where the asymptotic properties, in the norm of $\mathcal{L}(H)$, of the estimator of $\rho^{\ast}$ formulated in Mas (1999), were derived, such that the weaker condition of compactness of $\rho$ was assumed. \textbf{Assumptions A1} and \textbf{A3}, and conditions set in Mas (1999), were also required. This compactness condition, jointly with $\left\| C^{-1/2} \rho \right\|_{\mathcal{L}(H)} < \infty$ (i.e., $\rho$ should be, at least, as smooth as $C^{1/2}$), was also imposed in Mas (2007), where the weak-convergence of the estimator of $\rho^{\ast}$ was addressed, under the convexity of the spectrum of $C$, when $k_n = o \left(\frac{n^{1/4}}{\ln(n)} \right)$. \'Alvarez-Li\'ebana et al. (2017) recently established a weakly-consistent diagonal componentwise estimator of $\rho$, in the norm of $\mathcal{S}(H)$, when $C$ and $\rho$ admit a diagonal spectral decomposition in terms of $\left\lbrace \phi_j, \ j \geq 1 \right\rbrace$. The mean-square convergence of the following estimator of $\rho$, when eigenvectors of $C$ are assumed to be known and $\rho$ is a symmetric operator, was proved, for a well-suited truncation parameter $k_n$:
\begin{equation}
\widehat{\rho}_{k_n} = \displaystyle \sum_{j=1}^{k_n} \widehat{\rho}_{n,j} \phi_j \otimes \phi_j, \quad \widehat{\rho}_{n,j} = \frac{\widehat{D}_{n,j}}{\widehat{C}_{n,j}} = \frac{n}{n-1} \frac{\displaystyle \sum_{i=0}^{n-2} X_{i,j} X_{i+1,j}}{\displaystyle \sum_{i=0}^{n-1} X_{i,j}^{2}}, \quad \widehat{C}_{n,j} \neq 0~a.s., \quad n \geq 2, \label{eq20}
\end{equation}

\noindent under the strictly positiveness of $C$ and  extra mild assumptions. A similar diagonal scenario will be developed in Section 8, where strongly-consistent estimators are provided, when eigenvectors of $C$ are known and unknown.

Alternative ARH(1) estimation parametric techniques, based on a modified version of the functional principal component analysis (FPCA) framework above-referred, have been developed. A spline-smoothed-penalized FPCA, with rank constraint, was presented in Besse \& Cardot (1996) (and later applied by Besse et al., 2000, on the forecasting of climatic variations). In that work, the paths were previously smoothed solving the following non-parametric optimization problem:
\begin{equation}
\displaystyle \min_{\widehat{X}_{i}^{q,\ell} \in H_q} \left\lbrace \frac{1}{np} \displaystyle \sum_{i=0}^{n-1}  \displaystyle \sum_{j=1}^{p} \left(X_i (t_j) - \widehat{X}_{i}^{q,\ell} (t_j) \right)^2 + \ell \left\| D^2 \widehat{X}_{i}^{q,\ell} \right\|_{L^2\left([0,1] \right)}^{2}\right\rbrace,~H_q \subset \left\lbrace f:~ \left\| D^2 f \right\|_{L^2 \left([0,1] \right)}^{2} < c,~c > 0\right\rbrace,\label{eq_8}
\end{equation}

\noindent being $\ell$ the penalized parameter and $\left\lbrace t_j, \ j=1,\ldots,p \right\rbrace$ the set of knots. The $q$-dimensional subspace $H_q$ is spanned by $\left\lbrace A \left( \ell \right)v_j, \ j =1,\ldots,q \right\rbrace$, being $A \left( \ell \right)$ the smoothing hat-matrix and $\left\lbrace v_j, \ j =1,\ldots,q \right\rbrace$ the eigenvectors associated to the first $q$-largest eigenvalues of the  matrix $S = \frac{1}{n}A \left( \ell \right)^{1/2} X^{\prime}  I_n X A \left( \ell \right)^{1/2}$. Estimator of $\rho$ was then built in Besse \& Cardot (1996) by $\widehat{\rho}_{q,\ell} = \widehat{D}_{q,\ell} \widehat{C}_{q,\ell}^{-1}$, with $\widehat{C}_{q,\ell} = \frac{1}{n} \displaystyle \sum_{i=0}^{n-1} \widehat{X}_{i}^{q,\ell} \otimes \widehat{X}_{i}^{q,\ell}$ and $ \widehat{D}_{q,\ell} = \frac{1}{n-1} \displaystyle \sum_{i=0}^{n-2} \widehat{X}_{i}^{q,\ell} \otimes \widehat{X}_{i+1}^{q,\ell}$. See also Cardot (1998), where a spline-smoothed-penalized FPCA was achieved into the Sobolev space $W^{2,2} \left([0,1] \right)$, providing a consistent componentwise truncated estimator of $\rho$ of an ARH(p) process, keeping in mind the regularized paths. Condition $\displaystyle \max_{j=1,\ldots,p-1} \left( t_j - t_{j+1} \right) = \mathcal{O}\left( p^{-1} \right)$, as well as the strictly positiveness of the eigenvalues of the autocovariance operator of the regularized trajectories, was assumed under a suitable choice for the truncation parameter.

It is also worth noting the work by Mas \& Menneatau (2003a), in relation to asymptotic results for the empirical functional second-order moments. Based on the perturbation theory, Mas \& Menneateau (2003b) proved how the asymptotic behaviour of a self-adjoint random operator is equivalent to that of its associated eigenvectors and eigenvalues. The results derived  in Mas \& Menneateau (2003a) are completed by Menneateau (2005), focusing on the law of the iterated logarithm, under the above-referred ARH(1) framework. In a more general framework, the lack of dependence of the functional linear model $Y_n = \Psi \left( X_n \right) + \varepsilon_n$, for each $n \in \mathbb{Z}$, was tested in Kokoszka et al. (2008), under \textbf{Assumptions A1}, \textbf{A3} and the asymptotic properties of $C_n$ derived in Bosq (2000). As discussed in Kokoszka et al. (2008), their approach can be adapted to the ARH(1) framework, and therefore, the nullity of the autocorrelation operator can be tested. In the above-referred works, the null hypotheses of the constancy of $\rho$  and the stationarity condition have been implicitly assumed. Horv\'ath et al. (2010) suggested a testing method on the stability of the autocorrelation operator of an ARH(1) process (against change point alternative), based on the componentwise context above-mentioned, while Horv\'ath et al. (2014) derived testing methods on the stationarity of functional time series (against change point alternative and the so-called two alternatives integrated and deterministic trend). ARH(1) forecasting of the electricity consumption was addressed in Andersson \& Lillest\o l (2010) (see also Cavallini et al., 1994). The asymptotic normality of the empirical principal components of a wide class of functional stochastic processes (even non-linear weakly dependent functional time series) was derived in Kokoszka \& Reimherr (2013a). In the context of linear regression, when both explanatory and response variables are valued in a function space, consistent forecasting of an ARH(1) process was achieved in H\"ormann \& Kidzi\'nski (2015), when the explanatory variables are allowed to be dependent. In the case of $\rho$ depends on an unknown real parameter $\theta$, the Ornstein-Uhlenbeck process (O.U. process) was characterized  by \'Alvarez-Li\'ebana et al. (2016) as an stationary ARH(1) process $ X_n = \rho_{\theta} \left(X_{n-1} \right)  + \varepsilon_{n} $, for any $n \in \mathbb{Z}$ and $\theta > 0$ (see also Kara-Terki \& Mourid, 2016, where the asymptotic normality of an ARH(1) process, with $\rho = \rho_{\theta}$, is studied).

% ******************************

% ******************************

\section{Extensions of the classical ARH(1) model}
\label{sec:3}

Enhancements to the classical ARH(1) model have been developed during the last decades. A great amount of them will be detailed in this Section, arranging the references in chronicle by blocks.

From the previous asymptotic results developed by Mourid (1993, 1996) in the Banach-valued context (see more details in Section 6), the natural extension of ARH(1) to ARH(p) processes, with $p$ greater than one,
was presented in Bosq (2000) as $X_n = \displaystyle \sum_{k=1}^{p} \rho_k \left(X_{n-k} \right) + \varepsilon_n$, for each $n \in \mathbb{Z}$, and $ \rho_k \in \mathcal{L}(H)$, for any $k=1,\ldots,p$, being $\rho_p$ a non-null operator on $H$. By its Markovian properties, ARH(p) model was rewritten by Bosq
(2000) as the $H^p$-valued ARH(1) process $Y_n = \rho^{\prime} \left(Y_{n-1} \right) + \varepsilon_{n}^{\prime}$, with $Y_n = \left(X_n, \ldots, X_{n-p+1} \right) \in H^p$, $\varepsilon_{n}^{\prime} = \left( \varepsilon_n, 0, \ldots, 0 \right) \in H^p$ and

\[\rho^{\prime}=
\begin{blockarray}{cccccc}
\begin{block}{(ccccc)c} 
\rho_1 & \rho_2 & \ldots & \rho_{p-1} &  \rho_p \\ 
I_H & 0_H & \ldots & 0_H & 0_H\\  
\vdots & \vdots & \ddots & \vdots & \vdots\\
0_H & 0_H & \ldots & I_H &  0_H & \longrightarrow \text{\textbf{p-th row}}\\ 
\end{block}
\end{blockarray} \hspace{-2.1cm} \in \mathcal{L}(H^p), \quad \langle \left(X_1,\ldots,X_p \right) ,\left(Y_1,\ldots,Y_p \right)  \rangle_p = \displaystyle \sum_{k=1}^{p} \langle X_k, Y_k \rangle_H, 
 \]
 
\noindent where $H^p$ denotes the cartesian product of $p$ copies of $H$, being a Hilbert space  endowed with $\langle \cdot,\cdot \rangle_p$. In the equation above, $I_H$ and $0_H$ denote the identity and null operators on $H$, respectively. The crucial choice of the lag order $p$ was discussed in Kokoszka \& Reimherr (2013c), when $\rho_{k} \in \mathcal{S}(H)$, for any  $k=1,\ldots,p$ and $\left\| \rho^{\prime}  \right\|_{\mathcal{L}(H)} < 1$. The following multistage testing procedure was proposed in the mentioned work, based on the estimation of the operators $\rho_{k}$, for each $k=1,\ldots,p$:
\begin{eqnarray}
&H_0&:~X\text{ is an i.i.d. sequence} \quad vs \quad H_{p-1}:~X\text{ is an ARH(1) process}, \nonumber \\
&H_{p-1}&:~X\text{ is an ARH(p-1) process} \quad vs \quad H_{p}:~X\text{ is an ARH(p) process}, \nonumber
\end{eqnarray}

\noindent in a manner that the method continues while a null hypothesis is not be rejected.

Aimed to include exogenous information in the dependence structure, ARH(1) processes with exogenous variables (ARHX(1) processes) were introduced in Guillas (2000) and Damon \& Guillas (2002) as follows:
\begin{equation}
X_n = \rho \left(X_{n-1}\right) + \displaystyle \sum_{k=1}^{p} a_k \left(Z_{n,k} \right) +  \varepsilon_n, \quad n \in \mathbb{Z}, \quad a_k,\rho \in \mathcal{L}(H),~k=1,\ldots,p. \label{eq_12}
\end{equation}

\noindent being $Z = \left\lbrace Z_{n,k}, \ n \in \mathbb{Z},~k=1,\ldots,p \right\rbrace$ the exogenous variables. Guillas (2000) initally proposed an AR(1) dependence structure in $Z$ (i.e., $a_k = 0_H$, for any $k=2,\ldots,p$), while the ARH(p) structure displayed in (\ref{eq_12}) was subsequently established in Damon \& Guillas (2002). See also Damon \& Guillas (2005), where a Markovian representation of (\ref{eq_12}) is adopted to reformulate it as a $H^{p}$-valued ARHX(1) process.

The first derivatives of the random paths of an ARH(1) process were included by Marion \& Pumo (2004) as the exogenous variables (so-called ARHD(1) process), when the trajectories belong to the Sobolev space $W^{2,1} \left([0,1]\right)$.  The ARH(1) process was given by $X_n = \rho \left( X_{n-1} \right) + \Psi \left(X_{n-1}^{\prime} \right) + \varepsilon_n$ for each $ n \in \mathbb{Z}$, with $ \rho,~\Psi \in \mathcal{K}(H)$, and was reformulated by Mas \& Pumo (2007) as the ARH(1) process:
\begin{equation}
X_{n} = A \left(X_{n-1} \right) + \varepsilon_n, \quad A = \Phi + \Psi D \in \mathcal{K}(H),~\left\| A \right\|_{\mathcal{L}(H)} < 1, \quad {\rm E} \left\lbrace \left\| X \right\|_{W}^{4} \right\rbrace < \infty, \quad D\left( f \right) = f^{\prime},
\end{equation}

\noindent with $\langle f,g \rangle_W = \displaystyle \int_{0}^{1} f(t) g(t) dt + \displaystyle \int_{0}^{1} f^{\prime}(t) g^{\prime}(t) dt$, for any $f,g \in W^{2,1} \left([0,1]\right)$.

After pointing out some extensions, where exogenous information has been additively incorporated, Guillas (2002) proposed an i.i.d. sequence of Bernoulli variables $I = \left\lbrace I_n, \ n \in \mathbb{Z} \right\rbrace$ to condition  an ARH(1) process, in a non-additive way. A conditional autoregressive Hilbertian process of order one (CARH(1) process, also known as doubly stochastic Hilbertian process of order one) was then formulated:
\begin{equation}
X_n = \rho_{I_n} \left(X_{n-1} \right) + \varepsilon_n = \begin{cases}
      \rho_0 \left(X_{n-1} \right) + \varepsilon_n, & \text{if}\ I_n = 0 \\
      \rho_1 \left(X_{n-1} \right) + \varepsilon_n, & \text{otherwise}
    \end{cases}, \quad \rho_0,~\rho_1 \in  \mathcal{L}(H),~n \in \mathbb{Z}.
 \label{eq__14}
\end{equation}

An extension of (\ref{eq__14}), where the latent process was considered as a continuous multivariate process  $V = \left\lbrace V_n, \ n \in \mathbb{Z} \right\rbrace$, was established in Cugliari (2013). Mourid (2004) proposed to consider the randomness of $\rho$ by defining it from a basic probability space $\left( \Omega, \mathcal{A}, \mathcal{P} \right)$ into $\mathcal{L}(H)$; i.e., $\rho^{\omega}  \in \mathcal{L}(H)$, for each $\omega \in \Omega$.  ARH(p) processes with random coefficients (RARH(p) processes) were then introduced. Its asymptotic properties, in the Hilbert-Schmidt norm, were subsequently derived in Allam \& Mourid (2014). 

In addition, weakly-dependent functional time series, based on a moment-based $m$-dependence, as the most direct relaxation of independence, were studied in H\"ormann \& Kokoszka (2010). The estimation of an ARH(1) process, in which $\rho$ is periodically correlated (PCARH(1) processes), was addressed in Soltani \& Hashemi (2011). In that work, the model $X_n = \rho_n \left(X_{n-1} \right) + \varepsilon_n$, with $\rho_n = \rho_{n+T}$, for each $n \in \mathbb{Z}$, was assumed (periodically correlated with period $T > 0$).

A new branch in the field of functional time series, when the data is gathered on a grid assuming a spatial interaction, was firstly introduced by Ruiz-Medina (2011). In that work, a novel family of spatial stochastic processes (SARH(1) processes), which can be seen as the Hilbert-valued extension of spatial autoregressive processes of order one (SAR(1) processes), was defined as follows:
\begin{equation}
X_{i,j} = R + \rho_1 \left(X_{i-1,j} \right) + \rho_2 \left(X_{i,j-1} \right) + \rho_3 \left(X_{i-1,j-1} \right) + \varepsilon_{i,j},~\left(i,j \right) \in \mathbb{Z}^2, \quad R \in H, \quad \rho_h  \in \mathcal{L}(H),~h=1,2,3,  \label{eq_14}
\end{equation}

\noindent based on the so-called Markov property of the three points for a spatial stochastic process. In (\ref{eq_14}), $\rho_h$ is assumed to be decomposed in terms of the eigenvalues $\left\lbrace  \lambda_{k,h} , \ k \geq 1 \right\rbrace$ and the biorthonormal systems of left and right eigenvectors, $\left\lbrace  \psi_k, \ k \geq 1 \right\rbrace$ and $\left\lbrace  \phi_{k} , \ k \geq 1 \right\rbrace$, respectively, for each $ h =1,2,3$. The spatial innovation process $\left\lbrace  \varepsilon_{i,j}, \ \left(i,j \right) \in \mathbb{Z}^2 \right\rbrace$ is imposed to be a two-parameter martingale difference sequence, with ${\rm E} \left\lbrace \varepsilon_{i,j} \otimes \varepsilon_{i,j} \right\rbrace$ not depending on the coordinates $\left(i,j \right) \in \mathbb{Z}^2$. Ruiz-Medina (2011) derived an unique stationary solution to the SARH(1) state equation (\ref{eq_14}), providing its inversion. The definition of SARH(1) processes, from the tensorial product  of ARH(1) processes, is provided as well. Extended classes of models of functional spatial time series are also formulated in that paper. Moment-based estimators of the functional parameters involved in the SARH(1) equation were proposed in Ruiz-Medina (2012), where their performance is illustrated with a real data application, for spatial functional prediction of ocean surface temperature.

A functional version of ARCH model, given by $X_n = \varepsilon_n \sigma_n$, with $\sigma_{n}^{2} = \delta + \rho \left( X_{n-1}^{2} \right)$, for each $n \in \mathbb{Z}$, was analysed in H\"ormann et al. (2013). A new set of sufficient conditions was provided in Ruiz-Medina \& \'Alvarez-Li\'ebana (2017a) for the asymptotic efficiency of diagonal componentwise estimators of the autocorrelation operator of a stationary ARH(1) process under both, classical and Beta-prior-based Bayesian, scenarios. In particular, under \textbf{Assumption A1}, $\rho$ is assumed to be linear bounded and self-adjoint operator, while the usual Hilbert-Schmidt condition is not imposed. Stronger assumptions for the eigenvalues $\left\lbrace \sigma_{j}^{2}, \ j \geq 1 \right\rbrace$ of $C_{\varepsilon} = {\rm E} \left\lbrace \varepsilon_n \otimes \varepsilon_n \right\rbrace$ were considered, to offset the slower decay rate of the eigenvalues  $\left\lbrace \rho_j, \ j \geq 1 \right\rbrace$ of $\rho$. Specifically, if $\rho = \displaystyle \sum_{j=1}^{\infty} \rho_j \phi_j \otimes \phi_j$, conditions $\rho_{j} = \sqrt{1 - \frac{\sigma_{j}^{2}}{C_j}}$ and $\frac{\sigma_{j}^{2}}{C_j} \leq 1$, with $\frac{\sigma_{j}^{2}}{C_j} = \mathcal{O} \left(j^{-1 - \gamma} \right)$, for any $\gamma > 0$ and $ \sigma_{j}^{2} = {\rm E} \left\lbrace \langle \varepsilon_n, \phi_j \rangle_{H}^{2} \right\rbrace$, for each $j \geq 1$, were assumed. The asymptotic equivalence of the estimators  was also provided, as well as of the their associated  plug-in predictors. The  Beta-prior-based Bayesian estimator of $\rho$ was derived in Ruiz-Medina \& \'Alvarez-Li\'ebana (2017a) as follows:
\begin{equation}
\widetilde{\rho}_n = \displaystyle \sum_{j=1}^{\infty} \widetilde{\rho}_{n,j} \phi_j \otimes \phi_j,~\widetilde{\rho}_{n,j}  = \frac{1}{2 \beta_{n,j}} \left( \left( \alpha_{n,j} + \beta_{n,j} \right) - \sqrt{\left( \alpha_{n,j} - \beta_{n,j} \right)^{2} - 4\beta_{n,j} \sigma_{j}^{2} \left( 2 - \left(a_j + b_j \right)\right)  }\right),
\end{equation}

\noindent with $\alpha_{n,j}  = \displaystyle \sum_{i=0}^{n-1} X_{i,j} X_{i+1,j}$ and $\beta_{n,j}  = \displaystyle \sum_{i=0}^{n-1} X_{i,j}^{2}$, for each $j \geq 1$ and $n \in \mathbb{Z}$, being $\left(a_j, b_j \right)$ the Beta parameters such that $\rho_j \sim \mathbb{B}\left(a_j, b_j \right)$, for any $j \geq 1$. We may also cite Ruiz-Medina \& \'Alvarez-Li\'ebana (2017b), where sufficient conditions for the strong-consistency, in the trace norm, of the above-formulated diagonal componentwise estimator of the autocorrelation operator of an ARH(1) process, are provided. Note that, in that paper, $\rho$ is not assumed to admit a diagonal spectral decomposition with respect to the eigenvectors of the autocovariance operator $C$.

See also Kowal et al. (2017), where a two-level hierarchical model has been recently proposed on the forecasting of an ARH(p) process, by using a Gibbs sampling algorithm. Their purpose is applied to the forecasting of the U.S. Treasury nominal yield curve.

% *******************************

% *******************************

\section{ARH estimation approaches based on alternative bases}
\label{sec:4}

In this section, we pay attention to the ARH(1) estimation, based on the projection into alternative bases to the eigenvectors of $C$. The sieves method introduced by Grenander (1981) was adapted by Bensmain \& Mourid (2001) for the estimation of the autocorrelation operator of an ARH(1) process. A novel consistent estimator was derived under both scenarios, when $\rho$ is a bounded linear operator, and under the Hilbert-Schmidt condition. Specifically, $\rho$ was estimated considering different subsets (so-called sieves) $\left\lbrace \Theta_m, \ m \in \mathbb{N} \right\rbrace$ of the parametric space $\Theta$, where $\rho$ takes its values, equipped with a metric $d$, such that $\Theta_m$ is a compact set, with $\Theta_m \subset \Theta_{m+1}$ and $\bigcup_{m \in \mathbb{N}} \Theta_m$ is dense in $\Theta$. 

In particular, in the former case, when $\rho$ is assumed to be a bounded linear operator $\rho \left( f \right) \left( t \right) = \displaystyle \int_{0}^{1} K \left( t -x \right) f(x) dx$, depending on a  kernel $K \left( \cdot \right)$, then $X_n \left( t \right) = \displaystyle \int_{0}^{1} K \left( t -s \right) X_{n-1} \left( s \right) ds + \varepsilon_n \left( t \right)$. The  Fourier basis $\left\lbrace \phi_0 \left( t \right) = I_{\left[0,1 \right]},~\phi_{2k} (t) = \sqrt{2} \cos \left( 2 \pi k t \right),~\phi_{2k + 1} (t) = \sqrt{2} \sin \left( 2 \pi k t \right), \ k \geq 1 \right\rbrace$ was considered, being $I_{\left[0,1 \right]}$ the identity function over the interval $\left[0,1 \right]$. The ARH(1) state equation  was developed as
\begin{equation}
\begin{cases}
a_0 \left(X_n \right) = a_0 \left(K \right) a_0 \left(X_{n-1} \right) + a_0 \left( \varepsilon_n \right), \quad a_k \left(X_n \right) = \left(a_k \left(K \right) a_k \left(X_{n-1} \right) - b_k \left(K \right) b_k \left(X_{n-1} \right) \right)  / 2 + a_k \left( \varepsilon_n \right) \\
b_k \left(X_n \right) = \left(a_k \left(K \right) b_k \left(X_{n-1} \right) + b_k \left(K \right) a_k \left(X_{n-1} \right) \right)  / 2 +b_k \left( \varepsilon_n \right) 
\end{cases} \label{eq_16}
\end{equation}

\noindent for each $n \in \mathbb{Z}$ and $k \geq 1$, being $\left\lbrace a_k \left(X_n \right),~a_k \left( \varepsilon_n \right),~a_k \left( K \right), \ k \geq 1 \right\rbrace$ and $\left\lbrace b_k \left(X_n \right),~b_k \left( \varepsilon_n \right),~b_k \left( K \right), \ k \geq 1 \right\rbrace$ the Fourier coefficients respect to cosine and sine functions, respectively. Bensmain \& Mourid (2001) assumed that $b_k(t) = 0$, for each $t \in \left[0,1 \right]$ and $k \geq 0$, in a manner that (\ref{eq_16}) becomes
\begin{equation}
x_{n,0} = c_0 x_{n-1,0} + \varepsilon_{n,0}, \quad x_{n,k} = \frac{1}{2} c_k x_{n-1,k} + \varepsilon_{n,k}, \quad x_{n,k} = a_k \left(X_n \right),~ c_k = a_k \left(K \right), \quad k \geq 1,~n \in \mathbb{Z}.
\end{equation}

Estimation of $\rho$ was then reached by forecasting the Fourier coefficients $\left\lbrace c_k, \ k \geq 0 \right\rbrace$ in the sieve
\begin{equation}
\Theta_{m_n} = \left\lbrace K \in L^2 \left( \left[ 0,1 \right] \right):~K(t) = c_{0}I_{\left[0,1 \right]} + \displaystyle \sum_{k=1}^{m_n} c_k \sqrt{2} \cos \left( 2\pi k t \right),\quad  \displaystyle \sum_{k=1}^{m_n} k^2 c_{k}^{2} \leq m_n,~m_n \longrightarrow^{n \to \infty} \infty \right\rbrace,
\end{equation}

\noindent providing, with $m_n \longrightarrow \infty$ as $n \to \infty$,
\begin{equation}
\widehat{c}_0 = \frac{\displaystyle \sum_{i=1}^{n} x_{i,0} x_{i-1,0}}{\displaystyle \sum_{i=1}^{n} x_{i-1,0}^{2}}, \quad \widehat{c}_k = \frac{\displaystyle \sum_{i=1}^{n} x_{i,k} x_{i-1,k}}{\displaystyle \sum_{i=1}^{n} \frac{1}{2} x_{i-1,k}^{2} + n 2 \lambda k}, \quad \text{ under } \displaystyle \sum_{k=1}^{m_n} k^2 \left( \frac{\displaystyle \sum_{i=1}^{n} x_{i,k} x_{i-1,0}}{\displaystyle \sum_{i=1}^{n} \frac{1}{2} x_{i-1,k}^{2} + n 2 \lambda k}\right)^2 = m_n.
\end{equation}

The non-diagonal componentwise estimator formulated in Bosq (2000) was used in Laukaitis \& Ra\u{c}kauskas (2002), by considering regularized paths in terms of a B-spline basis. In that work, the forecasting of the intensity of both cash withdrawal in cash machines (so-called automatic teller machines or ATM) and transactions in points of sale (so-called POS), depending on Vilnius Bank, was achieved. Antoniadis \& Sapatinas (2003) discussed how the prediction of functional stochastic processes can be seen as a linear ill-posed inverse problem, providing a few approaches about the regularization techniques required. In the context of 1-year-ahead forecasting of the climatological ENSO time series, they also proposed three linear wavelet-basis-based ARH(1) predictors, one of which is based on the resolvent estimators of $\rho$ formulated in Mas (2000). From the componentwise estimation framework developed in Bosq (2000), they derived regularized wavelet estimators, by means of a previously wavelet-basis-based smoothing method:
\begin{equation}
\widetilde{Y}_{i,\widehat{\lambda}^{M}} = \widetilde{X}_{i,\widehat{\lambda}^{M}} -  \frac{1}{n} \displaystyle \sum_{i=0}^{n-1} \widetilde{X}_{i,\widehat{\lambda}^{M}}, \quad \widetilde{X}_{i,\widehat{\lambda}^{M}} = \displaystyle  \sum_{k=0}^{2^{j_0} - 1} \widehat{\alpha_{j_{0}k}^{i}} \phi_{j_{0} k} + \displaystyle \sum_{j=j_0}^{J-1} \displaystyle \sum_{k=0}^{2^j - 1} \widehat{\beta_{jk}^{i}} \psi_{jk}, \quad i \in \mathbb{Z}, \label{wavelets1a}
\end{equation}
\noindent with smoothing parameter $\widehat{\lambda}^{M} = \left( \displaystyle \sum_{j=1}^{M} \sigma_{j}^{2} \right)\left( \displaystyle \sum_{j=1}^{M} C_{j} \right) / N$. The plug-in predictor was given by
\begin{eqnarray}
\widetilde{\rho}_{n,\widehat{\lambda}^{M}} \left(X_{n-1} \right) &=& \displaystyle \sum_{j=1}^{k_n} \widetilde{\rho}_{n,\widehat{\lambda}^{M},j} \left(X_{n-1} \right) \widetilde{\phi}_{j}^{M}, \nonumber \\
 \widetilde{\rho}_{n,\widehat{\lambda}^{M},j} \left(X_{n-1} \right) &=& \frac{1}{n-1} \displaystyle \sum_{k=1}^{k_n}  \displaystyle \sum_{i=0}^{n-2} \frac{1}{\widetilde{C}_{n,\widehat{\lambda}^{M},k}}  \widetilde{X}_{n-1,\widehat{\lambda}^{M},k} \widetilde{Y}_{i,\widehat{\lambda}^{M},k} \widetilde{Y}_{i+1,\widehat{\lambda}^{M},j}, \label{a105}
\end{eqnarray}
\noindent with $\widetilde{X}_{n-1,\widehat{\lambda}^{M},j} = \langle \widetilde{\phi}_{j}^{M}, X_{n-1} \rangle_H$ and $\widetilde{Y}_{i+1,\widehat{\lambda}^{M},j} = \langle \widetilde{\phi}_{j}^{M}, \widetilde{Y}_{i+1,\widehat{\lambda}^{M}} \rangle_H$, for each $j=1,\ldots,k_n$ and $i=0,\ldots,n-1$, where $\left\lbrace \widetilde{C}_{n,\widehat{\lambda}^{M},j}, \ j \geq 1 \right\rbrace$ and  $\left\lbrace \widetilde{\phi}_{j}^{M}, \ j \geq 1 \right\rbrace$ denote the eigenvalues and eigenvectors, respectively, of the empirical estimator $\widetilde{C}_{n,\widehat{\lambda}^{M}} = \frac{1}{n} \displaystyle \sum_{i=0}^{n-1} \widetilde{Y}_{i,\widehat{\lambda}^{M}} \otimes \widetilde{Y}_{i,\widehat{\lambda}^{M}}$. 
Values $\left\lbrace \widehat{\alpha_{j_{0}k}^{i}},~\phi_{j_{0}k}, \ k=0,\ldots,2^{j_0}-1 \right\rbrace$ and $\left\lbrace \widehat{\beta_{jk}^{i}},~\psi_{jk}, \ j \geq j_0,~k=0,\ldots,2^{j}-1 \right\rbrace$, for $i=0,\ldots,n-1$, in equation (\ref{wavelets1a}), denote the scaling coefficients, at $J - j_0$ resolutions levels, for a primary resolution level $j_0 < J$. \textbf{Assumptions A1} and \textbf{A3} were imposed, along with
\begin{equation}
n C_{k_n}^{4} \to \infty, \quad \frac{1}{n} \displaystyle \sum_{j=1}^{k_n} \frac{b_j}{C_{j}^{2}} \to 0,~\text{as } n \to \infty, \quad b_j = \displaystyle \max \left( \left(C_{j-1} - C_j \right)^{-1}, \left(C_{j} - C_{j+1} \right)^{-1} \right). \label{wavelets1c}
\end{equation}

Hyndman \& Ullah (2007) detailed an alternative robust version of FPCA, to avoid the instability induced by outlying observations. Forecasting of mortality and fertility rates, as continuous curves, was performed in Hyndman \& Ullah (2007):
\begin{equation}
y_t(x_i) = f_t(x_i) + \sigma_t (x_i) \varepsilon_{t,i}, \quad \left\lbrace  \varepsilon_{t,i}, \ t=1,\ldots,n,~i=1,\ldots, p \right\rbrace~\text{ i.i.d. sequence of standard normal},
\end{equation}
\noindent where $\left\lbrace x_i, \ i=1,\ldots, p \right\rbrace$ denotes the ages covered, being $\left\lbrace y_t(x_i), \ i=1,\ldots, p,~t=1,\ldots,n \right\rbrace$ the log-rates of mortality (or fertility) for age $x_i$ in year $t$. Curves are decomposed in terms of
\begin{equation}
f_t (x) =  \displaystyle \sum_{k=1}^{K} \beta_{t,k} \psi_k (x) + e_t (x), \quad e_t(x) \sim \mathcal{N} \left(0,v(x) \right), \quad t=0,1,\ldots,n,
\end{equation}

\noindent being $\left\lbrace \psi_k, \ k=1,\ldots,K \right\rbrace$ an orthonormal basis, with the coefficients $\left\lbrace \beta_{t,k}, \ k=1,\ldots,K \right\rbrace$ being predicted by using a real-valued ARIMA process. A weighted version of the approach presented in Hyndman \& Ullah (2007) considering the largest weights for the most recent data (required in fields such as demography), was developed in Hyndman \& Shang (2009). Instead of the curve-by-curve forecasting established in Hyndman \& Ullah (2007) and Hyndman \& Shang (2009), a multivariate VARMA model was applied by Aue et al. (2015), to avoid the loss of information invoked by the uncorrelated assumption of the FPC scores, imposed in those works.

Kargin \& Onatski (2008) focused on the predictor of an ARH(1) process, instead of on the operators $\rho$ and $C$ themselves. They proposed to replace the functional principal components with directions more relevant to forecasting, by searching a reduced-rank approximation (see also Didericksen et al., 2012, where a comparative study, between approaches in Bosq, 2000, and Kargin \& Onatski, 2008, was undertaken). Their method, so-called predictive factor decomposition, is built under the searching of a minimal operator $\rho \in R_p$, aimed to minimize  ${\rm E} \left\lbrace \left\| X_n - \rho \left(X_{n-1} \right) \right\|_{H}^{2} \right\rbrace$, being $R_p$ the set of $p$-rank operator. The predictor was then given by
\begin{equation}
\widehat{X}_n = \displaystyle \sum_{l=1}^{p} \langle X_{n-1}, \widehat{b}_{l}^{\alpha}\rangle_H D_n \widehat{b}_{l}^{\alpha}, \quad \widehat{b}_{l}^{\alpha} = \alpha \widehat{x}_{l}^{\alpha}  + \displaystyle \sum_{j=1}^{K} \frac{\langle \widehat{x}_{l}^{\alpha}, \phi_{n,j} \rangle_H}{C_{n,j}^{1/2}}\phi_{n,j} , \quad l=1,\ldots,p,
\end{equation}

\noindent being $\left\lbrace \widehat{x}_{l}^{\alpha}, \ l=1,\ldots,p \right\rbrace$ a linear combination of the eigenvectors $\left\lbrace \phi_{n,j}, \ j=1,\ldots,p \right\rbrace$ of the empirical autocovariance operator. Kargin \& Onatski (2008) proposed to fix $\alpha = 0.75$. 

A dynamic functional principal components analysis (DFPCA) approach was formulated by Panaretos \& Tavakoli (2013), based on an harmonic decomposition (so-called Cram\'er-Karhunen-Lo\`eve decomposition) of the paths by using a set of spectral density operators $\left\lbrace F_{\omega}, \ \omega \in \left[-\pi, \pi \right] \right\rbrace$. These are defined as the discrete-time Fourier transform $F_{\omega} =  \frac{1}{2 \pi} \displaystyle \sum_{h \in \mathbb{Z}} e^{-i \omega h} C_h$ of the covariance operators $C_h = {\rm E} \left\lbrace X_n \otimes X_{n+h} \right\rbrace$, for each $h \in \mathbb{Z}$. The formulated predictor is given by a stochastic integral, depending on a finite sum of tensorial products of eigenvectors of the spectral density operators, as an extension of the Brillinger's information theory. See also H\"orman et al. (2015), where a DFPCA was also established, by replacing the usual FPC scores $\langle X_n, \phi_j \rangle_H$, for each $j \geq 1$ and $n \in \mathbb{Z}$, with an explicitly construction of dynamic FPC scores $Y_{n,j} = \displaystyle \sum_{l \in \mathbb{Z}} \langle X_{n-1} \left( \cdot \right), \psi_{j} \left( \cdot,l \right)\rangle_H$:
\begin{equation}
\psi_{j}(\cdot,l) = \frac{1}{2 \pi} \displaystyle \int_{-\pi}^{\pi} \varphi_{j} (u, \omega) e^{-i l \omega} d \omega, \quad F_{\omega} \left( \cdot \right) =  \frac{1}{2 \pi} \displaystyle \sum_{h \in \mathbb{Z}} e^{-i \omega h} C_h \left( \cdot \right)= \displaystyle \sum_{j=1}^{\infty} \lambda_{j} (\omega) \langle \cdot, \varphi_{j}(\cdot , \omega)  \rangle_{H} \varphi_{j}(\cdot, \omega) .
\end{equation}

% ******************************

% ******************************

\section{Hilbert-valued moving-average and general linear processes}

This section is devoted to describe the main contributions in the field of Hilbertian moving-average processes (MAH processes), including the general case of Hilbertian general linear processes (LPH). The case of ARMAH processes is considered as well. From the Wold decomposition $X_n = \varepsilon_n + \displaystyle \sum_{k=1}^{\infty} a_k \left(\varepsilon_{n-k} \right)$ of a LPH, for each $n \in \mathbb{Z}$ and $a_k \in \mathcal{L}(H)$, for any $k \geq 1$, the stationarity is held as long as $\varepsilon = \left\lbrace \varepsilon_n, \ n \in \mathbb{Z} \right\rbrace$ is a $H$-valued SWN and $\displaystyle \sum_{k=1}^{\infty} \left\| a_k \right\|_{\mathcal{L}(H)}^{2} < \infty$. Building on the early works by Bosq (1991) and Mourid (1993), the invertibility of LPH was proved in Merlev\`ede (1995) if and only if $1 - \displaystyle \sum_{j=1}^{\infty} z^j \left\| a_j \right\|_{\mathcal{L}(H)} \neq 0$, for $\left| z \right| < 1$. Asymptotic properties were subsequently derived in Merlev\`ede (1996). 

Merlev\`ede (1997) provided a Markovian representation of stationary and invertible LPH in a subspace $H_\beta = \left\lbrace X:~ \left\| X \right\|_{H_{\beta}} = \displaystyle \sum_{k=1}^{\infty} \beta_k \left\| X_k \right\|_{H}^{2} < \infty \right\rbrace$, being $\beta = \left\lbrace \beta_k, \ k \geq 1 \right\rbrace$ a strictly positive decreasing and summable sequence. Let us define the $H_{\beta}$-random variables $Y_n = \left(X_n, X_{n-1}, \ldots, X_{n-p-1}, X_{n-p}, \ldots \right)^{\prime}$ and $e_{n} = \left(\varepsilon_n, 0,0,\ldots \right)^{\prime}$, for each $n \in \mathbb{Z}$. A strongly-consistent plug-in predictor was derived in Merlev\`ede (1997), by estimating the operator

\[R=
\begin{blockarray}{cccccc}
\begin{block}{(ccccc)c} 
\rho_1 & \rho_2 & \ldots & \rho_p & \ldots \\ 
I_H & 0_H & \ldots & 0_H & \ldots \\ 
\vdots & \vdots & \ddots & \vdots & \ddots \\ 
0_H & 0_H & \ldots & I_H & \ldots & \longrightarrow \text{\textbf{p-th row}}\\ 
\vdots & \vdots & \ddots & \vdots & \ddots &  \\
\end{block}
\end{blockarray}
 \]

\noindent under $\left\| R \right\|_{\mathcal{L} (H_\beta)} < 1$ and ${\rm E} \left\lbrace \left\| Y_0 \right\|_{H_\beta}^{4} \right\rbrace < \infty$. Mas (2002) studied the weak-convergence for the empirical autocovariance and cross-covariance operators of LPH. In particular, under ${\rm E} \left\lbrace \left\| \varepsilon_0 \right\|_{H}^{4} \right\rbrace < \infty$ and $\displaystyle \sum_{k=1}^{\infty} \left\| a_k \right\|_{\mathcal{L}(H)} < \infty$,
\begin{equation}
\sqrt{n} \begin{pmatrix} C_{n,0} - C_0  \\C_{n,1} - C_1 \\ \vdots  \\ C_{n,h} - C_h \end{pmatrix} \longrightarrow^{w} \mathcal{N} \left(0, \Sigma \right), \quad C_h = {\rm E} \left\lbrace X_0 \otimes X_h \right\rbrace, \quad C_{n,h} = \frac{1}{n-h} \displaystyle \sum_{i=0}^{n-h-1} X_i \otimes X_{i+h}, \quad h \in \mathbb{N}.
\end{equation}

MAH(q) and ARHMAH(p,q) processes, with $p$ and $q$ greater than one, as a particular case of LPH, were defined in Bosq \& Blanke (2007) as
\begin{equation}
X_n =  \varepsilon_n + \displaystyle \sum_{k=1}^{q} l_k \left( \varepsilon_{n-k} \right), \quad l_k \in \mathcal{L}(H), \quad \left\| l_{k} \right\|_{\mathcal{L}(H)} < 1,
\end{equation}
and
\begin{equation}
X_n =  \varepsilon_n + \displaystyle \sum_{j=1}^{p} \rho_j \left( X_{n-j} \right) + \displaystyle \sum_{k=1}^{q} l_k \left( \varepsilon_{n-k} \right), \quad l_k,~\rho_j \in \mathcal{L}(H),~\left\| l_{k} \right\|_{\mathcal{L}(H)} < 1,~\left\| \rho_{j} \right\|_{\mathcal{L}(H)} < 1,
\end{equation}

\noindent respectively, for each $n \in \mathbb{Z}$ and $k=1,\ldots,q,~j=1,\ldots,p$. LPH in a wide sense, when $\left\lbrace a_j, \ j \geq 1 \right\rbrace$ are allowed to be unbounded, were studied in Bosq (2007) and Bosq \& Blanke (2007). In that framework, linear closed subspaces of $L_{H}^{2} \left(\Omega, \mathcal{A}, \mathcal{P} \right)$ (see Fortet, 1995) are crucial to extend aspects as Wold decomposition, orthogonality and Markovianess. Unlike the estimation of an ARH(1) process, troubles  in the estimation of the operator $l$ of a MAH(1) process  arise from the non-linear behaviour of the moment equation. We may cite Turbillon et al. (2008), where  the estimation of the MAH(1) model $X_n = \varepsilon_n + l \left( \varepsilon_{n-1} \right)$, being $l \in \mathcal{K}(H)$ under $\left\|D C^{-1} \right\|_{\mathcal{L}(H)} < 1/2$ and $\left\|D^{*} C^{-1} \right\|_{\mathcal{L}(H)} < 1/2$, was reached. A special framework was introduced in Wang (2008), where a real-valued
non-linear ARIMA(p,d,q) model was modified, in a manner that functional MA coefficients were included:
\begin{equation}
X_n + \displaystyle \sum_{j=1}^{p} \rho_j X_{n-j} = \varepsilon_n +  \displaystyle \sum_{k=1}^{q} f_k \left( X_{n-k -d} \right) \varepsilon_{n-k}, \quad n \in \mathbb{Z}, \label{eq_27}
\end{equation}

\noindent being $\left\lbrace f_k, \ k \geq 1 \right\rbrace$ a set of arbitrary univariate functions. Forecasting of the Chinese Consumer Price Index, which monthly collects prices paid by middle-class consumers for a standard basket of goods and services (e.g., fuel, oil, milk, drugs, etc), was achieved in Chen et al. (2016), adopting smooth functions as functional MA coefficients in equation (\ref{eq_27}).

Furthermore, a survey about the asymptotic properties of LPH, derived in the above-referred works by Merlev\`ede (1995, 1996, 1997), was achieved in Bosq (2000) and Bosq \& Blanke (2007). Dedecker \& Merlev\`ede (2003) extended the conditional central limit theorem (see Dedecker \& Merlev\`ede, 2002) to LPH.  Useful tools proposed by Hyndman \& Shang (2008), such as visualization and outlier detection, can be applied to observed ARMAH processes, obeying a functional linear model. Outlier detection in French male age-specific mortality data was also achieved in that work. Bosq (2010) addressed the structure of tensorial products for ARMAH models, when the innovations are assumed to be martingale difference functional sequences, by using the linear close subspace theory.

% ******************************

% ******************************

% ******************************

% ******************************

\section{Banach-valued autoregressive processes}

The study of functional time series, with values in a real separable Banach space $B$, is reviewed in this Section. The Kuelb's Lemma (see Kuelbs, 1976) plays a crucial role on the derivation of the estimation results in this framework. Given a real separable Banach space $\left( B, \left\| \cdot \right\|_{B} \right)$, the Kuelb's Lemma proves that there exists an inner product $\langle \cdot, \cdot \rangle_0$ on $B$, with its associated norm $\left\| \cdot \right\|_0$ weaker than $\left\| \cdot \right\|_{B}$, providing a dense and continuous embedding $B \hookrightarrow H$, where $H$ is the completation of $B$ under $\left\| \cdot \right\|_0$. Specifically, in the ARB(1) context, the componentwise estimation of the autocorrelation operator is achieved in Labbas \& Mourid (2002), by considering the corresponding orthonormal basis, in the Hilbert space $H$ associated with $B$ by a continuous embedding. Strong-consistency of the formulated estimator, in the norm of bounded linear operators on $H$, was also derived.

Simultaneously to the early work by Bosq (1991), Pumo (1992, 1998) considered a particular case of the referred framework, adopting the Banach space of continuous functions on the interval $\left[0,1 \right]$, so-called $\mathcal{C} = \mathcal{C} \left([0,1] \right))$. Particularly, the $AR\mathcal{C}(1)$ process was formulated as follows:
 \begin{equation}
 X_n = \rho \left(X_{n-1} \right) + \varepsilon_n, \quad X_n,~\varepsilon_n \in \mathcal{C}, \quad \rho (X) (t) = \displaystyle \int_{0}^{1} r(s,t) X(s) ds, ~X \in \mathcal{C}, \quad \left\| r \right\|_{\mathcal{C} \left([0,1]^2 \right)} < 1, \quad n \in \mathbb{Z}. \label{eq_34}
 \end{equation}

Pumo (1992, 1998) consider  the real separable Hilbert space $H=L^2 \left([0,1], \beta_{[0,1]}, \lambda \right)$, where $\lambda$ denotes the Lebesgue measure, and the corresponding continuous extension of $\rho$ to $\rho^{\prime}$, defined on $L^2 \left([0,1], \beta_{[0,1]}, \lambda \right)$, such that $\left\| \rho^{\prime} \right\|_{\mathcal{L}(L^2 \left([0,1], \beta_{[0,1]}, \lambda \right))} < 1$. Specifically, an ARH(1) process $X^{\prime}$, associated with the ARB(1) process  $X$, is defined, from projection into an orthonormal basis $\left\lbrace e_{j}, \ j\geq 1 \right\rbrace$ of  $L^2 \left([0,1], \beta_{[0,1]}, \lambda \right)$. The restriction to $\mathcal{C}$ of the componentwise estimator of $\rho^{\prime}$, computed in terms of the eigenvectors of the autocovariance operator of $X^{\prime}$, provides an estimator of $\rho$, under strongly-mixing and Cramer conditions. Strong-consistency of the formulated estimator, in the norm of bounded linear operators on $\mathcal{C}$, was derived in Bosq (2000). The natural extension of $AR\mathcal{C}(1)$ to $AR\mathcal{C}(p)$ processes, with $p$ greater than one, was firstly proposed in Mourid (1993, 1996). In that works, the characterization of some continuous time processes (such as Ornstein-Uhlenbeck processes) as $AR\mathcal{C}(1)$ processes, was also provided. As commented in Section 2, a continuous-time stochastic process $\xi = \left\lbrace \xi_t, \ t \geq 0 \right\rbrace$ is turned into a set of  functional  random variables $X = \left\lbrace X_n(t) = \xi_{n\delta + t}, \ n \in \mathbb{Z} \right\rbrace$, with $X_n(t)$ taking values in the interval $\left[ 0, \delta \right]$, for each $n \in \mathbb{Z}$ and $\delta > 0$. A method to estimate the periodicity $\delta$ for an $AR\mathcal{C}(p)$ process was developed by Benyelles \& Mourid (2001), by using the work by Martin (1982). 

In the particular Banach space $\mathcal{C}$, a non-plug-in predictor of a stationary $AR\mathcal{C}(1)$ process $X$, based on the projection into an orthonormal basis of the Reproducing  Kernel Hilbert space, associated with the autocovariance operator of $X^{\prime}$, was derived in Mokhtari \& Mourid (2003) via Parzen's approach (see also Parzen, 1961). In the case of the eigenvalues $\left\lbrace C_{j}^{\prime}, \ j \geq 1 \right\rbrace$ of the autocovariance operator of $X^{\prime}$ are unknown, they assumed that $C_{n,1} > \ldots > C_{n,k_n} > 0~a.s.$, being $\left\lbrace C_{n,j}^{\prime}, \ j \geq 1 \right\rbrace$ the eigenvalues of the empirical estimator of the autocovariance operator of $X^{\prime}$, for a suitable truncation parameter $k_n$. Under mild assumptions over the eigenvectors of the empirical estimator of the autocovariance operator of $X^{\prime}$,  Mokhtari \& Mourid (2003) derived a consistent non-plug-in predictor. As discussed in that work, the spectral decomposition of cross-covariance operator $D$ is not needed in the Parzen-approach-based framework, as required in Pumo (1992, 1998). The equivalence of the asymptotic behaviour of both approaches is also derived.

%\textbf{\textcolor{blue}{[*** FALTA REESCRIBIR]}} As follow-up to the analysis of heavy-tailed functional-valued distributions analysed by Haan and Lin (2001), the extreme value behaviour of space-time processes, valued in the space $\mathcal{D}^d = \mathcal{D} \left([0,1]^d \right)$, defined as the space of real-valued c\`adl\`ag functions (right continuous functions with left limits) on $\left[0,1\right]^d$, was studied in Davis and Mikosch (2008). In particular, the forecasting of the maximum ozone level was performed as an $AR\mathcal{D}^d(1)$ process.  \textbf{\textcolor{blue}{[***]}}
%\textcolor{blue}{
%\begin{itemize}
%\item \textbf{Raluca Ralan, Regular Variation of Infinite Series of Processes with Random Coefficients}. Dicen que es una serie temporal $\mathcal{D}$-valued, con los coeficientes del proceso lineal siendo funciones deterministicas, y sus resultados de heavy-tailed son validos para cadlag processes indexados en $[0,1]$. Ademas ellos lo extienden a cuando los coeficientes son  random processes.
%\item \textbf{Richard A. Davis et al., Max-stable processes for modeling extremes observed in
%space and time.} "We mention two approaches proposed in the literature concerning the analysis and quantification of the extremal behavior of
%processes observed both in space and time. One idea can be found in Davis and Mikosch (2008), who study the extremal
%properties of a moving average process of spatial fields, where the coefficients and the white-noise process depend on
%space and time."
%\end{itemize}}

When the  space $\mathcal{D}= \mathcal{D}\left(\left[0,1\right] \right)$ (defined as the space of right continuous functions with left limits on $\left[0,1\right]$) is adopted, El Hajj (2011, 2013) deeply addressed the estimation and prediction of $AR\mathcal{D}(1)$ processes, when $\mathcal{D}$ is equipped with the Skorokhod topology. While $\mathcal{D}$ is a non-separable Banach space under the supremum norm, the Skorokhod topology provides the separability property to the metric space. The asymptotic properties of this special class of functional autoregressive processes were provided in El Hajj (2011). When the autocorrelation operator takes values in $\mathcal{C}$, by considering the continuous embedding into $L^2 \left([0,1], \beta_{[0,1]}, \lambda \right)$ already commented, the estimation and prediction of $AR\mathcal{D}(1)$ and $MA\mathcal{D}(1)$ processes was addressed  by El Hajj (2013).

Additionally, some general results, useful for the development of the  theory  of linear processes in Banach spaces, are now noted. Guillas (2000) generalized some notions about non-causality to Banach-valued context. Extending the results by Bosq (2000), Dehling \& Sharipov (2005) provided the asymptotic properties of functional second-order moments of an ARB process, allowing weakly dependent innovation processes. The sieves method already detailed in Section 4, and initially proposed on the ARH(1) forecasting by Bensmain \& Mourid (2001), was applied in Rachedi (2005), on the estimation and prediction of an ARB(1) process, based on the dual space of $B$. Rates of convergences of the formulated estimator were provided in the mentioned work, when $\rho$ is assumed to be a $p$-summable operator; i.e., for each $X_0,\ldots,X_{n-1}$ on $B$, there exist $p \in \left(1,\infty \right)$ and constant $c > 0$ such that
 \begin{equation}
 \left( \displaystyle \sum_{i=0}^{n-1} \left\| \rho \left( X_i \right) \right\|^p \right)^{1/p} \leq c \displaystyle \sup_{\left\| X^{*} \right\|\leq 1} \left( \displaystyle \sum_{i=0}^{n-1} \left| \left(X^{*}, X_i \right) \right|^p \right)^{1/p}, \quad X^{*} \in B^{*},\label{eq__36}
 \end{equation}
 \noindent providing the $p$-integrable norm $\pi_p \left( \rho \right)$, as the minimum value of constant $c$ which verifies the equation (\ref{eq__36}). If $\Pi_p (B)$ is the set of  $p$-summable operators on $B$, when $H$ is a Hilbert space, then $\Pi_2 (B)$ coincides with the space of Hilbert-Schmidt operators on $H$, with the Hilbert-Schmidt norm  $\pi_2$. A decomposition of $\rho$ in terms of so-called Schauder and Markushevich bases, was also obtained. In the Skorokhod space $\mathcal{D}$ studied in El Hajj (2011, 2013), Blanke \& Bosq (2014) analysed the intensity of jumps of $\mathcal{D}$-valued linear processes, providing some limit theorems for $ARMA\mathcal{D}(1,1)$ processes, when both fixed and random number of jumps are regarded. As discussed in Blanke \& Bosq (2014), the estimation of these jumps can be used in the prediction of compound Poisson processes, which are used, e.g., for forecasting the payments, at fixed instants, refunded to the holders of an insurance policy. An extension of ARMA processes to general complex separable Banach spaces was proposed in Spangenberg (2013). Firstly, the stationarity of the ARMAB(1,q) process
\begin{equation}
X_n = \varepsilon_n + \rho_1 \left( X_{n-1} \right) + \displaystyle \sum_{k=1}^{q} l_k \left( \varepsilon_{n-k} \right), \quad \rho_1,~l_k \in \mathcal{L}(B),~k=1,\ldots,q,
\end{equation}
\noindent was proved, under the hyperbolic property over $\rho_1$ (i.e., $\sigma\left(\rho_1 \right) \cap \mathbb{S} = \left\lbrace \emptyset \right\rbrace$, being $\mathbb{S}$ the unit circle and $\sigma\left(\rho_1 \right)$ the spectrum of $\rho_1$) and $\displaystyle \log^{+}$-moment conditions. Stationarity of ARMAB(p,q) processes was subsequently derived by a $B^p$-valued ARMA(1,q) representation. From results in Soltani \& Hashemi (2011), where PCARH(1) processes were introduced, Parvardeh et al. (2017) derived the asymptotic properties for Banach-valued autoregressive periodically correlated processes of order one (PCARB(1) processes).

\section{Non-parametric functional time series framework}
\label{sec:6}

Let us see the main references in the context of non-parametric functional time series and functional linear regression, when both explanatory and response variables take values in a space of functions.

As a functional extension of the multivariate-based electricity consumption forecasting approach developed by Poggi (1994), a non-parametric kernel-based predictor was formulated in Besse et al. (2000):
\begin{equation}
\widehat{X}_{n}^{h_n} = \widehat{\rho}_{h_n} \left( X_{n-1} \right), \quad \widehat{\rho}_{h_n} \left( X_{n-1} \right) = \frac{\displaystyle \sum_{i=0}^{n-2} \widehat{X}_{i+1} K \left(\frac{\left\| \widehat{X}_i - X_{n-1} \right\|_{L^2 \left(\left[a,b \right] \right)}^{2}}{h_n} \right)}{\displaystyle \sum_{i=0}^{n-2} K \left(\frac{\left\| \widehat{X}_i - X_{n-1} \right\|_{L^2 \left(\left[a,b \right] \right)}^{2}}{h_n} \right)}, \quad \widehat{X}_i = argmin \left\| D \widehat{X}_i \right\|_{L^2 \left(\left[a,b \right] \right)}^{2},\label{a100}
\end{equation}

\noindent being $K$ the usual Gaussian kernel and $D$ a $d$-th order differential operator. Forecasting of climatological time series, so-called ENSO time series, was therein achieved. Cuevas et al. (2002) addressed the strong-consistency estimation of the underlying linear operator of a linear regression, when both explanatory and response variables are assumed to be $H$-valued random variables, with $H = L^2 \left(\left[a,b \right] \right)$. In particular, the design is given by the triangular array $\left\lbrace X_{i,n} \left(t \right), \ 1 \leq i \leq n \right\rbrace$, providing the model $Y_{i,n} = \Psi \left( X_{i,n} \right) + \varepsilon_{i,n}$, with $ X_{i,n} \in L^2 \left(\left[a,b \right] \right)$ and $Y_{i,n} \in L^2 \left(\left[c,d \right] \right)$, under $\Psi \in \mathcal{L} \left(L^2 \left(\left[a,b \right] \right),L^2 \left(\left[c,d \right] \right)\right)$.

Antoniadis et al. (2006) introduced the  two-steps prediction approach so-called kernel wavelet functional (KWF) method, where strongly-mixing conditions are imposed. An expansion of strictly stationary functional time series into a discrete wavelet basis $\left\lbrace \psi_{k}^{J}, \ k=0,\ldots,2^J - 1 \right\rbrace$, at scale $J$, is achieved, and the forecasting of $\widehat{X}_{n} = {\rm E} \left\lbrace X_n | X_{n-1}, \ldots, X_0 \right\rbrace$, for each $n \in \mathbb{Z}$, was then performed by
\begin{equation}
\widehat{X}_{n}^{J}  \left( \cdot \right) = \displaystyle \sum_ {k=0}^{2^J - 1} \widehat{\xi}_{n,k}^{J} \psi_{k}^{J} \left( \cdot \right), \quad \widehat{\Xi}_{n} = \frac{\displaystyle \sum_{i=0}^{n-2} K \left( D \left(P \left( \Xi_n \right),D \left(P \left( \Xi_i \right) \right)\right) / h_n \right) \Xi_{i+1}}{1/n +  \displaystyle \sum_{i=0}^{n-2} K \left( D \left(P \left( \Xi_n \right),D \left(P \left( \Xi_i \right) \right)\right) / h_n \right)}, \label{eq__37}
\end{equation}

\noindent where $\widehat{\Xi}_{n} = \left\lbrace \widehat{\xi}_{n,k}^{J}:~k=0,1,\ldots,2^J - 1 \right\rbrace$ denotes, for each $n \in \mathbb{Z}$, the set of predicted scaling coefficients,  at scale $J$, being $P \left( \Xi_i \right)$ the set of wavelet coefficients derived by the so-called pyramid algorithm (see Mallat, 1989), for any $i=0,1,\ldots,n-1$. Distance $D \left( \cdot, \cdot \right)$ in (\ref{eq__37}), for a two set of discrete wavelet coefficients $\left\lbrace \theta_{j,k}^{i}, \ i=1,2 \right\rbrace$, at scale $j=j_0,\ldots,J-1$ and location $k = 0,\ldots,2^j - 1$, is given by
\begin{equation}
D \left( \theta^{1}, \theta^2 \right) = \displaystyle  \sum_{j=j_0}^{J - 1} 2^{-j / 2} d_j \left(\theta^{1}, \theta^2 \right), \quad d_j \left(\theta^{1}, \theta^2 \right) = \left( \displaystyle \sum_{k=0}^{2^j -1} \left( \theta_{j,k}^{1} - \theta_{j,k}^{2} \right)^2 \right)^{1/2}.
\end{equation}

Pointwise prediction intervals were also built. See also Cugliari (2011), where a continuous wavelet transform (CWT) is also considered. Slightly modifications have been also proposed in the French electricity consumption forecasting addressed by  Antoniadis et al. (2012), when the hypothesis of stationarity is not held. We may also cite the work by Antoniadis et al. (2009), in which a method on the selecting the properly bandwidth $h_n$, for kernel-based forecasting of functional time series, was developed.

Functional versions of partial least-squares regression and principal component regression (FPLSR and FPCR, respectively) were formulated in Reiss \& Ogden (2007). In this work, a functional smoothing-based approach to signal regression was adopted, where decompositions in terms of B-spline bases and roughness penalties are involved. Let us consider a general functional linear regression model, when  Hilbert-valued response and $\mathcal{F}$-valued explanatory variables are considered, when $\mathcal{F}$ is defined as a general functions space, equipped with a semi-metric $d$ and its associated topology $\mathcal{T}_{\mathcal{F}} \left( X,t \right) = \left\lbrace X_1 \in \mathcal{F}:~d \left( X_1,X \right) \leq t \right\rbrace$. In this framework, a non-parametric kernel-based estimator of the underlying regression operator was derived in Ferraty et al. (2012) as follows, for each $i=0,1,\ldots,n-1$:
\begin{equation}
Y_i = \Psi \left( X_i \right) + \varepsilon, \quad \widehat{Y_n} = \widehat{\Psi}_{h_n}\left( X_n \right),\quad \widehat{\Psi}_{h_n}\left( X_n \right)= \frac{\displaystyle \sum_{i=0}^{n-2} X_{i+1} K \left(\frac{d \left(X_i, X_{n-1} \right)}{h_n} \right)}{\displaystyle \sum_{i=0}^{n-2} K \left(\frac{d \left(X_i, X_{n-1} \right)}{h_n} \right)},
\end{equation}

\noindent being $K$ a Gaussian kernel (see Ferraty \& Vieu, 2006, about the choice of a semi-metric $d$).

% ******************************

% ******************************

\section{ARH(1) strongly-consistent diagonal componentwise  parameter estimators}
\label{sec:8}

In this section, we restrict our attention to the case where $C$ and $\rho$ admit a diagonal spectral decomposition in terms of the common eigenvectors system $\left\lbrace \phi_j, \ j \geq 1 \right\rbrace$, since in that case, an important dimension reduction is achieved. This spectral diagonalization can be reached under a wide range of scenarios, leading to a sparse representation of kernels of the associated integral operators (see, e.g., Meyer \& Coifman, 1997, for the case of spectral diagonalization of Calder\'on-Zygmund operators in terms of wavelet bases, in Besov spaces). As discussed in \'Alvarez-Li\'ebana et al. (2017), this particular scenario is naturally obtained when $\rho$ and $C$ are linked by a continuous function. Both scenarios, when $\left\lbrace \phi_j, \ j \geq 1 \right\rbrace$ are known and unknown, are now covered, providing the corresponding strong-consistency results, for the formulated diagonal componentwise estimators of $\rho$, in the norm of bounded linear operators. Proof details are provided in Sections S.1-S.2 of the Supplementary Material. Under a non-diagonal scenario, an almost sure upper-bound for the $S(H)$ norm of the error associated with the diagonal componentwise estimator of $\rho$, when eigenvectors of C are unknown, is provided in Section S.3 of the Supplementary Material.

% ******************************

% ******************************

\subsection{ARH(1) model: diagonal framework}
\label{sec:31}

 Let  $H$ be a real separable Hilbert space, and let  $X = \left\lbrace X_n,\ n \in \mathbb{Z} \right\rbrace $
 be  a zero-mean stationary ARH(1) process on the basic probability space $(\Omega,\mathcal{A},P)$, satisfying:
\begin{equation}
X_{n } = \rho \left(X_{n-1} \right) + \varepsilon_{n}, \quad \rho \in \mathcal{L}(H), \quad \left\| \rho \right\|_{\mathcal{L}(H)} < 1,\quad n \in
\Z,\label{24bb}
\end{equation}

\noindent  when $H$-valued
innovation process $\varepsilon= \left \lbrace \varepsilon_{n}, \ n\in \mathbb{Z}\right\rbrace$ is assumed to be SWN, and to be  uncorrelated with $X_0$, with 
$\sigma^{2}_{\varepsilon}= {\rm E} \left\lbrace
\|\varepsilon_{n}\|_{H}^{2} \right\rbrace <\infty,$  for all $n\in
\mathbb{Z}.$ In addition, let us consider the following assumptions:

\vspace{0.3cm}

\noindent \textbf{Assumption A1.} The autocovariance operator
$C={\rm E} \left\lbrace X_{n}\otimes X_{n} \right\rbrace,$ for every
$n \in \mathbb{Z},$ is a strictly positive and self-adjoint operator, in the  trace class, with $Ker\left( C \right) = \left\lbrace \emptyset \right\rbrace$. Its eigenvalues $\left\lbrace
C_{j}, \ j\geq 1\right\rbrace $ then satisfy
\begin{equation}
 \displaystyle \sum_{j=1}^{\infty} C_j < \infty,\quad  C_1  > \ldots > C_j > C_{j+1} > \ldots > 0, \quad C (f)(g)=\displaystyle \sum_{j=1}^{\infty}C_{j}\left\langle \phi_{j},f\right\rangle_{H}\left\langle \phi_{j},g\right\rangle_{H},~ \forall f,g\in H.
\label{eqautocov2}
\end{equation}

\vspace{0.3cm}

\noindent \textbf{Assumption A2.} The autocorrelation operator
$\rho$ is a self-adjoint and Hilbert-Schmidt operator, admitting
the following diagonal spectral decomposition:
\begin{equation}
\rho (f)(g)= \displaystyle \sum_{j=1}^{\infty} \rho_j \left\langle \phi_{j},f\right\rangle_{H}\left\langle \phi_{j},g\right\rangle_{H}, \quad \displaystyle \sum_{j=1}^{\infty} \rho_{j}^{2} < \infty, \quad  \forall f,g\in H,
\label{12}
\end{equation}
\noindent where $\left\lbrace \rho_j,\  j \geq 1 \right\rbrace $ is
the system of eigenvalues of $\rho,$
with respect to the orthonormal system $\left\lbrace \phi_j,\  j \geq 1 \right\rbrace $.

\vspace{0.3cm}

Under \textbf{Assumptions A1-A2}, the cross-covariance operator $D = \rho C ={\rm E} \left\lbrace X_{n}\otimes X_{n+1} \right\rbrace$, for each $n \in \mathbb{Z}$, can be also diagonally decomposed, with regard to the eigenvectors of $C$, providing a set of eigenvalues $\left\lbrace D_j = \rho_j C_j, \ j \geq 1 \right\rbrace$. Moreover, projections of (\ref{24bb}) into  $\left\lbrace \phi_j,\ j \geq 1 \right\rbrace $ lead to the stationary zero-mean AR(1) representation, under $\left\| \rho \right\|_{\mathcal{L}(H)} = \displaystyle \sup_{j \geq 1} \left| \rho_j \right| < 1$:
\begin{equation}
X_{n,j}=\rho_{j}X_{n-1,j}+\varepsilon_{n,j},\quad X_{n,j} = \left\langle
X_{n},\phi_{j}\right\rangle_{H},~\varepsilon_{n,j}=\left\langle \varepsilon_{n},\phi_{j}\right\rangle_{H}, \quad \rho_j \in \mathbb{R},~\left| \rho_j \right| < 1, \quad j\geq 1, ~n \in \mathbb{Z}.
\label{14}
\end{equation}

\subsection{Diagonal strongly-consistent estimator when the eigenvectors of $C$ are known}
\label{sec:32}

In the following of this subsection, \textbf{Assumption A3} will be imposed, jointly with $\langle X_0, \phi_j \rangle_H \neq 0$ a.s. for any $j \geq 1$ (so-called \textbf{Assumption A4}). In the case of $\left\lbrace \phi_j,\ j \geq 1 \right\rbrace $  are assumed to be known, and under \textbf{Assumptions A1-A2}, the following estimators of $C$ and $D$, based on the empirical moment-based estimation of their eigenvalues, is developed, for any $n \geq 2$ and $j \geq 1$:
\begin{eqnarray}
\widehat{C}_n = \displaystyle \sum_{j=1}^{\infty} \widehat{C}_{n,j} \phi_j \otimes \phi_j,~\widehat{C}_{n,j} = \frac{1}{n} \displaystyle \sum_{i=0}^{n-1} X_{i,j}^{2},~\widehat{D}_n = \displaystyle \sum_{j=1}^{\infty} \widehat{D}_{n,j} \phi_j \otimes \phi_j,~\widehat{D}_{n,j} = \frac{1}{n-1} \displaystyle \sum_{i=0}^{n-2} X_{i,j}X_{i+1,j}. \label{65aa} 
\end{eqnarray}

From  \textbf{Assumption A4}, let us consider the diagonal  componentwise estimator of $\rho$
\begin{eqnarray}
\widehat{\rho}_{k_n} &=& \displaystyle \sum_{j=1}^{k_n} \widehat{\rho}_{n,j} \phi_j \otimes \phi_j, \quad \widehat{\rho}_{n,j} = \frac{\widehat{D}_{n,j}}{\widehat{C}_{n,j}}  = \frac{n}{n-1} \frac{\displaystyle \sum_{i=0}^{n-2} X_{i,j}X_{i+1,j}}{\displaystyle \sum_{i=0}^{n-1} X_{i,j}^{2}}, \quad j \geq 1,~n \geq 2. \label{48}
\end{eqnarray}

Under \textbf{Assumptions A1-A4}, the following proposition provides the strong-consistency, in the norm of $\mathcal{L}(H)$, of the estimator (\ref{48}), as well as of its associated ARH(1) plug-in predictor, in the underlying Hilbert space. Proof details are provided in Sections S.1-S.2 in the Supplementary Material.

\begin{proposition}
\label{pr2m}
\textit{Under \textbf{Assumptions A1--A4}, for a truncation parameter $k_n < n$, with $\displaystyle \lim_{n \to \infty} k_n= \infty$,
\begin{equation}
\frac{n^{1/4}}{\left(\ln(n) \right)^{\beta}}\left\| \widehat{\rho}_{k_n}- \rho\right\|_{\mathcal{L}(H)} \longrightarrow^{a.s.} 0, \quad 
\| \left( \widehat{\rho}_{k_n} -\rho \right) (X_{n-1})\|_{H} \longrightarrow^{a.s.} 0,\quad n\rightarrow \infty.
\end{equation}}
\end{proposition}

\subsection{Diagonal strongly-consistent estimator when the eigenvectors of $C$ are unknown}
\label{sec:33}

 In the case of $\{\phi_{j},\ j\geq 1\}$ are unknown, as is often in practice, $C_n = \frac{1}{n} \displaystyle \sum_{i=0}^{n-1} X_i \otimes X_i$ admits a diagonal spectral decomposition in terms of  $\{C_{n,j},\ j\geq 1\}$ and $\{\phi_{n,j},\ j\geq 1\}$, satisfying, for each $n \geq 2$:
 \begin{eqnarray}
C_{n} \left(\phi_{n,j} \right) &=& C_{n,j}\phi_{n,j},\quad j\geq 1,\quad C_{n,1}\geq \dots\geq C_{n,n}\geq 0=C_{n,n+1}=C_{n,n+2}=\dots \label{ev1}\\
C_n &=& \frac{1}{n} \displaystyle \sum_{i=0}^{n-1} X_i \otimes X_i=\displaystyle \sum_{j=1}^{\infty} C_{n,j} \phi_{n,j} \otimes \phi_{n,j}, \quad C_{n,j} =\frac{1}{n} \displaystyle \sum_{i=0}^{n-1} \widetilde{X}_{i,j}^{2},~j \geq 1.\label{66aa}
\end{eqnarray}

In the remainder, we will denote $\widetilde{X}_{i,j} = \langle X_i, \phi_{n,j} \rangle_H$ and $\phi_{n,j}^{\prime} = \mbox{sgn}\left\langle \phi_{n,j} , \phi_{j} \right\rangle_{H}\phi_{j}$, for each $ i\in \mathbb{Z},~j\geq 1$ and $n\geq 2$, where $\mbox{sgn} \langle \phi_{n,j}, \phi_j \rangle_H= \mathbf{1}_{\langle \phi_{n,j}, \phi_j \rangle_H\geq 0}-\mathbf{1}_{\langle \phi_{n,j}, \phi_j \rangle_H< 0}$. Since $\left\lbrace \phi_{n,j} , \ j \geq 1 \right\rbrace$ is a complete orthonormal system of eigenvectors, for each $n\geq 2$, operator $D_n = \frac{1}{n-1} \displaystyle \sum_{i=0}^{n-2} X_i \otimes X_{i+1}$ admits the following non-diagonal spectral representation:
\begin{equation}
D_n = \frac{1}{n-1}\displaystyle \sum_{i=0}^{n-2} X_i \otimes X_{i+1}= \displaystyle \sum_{j=1}^{\infty} \displaystyle
\sum_{l=1}^{\infty} D_{n,j,l}^{*} \phi_{n,j} \otimes \phi_{n,l} = \displaystyle \sum_{j=1}^{\infty}\sum_{l=1}^{\infty}\frac{1}{n-1}\displaystyle
\sum_{i=0}^{n-2} \widetilde{X}_{i,j}\widetilde{X}_{i+1,l}\phi_{n,j} \otimes \phi_{n,l}, \label{37aaa}
\end{equation}
\noindent where $D_{n,j,l}^{*} = \langle D_n \left(\phi_{n,j} \right), \phi_{n,l}  \rangle_H= \frac{1}{n-1}\displaystyle
\sum_{i=0}^{n-2}\widetilde{X}_{i,j}\widetilde{X}_{i+1,l},$ for each $j,l \geq 1$ and $n \geq 2$. In particular, we will use the notation $D_{n,j}=D_{n,j,j}^{*} =  \langle D_n \left(\phi_{n,j} \right), \phi_{n,j}  \rangle_H$. The following  assumption is here deemed:

\vspace{0.3cm}

\noindent \textbf{Assumption A5.} $C_{n,k_n} > 0~a.s,$ where $k_n$ is a suitable truncation parameter $k_n < n$, with $\displaystyle \lim_{n \to \infty} k_n = \infty$.

\vspace{0.3cm}

From \textbf{Assumption A5}, the following  diagonal componentwise estimator of  $\rho$ is outlined:
\begin{equation}
\widetilde{\rho}_{k_n} = \displaystyle \sum_{j=1}^{k_n} \widetilde{\rho}_{n,j} \phi_{n,j} \otimes \phi_{n,j}, \quad \widetilde{\rho}_{n,j} =  \frac{D_{n,j}}{C_{n,j}}  = \frac{n}{n-1} \frac{\displaystyle \sum_{i=0}^{n-2} \widetilde{X}_{i,j} \widetilde{X}_{i+1,j}}{\displaystyle \sum_{i=0}^{n-1} \widetilde{X}_{i,j}^{2}}, \quad j \geq 1,~n \geq 2. \label{141}
\end{equation}

 Under  \textbf{Assumptions A1-A3} and \textbf{A5},  the strong-consistency of the diagonal estimator  $\widetilde{\rho}_{k_n}$  is reached in Proposition \ref{proposition5}, and proved in Sections S.1-S.2 of the Supplementary Material. The large-sample behaviour of (\ref{141}) is illustrated in Section S.4 of the Supplementary Material, under diagonal scenarios.

\begin{proposition}
\label{proposition5}
\textit{Let $\left\lbrace k_{n}, \ n \in \mathbb{N} \right\rbrace$ a sequence of integers such that, for $\widetilde{n}_{0}$ sufficiently large, and $\beta > \frac{1}{2},$
\begin{eqnarray}
\Lambda_{k_{n}}&=&o\left(n^{1/4}(\ln(n))^{\beta -1/2} \right), \quad k_{n}C_{k_{n}}<1,~n\geq \widetilde{n}_{0}, \label{eq__46} \\
\frac{1}{C_{k_n}} \displaystyle \sum_{j=1}^{k_n} a_j& =& \mathcal{O} \left( n^{1/4} \left(\ln(n)\right)^{-\beta} \right),\label{eq__47}
\end{eqnarray}
\noindent where $\Lambda_{k_n}=\displaystyle \sup_{1\leq j\leq k_n}(C_{j}-C_{j+1})^{-1}$, $a_1 = 2 \sqrt{2} \frac{1}{C_1 - C_2}$ and $a_j = 2 \sqrt{2} \displaystyle \max \left( \frac{1}{C_{j-1} - C_j},\frac{1}{C_{j} - C_{j+1}} \right)$, for any $2 \leq j \leq k_n$. Then, under  \textbf{Assumptions A1-A3}  and \textbf{A5},
\begin{equation}
\left\| \widetilde{\rho}_{k_n}- \rho\right\|_{\mathcal{L}(H)} \longrightarrow^{a.s.} 0,\quad \| \left(\widetilde{\rho}_{k_n} - \rho \right)(X_{n-1})\|_{H} \longrightarrow^{a.s.} 0,\quad n\rightarrow \infty. 
\label{eqboundedoperators}
\end{equation}}
\end{proposition}

When $\rho$ does not admit a diagonal spectral representation, an almost sure upper bound for the error $\left\| \widetilde{\rho}_{k_n} - \rho \right\|_{\mathcal{S}(H)}^{2}$ is provided in Section S.3 of the Supplementary Material, under \textbf{Assumption A5} and conditions imposed in Lemma 2 (see Section S.2 of the Supplementary Material).

% **********************************

% **********************************

\section{Comparative study: an evaluation of the performance}
\label{sec:comparative}

A comparative study is undertaken to illustrate the performance of the ARH(1) predictor formulated in Section 8, and those ones given by Antoniadis \& Sapatinas (2003), Besse et al. (2000), Bosq (2000) and Guillas (2001), under different diagonal, pseudo-diagonal and non-diagonal scenarios, when $\left\lbrace \phi_j, \ j \geq 1 \right\rbrace$ are unknown. Additionally to the figures displayed in this Section, more details of the numerical results obtained can  also be found in the tables included in Section S.5 of the Supplementary Material. In all of the scenarios considered, autocovariance operator $C$ is given by
\begin{equation}
 C(f) (g) = \sum_{j=1}^{\infty} C_j \langle \phi_j, f \rangle_H \langle \phi_j, g \rangle_H, \quad \phi_{j} \left( x\right)= \sqrt{\frac{2}{b-a}} \sin \left(
\frac{\pi j x}{b-a} \right), \quad f,g\in H = L^{2}((a,b)),~x\in (a,b).\label{aco}
\end{equation}

In the remaining, we fix $(a,b)=(0,4)$ and, under \textbf{Assumption A1}, $C_j = c_1 j^{-\delta_1}$, for any $j \geq 1$, being $c_1$ a positive constant. Different rates of convergence to zero of the eigenvalues of $C$ are studied, corresponding to the values of the shape parameter $\delta_{1}\in (1,2)$. Table \ref{tab:Table_Scen} and Figure \ref{fig:F0} below show the scenarios considered in the illustration of the performance of the componentwise estimators of $\rho$ compared (see equations (\ref{eq_nueva_5})-(\ref{eq_nueva_6}) above).

\setlength{\heavyrulewidth}{0.3em}
\setlength{\lightrulewidth}{0.15em}

\vspace{0.25cm}
\begin{table}[h!]
 \caption{{\small Parameters involved in the definition of $\rho$ and $C_{\varepsilon}$ (see equations (\ref{eq_nueva_5})-(\ref{eq_nueva_6}) above), under different diagonal, pseudo-diagonal and non-diagonal Gaussian generations, for any $j,h \geq 1$, $M = 50$, $\delta_2 = 11/10$ and different shape parameters $\delta_{1}$  detailed below equation (\ref{aco}) (see Tables \ref{tab:Table_Scen_A}-\ref{tab:Table_Scen_D}).}}
\centering
\vspace{-0.17cm}
\begin{tabular}{ccccc}
\toprule
 Framework & $\rho_{j,j}$ &  $\rho_{j,h}$, with $j \neq h$ & $\sigma_{j,j}^{2}$ &  $\sigma_{j,h}^{2}$, with $j \neq h$ \\
 \midrule
Diagonal & $c_2 j^{-\delta_2}$  & 0 & $C_j \left( 1 - \rho_{j,j}^{2} \right)$ & 0  \\
Pseudo-diagonal & $c_2 j^{-\delta_2}$ &  $e^{- \left| j - h \right| / W},~W = 0.2$ & $C_j \left( 1 - \rho_{j,j}^{2} \right)$ & $e^{- \left| j - h \right|^{2} / W},~W = 0.2$  \\
Non-diagonal &  $c_2 j^{-\delta_2}$ & $ \frac{1}{K} \frac{1}{\left| j - h \right|^2 + 1},~W=0.2,~\frac{1}{K} = 0.275$ & $C_j \left( 1 - \rho_{j,j}^{2} \right)$ & $e^{- \left| j - h \right|^{2} / W},~W = 0.2$   \\
\bottomrule
\end{tabular}
  \label{tab:Table_Scen}
\end{table}

In Table \ref{tab:Table_Scen} above, $c_2$ is a constant which belongs to the interval $\left(0, 1 \right)$, such that $\rho$ is a Hilbert-Schmidt operator.

\vspace{-0.1cm}
\begin{figure}[H]
  \hspace{-1.04cm} \includegraphics[width=9.42cm,height=6.1cm]{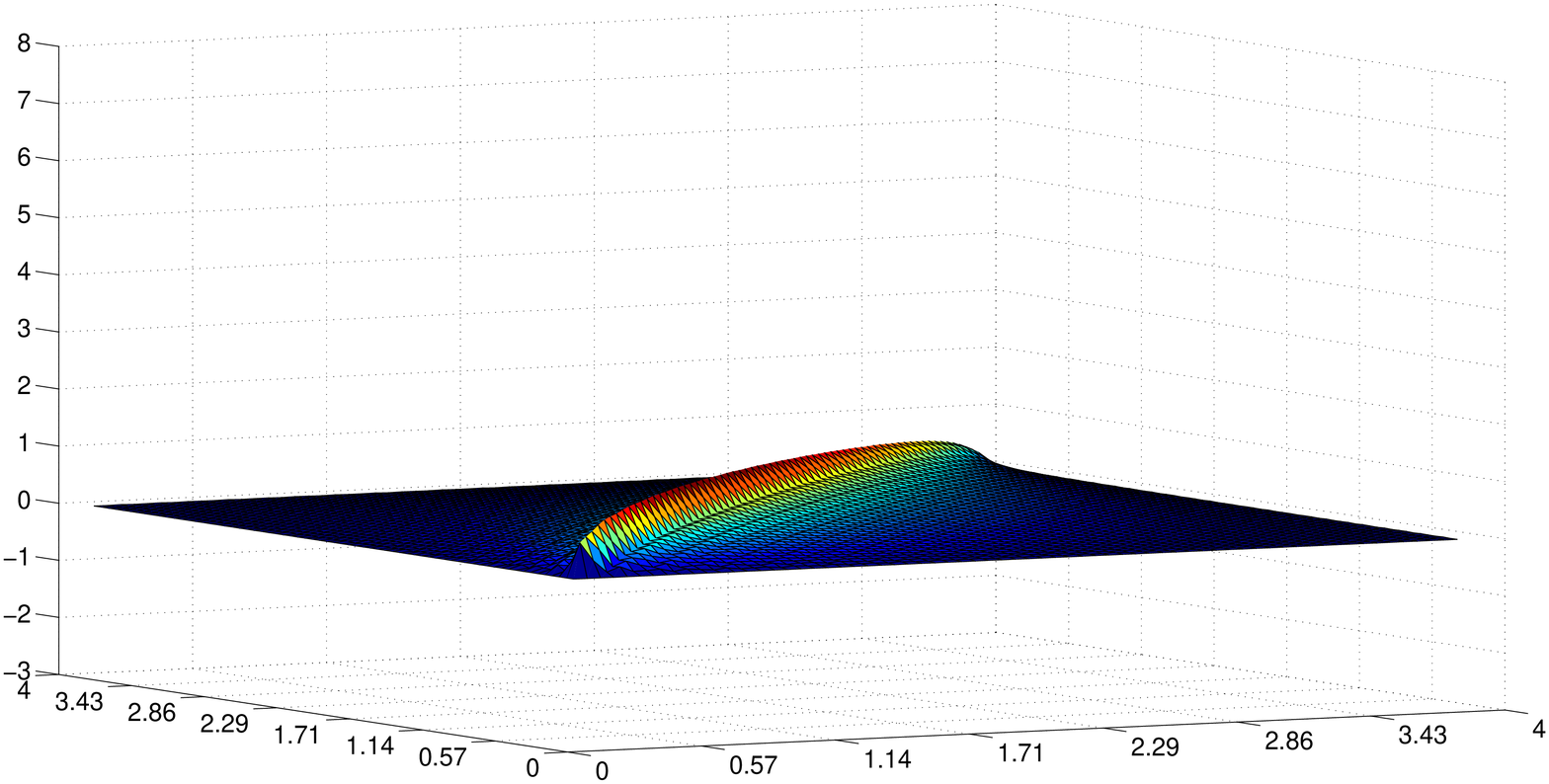}
 \hspace{-1.03cm}  \includegraphics[width=9.42cm,height=6.1cm]{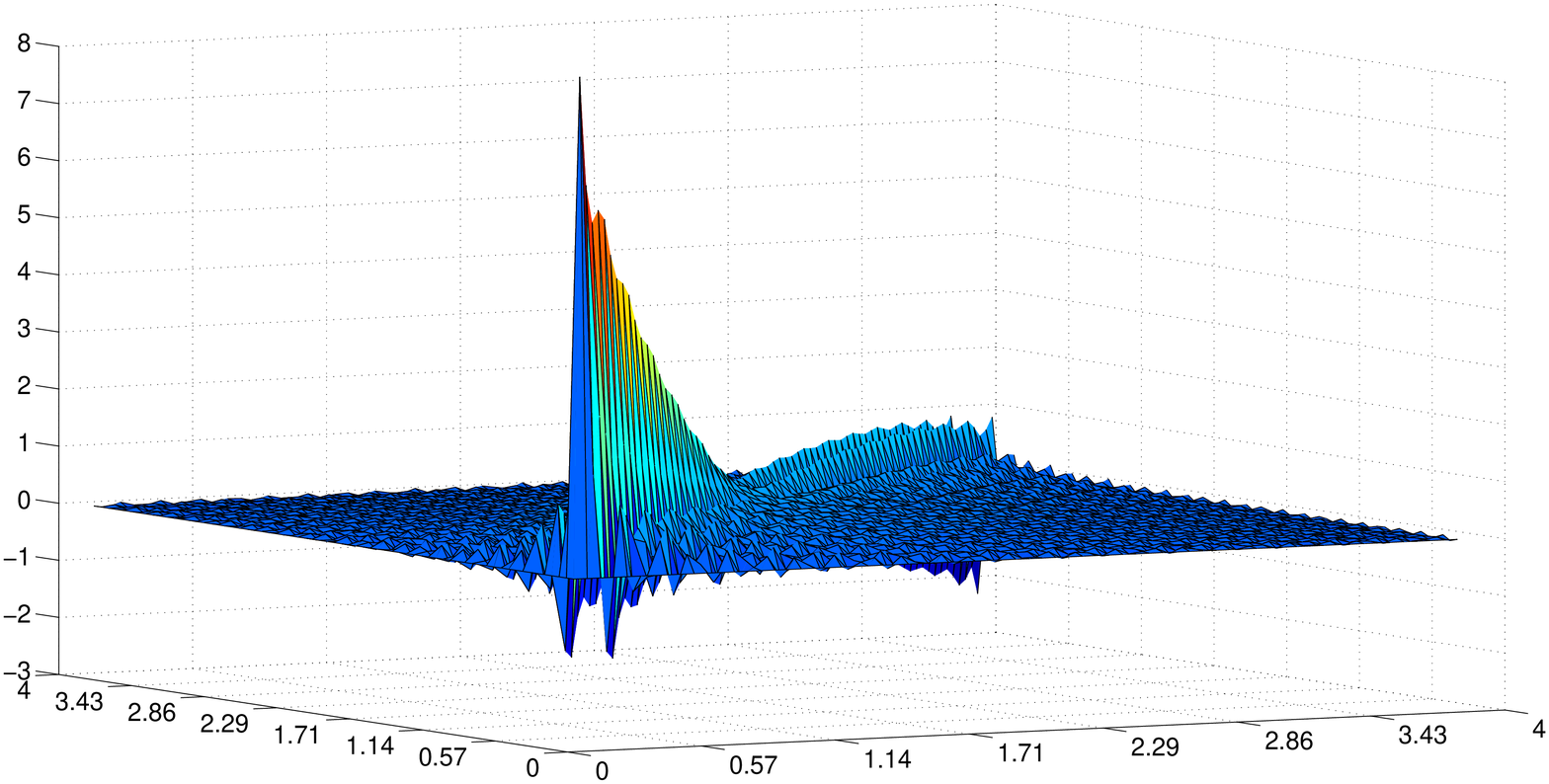}
\vspace{-1.2cm}
\caption{Operator $\psi \left( \cdot, \cdot \right)$ associated with the autocorrelation operator $\rho$, valued at the grid $\left[a,b \right] \times \left[a,b \right] $, for pseudo-diagonal (on left) and non-diagonal  (on right) scenarios (see Table \ref{tab:Table_Scen}). Discretization step $\Delta t = 0.06$ and shape parameter $\delta_1 = 3/2$ are adopted.} \label{fig:F0}
\end{figure}

\subsection{Large-sample behaviour of the ARH(1) plug-in predictors}

Large-sample behaviour of the ARH(1) plug-in predictor formulated in Section 8.3, as well as those ones in Bosq (2000) and Guillas (2001) (see equations  (\ref{eq23}) and (\ref{eq25}) above, respectively), will be illustrated. Remark that the ARH(1) plug-in predictor established in Section 8.3  will be only considered under diagonal scenarios. Strong-consistency results for the estimator $\widetilde{\rho}_{k_n}$, in the trace norm, when $\rho$ is a positive and trace operator, which does not admit a diagonalization in terms of the eigenvectors of $C$, have been recently provided in Ruiz-Medina \& \'Alvarez-Li\'ebana (2017b). See also Section S.3 of the Supplementary Material, where an almost sure upper bound for $\left\| \widetilde{\rho}_{k_n} - \rho \right\|_{\mathcal{S}(H)}^{2}$ is theoretically derived, when $\rho$ is not assumed to admit a diagonal spectral representation, nor to be a trace operator, but it is assumed to be a Hilbert-Schmidt operator.

As commented earlier (see equation (\ref{eq23}) above), \textbf{Assumptions A1}, \textbf{A3} and \textbf{A5}, jointly with the boundedness of $X_0$ and the Hilbert-Schmidt assumption of $\rho$, are required in the strong-consistency results by Bosq (2000). Condition (\ref{eq__47}) was also imposed in that work. From Bosq (2000, Example 8.6), conditions therein considered are held under any scenario in which the truncation parameter $k_{n} = \lceil \log(n) \rceil$ is adopted, under \textbf{Assumptions A1-A3} and \textbf{A5} (it can be proved as condition (\ref{eq__46}) is also verified when $k_{n} = \lceil \log(n) \rceil$). In the formulation of mean-square convergence, Guillas also considered \textbf{Assumptions A1}, \textbf{A3} and \textbf{A5}. From Guillas (2001, Theorem 2, Example 4), if the regularization sequence above-referred (see equation (\ref{eq25})) verifies $\alpha \frac{C_{k_n}^{\gamma}}{n^{\epsilon}} \leq u_n \leq \beta C_{k_n}$, for $0 < \beta < 1$ and $\alpha > 0$, with $\gamma = 1$ and $\epsilon = 0$, then the mean-square consistency is achieved under a suitable truncation parameter. In particular, if $k_{n} = \lceil e^{\prime} n^{1/\left(8 \delta_1 + 2 \right)} \rceil$, with $e^{\prime} = 17/10$, the rate of convergence in quadratic mean is of order of $n^{-\delta_1 / \left( 4 \delta_1 + 1 \right)}$. Remark that, since  $ \lceil \left( 17 / 10 \right) n^{1/\left(8 \delta_1 + 2 \right)} \rceil < \lceil \ln(n) \rceil$ for $n$ large enough, conditions (\ref{eq__46})-(\ref{eq__47}) are also verified when $k_{n} = \lceil e^{\prime} n^{1/\left(8 \delta_1 + 2 \right)} \rceil$, with $e^{\prime} = 17/10$, is fixed.

For sample sizes $n_t = 35000 + 40000 \left(t -1 \right),~t=1,\ldots,10$, and a suitable truncation parameter $k_n$,
\begin{equation}
F \left(k_n, n_t, \beta \right) = \left( \displaystyle \sum_{l=1}^{N} \mathbf{1}_{\left(\xi_{n_t,\beta},\infty \right)} \left( \left\| \left(\rho - \overline{\rho}_{k_n}^{l} \right) \left(X_{n-1}^{l} \right) \right\|_{H}^{k_n} \right) \right) /N, \label{a955}
\end{equation}
\noindent will be displayed (see Figures \ref{fig:F1}-\ref{fig:F2a} below), being $\mathbf{1}_{\left(\xi_{n_t,\beta},\infty \right)}$ the indicator function over the interval $\left(\xi_{n_t,\beta},\infty \right)$, where $\xi_{n_t,\beta}$ numerically fits the almost sure rate of convergence of $\left\| \left(\rho - \overline{\rho}_{k_n}^{l} \right) \left(X_{n-1}^{l} \right) \right\|_{H}^{k_n}$ (see Tables \ref{tab:Table_Scen_A}-\ref{tab:Table_Scen_B} below). When the data is generated in the diagonal framework,

\begin{equation}
\left\| \left(\rho - \overline{\rho}_{k_n}^{l} \right) \left(X_{n-1}^{l} \right) \right\|_{H}^{k_n} = \sqrt{\displaystyle \int_{a}^{b} \left(\displaystyle \sum_{j=1}^{k_n} \rho_j X_{n-1,j}^{l} \phi_j (t)  - \displaystyle \sum_{j=1}^{k_n} \overline{\rho}_{n,j}^{l} \left( X_{n-1}^{l} \right) \phi_{n,j}^{l}(t) \right)^2 dt}, \label{a755}
\end{equation}
\noindent is computed, being $\overline{\rho}_{k_n}^{l} \left( X_{n-1}^{l} \right)$ the  predictors defined in (\ref{ev1})-(\ref{141}), (\ref{eq23}) and (\ref{eq25}), respectively, for any $j=1,\ldots,k_n$, and based on the $l$th generation of the values  $\widetilde{X}_{i,j}^{l} = \langle X_{i}^{l}, \phi_{n,j}^{l} \rangle_H$, for $l=1, \ldots, N$, with $N=500$ simulations. The following parameter values will be considered, when the diagonal data generation is assumed:

\setlength{\heavyrulewidth}{0.3em}
\setlength{\lightrulewidth}{0.15em}

\begin{table}[h!]
 \caption{{\small Diagonal scenarios considered (see Figure \ref{fig:F1}  below, and Table 4 in the Supplementary Material), with $\delta_2 = 11/10$,  $n_t = 35000 + 40000(t-1),~t=1,\ldots,10$, and $\xi_{n_t,\beta} = \frac{\left(\ln(n_t) \right)^{\beta}}{n_{t}^{1/2}}$, for $\beta = 65/100$.}}
\centering
\vspace{-0.17cm}
\begin{tabular}{ccc}
\toprule
  Scenario & $\delta_1$ &  $k_n$ \\
 \midrule
1 & $3/2$  & $ \lceil \ln(n) \rceil$  \\
2 & $24/10$  & $\lceil \ln(n) \rceil$ \\
3 &  $3/2$  & $ \lceil e^{\prime} n^{1/\left(8 \delta_1 + 2 \right)} \rceil,~e^{\prime} = 17/10$  \\
4& $24/10$  & $\lceil e^{\prime} n^{1/\left(8 \delta_1 + 2 \right)} \rceil,~e^{\prime} = 17/10$  \\
\bottomrule
\end{tabular}
  \label{tab:Table_Scen_A}
\end{table}

As discussed, conditions formulated in Bosq (2000) and Proposition \ref{proposition5} of the current paper are held for scenarios 1-2 in Table \ref{tab:Table_Scen_A}, while in scenarios 3-4, the conditions assumed in Proposition \ref{proposition5} and the approaches by Bosq (2000) and Guillas (2001), are verified (but not in an optimal sense).

\vspace{-0.4cm}
\begin{figure}[H]
  \hspace{-1.06cm} \includegraphics[width=9.3cm,height=5.9cm]{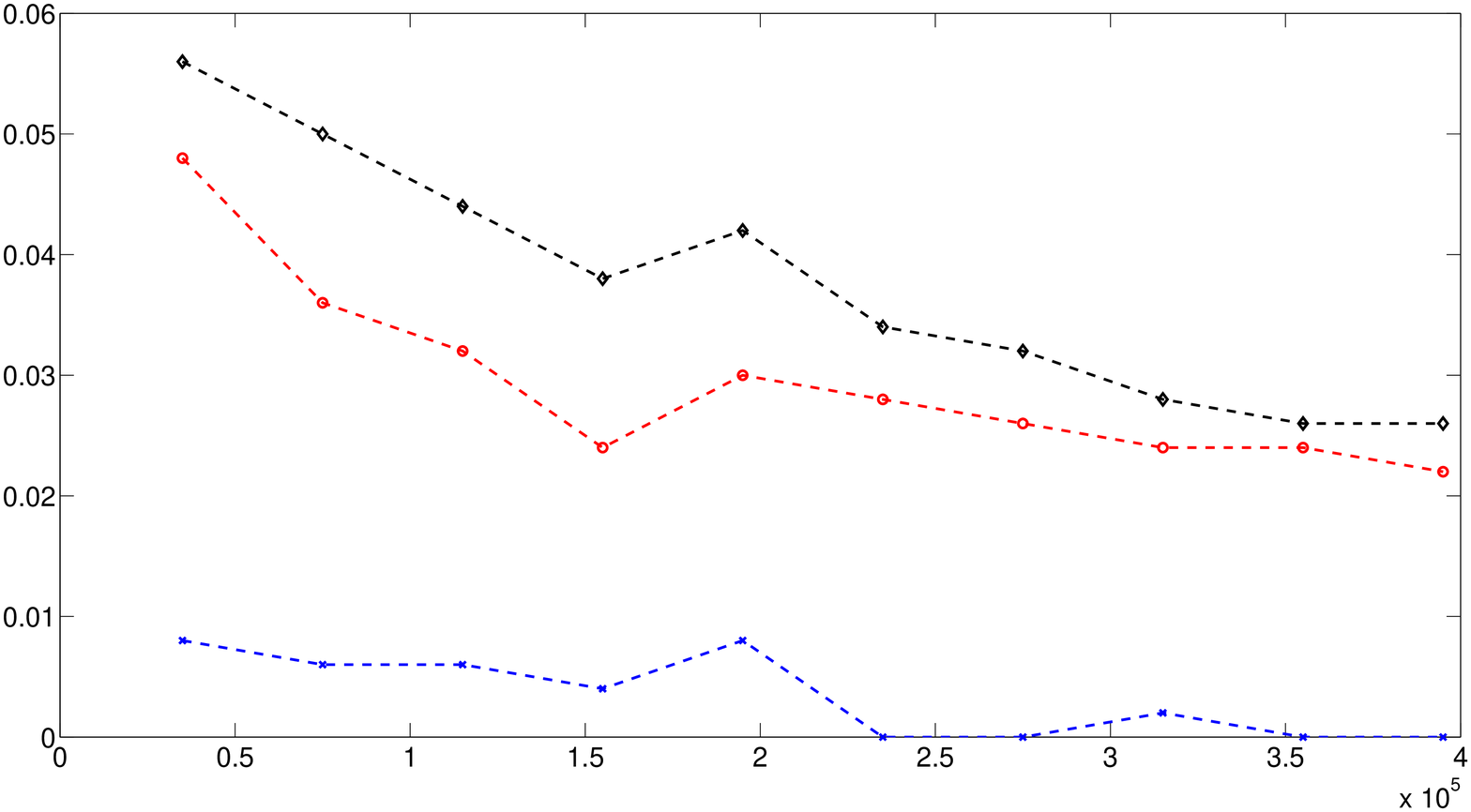}
 \hspace{-1cm}  \includegraphics[width=9.3cm,height=5.9cm]{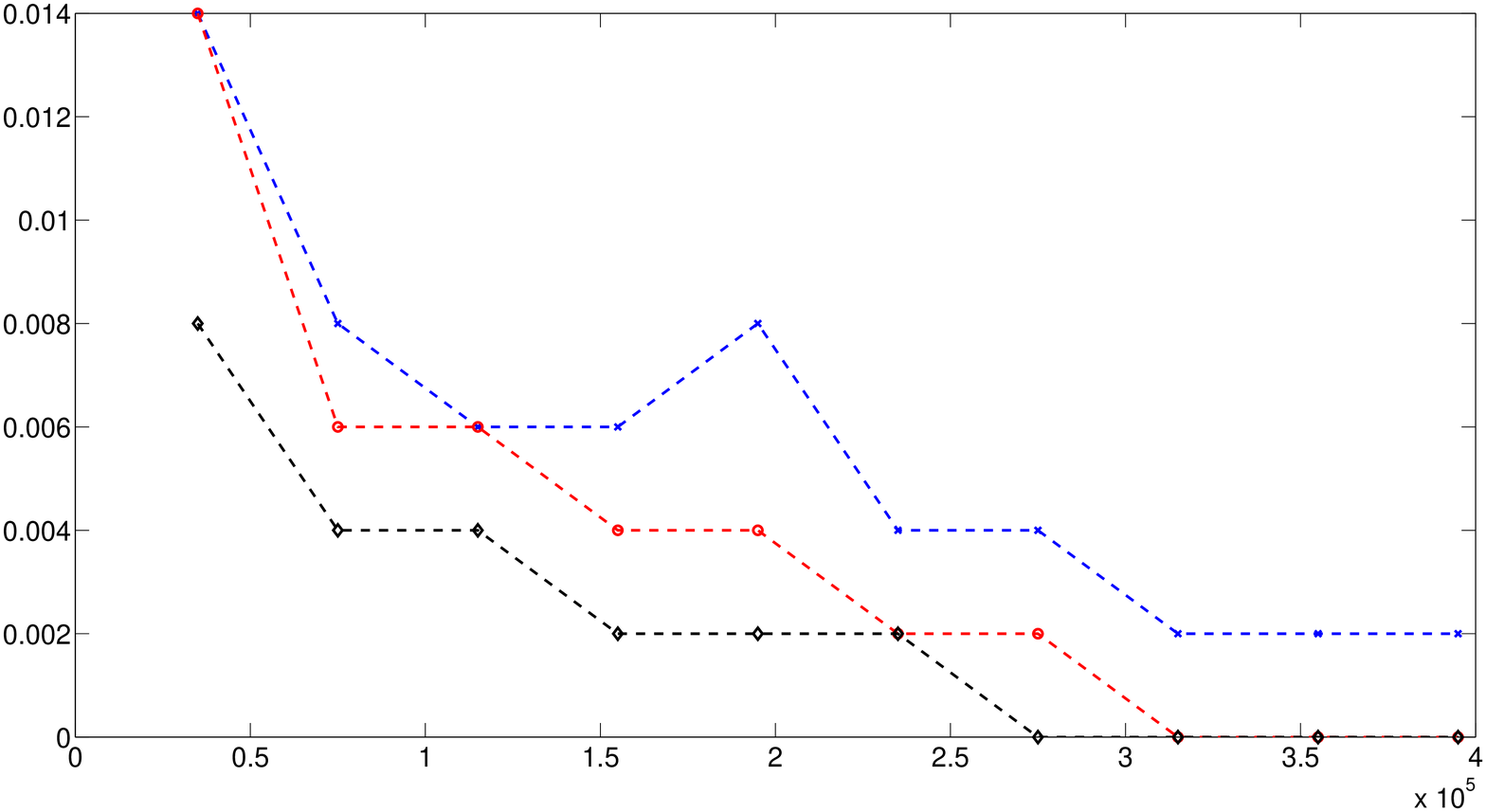}
\vspace{-1.22cm}
\caption{{\small $F \left(k_n, n_t, \beta \right)$ values, for scenario 2 (on left) and scenario 4 (on right), for our approach (blue star dotted line) and those one presented in Bosq (2000) (red circle dotted line) and Guillas (2001) (black diamond dotted line). The curve $\xi_{n_t,\beta} = \frac{ \left(\ln(n_t) \right)^{\beta} }{n_{t}^{1/2}}$, with $\beta = 65/100$, is adopted.}} \label{fig:F1}
\end{figure}

Under pseudo-diagonal and non-diagonal frameworks, the following truncated norm is then computed, instead of (\ref{a755}):
\begin{equation}
\left\| \left(\rho - \overline{\rho}_{k_n}^{l} \right) \left(X_{n-1}^{l} \right) \right\|_{H}^{k_n} = \sqrt{\displaystyle \int_{a}^{b} \left( \displaystyle \int_{a}^{b} \left(\displaystyle \sum_{j,k=1}^{k_n} \rho_{j,k} \phi_j (t) \phi_k (s) \right)ds  - \displaystyle \sum_{j=1}^{k_n} \overline{\rho}_{n,j}^{l} \left( X_{n-1}^{l} \right) \phi_{n,j}^{l}(t) \right)^2 dt}. \label{a766}
\end{equation}

Pseudo-diagonal  and non-diagonal scenarios, for approaches formulated in Bosq (2000) and Guillas (2001), are outlined in Table \ref{tab:Table_Scen_B}.

\begin{table}[h!]
 \caption{{\small Pseudo-diagonal and non-diagonal scenarios considered (see Figures \ref{fig:F2}-\ref{fig:F2a}  below, and Tables 5-6 in the Supplementary Material), with $\delta_2 = 11/10$,  $n_t = 35000 + 40000(t-1),~t=1,\ldots,10$, and $\xi_{n_t,\beta} = \frac{\left(\ln(n_t) \right)^{\beta}}{n_{t}^{1/3}}$.}}
\centering
\vspace{-0.17cm}
\begin{tabular}{cccc||cccc}
\toprule
\multicolumn{4}{c||}{Pseudo-diagonal scenarios} & \multicolumn{4}{c}{Non-diagonal scenarios} \\ 
\hline
  Scenario & $\delta_1$ &  $k_n$ & $\beta$ & Scenario & $\delta_1$ &  $k_n$ & $\beta$ \\
 \midrule
5 & $3/2$  & $ \lceil \ln(n) \rceil$ & $3/10$ & 9 & $3/2$  & $ \lceil \ln(n) \rceil$ & $125/100$\\
6 & $24/10$  & $\lceil \ln(n) \rceil$ & $3/10$ & 10 &  $24/10$  & $ \lceil \ln(n) \rceil$ & $125/100$\\
7&  $3/2$  & $ \lceil (17/10) n^{1/\left(8 \delta_1 + 2 \right)} \rceil$  & $3/10$ & 11 & $3/2$  & $ \lceil(17/10)  n^{1/\left(8 \delta_1 + 2 \right)} \rceil$ & $125/100$\\
8& $24/10$  & $\lceil (17/10) n^{1/\left(8 \delta_1 + 2 \right)} \rceil$  & $3/10$ & 12 &  $24/10$  & $ \lceil (17/10) n^{1/\left(8 \delta_1 + 2 \right)} \rceil$ & $125/100$ \\
\bottomrule
\end{tabular}
  \label{tab:Table_Scen_B}
\end{table}

Scenarios 5-6 and 9-10 verify conditions required in Bosq (2000), while scenarios 7-8 and 11-12 are included in both setting of conditions, proposed in Bosq (2000) and Guillas (2001).

\begin{figure}[H]
 \hspace{-1.06cm} \includegraphics[width=9.3cm,height=5.9cm]{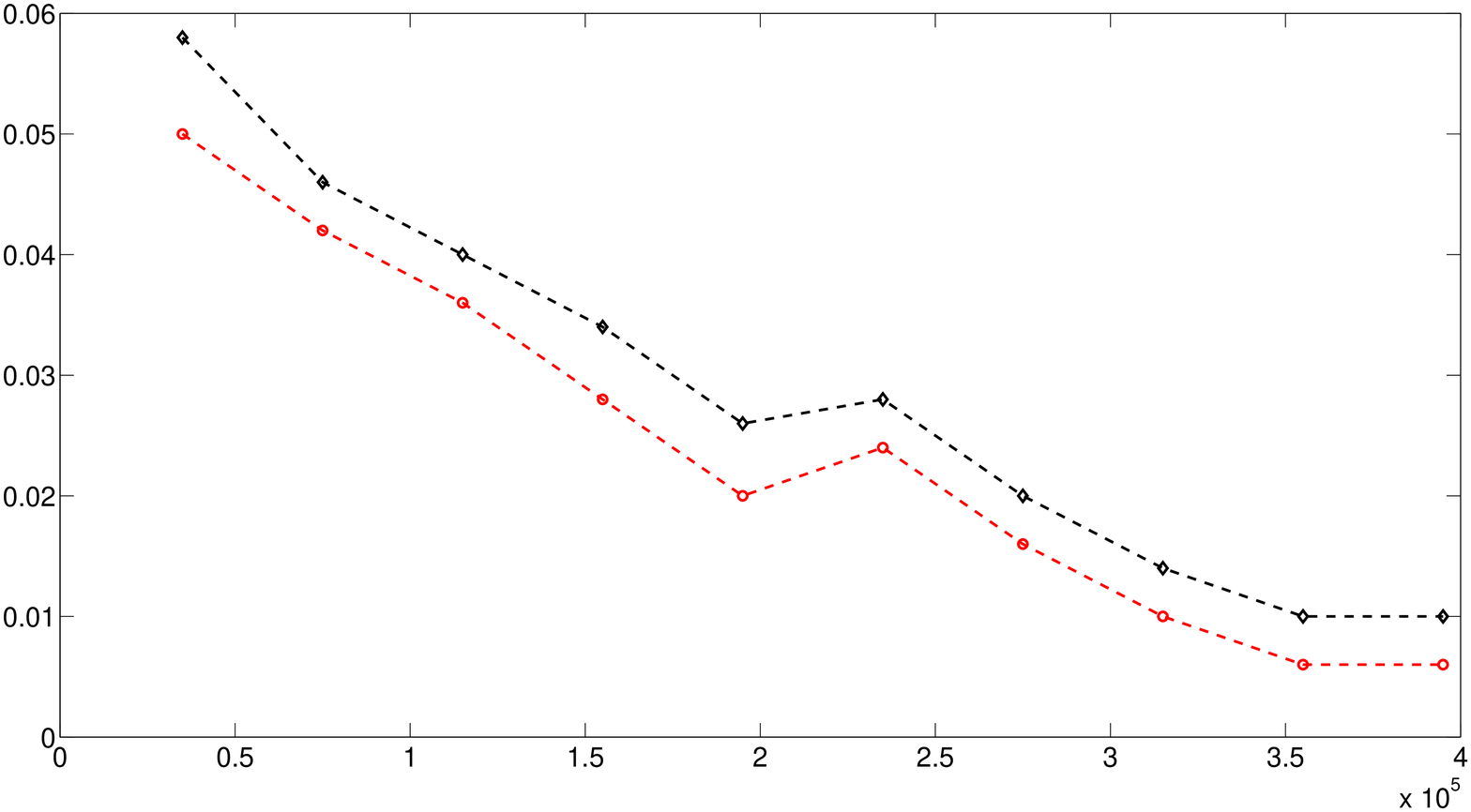}
 \hspace{-1cm}  \includegraphics[width=9.3cm,height=5.9cm]{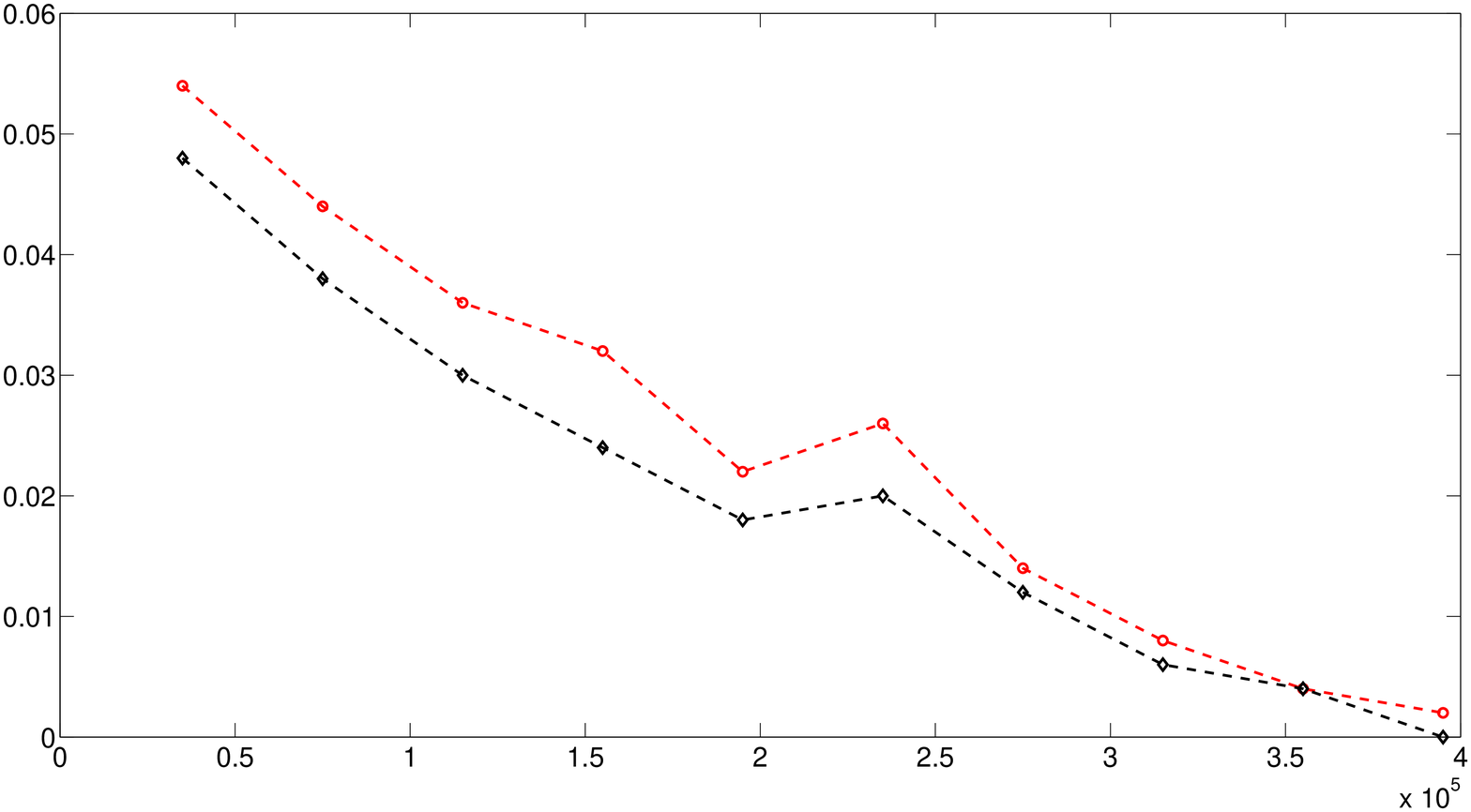}
\vspace{-1.22cm}
\caption{{\small $F \left(k_n, n_t, \beta \right)$ values, for scenario 6 (on left) and scenario 8 (on right), for approaches presented in Bosq (2000) (red circle dotted line) and Guillas (2001) (black diamond dotted line). The curve $\xi_{n_t,\beta} = \frac{\left(\ln(n_t) \right)^{\beta} }{n_{t}^{1/3} }$, with $\beta = 3/10$, is adopted.}} \label{fig:F2}
\end{figure}

\begin{figure}[H]
 \hspace{-1.06cm} \includegraphics[width=9.3cm,height=5.9cm]{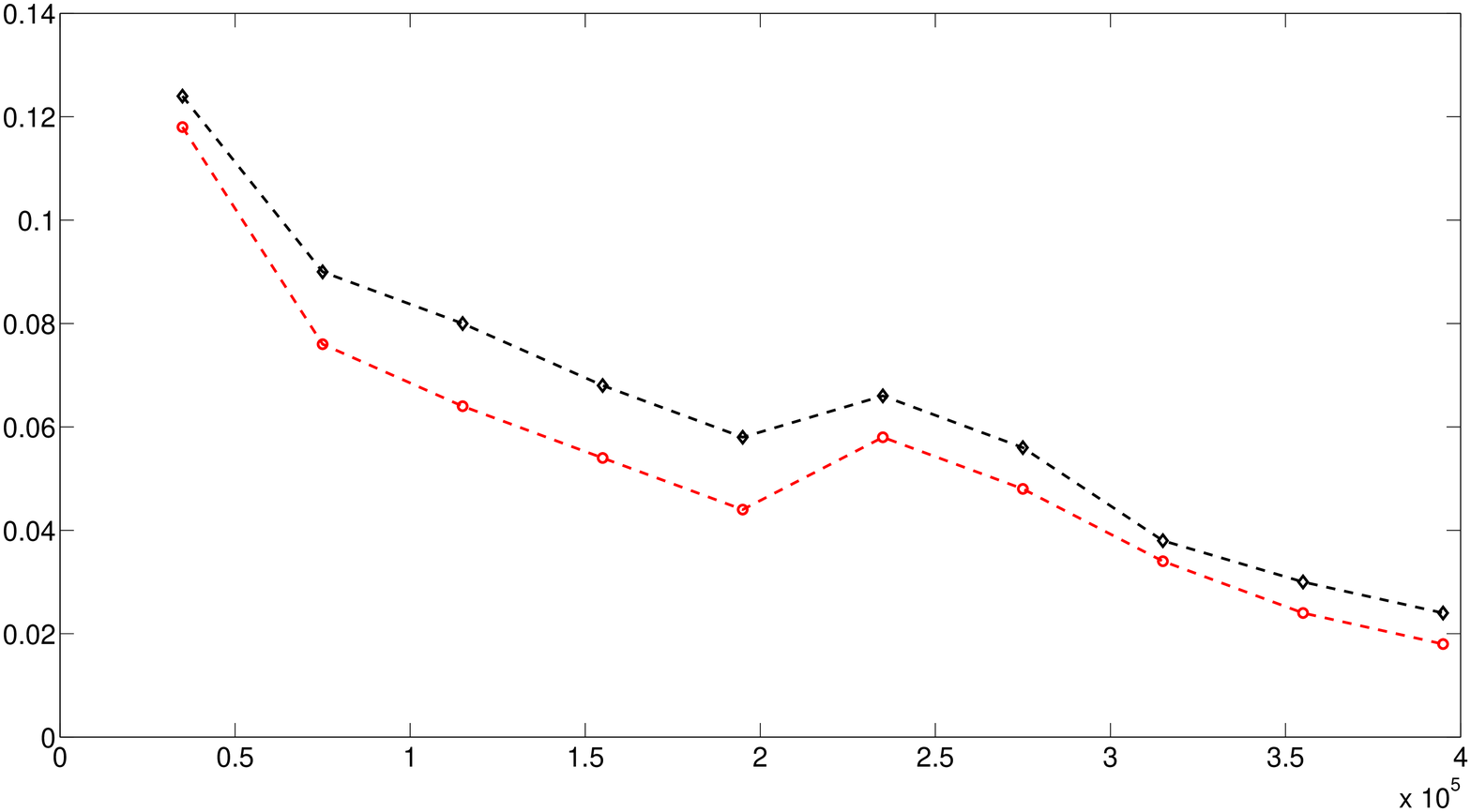}
 \hspace{-1cm}  \includegraphics[width=9.3cm,height=5.9cm]{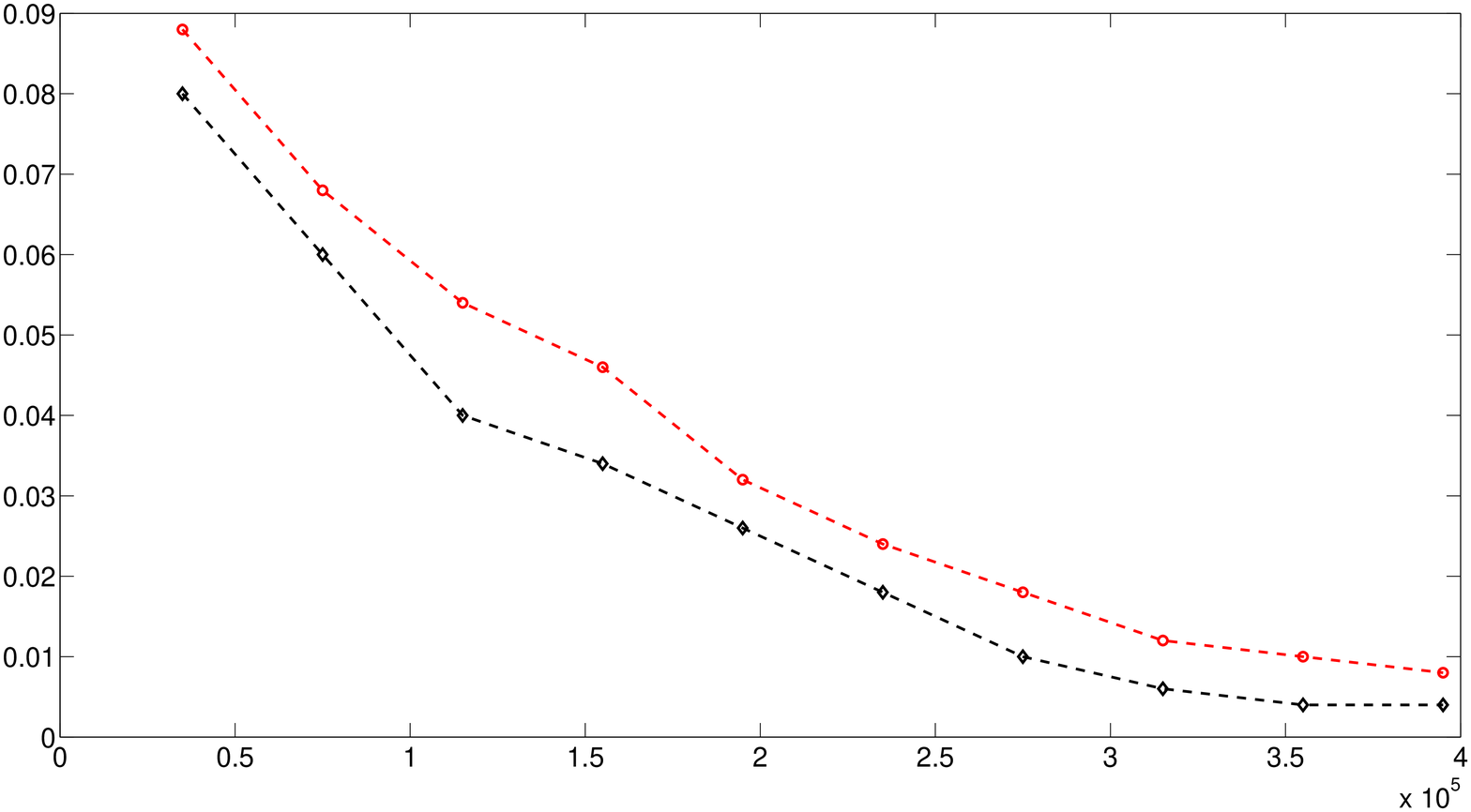}
\vspace{-1.22cm}
\caption{{\small $F \left(k_n, n_t, \beta \right)$ values, for scenario 10 (on left) and scenario 12 (on right), for approaches presented in Bosq (2000) (red circle dotted line) and Guillas (2001) (black diamond dotted line). The curve $\xi_{n_t,\beta}= \frac{\left(\ln(n_t) \right)^{\beta} }{n_{t}^{1/3} }$, with $\beta = 125/100$, is adopted.}} \label{fig:F2a}
\end{figure}

Results obtained in the diagonal scenarios 1-4, which are reflected in Table \ref{tab:Table_Scen_A}, have been applied to the three componentwise ARH(1) plug-in predictor approaches. As expected, the amount of values $\left\| \left(\rho - \overline{\rho}_{k_n}^{l} \right) \left(X_{n-1}^{l} \right) \right\|_{H}^{k_n}$, which lie within the band $\left[0, \xi_{n_t,\beta} \right)$, is greater as long as the decay rate of the eigenvalues of $C$ is faster. Since a diagonal framework is considered in scenarios 1-2, a better performance of the approach formulated in Section 8.3 can be noticed in comparison with those ones by Bosq (2000) and Guillas (2001), where errors appear,  when sample sizes are not sufficiently large, in the estimation of the non-diagonal componentwise coefficients of $\rho$, which must be zero. This possible effect of the non-diagonal design, under a diagonal scenario, is not observed, for  truncation rules selecting a very small number of terms, in relation to the sample size. This fact occurs  in the truncation rule adopted in Guillas (2001).   

In the pseudo-diagonal and non-diagonal scenarios outlined in Table \ref{tab:Table_Scen_B}, methodologies in Bosq (2000) and Guillas (2001) are compared, such that the curve $\xi_{n_t,\beta} = \frac{\left(\ln(n_t) \right)^{\beta} }{n_{t}^{1/3} }$, for $\beta = 3/10$ and $\beta = 125/100$, numerically fits their almost sure rate of convergence. As observed (see also Tables 4-6 in Section S.5 of the Supplementary material), the sample-size dependent truncation rule, according to the rate of convergence to zero of the eigenvalues of $C$, plays a crucial role in the observed performance of both approaches.

%\vspace{-0.6cm}
%\begin{figure}[H]
% \hspace{-1.06cm} \includegraphics[width=9.3cm,height=5.9cm]{Comp_CMedia_TruncNuestro_NoNoDiag.eps}
% \hspace{-1.02cm}  \includegraphics[width=9.3cm,height=5.9cm]{Comp_CMedia_TruncSuyo_NoNoDiag.eps}
%\vspace{-1.22cm}
%\caption{{\small Empirical mean of $\left\| \left(\rho - \overline{\rho}_{k_n}^{l} \right) \left(X_{n-1}^{l} \right) \right\|_{H}^{k_n}$, for scenario 9 (on left) and scenario 11 (on right), for our approach (blue star dotted line) and those one presented in Bosq (2000) (red circle dotted line) and Guillas (2001) (black diamond dotted line). The curve $\xi_{n_t}^{\beta} = \frac{\left(\ln(n_t) \right)^{\beta} }{n_{t}^{1/3} },~\beta = 125/100$, is drawn (green dotted line).}} \label{fig:F3}
%\end{figure}

\subsection{Small-sample behaviour of the ARH(1) plug-in and non-plug-in predictors}

Smaller sample sizes must be adopted in this subsection, since computational limitations arise when regularized wavelet plug-in predictor formulated in Antoniadis  \& Sapatinas (2003) (see equations (\ref{wavelets1a})-(\ref{wavelets1c}) above), as well as penalized predictor and non-parametric kernel-based predictor applied in Besse et al. (2000) (see equations (\ref{eq_8}) and (\ref{a100}), respectively), are included in the comparative study. See also Section S.5 of the Supplementary Material, where extra numerical results are provided. As above, the diagonal componentwise estimator here formulated will be only considered under diagonal scenarios.

On the one hand, \textbf{Assumptions A1} and \textbf{A3}, and conditions in (\ref{wavelets1c}), are required when regularized-wavelet-based prediction approach is applied. In particular, since $C_j = c_1 j^{-\delta_1}$, for any $j \geq 1$, if $k_n = \lceil n^{1/\alpha} \rceil$ is adopted, then $1- \frac{4 \delta_1}{\alpha} > 0$, leading to $\alpha > 4 \delta_1$. Additionally to $k_n = \lceil \ln(n) \rceil$, the truncation parameter $k_n = \lceil n^{1/\alpha} \rceil$ will be adopted (see Table \ref{tab:Table_Scen_C}-\ref{tab:Table_Scen_D}), with $\alpha=6.5$ and $\alpha = 10$, for $\delta_1 = 3/2$ and $\delta_1=24/10$, respectively. Furthermore, $F \left(k_n,n_t, \beta \right)$ values defined in (\ref{a955})-(\ref{a766}) are computed for the wavelet-based approach just replacing $\left\lbrace \phi_{n,j}, \ j \geq 1 \right\rbrace$ by $\left\lbrace \widetilde{\phi}_{j}^{M},\ j \geq 1 \right\rbrace$ (see equations (\ref{wavelets1a})-(\ref{wavelets1c})). As before, since $\lceil n^{1/\alpha} \rceil < \lceil \ln(n) \rceil$ and $\lceil n^{1/\alpha} \rceil <  \lceil (17/10) n^{1/\left(8 \delta_1 + 2 \right)} \rceil$, for $\alpha=6.5$ and $\alpha = 10$, conditions imposed for the estimator formulated in Section 8, as well as in Bosq (2000) and Guillas (2001), are verified when the truncation parameter $k_n = \lceil n^{1/\alpha} \rceil$, with $\alpha=6.5$ and $\alpha = 10$, is studied.

On the other hand, let us also compare with the  techniques presented in Besse et al. (2000), and detailed in equations (\ref{eq_8}) and (\ref{a100}), based on penalized prediction and non-parametric kernel-based prediction, respectively. In those techniques, they assume that the functional values of the stationary process are in the Sobolev space $W^{2,2} \left( \left[0, 1 \right] \right)$. When the referred methodologies in Besse et al. (2000) are implemented, the following alternative norm replaces the norm reflected in (\ref{a755})-(\ref{a766}), for values $F \left( k_n, n_t, \beta \right)$:
\begin{eqnarray}
\left\| \left(\rho - \overline{\rho}_{k_n}^{l} \right) \left(X_{n-1}^{l} \right) \right\|_{H} &=& \sqrt{\displaystyle \int_{a}^{b} \left( \rho \left(X_{n-1}^{l} \right) (t) - \overline{\rho}_{k_n}^{l} \left(X_{n-1}^{l} \right) (t) \right)^2 dt }, \quad l=1,\ldots,N.\label{new2}
\end{eqnarray}

\vspace{-0.3cm}

Hence, the following diagonal scenarios are regarded:

\begin{table}[h!]
 \caption{{\small Diagonal scenarios considered (see Figure \ref{fig:F5} below, and Tables 9-10 in the Supplementary Material), with $M=50$, $q=10$, $\delta_2 = 11/10$,  $n_t = 750 + 500(t-1),~t=1,\ldots,13$, and $\xi_{n_t,\beta} = \frac{\left(\ln(n_t) \right)^{\beta}}{n_{t}^{1/2}},~\beta = 65/100$.}}
\centering
\vspace{-0.17cm}
\begin{small}
\begin{tabular}{cccc}
\toprule
  Scenario & $\delta_1$ &  $k_n$ & $h_n$ \\
 \midrule
13 & $3/2$  & $ \lceil \ln(n) \rceil$ & $0.15,~0.25$ \\
14 & $24/10$  & $\lceil \ln(n) \rceil$ & $0.15,~0.25$ \\
15 &  $3/2$  & $ \lceil n^{1/\alpha} \rceil,~\alpha = 6.5$  \\
16 & $24/10$  & $ \lceil n^{1/\alpha} \rceil,~\alpha = 10$ \\
\bottomrule
\end{tabular}
\end{small}
  \label{tab:Table_Scen_C}
\end{table}

Remark that, since both approaches formulated in Besse et al. (2000) not depend on the truncation parameter $k_n$ adopted, we only perform them for scenarios 13-14, where different rates of convergence to zero of the eigenvalues of $C$  are considered, and conditions imposed in that paper are verified. Conditions formulated in Bosq (2000) and Proposition \ref{proposition5} of the current paper are held for all scenarios, while the conditions assumed in Antoniadis \& Sapatinas (2003) and Guillas (2001) are only verified under scenarios 15-16. 

\vspace{-0.4cm}
\begin{figure}[H]
 \hspace{-1.06cm} \includegraphics[width=9.2cm,height=5.9cm]{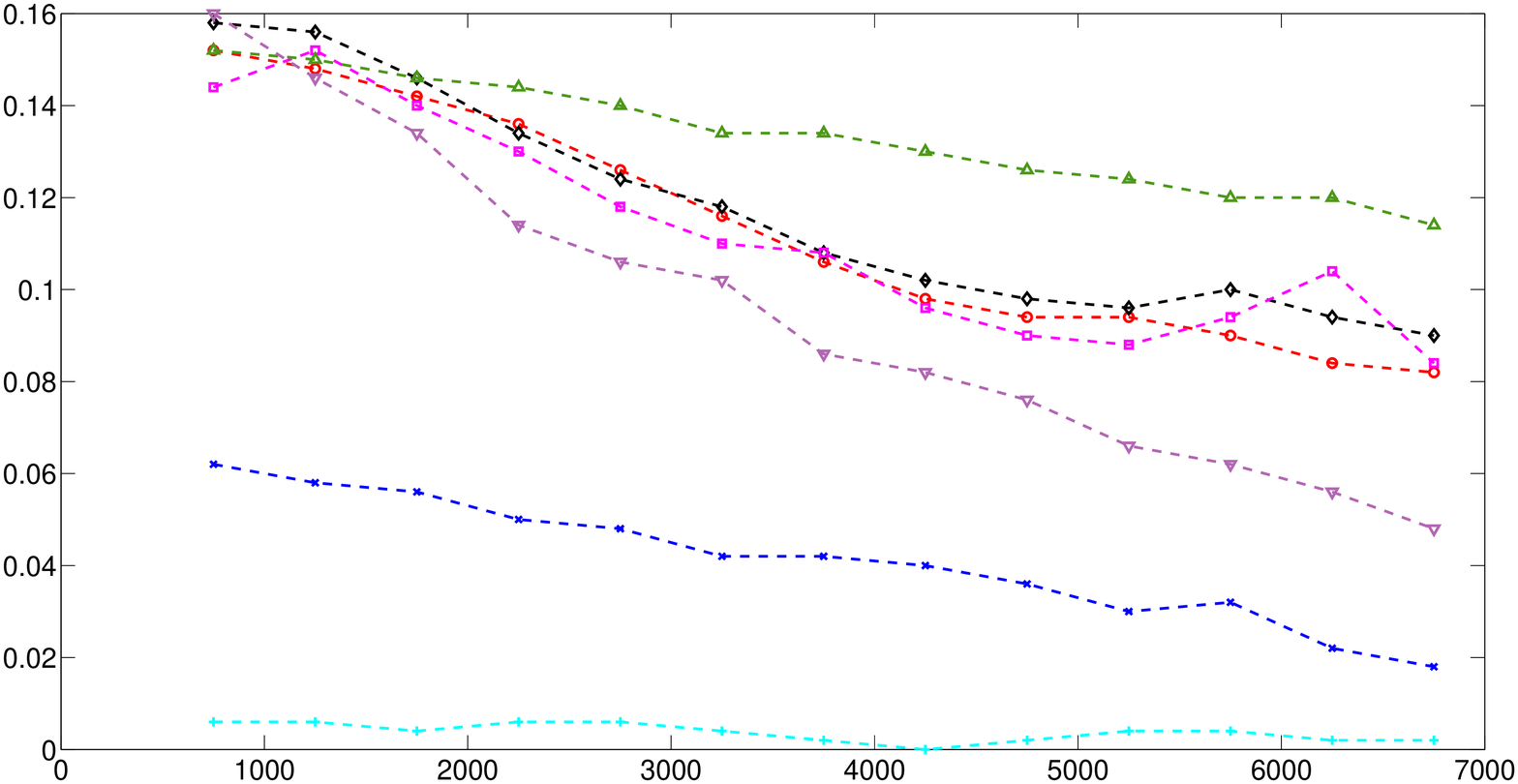}
 \hspace{-1.02cm}  \includegraphics[width=9.2cm,height=5.9cm]{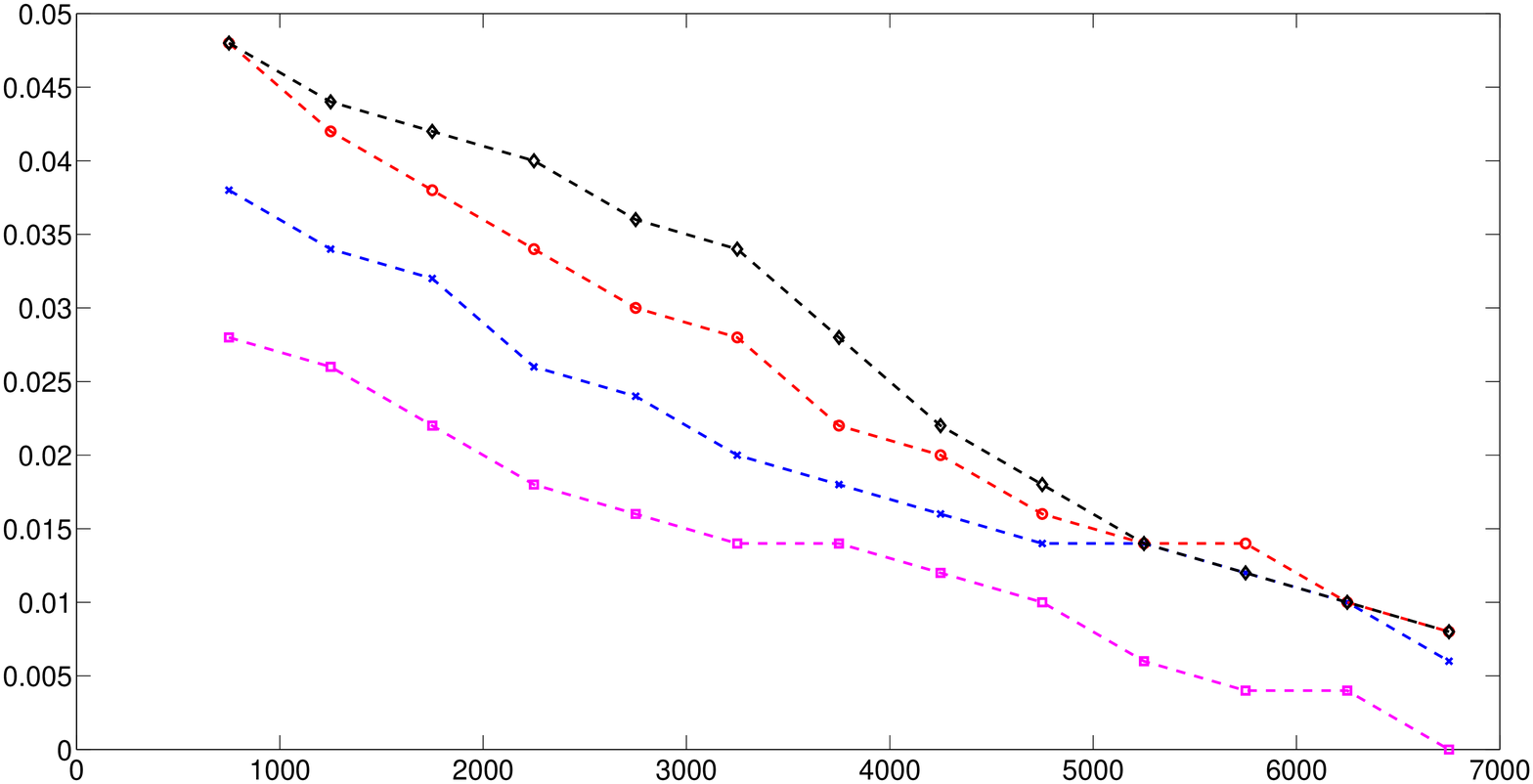}
\vspace{-1.22cm}
\caption{{\small $F \left( k_n, n_t, \beta \right)$ values, for scenario 14 (on left) and scenario 16 (on right), for our approach (blue star dotted line) and those one presented in Antoniadis and Sapatinas (2001) (pink square dotted line), Besse et al. (2000) (cyan blue plus dotted line for penalized prediction; dark green upward-pointing triangle  and purple downward-pointing triangle dotted lines, for kernel-based prediction, for $h_n = 0.15$ and $h_n = 0.25$, respectively), Bosq (2000) (red circle dotted line) and Guillas (2001) (black diamond dotted line). The curve $\xi_{n_t,\beta} = \frac{\left(\ln(n_t) \right)^{\beta} }{n_{t}^{1/2} }$, with $\beta = 65/100$, is drawn (light green dotted line).}} \label{fig:F5}
\end{figure}

 Pseudo-diagonal and non-diagonal scenarios are detailed in Table \ref{tab:Table_Scen_D}.
\begin{table}[h!]
 \caption{{\small Pseudo-diagonal and non-diagonal scenarios considered (see Figures \ref{fig:F6}-\ref{fig:F7} below, and Tables 11-14 in the Supplementary Material), with $\delta_2 = 11/10$,  $n_t = 750 + 500(t-1),~t=1,\ldots,13$, and $\xi_{n_t,\beta} = \frac{\left(\ln(n_t) \right)^{\beta}}{n_{t}^{1/3}}$.}}
\centering
\vspace{-0.17cm}
\begin{small}
\begin{tabular}{ccccc||ccccc}
\toprule
\multicolumn{5}{c||}{Pseudo-diagonal scenarios} & \multicolumn{5}{c}{Non-diagonal scenarios} \\ 
\hline
  Scenario & $\delta_1$ &  $k_n$ & $\beta$ & $h_n$ & Scenario & $\delta_1$ &  $k_n$ & $\beta$ & $h_n$  \\
 \midrule
17 & $3/2$  & $ \lceil \ln(n) \rceil$ & $3/10$ & $1.2,~1.7$ & 21 & $3/2$  & $ \lceil \ln(n) \rceil$ & $125/100$ & $1.2,~1.7$\\
18 & $24/10$  & $\lceil \ln(n) \rceil$ & $3/10$ & $1.2,~1.7$ & 22 &  $24/10$  & $ \lceil \ln(n) \rceil$ & $125/100$& $1.2,~1.7$\\
19&  $3/2$  & $ \lceil n^{1/\alpha} \rceil,~\alpha = 6.5$   & $3/10$ & & 23 & $3/2$  & $ \lceil n^{1/\alpha} \rceil$ & $125/100$ & \\
20& $24/10$  & $ \lceil n^{1/\alpha} \rceil,~\alpha = 10$  & $3/10$ & & 24 &  $24/10$  & $ \lceil n^{1/\alpha} \rceil$ & $125/100$&  \\
\bottomrule
\end{tabular}
\end{small}
  \label{tab:Table_Scen_D}
\end{table}

As noted before, approaches formulated in Besse et al. (2000) are only tested for scenarios 17-18 and 21-22, and conditions in Bosq (2000) are verified for all scenarios. Scenarios developed by Antoniadis \& Sapatinas (2003) and Guillas (2001) are only held when the truncation parameter proposed in Antoniadis \& Sapatinas (2003)  is adopted.

\vspace{-0.6cm}
\begin{figure}[H]
 \hspace{-1.06cm} \includegraphics[width=9.2cm,height=5.9cm]{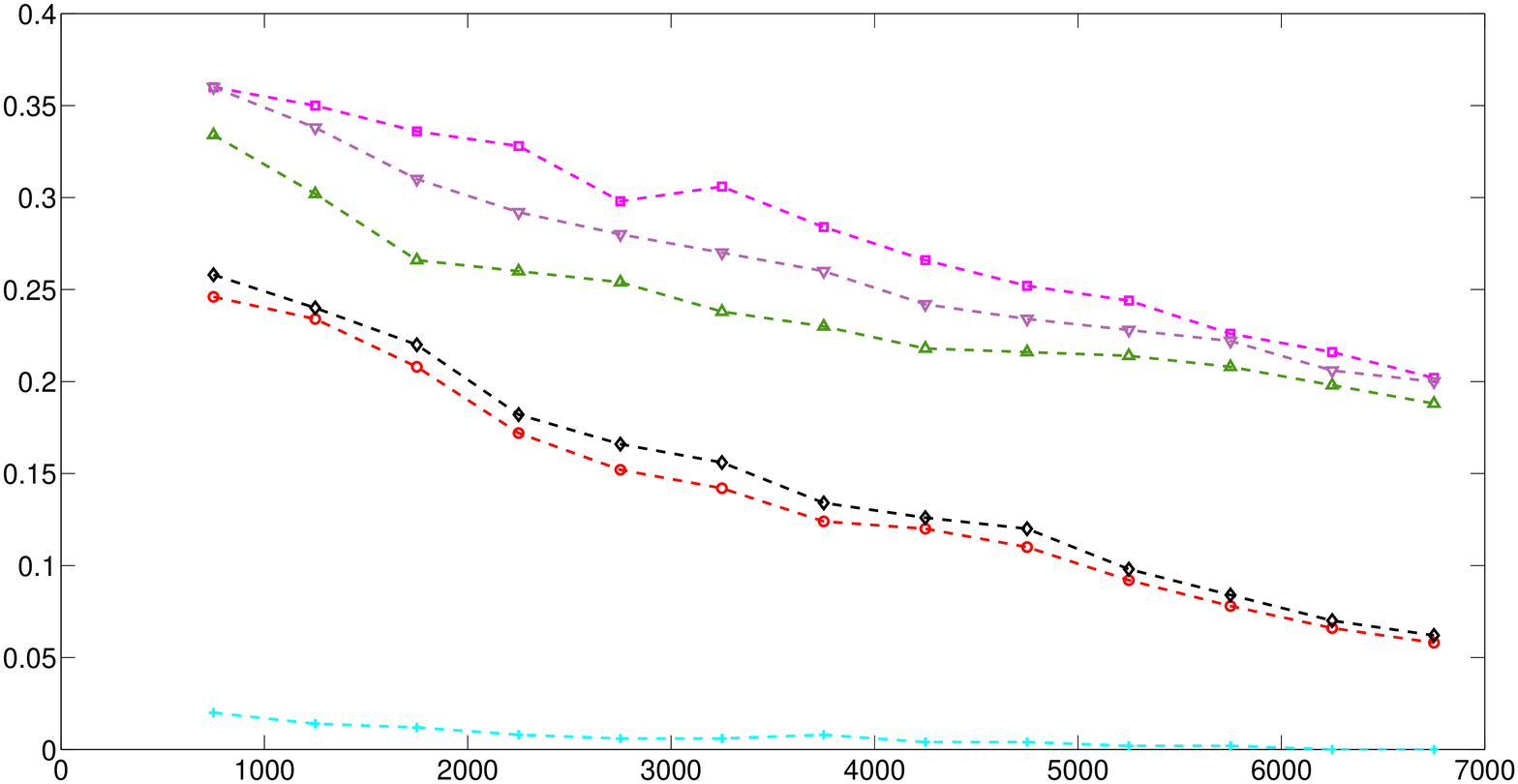}
 \hspace{-1.02cm}  \includegraphics[width=9.2cm,height=5.9cm]{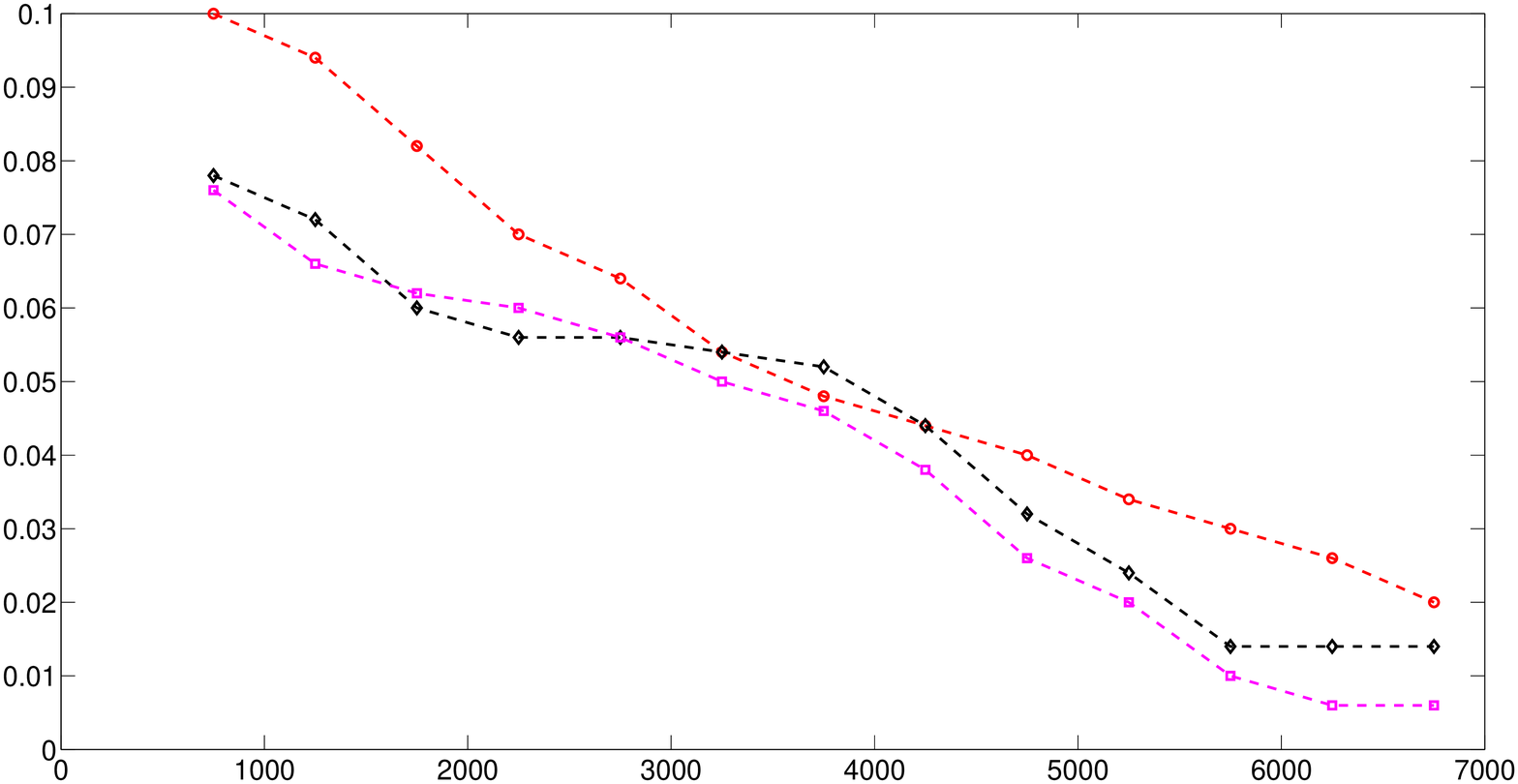}
\vspace{-1.22cm}
\caption{{\small $F \left( k_n, n_t, \beta \right)$ values, for scenario 18 (on left) and scenario 20 (on right), for approaches presented in Antoniadis \& Sapatinas (2003) (pink square dotted line), Besse et al. (2000) (cyan blue plus dotted line for penalized prediction; dark green upward-pointing triangle  and purple downward-pointing triangle dotted lines, for kernel-based prediction, for $h_n = 1.2$ and $h_n = 1.7$, respectively), Bosq (2000) (red circle dotted line) and Guillas (2001) (black diamond dotted line). The curve $\xi_{n_t,\beta} = \frac{\left(\ln(n_t) \right)^{\beta} }{n_{t}^{1/3} }$, with $\beta = 3/10$, is drawn (light green dotted line).}} \label{fig:F6}
\end{figure} 

\vspace{-0.47cm}
\begin{figure}[H]
 \hspace{-1.06cm} \includegraphics[width=9.2cm,height=5.9cm]{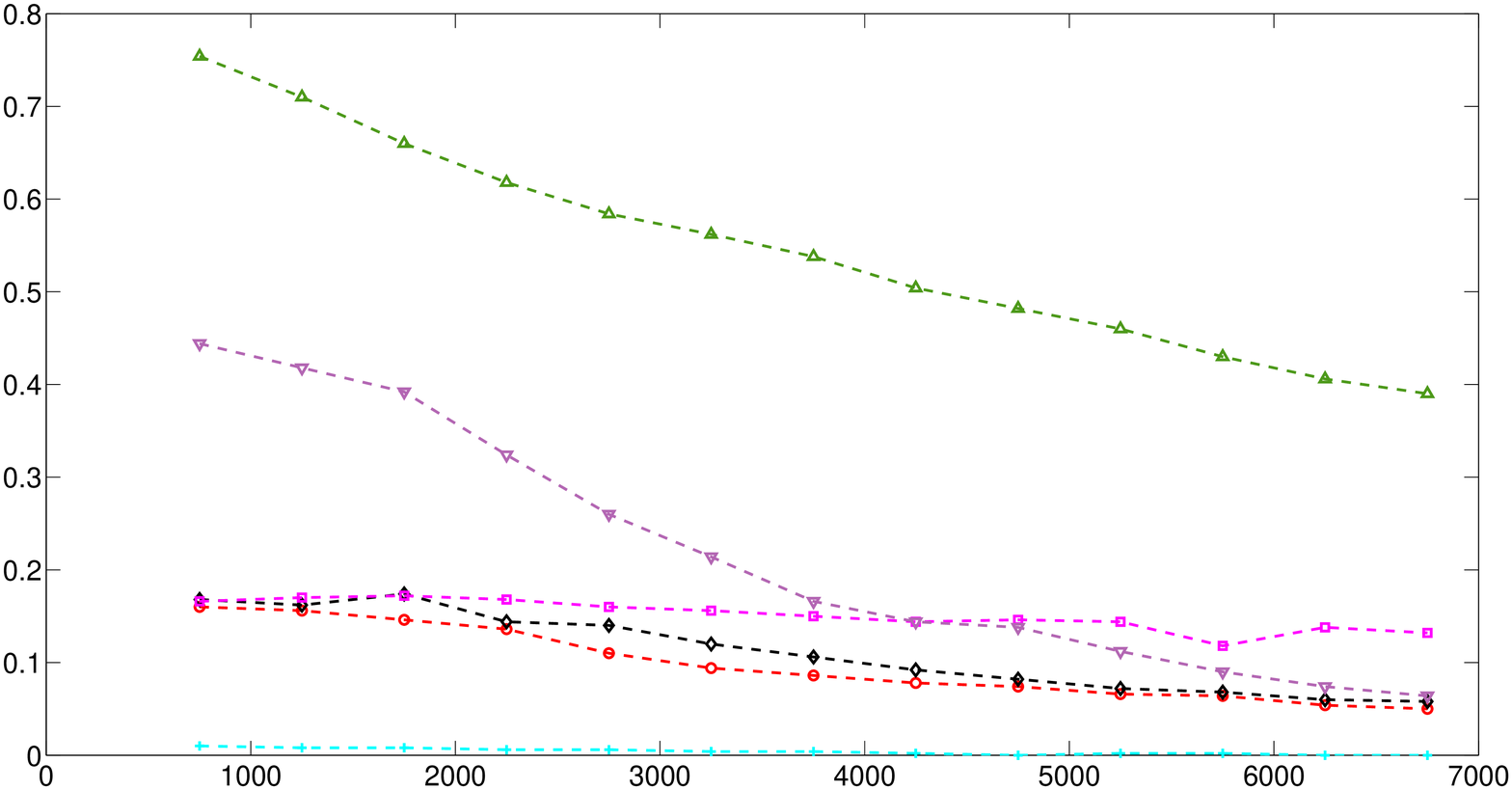}
 \hspace{-1.02cm}  \includegraphics[width=9.2cm,height=5.9cm]{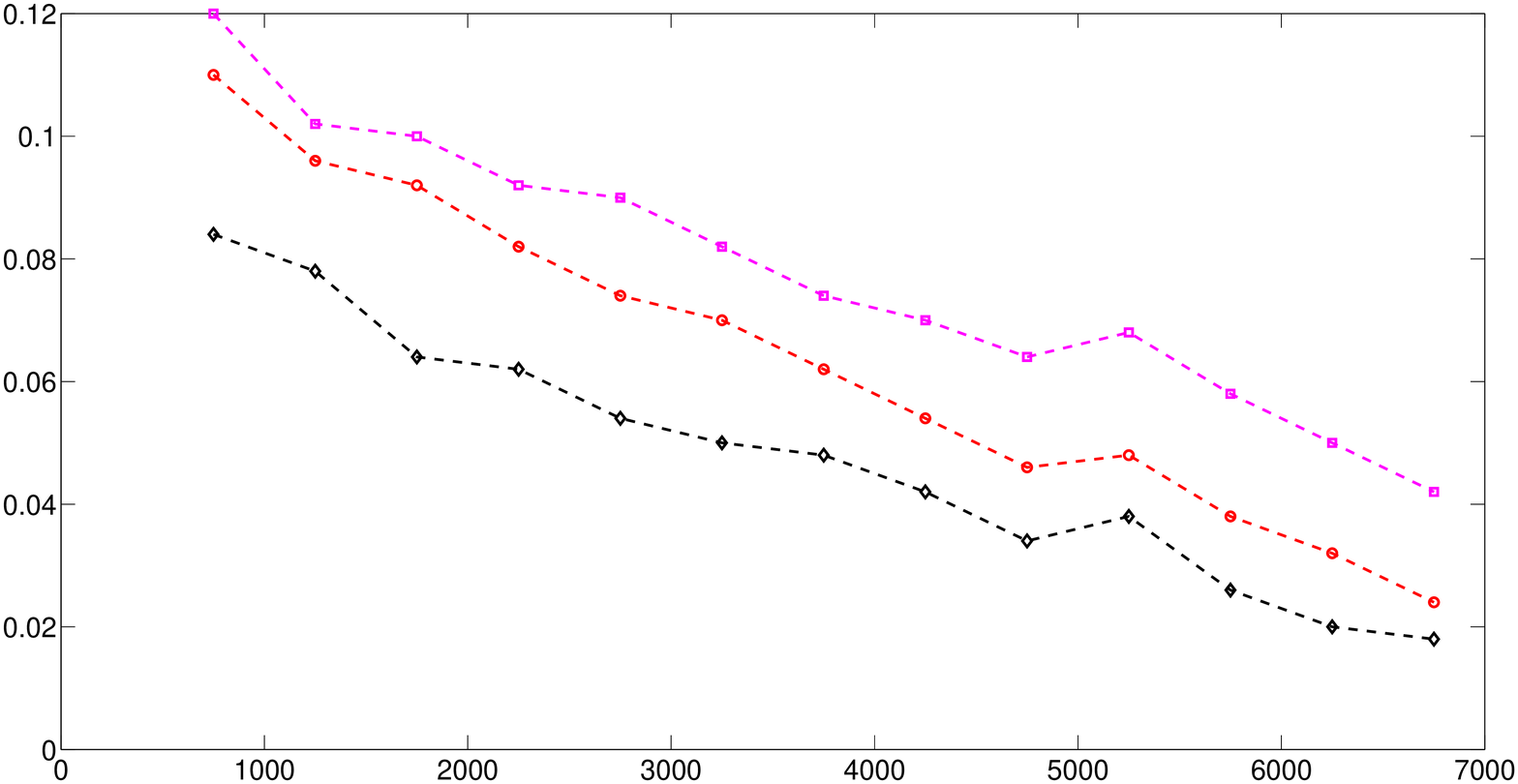}
\vspace{-1.22cm}
\caption{{\small $F \left( k_n, n_t, \beta \right)$ values, for scenario 22 (on left) and scenario 24 (on right), for approaches presented in Antoniadis \& Sapatinas (2003) (pink square dotted line), Besse et al. (2000) (cyan blue plus dotted line for penalized prediction; dark green upward-pointing triangle  and purple downward-pointing triangle dotted lines, for kernel-based prediction, for $h_n = 1.2$ and $h_n = 1.7$, respectively), Bosq (2000) (red circle dotted line) and Guillas (2001) (black diamond dotted line). The curve $\xi_{n_t,\beta} = \frac{\left(\ln(n_t) \right)^{\beta} }{n_{t}^{1/3} }$, with $\beta = 125/100$, is drawn (light green dotted line).}} \label{fig:F7}
\end{figure}

When smaller sample sizes are adopted, and approaches formulated in Antoniadis \& Sapatinas (2003) and Besse et al. (2000) are included in the comparative study, scenarios 13-24 have been considered and reflected in Tables \ref{tab:Table_Scen_C}-\ref{tab:Table_Scen_D}. As expected, the larger sample are used, the better performance is obtained for the empirical-eigenvector-based componentwise approaches tested in the previous subsection. Note that even when small sample sizes are studied, a good performance of the ARH(1) plug-in predictor given in equations (\ref{ev1})-(\ref{141}) is observed. As well as the regularized wavelet-based approach detailed in Antoniadis \& Sapatinas (2003) becomes the best methodology for small sample sizes, in comparision with the componentwise techniques above mentioned. 

Note that the good performance observed corresponds to the truncation rule proposed by these authors, with a small number of terms. While, when a larger number of terms is considered, according to the alternative truncation rules tested, the observed outperformance  does not hold.  While the penalized prediction approach proposed in Besse et al. (2000) has been shown as the more accurate, is, however, less affected by the regularity conditions imposed on the autocovariance kernel (see Tables 9-14 in Section  S.5 of the supplementary material). The non-parametric kernel-based purpose by Besse et al. (2000) requires to solve the selection problem  associated with  the bandwidth parameter.  Furthermore, a drawback of both approaches in Antoniadis \& Sapatinas (2003) and Besse et al. (2000) is that they require large computational times in their implementations. The underlying dependence structure, given by the covariance operators and their spectral decompositions, cannot be provided in those approaches.

% *******************************

% *******************************

\vspace*{0.5cm}

 \noindent {\bfseries {\large Acknowledgments}}

This work has been supported in part by project MTM2015--71839--P (co-funded by Feder funds), of the DGI, MINECO, Spain.

% **********************************
 
% **********************************


\vspace{-0.27cm}
\begin{thebibliography}{}
\vspace{-0.05cm}
\begin{small}

\vspace{-0.27cm}
\bibitem{Allam14}
Allan, A. \& Mourid, T. (2014). Covariance operator estimation of a functional autoregressive process with random coefficients. \emph{Statist. Probab. Lett.}, \textbf{84}, 1--8.

\vspace{-0.27cm}
\bibitem{Alvarez16}
\'Alvarez-Li\'ebana, J., Bosq, D. \& Ruiz-Medina, M. D. (2016). Consistency of the plug-in functional predictor of the Ornstein-Uhlenbeck process in
Hilbert and Banach spaces. \emph{Statist. Probab. Lett.}, \textbf{117}, 12--22.

\vspace{-0.27cm}
\bibitem{Alvarez17}
 \'Alvarez-Li\'ebana, J., Bosq, D. \& Ruiz-Medina, M. D. (2017). Asymptotic properties of a componentwise ARH(1) plug-in  predictor. \emph{J. Multivariate Anal.}, \textbf{155}, 12--34.



%\vspace{-0.27cm}
%\bibitem{Alvarez17}
%\'Alvarez-Li\'ebana, J. \& Ruiz-Medina, M. D. (2017). The effect of the spatial domain in FANOVA models with ARH(1) term. \emph{Stat. Interface}, in press.


\vspace{-0.27cm}
\bibitem{Andersson10}
Andersson, J. \& Lillest\o l, J. (2010). Modeling and forecasting electricity consumption by functional data analysis. \emph{The Journal of Energy Markets}, \textbf{3}, 3--15.

%\vspace{-0.27cm}
%\bibitem{Antoch08}
%Antoch, J., Prchal, L., De Rosa, M.R. \& Sarda, P. (2008). \emph{Functional linear regression with functional response: application to prediction of electricity consumption}. In: Functional and Operatorial Statistics. Contributions to Statistics. Physica-Verlag HD.


\vspace{-0.27cm}
\bibitem{Antoniadis12}
Antoniadis, A., Brosat, X., Cugliari, J. \& Poggi, J. M. (2012). Pr\'evision d'un processus \`a valeurs fonctionnelles en pr\'esence de non stationnarit\'es. Application \`a la consommation d'\'electricit\'e. \emph{J. SFdS}, \textbf{153}, 52--78.

\vspace{-0.27cm}
\bibitem{Antoniadis06}
Antoniadis, A., Paparoditis, E. \& Sapatinas, T. (2006). A functional wavelet-kernel approach for time series prediction. \emph{J. R. Stat. Soc. Ser. B. Stat. Methodol.}, \textbf{68}, 837--857.


\vspace{-0.27cm}
\bibitem{Antoniadis09}
Antoniadis, A., Paparoditis, E. \& Sapatinas, T. (2009). Bandwidth selection for functional time series prediction. \emph{Statist. Probab. Lett.}, \textbf{79}, 733--740.

\vspace{-0.27cm}
\bibitem{Antoniadis03}
Antoniadis, A. \& Sapatinas, T. (2003). Wavelet methods for continuous-time
prediction using Hilbert-valued autoregressive processes. \emph{J.
Multivariate Anal.}, \textbf{87}, 133--158.


\vspace{-0.27cm}
\bibitem{Aueetal15}
Aue, A., Norinho, D. \&  H\"ormann, S. (2015). On the prediction of stationary functional time series. \emph{J. Amer. Statist. Assoc.}, \textbf{110}, 378--392 .



%\vspace{-0.27cm}
%\bibitem{Benko09}
%Benko, M., H\"ardle, W. \& Kneip, A. (2009). Common functional principal components. \emph{Ann. Statist.}, \textbf{37}, 1--34.

\vspace{-0.27cm}
\bibitem{Bensmain01}
Bensmain, N. \& Mourid, T. (2001). Estimateur ``sieve" de l'op\'erateur d'un
processus ARH$(1)$. \emph{C. R. Acad. Sci. Paris S\'er. I Math.}, \textbf{332},
1015--1018.

\vspace{-0.27cm}
\bibitem{Benyelles01}
Benyelles, W. \& Mourid, T. (2001). Estimation de la p\'eriode d'un processus \`a temps continu \`a repr\'esentation autor\'egressive. \emph{C. R. Acad. Sci. Paris S\'er. I Math.}, \textbf{333}, 245--248.

%\vspace{-0.27cm}
%\bibitem{Berkes09}
%Berkes, I., Gabrys, R., Horv\'ath, L. \& Kokoszka, P. (2009). Detecting changes in the mean of functional observations. \emph{J. R. Stat. Soc. Ser. B. Stat. Methodol.}, \textbf{71}, 927--946.
%
%\vspace{-0.27cm}
%\bibitem{Bernard97}
%Bernard, P. (1997). \emph{Analyse de signaux physiologiques}. Mémoire Univ. Cathol., Angers.

\vspace{-0.27cm}
\bibitem{Besse96}
Besse, P.C. \& Cardot, H. (1996). Approximation spline de la pr\'evision d'un processu fonctionnel autoregr\'essif d'ordre 1. \emph{Canad. J. Statist.}, \textbf{24}, 467--487.

%\vspace{-0.27cm}
%\bibitem{Besse97}
%Besse, P.C., Cardot, H. \& Ferraty, F. (1997). Simultaneous non-parametric
%regressions of unbalanced longitudinal data. \emph{Comput. Statist. Data Anal.}, \textbf{24}, 255--270.


\vspace{-0.27cm}
\bibitem{Besse00}
Besse, P. C., Cardot, H. \& Stephenson, D. B. (2000). Autoregressive
forecasting of some functional climatic variations. \emph{Scand. J.
Stat.}, \textbf{27}, 673--687.



\vspace{-0.27cm}
\bibitem{Blanke14}
Blanke, D. \& Bosq, D. (2014). Exponential bounds for intensity of jumps. \emph{Math. Methods Statist.}, \textbf{23}, 239--255.


\vspace{-0.27cm}
\bibitem{Bosq91}
Bosq, D. (1991). Modelization, non-parametric estimation and prediction for continuous time processes. \emph{Nonparametric functional estimation and related topics, NATO, ASI Series},  \textbf{335}, 509--529.

%\vspace{-0.27cm}
%\bibitem{Bosq96}
%Bosq, D. (1996). Limit theorems for Banach-valued autoregressive processes. Applications to real continuous time processes. \emph{Bull. Belg. Math. Soc. Simon Stevin}, \textbf{9}, 537--555.

\vspace{-0.27cm}
\bibitem{Bosq99a}
Bosq, D. (1999a). Autoregressive representation for the empirical covariance operator of an ARH$(1)$. \emph{C. R.
Acad. Sci. Paris S\'er. I Math.}, \textbf{329}, 531--534.

\vspace{-0.27cm}
\bibitem{Bosq99b}
Bosq, D. (1999b). Autoregressive Hilbertian Processes. \emph{Annales de l'ISUP}, \textbf{43}, 25--55.

\vspace{-0.27cm}
\bibitem{Bosq00}
Bosq, D. (2000). \emph{Linear Processes in Function Spaces}. Springer, New York.

%
%\vspace{-0.27cm}
%\bibitem{Bosq03}
%Bosq, D. (2003a). Berry-Esseen inequality for linear processes in Hilbert spaces. \emph{Statist. Probab. Lett.}, \textbf{63}, 243--247.

%\vspace{-0.27cm}
%\bibitem{Bosq03}
%Bosq, D. (2003b). Estimation of autocorrelation operator and prediction for infinite dimensional autoregressive processes. \emph{Math. Methods Statist.}, \textbf{11}, 381--401.

\vspace{-0.27cm}
\bibitem{Bosq07}
Bosq, D. (2007). General linear processes in Hilbert spaces and prediction.
\emph{J. Stat. Planning and Inference}, \textbf{137}, 879--894.

%\vspace{-0.27cm}
%\bibitem{Bosq09}
%Bosq, D. (2009). A note on asymptotic parametric prediction. \emph{J. Statist. Plann. Inference}, \textbf{139}, 1506--1513.

\vspace{-0.27cm}
\bibitem{Bosq10}
Bosq, D. (2010). Tensorial products of functional ARMA processes. \emph{J. Multivariate Anal.}, \textbf{101}, 1352--1363.

\vspace{-0.27cm}
\bibitem{Bosq07}
Bosq, D. \& Blanke, D. (2007). \emph{Inference and Predictions in Large
Dimensions}. Wiley, Chichester.

%\vspace{-0.27cm}
%\bibitem{Bosq99}
%Bosq, D. \& Shen, J. (1999). Estimation of an autoregressive semiparametric model with exogenous variables. \emph{J. Statist. Plann. Inference}, \textbf{68}, 105--127.

%\vspace{-0.27cm}
%\bibitem{Bradley97}
%Bradley, R. C. (1997). On quantiles and the central limit theorem question for strongly mixing sequences. \emph{J. Theoret. Probab.}, \textbf{10}, 507--555.
%
%\vspace{-0.27cm}
%\bibitem{Burfield15}
%Burfield, R., Neumann, C. \& Saunders, C. P. (2015). Review and application of functional data analysis to chemical data - The example of the comparison, classification, and database search of forensic ink chromatograms. \emph{Chemometrics and Intelligent Laboratory Systems}, \textbf{149}, 97--106.


%\vspace{-0.27cm}
%\bibitem{Cai06}
%Cai, T. T. \& Hall, P. (2006). Prediction in functional linear regression. \emph{Ann. Statist.}, \textbf{34}, 2159--2179.

%\vspace{-0.27cm}
%\bibitem{Canale16}
%Canale, A. \& Vantini, S. (2016). Constrained functional time series: applications to the Italian gas market. \emph{Int. J. Forecasting}, \emph{32}, 1340--1351.

%
%\vspace{-0.27cm}
%\bibitem{Cardot99}
%Cardot, H., Ferraty, F. \& Sarda, P. (1999). Functional linear model. \emph{Statist. Probab. Lett.}, \textbf{45}, 11--22.

\vspace{-0.27cm}
\bibitem{Cardot98}
Cardot, H. (1998). Convergence du lissage spline de la pr\'evision
des processus autor\'egressifs fonctionnels. \emph{C. R.
Acad. Sci. Paris S\'er. I Math.}, \textbf{326}, 755--758.

%\vspace{-0.27cm}
%\bibitem{Cardot03}
%Cardot, H., Ferraty, F. \& Sarda, P. (2003). Spline estimators for the functional linear model. \emph{Statist. Sinica}, \textbf{13}, 571--591.
%
%
%\vspace{-0.27cm}
%\bibitem{Cardot07}
%Cardot, H., Mas, A. \& Sarda, P. (2007). CLT in functional linear regression models. \emph{Probab. Theory Related Fields}, \textbf{138}, 325--361.
%
%\vspace{-0.27cm}
%\bibitem{Cardot05}
%Cardot, H. \& Sarda, P. (2005). Estimation in generalized linear models for functional data via penalized likelihood. \emph{J. Multivariate Anal.}, \textbf{92}, 24--41.



\vspace{-0.27cm}
\bibitem{Cavallini94}
Cavallini, A., Montanari, G. C., Loggini, M., Lessi, M. \& Cacciari, M. (1994). Nonparametric prediction of harmonic levels in electrical networks. \emph{Proceed. IEEE ICHPS VI Bologna}, 165--171.


%\vspace{-0.27cm}
%\bibitem{Crambes09}
%Crambes, C., Kneip, A. \& Sarda, P. (2009). Smoothing splines estimators for functional linear regression. \emph{Ann. Statist.}, \textbf{37}, 35--72.
%
%\vspace{-0.27cm}
%\bibitem{Cuevas14}
%Cuevas, A. (2014). A partial overview of the theory of statistics with functional data. \emph{J. Statist. Plann. Inference}, \textbf{147}, 1--23.




\vspace{-0.27cm}
\bibitem{Chen16}
Chen, S. X., Lei, L. \& Tu, Y. (2016). Functional coefficient moving average model with applications to forecasting chinese CPI. \emph{Statist. Sinica}, \textbf{26}, 1649--1672.


\vspace{-0.27cm}
\bibitem{Cuevas02}
Cuevas, A., Febrero, M. \& Fraiman, R. (2002). Linear functional regression: the case of fixed design and functional response. \emph{Canad. J. Statist.}, \textbf{30}, 285--300.



\vspace{-0.27cm}
\bibitem{Cugliari11}
Cugliari, J. (2011). \emph{Pr\'evision non param\'etrique de processus \`a
valeurs fonctionnelles. Application \`a la consommation
d'\'electricit\'e}. PhD thesis, University of Paris-Sud 11, Paris.




\vspace{-0.27cm}
\bibitem{Cugliari13}
Cugliari, J. (2013). Conditional Autoregressive Hilbertian processes. arXiv:1302.3488.

\vspace{-0.27cm}
\bibitem{Damon02}
Damon, J. \& Guillas, S. (2002). The inclusion of exogenous variables in functional autoregressive ozone forecasting. \emph{Environmetrics}, \textbf{13}, 759--774.

\vspace{-0.27cm}
\bibitem{Damon05}
Damon, J. \& Guillas, S. (2005). Estimation and simulation of autoregressive
processes with exogenous variables. \emph{Stat. Inference Stoch. Process.}, \textbf{8}, 185--204.



%\vspace{-0.27cm}
%\bibitem{Dautray90}
%Dautray, R. \& Lions, J.-L. (1990). \emph{Mathematical Analysis and
%Numerical Methods for Science and Technology Volume 3: Spectral
%Theory and Applications}. Springer, New York.

%
%
%\vspace{-0.27cm}
%\bibitem{Davis08}
%Davis, R. A. \& Mikosch, T. (2008). Extreme value theory for space–time processes with heavy-tailed distributions. \emph{Stochastic Process. Appl.}, \textbf{118}, 560--584.

\vspace{-0.27cm}
\bibitem{Dedecker02}
Dedecker, J. \& Merlev\`ede, F. (2002). Necessary and sufficient conditions for the conditional central limit
theorem. \emph{Ann. Probab.}, \textbf{30}, 1044--1081.

\vspace{-0.27cm}
\bibitem{Dedecker03}
Dedecker, J. \& Merlev\`ede, F. (2003). The conditional central limit theorem in Hilbert spaces. \emph{Stochastic Process. Appl.}, \textbf{108}, 229--262.

%%
%%\vspace{-0.27cm}
%%\bibitem{Dede09}
%%Dede, S. (2009). Moderate deviations for stationary sequences of Hilbert-valued bounded random variables. \emph{J. Math. Anal. Appl.}, \textbf{349}, 374--394.
%%
%%\vspace{-0.27cm}
%%\bibitem{Dehling86}
%%Dehling, H., Denker, M. \& Philipp, W. (1986). A bounded law of the iterated logarithm for Hilbert space valued martingales and its application to U-statistics. \emph{Probab. Theory Related Fields}, \textbf{72}, 111--131.
%

\vspace{-0.27cm}
\bibitem{Dehling05}
Dehling, H. \& Sharipov, O. S. (2005). Estimation of mean and covariance
operator for Banach space valued autoregressive processes with
dependent innovations. \emph{Stat. Inference Stoch. Process.}, \textbf{8}, 137--149.


\vspace{-0.27cm}
\bibitem{Didericksen12}
Didericksen, D., Kokoszka, P. \& Zhang, X. (2012). Empirical properties of forecast with the functional autoregressive model. \emph{Comput. Stat.}, \textbf{27}, 285--298.

%%\vspace{-0.27cm}
%%\bibitem{Diestel95}
%%Diestel, J., Jarchow, H. \& Tonge, A. (1995). \emph{Absolutely summing operators}.
%%Cambridge Stud. Adv. Math., volume 43. Cambridge University Press, Cambridge.

%\vspace{-0.27cm}
%\bibitem{Elezovic09}
%Elezovi\'c, S. (2009). Functional modelling of volatility in the Swedish limit order book. \emph{Comput. Statist. Data Anal.}, \textbf{53}, 2107--2118.

%\vspace{-0.27cm}
%\bibitem{Fernandez05}
%Fern\'andez de Castro, B.M., Gonz\'alez-Manteiga \& Guillas, S. (2005). Functional Samples and Bootstrap for Predicting Sulfur Dioxide Levels. \emph{Technometrics}, \textbf{2}, 212--222.

%\vspace{-0.27cm}
%\bibitem{Ferraty02}
%Ferraty, F., Goia, A. \& Vieu, P. (2002). Functional nonparametric model for time series: a fractal approach for dimension reduction. \emph{Test}, \textbf{11}, 317--344.

\vspace{-0.27cm}
\bibitem{Ferraty12}
Ferraty, F., Van Keilegom, I. \& Vieu, P. (2012). Regression when both response and predictor are functions. \emph{J. Multivariate Anal.}, \textbf{109}, 10--28.


%
%%\vspace{-0.27cm}
%%\bibitem{Ferraty02}
%%Ferraty, F.  \& Vieu, P. (2002). The functional nonparametric model and its applications to spectrometric data. \emph{Comput. Statist. Data Anal.}, \textbf{17}, 545--564.

\vspace{-0.27cm}
\bibitem{Ferraty06}
Ferraty, F.  \& Vieu, P. (2006).  \emph{Nonparametric functional data analysis: Theory and practice}. Springer, Berlin.

\vspace{-0.27cm}
\bibitem{Fortet95}
Fortet, R. (1995). \emph{Vecteurs, fonctions et distributions aleatoires dans les espaces de Hilbert}. Hermes, Paris.
%
%\vspace{-0.27cm}
%\bibitem{Fraiman01}
%Fraiman, R. \& Muniz, G. (2001). Trimmed means for functional data. \emph{Test}, \textbf{10}, 419--440.

%%\vspace{-0.27cm}
%%\bibitem{Goia16}
%%Goia, A. \& Vieu, P. (2016). An introduction to recent advances in high/infinite  dimensional statistics. \emph{J. Multivariate Anal.}, \textbf{146}, 1--6.


%
%\vspace{-0.27cm}
%\bibitem{Grebenkov13}
%Grebenkov, D. S. \& Nguyen, B. T. (2013). Geometrical structure of Laplacian eigenfunctions. \emph{SIAM Rev.}, \textbf{55}, pp. 601--667.


\vspace{-0.27cm}
\bibitem{Grenander81}
Grenander, U. (1981). \emph{Abstract Inference}. Wiley, New York.

\vspace{-0.27cm}
\bibitem{Guillas00}
Guillas, S. (2000). Noncausality and functional discretization, limit
theorems for an ARHX$(1)$ process. \emph{C. R. Acad. Sci. Paris S\'er.
I Math.}, \textbf{331}, 91--94.

\vspace{-0.27cm}
\bibitem{Guillas01}
Guillas, S. (2001). Rates of convergence of autocorrelation estimates for
autoregressive Hilbertian processes. \emph{Statist. Probab. Lett.}, \textbf{55}, 281--291.

\vspace{-0.27cm}
\bibitem{Guillas02}
Guillas, S. (2002). Doubly stochastic Hilbertian processes. \emph{J. Appl. Probab.}, \textbf{39}, 566--580.


%\vspace{-0.27cm}
%\bibitem{Haan01}
%de Haan, L. \& Lin, T. (2001). On convergence toward an extreme value distribution in $\mathcal{C}\left([0,1]\right)$. \emph{Ann. Probab.}, \textbf{29}, 467--483.


\vspace{-0.27cm}
\bibitem{Hajj11}
El Hajj, L. (2011). Limit theorems for $D[0,1]$-valued autoregressive processes. \emph{C. R. Acad. Sci. Paris S\'er. I Math.}, \textbf{349}, 821--825.

\vspace{-0.27cm}
\bibitem{Hajj13}
El Hajj, L. (2013). \emph{Inf\'erence statistique pour des variables
fonctionnelles \`a Sauts}.  PhD thesis, University Paris VI, Paris.


%%\vspace{-0.27cm}
%%\bibitem{Hoffmann76}
%%Hoffmann-J\o rgensen, J. \& Pisier, G. (1976). The law of large numbers and the central limit theorem in Banach spaces. \emph{Ann. Probab.}, \textbf{4}, 587--599.

\vspace{-0.27cm}
\bibitem{Hormann13}
H\"ormann, S., Horv\'ath, L. \& Reeder, R. (2013). A functional version of the ARCH model. \emph{Econometric Theory}, \textbf{29}, 267--288.

\vspace{-0.27cm}
\bibitem{Hormann15}
H\"ormann, S., Kidzi\'nski, L. \& Hallin, M. (2015). Dynamic functional principal components. \emph{J. R. Stat. Soc. Ser. B. Stat. Methodol.}, \textbf{77}, 319--348.


\vspace{-0.27cm}
\bibitem{Hormann14}
H\"ormann, S. \& Kidzi\'nski, L. (2015). A Note on Estimation in Hilbertian Linear Models. \emph{Scand. J. Stat.}, \textbf{42}, 43--62.



\vspace{-0.27cm}
\bibitem{Hormann10}
H\"ormann, S. \& Kokoszka, P. (2010). Weakly dependent functional data. \emph{ Ann. Statist.}, \textbf{38}, 1845--1884.


\vspace{-0.27cm}
\bibitem{Horvath10}
Horv\'ath, L., Huskov\'a, M. \& Kokoszka, P. (2010). Testing the stability of the functional autoregressive process. \emph{J. Multivariate Anal.}, \textbf{101}, 352--367.

%%\vspace{-0.27cm}
%%\bibitem{Horvath12}
%%Horv\'ath, L. \& Kokoszka, P. (2012). \emph{Inference for functional data with applications}. Springer, New York.
%%
\vspace{-0.27cm}
\bibitem{Horvath14}
Horv\'ath, L., Kokoszka, P. \& Rice, G. (2014). Testing stationarity of functional time series. \emph{J. Econometrics}, \textbf{179}, 66--82.


\vspace{-0.27cm}
\bibitem{Hyndman}
Hyndman, R. J. \&  Shang, H. L. (2008). \emph{Bagplots, Boxplots and Outlier Detection for Functional Data}. In: Functional and Operatorial Statistics. Contributions to Statistics. Physica-Verlag HD.

\vspace{-0.27cm}
\bibitem{Hyndman}
Hyndman, R. J. \&  Shang, H. L. (2009). Forecasting functional time series.  \emph{J. Korean Statist. Soc.}, \textbf{38}, 199--211.


\vspace{-0.27cm}
\bibitem{Hyndman07}
Hyndman, R. J. \& Ullah, Md. S (2007). Robust forecasting of mortality and fertility rates: A functional data approach. \emph{Comput. Statist. Data Anal.}, \textbf{51}, 4942--4956.


\vspace{-0.27cm}
\bibitem{Kara16}
Kara-Terki, N. \& Mourid, T. (2016). Local asymptotic normality of Hilbertian autoregressive processes.  \emph{C. R. Acad. Sci. Paris S\'er. I}, \textbf{354}, 634--638.


\vspace{-0.27cm}
\bibitem{Kargin08}
Kargin, V. \& Onatski, A. (2008). Curve forecasting by functional
autoregression. \emph{J. Multivariate Anal.}, \textbf{99}, 2508--2526.

%\vspace{-0.27cm}
%\bibitem{Klepsch17}
%Klepsch, J. \& Kl\"uppelberg, C. (2017). An innovations algorithm for the prediction of functional linear processes. \emph{J. Multivariate Anal.}, \textbf{155}, 252--271.
%%
%
\vspace{-0.27cm}
\bibitem{Kokoszka08}
Kokoszka, P., Maslova, I., Sojka, J. \& Zhu, L. (2008). Testing for lack of dependence in the functional linear model. \emph{Canad. J. Statist.}, \textbf{36}, 207--222.


\vspace{-0.27cm}
\bibitem{Kokoszka13a}
Kokoszka, P. \& Reimherr, M. (2013a). Asymptotic normality of the principal components of functional time series. \emph{Stochastic Process. Appl.}, \textbf{123}, 1546--1562.

%\vspace{-0.27cm}
%\bibitem{Kokoszka13b}
%Kokoszka, P. \& Riemherr, M. (2013b). Predictability of shapes of intraday price curves. \emph{Econom. J.}, \textbf{16}, 285--308.

\vspace{-0.27cm}
\bibitem{Kokoszka13b}
Kokoszka, P. \& Riemherr, M. (2013b). Determining the order of the functional autoregressive model. \emph{J. Time Series Anal.}, \textbf{34}, 116--129.

\vspace{-0.27cm}
\bibitem{Kowal17}
Kowal, D. R., Matteson, D. S. \& Ruppert, D. (2017). Functional autoregression for sparsely sample data. \emph{J. Bus. Econom. Statist.}, to appear.

%\vspace{-0.27cm}
%\bibitem{Kuelbs70}
%Kuelbs, J. (1970). Gaussian measures on a Banach spaces. \emph{J. Funct. Anal.}, \textbf{5}, 354--367.

\vspace{-0.27cm}
\bibitem{Kuelbs76}
Kuelbs, J. (1976). A strong convergence theorem for Banach valued random variables. \emph{Ann. Prob.}, \textbf{4}, 744--771.

\vspace{-0.27cm}
\bibitem{Labbas02}
Labbas, A. \& Mourid, T. (2002). Estimation et pr\'evision d'un processus autor\'egressif Banach. \emph{C. R. Acad. Sci. Paris S\'er. I}, \textbf{335}, 767--772.

\vspace{-0.27cm}
\bibitem{Laukaitis02}
Laukaitis, A. \& Ra\u{c}kauskas, A. (2002). Functional data analysis of payment systems. \emph{Nonlinear Anal. Model. Control}, \textbf{7}, 53--68.


%\vspace{-0.27cm}
%\bibitem{Liu16}
%Liu, X., Xiao, H. \& Chen, R. (2016). Convolutional autoregressive models for functional time series. \emph{J. Econometrics}, \textbf{194}, 263--282.

%%\vspace{-0.27cm}
%%\bibitem{Maas97}
%%Maa\ss, P. \& Rieder, A. (1997). \emph{Wavelet-accelerated Tikhonov-Phillips regularization with applications}. in: Inverse Problems in Medical Imaging and Nondestructive Testing. Springer-Verlag, Wien, 134--159.
%
%%\vspace{-0.27cm}
%%\bibitem{Mair96}
%%Mair, B. A. \& Ruymgaart, F. H. (1996). Statistical inverse estimation in Hilbert scales. \emph{SIAM J. Appl. Math.}, \textbf{56}, 1424--1444.
%
%\vspace{-0.27cm}
%\bibitem{Malfait03}
%Malfait, N. \& Ramsay, J. O. (2003). The historical functional linear model. \emph{Canad. J. Statist.}, \textbf{31},  115--128.
%
%%\vspace{-0.27cm}
%%\bibitem{Menneteau05}
%%Menneteau, L. (2005). Some laws of the iterated logarithm in Hilbertian
%%autoregressive models. \emph{J. Multivariate Anal.}, \textbf{92}, 405--425.
%
%

\vspace{-0.27cm}
\bibitem{Mallat89}
Mallat, S. G. (1989). A theory of multiresolution signal decomposition: the wavelet representation. \emph{IEEE Trans. Pattn. Anal. Mach. Intell.}, \textbf{11}, 674--693.


\vspace{-0.27cm}
\bibitem{Marion04}
Marion, J.-M. \& Pumo, B. (2004). Comparison of ARH$(1)$ and ARHD$(1)$ models
on physiological data. \emph{Ann. I.S.U.P.}, \textbf{48},  29--38.

\vspace{-0.27cm}
\bibitem{Martin82}
Martin, D. E. K. (1982). \emph{Estimation of the minimal period of periodically correlated sequences}. Ph. D. Thesis. University of Maryland, Maryland.

\vspace{-0.27cm}
\bibitem{Mas99}
Mas, A. (1999). Normalit\'e asymptotique de l'estimateur empirique de
l'op\'erateur d'autocorr\'elation d'un processus ARH$(1)$. \emph{C. R.
Acad. Sci. Paris S\'er. I Math.}, \textbf{329}, 899--902.

\vspace{-0.27cm}
\bibitem{Mas00}
Mas, A. (2000). \emph{Estimation d'op\'erateurs de corr\'elation de processus fonctionnels: lois limites, tests, d\'eviations mod\'er\'ees}. Ph. D. Thesis. University of Paris VI, France. 


\vspace{-0.27cm}
\bibitem{Mas02}
Mas, A. (2002). Weak convergence for the covariance operators of a Hilbertian linear process. \emph{Stochastic Process. Appl.}, \textbf{99}, 117--135.

\vspace{-0.27cm}
\bibitem{Mas04}
Mas, A. (2004). Consistance du pr\'edicteur dans le mod\`ele ARH$(1)$: le cas
compact. \emph{Ann. I.S.U.P.}, \textbf{48}, 39--48.

\vspace{-0.27cm}
\bibitem{Mas07}
Mas, A. (2007). Weak-convergence in the functional autoregressive model.
\emph{J. Multivariate Anal.}, \textbf{98}, 1231--1261.

\vspace{-0.27cm}
\bibitem{Mas03a}
Mas, A. \& Menneteau, L. (2003a). Large and moderate deviations for infinite
dimensional autoregressive processes. \emph{J. Multivariate Anal.}, \textbf{87}, 241--260.

\vspace{-0.27cm}
\bibitem{Mas03b}
Mas, A. \& Menneteau, L. (2003b). Perturbation approach applied to the asymptotic study of random operators. \emph{Progress in Probability}, \textbf{55}, 127--134.

\vspace{-0.27cm}
\bibitem{MasPumo07}
Mas, A. \& Pumo, B. (2007). The ARHD model. \emph{J. Statist. Plann. Inference}, \textbf{137}, 538--553.


\vspace{-0.27cm}
\bibitem{Menneteau05}
Menneteau, L. (2005). Some laws of the iterated logarithm in Hilbertian autoregressive models. \emph{J. Multivariate Anal.}, \textbf{92}, 405--425.

%
%\vspace{-0.27cm}
%\bibitem{Meinguet10}
%Meinguet, T. \& Segers, J. (2010). Regularly varying time series in Banach
%spaces. arXiv:1001.3262


\vspace{-0.27cm}
\bibitem{Merlevede95}
Merlev\`ede, F. (1995). Sur l'inversibilit\'e des processus lin\'eaires \`a valeurs dans un espace de Hilbert. \emph{C. R. Acad. Sci. Paris, S\'er. I}, \textbf{321}, 477--480.

\vspace{-0.27cm}
\bibitem{Merlevede96}
Merlev\`ede, F. (1996). Central limit theorem for linear processes with values in a Hilbert space. \emph{Stochastic Process. Appl.}, \textbf{65}, 103--114.


\vspace{-0.27cm}
\bibitem{Merlevede97}
Merlev\`ede, F. (1997). R\'esultats de convergence presque s\^ure pour l'estimation et la pr\'evision des processus lin\'eaires hilbertiens. \emph{C. R. Acad. Sci. Paris, S\'er. I}, \textbf{324}, 573--576.


\vspace{-0.27cm}
\bibitem{Merlevedeetal97}
Merlev\`ede, F., Peligrad, M. \& Utev, S. (1997). Sharp Conditions for the CLT of Linear Processes in a Hilbert Space. \emph{J. Theoret. Probab. }, \textbf{10}, 681--693.

\vspace{-0.27cm}
\bibitem{Meyer97}
 Meyer, Y. \& Coifman, R. (1997). \emph{Wavelets, Calder\'on-Zygmund and Multilinear Operators}. Cambridge University Press, Cambridge.

\vspace{-0.27cm}
\bibitem{Mokhtari03}
Mokhtari, F. \& Mourid, T. (2003). Prediction of continuous time autoregressive processes via the Reproducing Kernel Spaces. \emph{Stat. Inference Stoch. Process.}, \textbf{6}, 247--266.

\vspace{-0.27cm}
\bibitem{Mourid93}
 Mourid, T. (1993). Processus autor\'egressifs banachiques d'ordre
sup\'erieur. \emph{C. R. Acad. Sci. Paris S\'er. I Math.}, \textbf{317}, 1167--1172.

\vspace{-0.27cm}
\bibitem{Mourid96}
Mourid, T. (1996). Repr\'esentation autor\'egressive dans un espace de Banach de processus r\'eels \`a temps continu et \'equivalence des lois. \emph{C. R. Acad. Sci. Paris S\'er. I Math.}, \textbf{322} , 1219-1224.

\vspace{-0.27cm}
\bibitem{Mourid04}
Mourid, T. (2004). Processus autor\'egressifs hilbertiens \`a coefficients
al\'eatoires. \emph{Ann. I.S.U.P.}, \textbf{48}, 79--85.



\vspace{-0.27cm}
\bibitem{Panaretos13}
Panaretos, V. M. \& Tavakoli, S. (2013). Cram\'er-Karhunen-Lo\`eve representation and harmonic principal component analysis of functional time series. \emph{Stoch. Processes Appl.}, \textbf{123}, 2779--2807.


\vspace{-0.27cm}
\bibitem{Parvardeh17}
Parvardeh, A., Mohammadi, N., Mahmoodi, S. \& Soltani, A. R. (2017). First order autoregressive periodically correlated model in Banach spaces: existence and central limit theorem. \emph{J. Multivariate Anal.}, \textbf{449}, 756--768.


\vspace{-0.27cm}
\bibitem{Parzen61}
Parzen, E. (1961). An approach to time series analysis. \emph{Ann. Math. Statist.}, \textbf{32}, 951--989.
%
%%\vspace{-0.27cm}
%%\bibitem{Parzen63}
%%Parzen, E. (1963). \emph{A new approach to synthesis of optimal smoothing and prediction systems}. Mathematical Optimization Techniques, University of California Press, Berkeley, 75--108.
%
%
%%\vspace{-0.27cm}
%%\bibitem{Peligrad97}
%%Peligrad, M. \& Utev, S. (1997). Central limit theorem for linear processes. \emph{Ann. Probab.}, \textbf{25}, 443--456.
%%
%%\vspace{-0.27cm}
%%\bibitem{Plato90}
%%Plato, R. \& Vainikko, G. (1990). On the regularization of projection method for solving ill-posed problems. \emph{Numer. Math.}, \textbf{57}, 63--79.

\vspace{-0.27cm}
\bibitem{Poggi94}
Poggi, J. M. (1994). Pr\'evision non param\'etrique de la consommation \'electrique. \emph{Revue de Statistique Appliqu\'ee}, \textbf{42}, 93--98.


%%\vspace{-0.27cm}
%%\bibitem{Preda05}
%%Preda, C. \& Saporta, G. (2005). PLS regression on a stochastic processes. \emph{Comput. Statist. Data Anal.}, \textbf{48}, 149--158.

%\vspace{-0.27cm}
%\bibitem{Prchal07}
%Prchal, L. \& Sarda, P. (2007). \emph{Spline estimator for functional linear regression with functional response}. Preprint.

\vspace{-0.27cm}
\bibitem{Pumo92}
Pumo, B. (1992). \emph{Estimation et pr\'evision de processus
autor\'egressifs fonctionnels}. PhD thesis, Universit\'e de Paris
6, Paris.

\vspace{-0.27cm}
\bibitem{Pumo98}
Pumo, B. (1998). Prediction of continuous time processes by
$C_{[0,1]}$-valued autoregressive process. \emph{Stat. Inference Stoch. Process.} , \textbf{1}, 297--309.


\vspace{-0.27cm}
\bibitem{Rachedi05}
Rachedi, F. (2005). Estimateurs cribles des processus autor\'egressifs Banachiques. PhD Thesis, Universit\'e Pierre et Marie Curie - Paris VI, Paris.


%\vspace{-0.27cm}
%\bibitem{Ramsay91}
%Ramsay, J. O. \& Dalzell, C. J. (1991). Sme tools for functional data analysis. \emph{J. R. Stat. Soc. Ser. B. Stat. Methodol.}, \textbf{53}, 539--572.
%
%
%%
%%\vspace{-0.27cm}
%%\bibitem{Ramsay05}
%%Ramsay, J. O. \& Silverman, B. W. (2005). \emph{Functional data analysis}. Springer, New York.
%%
\vspace{-0.27cm}
\bibitem{Reiss07} 
Reiss, P. T. \& Ogden, T. (2007). Functional principal component regression and functional partial least-squares. \emph{J. Amer. Statist. Assoc.}, \textbf{102}, 984--996.

\vspace{-0.27cm}
\bibitem{Ruiz} 
Ruiz-Medina, M. D. (2011). Spatial autoregressive and   moving
average Hilbertian processes. \emph{J. Multivariate
Anal.}, \textbf{102}, 292--305.

\vspace{-0.27cm}
\bibitem{Ruiz12}
Ruiz-Medina, M. D. (2012). Spatial functional prediction from spatial
autoregressive Hilbertian processes. \emph{Environmetrics} \textbf{23}, 119--128.


\vspace{-0.27cm}
\bibitem{Ruiz17a}
Ruiz-Medina, M. D. \& \'Alvarez-Li\'ebana, J. (2017a). Classical and Bayesian
componentwise predictors for non-ergodic ARH(1) processes. \emph{REVSTAT}, in press.

\vspace{-0.27cm}
\bibitem{Ruiz17b}
Ruiz-Medina, M. D. \& \'Alvarez-Li\'ebana, J. (2017b). \emph{Consistent diagonal componentwise ARH(1) prediction}. Submitted.


%
%\vspace{-0.27cm}
%\bibitem{Shao07}
%Shao, X. \& Wu, W. B. (2007). Asymptotic spectral theory for non-linear time series. \emph{Ann. Statist.}, \textbf{35}, 1773--1801.

%%
%%%\vspace{-0.2cm}
%%%\bibitem{Saadatmand17}
%%%Saadatmand, A., Nematollahi, A. R. \& Sadooghi-Alvandi, M. (2017). On the estimation of missing values in AR(1) model with exponential innovations. \emph{Journal Communications in Statistics - Theory and Methods},  \textbf{46}, 3393--3400.
%
%
%%\vspace{-0.27cm}
%%\bibitem{Shang16}
%%Shang, H. L., Smith, P. W. F., Bijak, J. \& Wi\'sniowski, A. (2016). A multilevel functional data method for forecasting population, with an application to the United Kingdom. \emph{Int. J. Forecasting}, \emph{32}, 629--649.
%
%%%\vspace{-0.2cm}
%%%\bibitem{Shang}
%%%Shang, H. L. (2013). Bayesian bandwidth estimation for a nonparametric functional regression model with unknown error density. \emph{Comput. Statist. Data Anal.}, \textbf{67}, 185--198.
%
%\vspace{-0.27cm}
%\bibitem{Sorensen13}
%S\o rensen, H., Goldsmith, J. \& Sangalli, L. M. (2013). An introduction with medical applications to functional data analysis. \emph{Stat. Med.}, \textbf{32}, 5222--5240.

\vspace{-0.27cm}
\bibitem{Soltani11}
Soltani, A. R. \& Hashemi, M. (2011). Periodically correlated autoregressive Hilbertian processes. \emph{Stat. Inference Stoch. Process.}, \textbf{14}, 177--188.


\vspace{-0.27cm}
\bibitem{Spangenberg13}
Spangenberg, F. (2013). Strictly stationary solutions of ARMA equations in Banach spaces. \emph{J. Multivariate Anal.}, \textbf{121}, 127--138.



\vspace{-0.27cm}
\bibitem{Turbillon08}
Turbillon, C., Bosq, D., Marion, J. M. \& Pumo, B. (2008). Estimation du param\`etre des moyennes mobiles hilbertiennes. \emph{C. R. Acad. Sci. Paris S\'er. I Math.}, \textbf{346}, 347--350.


%\vspace{-0.27cm}
%\bibitem{Varadhan62}
%Varadhan, S. R. S. (1962). Limit theorems for sums of independent random
%variables with values in a Hilbert space. \emph{Sankhya A}, \textbf{24}, 213--238.
%
\vspace{-0.27cm}
\bibitem{Wang08}
Wang, H.B. (2008). Nonlinear ARMA models with functional MA coefficients. \emph{J. Time Series Anal.}, \textbf{29}, 1032--1056.



\end{small}
\end{thebibliography}
\end{document}

% ********************************

% ********************************

\title{Supplementary Material: A review and comparative study on functional time series techniques}

\author{Javier \'Alvarez-Li\'ebana$^1$}
\maketitle
\begin{flushleft}
$^1$ Department of Statistics and O.R., University of Granada, Spain. 

\textit{E-mail: javialvaliebana@ugr.es}
\end{flushleft}

\doublespacing

% ********************************

% ********************************

\renewcommand{\absnamepos}{flushleft}
\setlength{\absleftindent}{0pt}
\setlength{\absrightindent}{0pt}
\renewcommand{\abstractname}{Summary}
\begin{abstract}
This document provides, as Supplementary Material to the paper entitled \textit{A review and comparative study on functional time series techniques}, the proof details of the propositions stated in Section 8 of the mentioned paper, as well as the auxiliary results required (see Sections S.1-S.2). Under a non-diagonal framework, Section S.3. provides a theoretical almost sure upper bound for the error in the norm of $\mathcal{S}(H)$ associated with the diagonal componentwise estimator of the autocorrelation operator considered in Section 8.3, when the eigenvectors of the autocovariance operator are unknown. A simulation study is undertaken in Section S.4 to illustrate the large sample behaviour of the formulated estimator in Section 8.3. Tables displaying more detailed  numerical results, corresponding to Figures in Section 9, are provided in Section S.5 of the current document.
\end{abstract}

\section*{S.1. Preliminaries}

 Let  $H$ be a real separable Hilbert space, and let  $X = \left\lbrace X_n,\ n \in \mathbb{Z} \right\rbrace $
 be  a zero-mean stationary ARH(1) process on the basic probability space $(\Omega,\mathcal{A},P)$, satisfying:
\begin{equation}
X_{n } = \rho \left(X_{n-1} \right) + \varepsilon_{n}, \quad \rho \in \mathcal{L}(H), \quad \left\| \rho \right\|_{\mathcal{L}(H)} < 1,\quad n \in
\Z,\label{24bb}
\end{equation}

\noindent  when $H$-valued
innovation process $\varepsilon= \left \lbrace \varepsilon_{n}, \ n\in \mathbb{Z}\right\rbrace$ is assumed to be SWN, and to be  uncorrelated with $X_0$, with 
$\sigma^{2}_{\varepsilon}= {\rm E} \left\lbrace
\|\varepsilon_{n}\|_{H}^{2} \right\rbrace <\infty,$  for all $n\in
\mathbb{Z}$. Under the above-setting of conditions, involved  in the introduction of equation (\ref{24bb}), $X$  admits the following MAH($\infty$) representation (see Bosq, 2000):
\begin{equation}
X_n = \displaystyle \sum_{k=0}^{\infty} \rho^{k} \left( \varepsilon_{n-k} \right), \quad n \in \mathbb{Z}, \label{2}
\end{equation}
\noindent  that provides  the unique stationary solution to equation (\ref{24bb}). In addition, let us consider the following assumptions:

\vspace{0.3cm}

\noindent \textbf{Assumption A1.} The autocovariance operator
$C={\rm E} \left\lbrace X_{n}\otimes X_{n} \right\rbrace,$ for every
$n \in \mathbb{Z},$ is a strictly positive and self-adjoint operator, in the  trace class, with $Ker\left( C \right) = \left\lbrace \emptyset \right\rbrace$. Its eigenvalues $\left\lbrace
C_{j}, \ j\geq 1\right\rbrace $ then satisfy
\begin{equation}
 \displaystyle \sum_{j=1}^{\infty} C_j < \infty,\quad  C_1  > \ldots > C_j > C_{j+1} > \ldots > 0, \quad C (f)(g)=\displaystyle \sum_{j=1}^{\infty}C_{j}\left\langle \phi_{j},f\right\rangle_{H}\left\langle \phi_{j},g\right\rangle_{H},~ \forall f,g\in H.
\label{eqautocov2}
\end{equation}

\vspace{0.3cm}

\noindent \textbf{Assumption A2.} The autocorrelation operator
$\rho$ is a self-adjoint and Hilbert-Schmidt operator, admitting
the following diagonal spectral decomposition:
\begin{equation}
\rho (f)(g)= \displaystyle \sum_{j=1}^{\infty} \rho_j \left\langle \phi_{j},f\right\rangle_{H}\left\langle \phi_{j},g\right\rangle_{H}, \quad \displaystyle \sum_{j=1}^{\infty} \rho_{j}^{2} < \infty, \quad  \forall f,g\in H,
\label{12}
\end{equation}
\noindent where $\left\lbrace \rho_j,\  j \geq 1 \right\rbrace $ is
the system of eigenvalues of $\rho,$
with respect to the orthonormal system $\left\lbrace \phi_j,\  j \geq 1 \right\rbrace $.

\vspace{0.3cm}

\noindent \textbf{Assumption A3.}  The random initial condition in (\ref{24bb}), $X_{0},$ satisfies ${\rm E} \left\lbrace \left\| X_0 \right\|_{H}^{4} \right\rbrace < \infty.$

\vspace{0.3cm}

Under \textbf{Assumptions A1-A2}, the cross-covariance operator $D = \rho C ={\rm E} \left\lbrace X_{n}\otimes X_{n+1} \right\rbrace$, for each $n \in \mathbb{Z}$, can be also diagonally decomposed, respect to the eigenvectors of $C$, providing the set of eigenvalues $\left\lbrace D_j = \rho_j C_j, \ j \geq 1 \right\rbrace$:
\begin{equation}
D(f)(g) = \rho C(f)(g) = \displaystyle \sum_{j=1}^{\infty} \rho_j C_j \langle \phi_j, f \rangle_H \langle \phi_j, g \rangle_H =  \displaystyle \sum_{j=1}^{\infty} D_j \langle \phi_j, f \rangle_H \langle \phi_j, g \rangle_H, \quad f,g \in H.
\end{equation}

 Moreover, projections of (\ref{24bb}) into  $\left\lbrace \phi_j,\ j \geq 1 \right\rbrace $ lead to the stationary zero-mean AR(1) representation, under $\left\| \rho \right\|_{\mathcal{L}(H)} = \displaystyle \sup_{j \geq 1} \left| \rho_j \right| < 1$:
\begin{equation}
X_{n,j}=\rho_{j}X_{n-1,j}+\varepsilon_{n,j},\quad X_{n,j} = \left\langle
X_{n},\phi_{j}\right\rangle_{H},~\varepsilon_{n,j}=\left\langle \varepsilon_{n},\phi_{j}\right\rangle_{H}, \quad \rho_j \in \mathbb{R},~\left| \rho_j \right| < 1, \quad j\geq 1, ~n \in \mathbb{Z}.
\label{14}
\end{equation}

\subsection*{S.1.1. Diagonal strongly-consistent estimator when the eigenvectors of $C$ are unknown}

When the  eigenvectors  $\left\lbrace \phi_j,\ j \geq 1 \right\rbrace $ of $C$ are known, and under \textbf{Assumption A2}, the following estimators of the covariance  operators, based on the estimation of the eigenvalues of their spectral decomposition, will be considered, for each $n \geq 2$:
\begin{eqnarray}
\widehat{C}_n &=& \displaystyle \sum_{j=1}^{\infty} \widehat{C}_{n,j} \phi_j \otimes \phi_j, \quad \widehat{C}_{n,j} = \frac{1}{n} \displaystyle \sum_{i=0}^{n-1} X_{i,j}^{2}, \quad j \geq 1, \label{65aa_a} \\
\widehat{D}_n &=& \displaystyle \sum_{j=1}^{\infty} \widehat{D}_{n,j} \phi_j \otimes \phi_j,\quad \widehat{D}_{n,j} = \frac{1}{n-1} \displaystyle \sum_{i=0}^{n-2} X_{i,j}X_{i+1,j}, \quad j \geq 1, \label{65aa}
\end{eqnarray}
\noindent where $\left\lbrace \phi_j,\ j \geq 1 \right\rbrace$ is the complete orthonormal eigenvectors system of $C$, with $\left\lbrace \widehat{C}_{n,j},\ j \geq 1 \right\rbrace$ and $\left\lbrace \widehat{D}_{n,j} ,\ j \geq 1 \right\rbrace$ being the eigenvalues of operators $\widehat{C}_n$ and $\widehat{D}_n$, respectively, for each $n \geq 2$.

\vspace{0.3cm}
\begin{remark}\label{rem1}  
Under definitions in equations (\ref{65aa_a})-(\ref{65aa}), the diagonal componentwise estimator, introduced in equation  (\ref{48}) below, for the autocorrelation operator $\rho ,$ naturally arises, which is different from the componentwise estimator approaches based on the projection of the natural  empirical covariance  operators $C_n$ and $D_n,$  given by
\begin{equation}
C_n = \frac{1}{n}\displaystyle \sum_{i=0}^{n-1} X_i \otimes X_i, \quad D_n = \frac{1}{n-1}\displaystyle
\sum_{i=0}^{n-2} X_i \otimes X_{i+1},\quad n\geq 2,\label{62aa}
\end{equation}

\noindent into  the eigenvectors $\{\phi_{j},\ j\geq 1\},$ in the case where they are known.
\end{remark}

\vspace{0.3cm}
For the derivation of the subsequently results,  we  will also need the following assumption:

\vspace{0.3cm}

\noindent \textbf{Assumption A4.}  $X_{0,j}^{2}= \left\langle X_{0}, \phi_{j}\right\rangle_{H}^{2}>0,$ a.s., for every $j\geq 1.$

\vspace{0.3cm}

\begin{remark}
\label{remark1aa}

From Cauchy-Schwarz's inequality, for any $j \geq 1$ and $n \geq 2,$
\begin{eqnarray}
\left| \frac{1}{n-1} \displaystyle \sum_{i=0}^{n-2} X_{i,j} X_{i+1,j} \right| &\leq &
 2 \left( \frac{1}{n}\displaystyle \sum_{i=0}^{n-1} X_{i,j} ^{2} \right), \label{12bb}
\end{eqnarray}

\noindent which implies, for each $j \geq 1$, under  \textbf{Assumption A4},
\begin{eqnarray}
\left| \frac{\frac{1}{n-1} \displaystyle \sum_{i=0}^{n-2} X_{i,j} X_{i+1,j}}{\frac{1}{n} \displaystyle \sum_{i=0}^{n-1} X_{i,j}^{2}} \right| \leq \frac{2 \left( \frac{1}{n}\displaystyle \sum_{i=0}^{n-1} X_{i,j} ^{2} \right)}{\frac{1}{n} \displaystyle \sum_{i=0}^{n-1} X_{i,j}^{2}} =2~a.s. \label{13bb}
\end{eqnarray}

\end{remark}

\vspace{0.3cm}

From  \textbf{Assumption A4}, let us consider the diagonal  componentwise estimator of $\rho$
\begin{eqnarray}
\widehat{\rho}_{k_n} &=& \displaystyle \sum_{j=1}^{k_n} \widehat{\rho}_{n,j} \phi_j \otimes \phi_j, \quad \widehat{\rho}_{n,j} = \frac{\widehat{D}_{n,j}}{\widehat{C}_{n,j}}  = \frac{n}{n-1} \frac{\displaystyle \sum_{i=0}^{n-2} X_{i,j}X_{i+1,j}}{\displaystyle \sum_{i=0}^{n-1} X_{i,j}^{2}}, \quad j \geq 1,~n \geq 2. \label{48}
\end{eqnarray}

\begin{remark}
\label{rem__3}

Note that,  under  \textbf{Assumption A1}, the  eigenvalues of $C$ are strictly positive, with multiplicity one, and  $C(H)=H,$
where $C(H)$  denotes the  range   of $C.$
For $f,~g \in C(H)$, there  exist $\varphi, \psi \in H$ such that $f = C\left( \varphi \right)$ and $g = C \left( \psi \right),$ and the following identities hold:
\begin{eqnarray}
 \langle f, g \rangle_{C(H)}&=&   \langle C^{-1}
C \left( \varphi \right), C^{-1} C \left( \psi \right) \rangle_{H} = \langle  \varphi ,  \psi  \rangle_{H}<\infty , \nonumber \\
 \left\| f \right\|_{C(H)}^{2} &=&  \langle C^{-1} C \left( \varphi \right), C^{-1} C \left( \varphi \right) \rangle_{H} =  \left\| \varphi \right\|_{H}^{2}<\infty. \label{15a}
\end{eqnarray}

\noindent Note that, from  Parseval's identity, for every $x \in C(H)$, $\left\| x \right\|_{C(H)}^{2} = \displaystyle \sum_{j=1}^{\infty} \frac{\left[\langle x, \phi_j \rangle_{H}\right]^{2}}{C_{j}^{2}} < \infty$. Thus, the range of $C$ can be also defined by
\begin{equation}
C(H) = \left\lbrace x \in H:\quad \displaystyle \sum_{j=1}^{\infty} \frac{\langle x, \phi_j \rangle_{H}^{2}}{C_{j}^{2}} < \infty \right\rbrace. \label{14a}
\end{equation}

\end{remark}

Under \textbf{Assumptions A1-A4}, the following proposition provides the strong-consistency, in the norm of $\mathcal{L}(H)$, of the estimator (\ref{48}) of the autocorrelation operator, as well as of its associated ARH(1) plug-in predictor, in the underlying Hilbert space. Asymptotic properties derived in Section S.2 below are required.

\begin{proposition}
\label{pr2m}
\textit{Under \textbf{Assumptions A1--A4}, for a truncation parameter $k_n < n$, with $\displaystyle \lim_{n \to \infty} k_n= \infty$,
\begin{equation}
\frac{n^{1/4}}{\left(\ln(n) \right)^{\beta}}\left\| \widehat{\rho}_{k_n}- \rho\right\|_{\mathcal{L}(H)} \longrightarrow^{a.s.} 0, \quad 
\| \left( \widehat{\rho}_{k_n} -\rho \right) (X_{n-1})\|_{H} \longrightarrow^{a.s.} 0,\quad n\rightarrow \infty.
\end{equation}}
\end{proposition}

\begin{proof}

Under \textbf{Assumption A1}, $C(H)=H,$ as a set of functions.
Then, for every $x\in C(H)=H,$ under \textbf{Assumptions A2} and \textbf{A4}, from Parseval's identity and Remark \ref{remark1aa},
\begin{eqnarray}
\|(\widehat{\rho}_{k_{n}}-\rho ) (x)\|_{H}^{2}&= &\sum_{j=1}^{k_{n}}\left[(\widehat{\rho}_{n,j}-\rho_{j})\left\langle x,\phi_{j}
\right\rangle_{H}\right]^{2}+\sum_{j=k_{n}}^{\infty}\left[\rho_{j}\left\langle x,\phi_{j}
\right\rangle_{H}\right]^{2} \leq \sum_{j=1}^{k_{n}}\left[\frac{D_{j}-\widehat{D}_{n,j}}{C_{j}}\left\langle x,\phi_{j}
\right\rangle_{H}\right]^{2}
\nonumber\\
&+&\sum_{j=1}^{k_{n}}\left[\frac{C_{j}-\widehat{C}_{n,j}}{C_{j}}\widehat{\rho}_{n,j}\left\langle x,\phi_{j}
\right\rangle_{H}\right]^{2}+\sum_{j=k_{n}}^{\infty}\left[\rho_{j}\left\langle x,\phi_{j}
\right\rangle_{H}\right]^{2}\nonumber\\
&\leq &\left[\sup_{1\leq j\leq k_{n}}\left|D_{j}-\widehat{D}_{n,j}\right|^{2}
+2\sup_{1\leq j\leq k_{n}}\left|C_{j}-\widehat{C}_{n,j}\right|^{2}\right]\nonumber\\
& \times &\sum_{j=1}^{k_{n}}\left[\frac{\left\langle x,\phi_{j}
\right\rangle_{H}}{C_{j}}\right]^{2}+\sum_{j=k_{n}}^{\infty}\left[\rho_{j}\left\langle x,\phi_{j}
\right\rangle_{H}\right]^{2},\quad \mbox{a.s.}
\label{eqrange}
\end{eqnarray}
Thus, taking the square root in booth sides of (\ref{eqrange}), and the supremum in $x\in H=C(H),$ with $\|x\|_{H}=1,$ at the left-hand side, we obtain
\begin{eqnarray}
 \|\widehat{\rho}_{k_{n}}-\rho \|_{\mathcal{L}(H)}& \leq & \sup_{x\in H,~\|x\|_{H}=1}
\left(\left[\sup_{1\leq j\leq k_{n}}\left|D_{j}-\widehat{D}_{n,j}\right|^{2}
+2\sup_{1\leq j\leq k_{n}}\left|C_{j}-\widehat{C}_{n,j}\right|^{2}\right]\right.\nonumber\\
& & \left. \times \sum_{j=1}^{k_{n}}\left[\frac{\left\langle x,\phi_{j}
\right\rangle_{H}}{C_{j}}\right]^{2}+\sum_{j=k_{n}}^{\infty}\left[\rho_{j}\left\langle x,\phi_{j}
\right\rangle_{H}\right]^{2}\right)^{1/2}\quad \mbox{a.s.}\label{lim0}
\end{eqnarray}

Furthermore, from \textbf{Assumptions A1-A2} and Remark \ref{rem__3}, for every $x\in C(H)=H,$
\begin{eqnarray}
\lim_{n\rightarrow \infty }\sum_{j=1}^{k_{n}}\left[\frac{\left\langle x,\phi_{j}
\right\rangle_{H}}{C_{j}}\right]^{2}=\|x\|^{2}_{C(H)}<\infty, \quad \lim_{n\rightarrow \infty}\sum_{j=k_{n}}^{\infty}\left[\rho_{j}\left\langle x,\phi_{j}
\right\rangle_{H}\right]^{2}=0.\label{lim2}
\end{eqnarray}

Under \textbf{Assumptions A1-A3}, from Corollary \ref{corollary2}  (see equations (\ref{c110})-(\ref{c11a}) in the Section S.2 below), as $n\rightarrow \infty,$
\begin{eqnarray}
\frac{n^{1/4}}{\left(\ln(n) \right)^{\beta}}\sup_{1\leq j\leq k_{n}}\left|C_{j}-\widehat{C}_{n,j}\right| \to^{a.s.}0, \quad \frac{n^{1/4}}{\left(\ln(n) \right)^{\beta}}\sup_{1\leq j\leq k_{n}}\left|D_{j}-\widehat{D}_{n,j}\right| \to^{a.s.}0.
\label{lim4}
\end{eqnarray}

Finally, from equations (\ref{lim0})--(\ref{lim4}), as $n\rightarrow \infty,$
\begin{equation}
\frac{n^{1/4}}{\left(\ln(n) \right)^{\beta}}\|\widehat{\rho}_{k_{n}}-\rho \|_{\mathcal{L}(H)}\to^{a.s.}0.
\end{equation}

Strong-consistency of the associated plug-in predictor is directly derived keeping in mind that
\begin{equation}
\left\| \left( \widehat{\rho}_{k_n} - \rho \right) \left( X_{n-1} \right) \right\|_H \leq \left\| \widehat{\rho}_{k_n} - \rho  \right\|_{\mathcal{L}(H)}  \left\| X_{n-1}  \right\|_H, \quad \left\| X_{n-1}  \right\|_H < \infty~a.s.
\end{equation}

\hfill \hfill $\blacksquare$
\end{proof}

\subsection*{S.1.2. Diagonal strongly-consistent estimator when the eigenvectors of $C$ are unknown}

 In the case of $\{\phi_{j},\ j\geq 1\}$ are unknown, as is often in practice, $C_n = \frac{1}{n} \displaystyle \sum_{i=0}^{n-1} X_i \otimes X_i$ admits a diagonal spectral decomposition in terms of  $\{C_{n,j},\ j\geq 1\}$ and $\{\phi_{n,j},\ j\geq 1\}$, satisfying, for each $n \geq 2$:
 \begin{eqnarray}
C_{n} \left(\phi_{n,j} \right) &=& C_{n,j}\phi_{n,j},\quad j\geq 1,\quad C_{n,1}\geq \dots\geq C_{n,n}\geq 0=C_{n,n+1}=C_{n,n+2}=\dots \label{ev1}\\
C_n &=& \frac{1}{n} \displaystyle \sum_{i=0}^{n-1} X_i \otimes X_i=\displaystyle \sum_{j=1}^{\infty} C_{n,j} \phi_{n,j} \otimes \phi_{n,j}, \quad C_{n,j} =\frac{1}{n} \displaystyle \sum_{i=0}^{n-1} \widetilde{X}_{i,j}^{2},~j \geq 1.\label{66aa}
\end{eqnarray}

In the remainder, we will denote $\widetilde{X}_{i,j} = \langle X_i, \phi_{n,j} \rangle_H$ and $\phi_{n,j}^{\prime} = \mbox{sgn}\left\langle \phi_{n,j} , \phi_{j} \right\rangle_{H}\phi_{j}$, for each $ i\in \mathbb{Z},~j\geq 1$ and $n\geq 2$, where $\mbox{sgn} \langle \phi_{n,j}, \phi_j \rangle_H= \mathbf{1}_{\langle \phi_{n,j}, \phi_j \rangle_H\geq 0}-\mathbf{1}_{\langle \phi_{n,j}, \phi_j \rangle_H< 0}$. Since $\left\lbrace \phi_{n,j} , \ j \geq 1 \right\rbrace$ is a complete orthonormal system of eigenvectors, for each $n\geq 2$, operator $D_n = \frac{1}{n-1} \displaystyle \sum_{i=0}^{n-2} X_i \otimes X_{i+1}$ admits the following non-diagonal spectral representation:
\begin{eqnarray}
D_n = \frac{1}{n-1}\displaystyle \sum_{i=0}^{n-2} X_i \otimes X_{i+1}= \displaystyle \sum_{j=1}^{\infty} \displaystyle
\sum_{l=1}^{\infty} D_{n,j,l}^{*} \phi_{n,j} \otimes \phi_{n,l} = \displaystyle \sum_{j=1}^{\infty}\sum_{l=1}^{\infty}\frac{1}{n-1}\displaystyle
\sum_{i=0}^{n-2} \widetilde{X}_{i,j}\widetilde{X}_{i+1,l}\phi_{n,j} \otimes \phi_{n,l}, \label{37aaa}
\end{eqnarray}
\noindent where $D_{n,j,l}^{*} = \langle D_n \left(\phi_{n,j} \right), \phi_{n,l}  \rangle_H= \frac{1}{n-1}\displaystyle
\sum_{i=0}^{n-2}\widetilde{X}_{i,j}\widetilde{X}_{i+1,l},$ for each $j,l \geq 1$ and $n \geq 2$. In particular, we will use the notation $D_{n,j}=D_{n,j,j}^{*} =  \langle D_n \left(\phi_{n,j} \right), \phi_{n,j}  \rangle_H$. The following  assumption is here deemed:

\vspace{0.3cm}

\noindent \textbf{Assumption A5.} $C_{n,k_n} > 0~a.s,$ where $k_n$ is a suitable truncation parameter $k_n < n$, with $\displaystyle \lim_{n \to \infty} k_n = \infty$.

\vspace{0.3cm}

From \textbf{Assumption A5}, the following  diagonal componentwise estimator of  $\rho$ is outlined:
\begin{eqnarray}
\widetilde{\rho}_{k_n} &=& \displaystyle \sum_{j=1}^{k_n} \widetilde{\rho}_{n,j} \phi_{n,j} \otimes \phi_{n,j}, \quad \widetilde{\rho}_{n,j} =  \frac{D_{n,j}}{C_{n,j}}  = \frac{n}{n-1} \frac{\displaystyle \sum_{i=0}^{n-2} \widetilde{X}_{i,j} \widetilde{X}_{i+1,j}}{\displaystyle \sum_{i=0}^{n-1} \widetilde{X}_{i,j}^{2}}, \quad j \geq 1,~n \geq 2. \label{141}
\end{eqnarray}

 Under  \textbf{Assumptions A1-A3} and \textbf{A5},  the strong-consistency of the diagonal  componentwise estimator  $\widetilde{\rho}_{k_n}$ of $\rho$ is reached in Proposition \ref{proposition5} (see auxiliary results in Section S.2 below).

\begin{proposition}
\label{proposition5}
\textit{Let $k_{n}$ a sequence of integers such that, for certain $\widetilde{n}_{0}$ sufficiently large, and $\beta > \frac{1}{2},$
\begin{eqnarray}
&&\Lambda_{k_{n}}=o\left(n^{1/4}(\ln(n))^{\beta -1/2} \right), \quad \frac{1}{C_{k_n}} \displaystyle \sum_{j=1}^{k_n} a_j = \mathcal{O} \left( n^{1/4} \left(\ln(n)\right)^{-\beta} \right), \quad k_{n}C_{k_{n}}<1,~n\geq \widetilde{n}_{0}, \label{1444}
\end{eqnarray}
\noindent where $\Lambda_{k_n}=\displaystyle \sup_{1\leq j\leq k_n}(C_{j}-C_{j+1})^{-1}$, $a_1 = 2 \sqrt{2} \frac{1}{C_1 - C_2}$ and $a_j = 2 \sqrt{2} \displaystyle \max \left( \frac{1}{C_{j-1} - C_j},\frac{1}{C_{j} - C_{j+1}} \right)$, for any $2 \leq j \leq k_n$. Then, under  \textbf{Assumptions A1-A3}  and \textbf{A5},
\begin{equation}
\left\| \widetilde{\rho}_{k_n}- \rho\right\|_{\mathcal{L}(H)} \longrightarrow^{a.s.} 0,\quad \| \left(\widetilde{\rho}_{k_n} - \rho \right)(X_{n-1})\|_{H} \longrightarrow^{a.s.} 0,\quad n\rightarrow \infty. 
\label{eqboundedoperators}
\end{equation}}
\textit{In particular, the following upper bound can be derived:
\begin{equation}
\left\| \widetilde{\rho}_{k_n}- \rho \right\|_{\mathcal{L}(H)} \leq \displaystyle \sup_{1 \leq j \leq k_n} \left| \widetilde{\rho}_{n,j} - \frac{D_{n,j}}{C_j} \right| + \displaystyle \sup_{1 \leq j \leq k_n} \left|  \frac{D_{n,j}}{C_j} - \rho_j \right| + 2\displaystyle \sum_{j=1}^{k_n} \frac{\left| D_{n,j} \right|}{C_j} \left\| \phi_{n,j} - \phi_{n,j}^{\prime} \right\|_{H} +\displaystyle \sup_{j > k_n} \left|  \rho_j \right|. \label{uppbound}
% \left\| \displaystyle \sum_{j=1}^{k_n} \left(\widetilde{\rho}_{n,j} - \frac{D_{n,j}}{C_j} \right) \phi_{n,j} \otimes \phi_{n,j} \right\|_{\mathcal{L}(H)} + \left\| \displaystyle \sum_{j=1}^{k_n} \left(\frac{D_{n,j}}{C_{j}} - \rho_j \right) \phi_{n,j}^{\prime} \otimes \phi_{n,j}^{\prime}  \right\|_{\mathcal{L}(H)} \nonumber \\
%&+& 2\displaystyle \sum_{j=1}^{k_n} \frac{D_{n,j}}{C_j} \left\| \phi_{n,j} 
\end{equation}}
\end{proposition}

\begin{proof}

Under \textbf{Assumptions A1-A2} and equation (\ref{141}), for every $x\in H,$
\begin{eqnarray}
\left\| \widetilde{\rho}_{k_n} (x)  - \rho(x) \right\|_{H} &\leq & \left\| \displaystyle \sum_{j=1}^{k_n} \widetilde{\rho}_{n,j} \langle \phi_{n,j}, x \rangle_H \phi_{n,j} - \displaystyle \sum_{j=1}^{k_n} \rho_j \langle \phi_{j}, x \rangle_H \phi_{j} \right\|_{H} \nonumber \\
&+& \left\|  \displaystyle \sum_{j=1}^{k_n} \rho_j \langle \phi_{j}, x \rangle_H \phi_{j} - \displaystyle \sum_{j=1}^{\infty} \rho_j \langle \phi_{j}, x \rangle_H \phi_{j} \right\|_{H}  = a_{k_n} (x) + b_{k_n} (x). \label{146}
\end{eqnarray}

Clearly, under \textbf{Assumption A2}, $\displaystyle \lim_{n \to \infty} b_{k_n} (x)  = 0$. Let us now study the behavior of the term  $a_{k_n} (x).$ From equations (\ref{ev1})-(\ref{141}), and under \textbf{Assumption A5},
\begin{eqnarray}
a_{k_n} (x) & \leq & \left\| \displaystyle \sum_{j=1}^{k_n} \left(\frac{D_{n,j}}{C_{n,j}} - \frac{D_{n,j}}{C_{j}}  \right) \langle \phi_{n,j}, x \rangle_H \phi_{n,j} \right\|_{H}  \nonumber \\
& + &  \left\| \displaystyle \sum_{j=1}^{k_n} \frac{D_{n,j}}{C_{j}}  \langle \phi_{n,j}, x \rangle_H \phi_{n,j} - \displaystyle \sum_{j=1}^{k_n} \rho_j \langle \phi_{n,j}^{'}, x \rangle_H \phi_{n,j}^{'} \right\|_{H} = a_{k_n,1} (x) + a_{k_n,2} (x), \label{149}
\end{eqnarray}
\noindent where $ \langle \phi_{j}, x \rangle_H \phi_{j} = \langle \phi_{n,j}^{'}, x \rangle_H \phi_{n,j}^{'} $,  with, as before, $\phi_{n,j}^{'} = \mbox{sgn} \langle \phi_{n,j}, \phi_j \rangle_H \phi_j,$ and  $\mbox{sgn} \langle \phi_{n,j}, \phi_j \rangle_H= \mathbf{1}_{\langle \phi_{n,j}, \phi_j \rangle_H\geq 0}-\mathbf{1}_{\langle \phi_{n,j}, \phi_j \rangle_H< 0},$ for each $j \geq 1$
and $n \geq 2.$

From equation (\ref{149}),
\begin{eqnarray}
a_{k_n,1} (x) \leq  \displaystyle \sum_{j=1}^{k_n} \left| D_{n,j} \right|  \frac{\left| C_{j} - C_{n,j} \right|}{C_{n,j} C_j} \left| \langle \phi_{n,j}, x \rangle_H  \right| \left\| \phi_{n,j} \right\|_{H}  \leq \left\| C - C_n \right\|_{\mathcal{L}(H)} \frac{1}{C_{k_n}} \displaystyle \sum_{j=1}^{k_n} \left| \frac{ D_{n,j}}{C_{n,j}} \right|  \left| \langle \phi_{n,j}, x \rangle_H  \right|.
\end{eqnarray}
Thus,  from Cauchy-Schwarz's inequality, in a similar way to  Remark \ref{remark1aa},
\begin{eqnarray}
a_{k_n,1} (x) &\leq &  \left\| C - C_n \right\|_{\mathcal{L}(H)} \frac{1}{C_{k_n}} \left(\displaystyle \sum_{j=1}^{k_n}   \frac{D_{n,j}^{2}}{C_{n,j}^{2}}  \right)^{1/2} \left( \displaystyle \sum_{j=1}^{\infty} \langle \phi_{n,j}, x \rangle_{H}^{2} \right)^{1/2}  \nonumber \\
&\leq &  2 \left\| C - C_n \right\|_{\mathcal{L}(H)} \frac{1}{C_{k_n}} k_{n}^{1/2} \left\| x \right\|_{H}~a.s. \label{157a}
\end{eqnarray}

From equation (\ref{1444}), for $n\geq \widetilde{n}_{0}$, $k_n < \frac{1}{C_{k_n}}$, which implies that, from Remark \ref{remark__4}  (see Section S.2 below),
\begin{eqnarray}
  a_{k_n,1} (x) &\leq  & 2 \left\| C - C_n \right\|_{\mathcal{L}(H)} C_{k_n}^{-3/2} \left\| x \right\|_{H} <  2 \left\| C - C_n \right\|_{\mathcal{L}(H)}\left\| x \right\|_{H} C_{k_n}^{-1/2} \displaystyle \sum_{j=1}^{k_n} a_j ~a.s. \label{157a}
\end{eqnarray}

From condition (\ref{1444}),  there also exists  a positive real number $M < \infty$ and an integer $n_0$ such that, for certain $\beta > \frac{1}{2}$ and $n \geq n_0$, with  $n_0$ large enough,
\begin{equation}
  C_{k_n}^{-1/2} \displaystyle \sum_{j=1}^{k_n} a_j < C_{k_n}^{-1} \displaystyle \sum_{j=1}^{k_n} a_j \leq M n^{1/4} \left(\ln(n) \right)^{-\beta}. \label{105aaa}
\end{equation}

From equations (\ref{157a})-(\ref{105aaa}), for $n\geq \max(\widetilde{n}_{0},n_{0})$, $a_{k_n,1} (x) < 2 M \frac{n^{1/4}}{\left(\ln(n) \right)^{\beta}} \left\| C - C_n \right\|_{\mathcal{L}(H)}    \left\| x \right\|_{H}$, with $\left\| x \right\|_{H} < \infty$, since $x \in H$. Hence, under \textbf{Assumption A3}, from Theorem \ref{theorem2} (see Section S.2 below),
\begin{equation}
a_{k_n,1} (x)  \to^{a.s.} 0,\quad n\rightarrow \infty. \label{158}
\end{equation}

Let us see now a bound for $a_{k_n,2}(x)$ in (\ref{149}):
\begin{eqnarray}
a_{k_n,2}(x) & \leq &  \left\| \displaystyle \sum_{j=1}^{k_n} \frac{D_{n,j}}{C_{j}}  \left( \langle \phi_{n,j}, x \rangle_H - \langle \phi_{n,j}^{'}, x \rangle_H \right) \phi_{n,j}  \right\|_{H}  \nonumber \\
&+& \left\|\displaystyle \sum_{j=1}^{k_n}\frac{D_{n,j}}{C_{j}}  \langle \phi_{n,j}^{'}, x \rangle_H \left(\phi_{n,j} - \phi_{n,j}^{'} \right) \right\|_{H} + \left\| \displaystyle \sum_{j=1}^{k_n} \left( \frac{D_{n,j}}{C_{j}} - \rho_j \right) \langle \phi_{n,j}^{'}, x \rangle_H \phi_{n,j}^{'} \right\|_{H}. \label{eq___36}
\end{eqnarray}

In a similar way to Remark \ref{remark1aa}, under \textbf{Assumptions A1-A2}, we get
\begin{eqnarray}
a_{k_n,2}(x) & \leq & 2 \displaystyle \sup_{j \geq 1} \left| C_{n,j} \right| C_{k_n}^{-1} \displaystyle \sum_{j=1}^{k_n} \left| \langle \phi_{n,j} -  \phi_{n,j}^{'}, x \rangle_H \right| \left\| \phi_{n,j}   \right\|_{H} \nonumber \\
&+& 2 \displaystyle \sup_{j \geq 1} \left| C_{n,j} \right| C_{k_n}^{-1} \displaystyle \sum_{j=1}^{k_n} \left| \langle \phi_{n,j}^{'}, x \rangle_H \right| \left\| \phi_{n,j} - \phi_{n,j}^{'} \right\|_{H} \nonumber \\
&+&  \displaystyle \sup_{j \geq 1} \left| D_{n,j} - D_j \right| C_{k_n}^{-1} \left\| \displaystyle \sum_{j=1}^{k_n}  \langle \phi_{n,j}^{'}, x \rangle_H \phi_{n,j}^{'} \right\|_{H}~a.s.
\end{eqnarray}
\noindent Hence, from Cauchy-Schwarz's inequality,
\begin{eqnarray}
a_{k_n,2}(x) & \leq &  2 \displaystyle \sup_{j \geq 1} \left| C_{n,j} \right|  \left\| x \right\|_{H} C_{k_n}^{-1} \displaystyle \sum_{j=1}^{k_n} \left\| \phi_{n,j} -  \phi_{n,j}^{'} \right\|_{H}  \nonumber \\
&+& 2 \displaystyle \sup_{j \geq 1} \left| C_{n,j} \right|  \left\| x \right\|_{H} C_{k_n}^{-1} \displaystyle \sum_{j=1}^{k_n} \left\|  \phi_{n,j}^{'}\right\|_{H}  \left\| \phi_{n,j} - \phi_{n,j}^{'} \right\|_{H} +  \sup_{j \geq 1} \left| D_{n,j} - D_j \right| \left\| x \right\|_{H} C_{k_n}^{-1} .
\end{eqnarray}

Since, for $n$ sufficiently large, from Theorem \ref{theorem2} under \textbf{Assumption A3}, $C_n$ admits a diagonal decomposition in terms of $\left\lbrace C_{n,j}, \ j \geq 1 \right\rbrace$,
\begin{eqnarray}
a_{k_n,2}(x) &\leq& 4 \left\| C_n  \right\|_{\mathcal{L}(H)}  \left\| x \right\|_{H} C_{k_n}^{-1}\displaystyle \sum_{j=1}^{k_n}  \left\| \phi_{n,j} - \phi_{n,j}^{'} \right\|_{H} + \sup_{j \geq 1} \left| D_{n,j} - D_j \right| \left\| x \right\|_{H}C_{k_n}^{-1} ~a.s.
\end{eqnarray}

From results in Bosq (2000, Lemma 4.3),
\begin{eqnarray}
a_{k_n,2}(x) & \leq &  4 \left\| C_n  \right\|_{\mathcal{L}(H)}  \left\| x \right\|_{H} \left\| C_n - C \right\|_{\mathcal{L}(H)} C_{k_n}^{-1} \displaystyle \sum_{j=1}^{k_n}  a_j  +  \sup_{j \geq 1} \left| D_{n,j} - D_j \right|  \left\| x \right\|_{H} C_{k_n}^{-1}  ~a.s.
\end{eqnarray}

On the one hand, from equation (\ref{1444}), there exists  a positive real number $M < \infty$ and an integer $n_0$ large enough such that, for certain  $\beta > \frac{1}{2}$ and $ n \geq n_0$,
\begin{eqnarray}
a_{k_n,2}(x) &\leq& 4 M \left\| C_n  \right\|_{\mathcal{L}(H)}  \left\| x \right\|_{H} \left\| C_n - C \right\|_{\mathcal{L}(H)} \frac{n^{1/4}}{\left(\ln(n) \right)^{\beta}} +  \sup_{j \geq 1} \left| D_{n,j} - D_j \right| \left\| x \right\|_{H} C_{k_n}^{-1} ~a.s.
\end{eqnarray}

On the other hand, applying Remark \ref{remark__4} below, with $n_1$ large enough, for certain $\beta > \frac{1}{2}$ and $ n \geq n_1$,
\begin{equation}
C_{k_n}^{-1}  < C_{k_n}^{-1} \displaystyle \sum_{j=1}^{k_n}  a_j,
\end{equation}
\noindent leading to
\begin{eqnarray}
a_{k_n,2}(x) & < &  M \left\| x \right\|_{H} \left(4 \left\| C_n  \right\|_{\mathcal{L}(H)}  \left\| C_n - C \right\|_{\mathcal{L}(H)}  + \sup_{j \geq 1} \left| D_{n,j} - D_j \right|  \right) \frac{n^{1/4}}{\left(\ln(n) \right)^{\beta}}~a.s., \nonumber \\ \label{AA}
\end{eqnarray}
for certain  $0<M < \infty.$

Hence, since $\left\| C_n \right\|_{\mathcal{L}(H)}< \infty$ and $\left\| x \right\|_{H} < \infty$, from Theorem \ref{theorem2}
and Corollary \ref{prr1} below, from conditions (\ref{1444}) and under \textbf{Assumptions A1-A3},
\begin{eqnarray}
a_{k_n,2}(x)\to^{a.s.} 0,\quad n\rightarrow \infty.   \label{163}
\end{eqnarray}

Taking supremum in $x\in H,$ with $\|x\|_{H}=1,$ at the left-hand side of equation (\ref{146}), from equations (\ref{149})--(\ref{163}), we obtain  the desired result on the  almost surely convergence to zero of $\left\| \widetilde{\rho}_{k_n}- \rho\right\|_{\mathcal{L}(H)} ,$ as $n\rightarrow \infty$. Strong-consistency of the associated plug-in predictor is obtained in analogous way as done in Proposition \ref{pr2m}.

\vspace{0.3cm}

The upper-bound in (\ref{uppbound}) can be directly obtained from $b_{k_n} (x)$, $a_{k_n, 1}(x)$ and $a_{k_n,2}(x)$, reflected in equations (\ref{146})-(\ref{149}) and  (\ref{eq___36}).

\hfill \hfill $\blacksquare$
\end{proof}

% *****************************

% *****************************

\section*{S.2. Asymptotic properties of the empirical eigenvalues and
eigenvectors.}

This section presents the auxiliary results needed on the formulation of the theoretical results derived in Section S.1. The asymptotic properties of the eigenvalues involved in the spectral decomposition of $\widehat{C}_n,~\widehat{D}_n,~C_n$ and $D_n$  will be obtained in Corollary 1 below. Corollary 2 provides the asymptotic properties of the diagonal coefficients of $D_n$,with respect to the eigenvectors of $C_n$. In the derivation of these results, the following theorem plays a crucial role (see Theorem 4.1, Corollary 4.1 and Theorem 4.8 in Bosq, 2000).

\begin{theorem}
\label{theorem2}
Under \textbf{Assumption A3}, for any $\beta > \frac{1}{2},$ as $n\rightarrow \infty,$
\begin{equation}
\frac{n^{1/4}}{\left(\ln(n) \right)^{\beta}}  \left\| C_n - C \right\|_{\mathcal{S}(H)} \to^{a.s.}0, \quad \frac{n^{1/4}}{\left(\ln(n) \right)^{\beta}} \left\| D_n - D \right\|_{\mathcal{S}(H)} \to^{a.s.}0,
\end{equation}
\noindent and, if $\left\| X_0 \right\|_{H}$ is bounded,
\begin{eqnarray}
\left\| C_n - C \right\|_{\mathcal{S}(H)} = \mathcal{O} \left(\left(\frac{\ln(n) }{n} \right)^{1/2} \right)~a.s., \quad \left\| D_n - D \right\|_{\mathcal{S}(H)} = \mathcal{O} \left(\left(\frac{\ln(n) }{n} \right)^{1/2} \right)~a.s.,
\end{eqnarray}
\noindent where $\left\| \cdot \right\|_{\mathcal{S}(H)}$ denotes the norm of the  Hilbert-Schmidt operators  on $H.$

\end{theorem}

From Theorem \ref{theorem2}, we obtain the following corollary on the asymptotic properties of the
eigenvalues $\left\lbrace \widehat{C}_{n,j},\ j \geq 1 \right\rbrace$ and $\left\lbrace \widehat{D}_{n,j} ,\ j \geq 1 \right\rbrace$
of $\widehat{C}_n$ and $\widehat{D}_n,$ respectively, as well as of the eigenvalues $\left\lbrace C_{n,j}, \ j \geq 1 \right\rbrace$ of the empirical estimator $C_n$ and the diagonal coefficients $\widetilde{D}_{n,j}=D_{n}(\widetilde{\phi}_{n,j})(\widetilde{\phi}_{n,j}),$ $j\geq 1,$ with, for $n$ sufficiently large, 
\begin{equation}
D_{n}(\widetilde{\phi}_{n,j})=\widetilde{D}_{n,j}\widetilde{\phi}_{n,j},\quad j\geq 1.
\label{sdcco}
\end{equation}

\begin{corollary}
\label{corollary2}

Under \textbf{Assumptions A1-A3}, the following identities hold, for any $\beta > \frac{1}{2}$:
\begin{eqnarray}
\frac{n^{1/4}}{\left(\ln(n) \right)^{\beta}} \displaystyle \sup_{j \geq 1} \left| \widehat{C}_{n,j} - C_j \right| & \leq & \frac{n^{1/4}}{\left(\ln(n) \right)^{\beta}} \left\| C_n - C \right\|_{\mathcal{S}(H)} \to^{a.s.} 0, \label{c110} \\
\frac{n^{1/4}}{\left(\ln(n) \right)^{\beta}} \displaystyle \sup_{j \geq 1} \left| \widehat{D}_{n,j} - D_j \right|  & \leq & \frac{n^{1/4}}{\left(\ln(n) \right)^{\beta}} \left\| D_n - D \right\|_{\mathcal{S}(H)} \to^{a.s.} 0, \label{c11a}
\end{eqnarray}

\noindent where, as before,  $\left\lbrace C_j, \ j \geq 1 \right\rbrace$ and $\left\lbrace D_j,  \ j \geq 1 \right\rbrace$ are the systems of eigenvalues of $C$ and $D$, respectively;  $\left\lbrace \widehat{C}_{n,j},\ j \geq 1 \right\rbrace$ and $\left\lbrace \widehat{D}_{n,j} ,\ j \geq 1 \right\rbrace$ are given in (\ref{65aa}).

In addition, for $n$ sufficiently large,
\begin{eqnarray}
\frac{n^{1/4}}{\left(\ln(n) \right)^{\beta}} \displaystyle \sup_{j \geq 1} \left| C_{n,j}- C_j \right|  & \leq & \frac{n^{1/4}}{\left(\ln(n) \right)^{\beta}} \left\| C_n - C \right\|_{\mathcal{S}(H)} \to^{a.s.} 0,\label{c11b}\\
\frac{n^{1/4}}{\left(\ln(n) \right)^{\beta}} \displaystyle \sup_{j \geq 1} \left| D_{n}(\widetilde{\phi}_{n,j})(\widetilde{\phi}_{n,j})- D_j \right|  & \leq & \frac{n^{1/4}}{\left(\ln(n) \right)^{\beta}} \left\| D_n - D \right\|_{\mathcal{S}(H)} \to^{a.s.} 0, \label{c11bhh}
\end{eqnarray}

\noindent where $\left\lbrace C_{n,j} ,\ j \geq 1 \right\rbrace$ are introduced in  (\ref{66aa}), and $\left\lbrace\widetilde{D}_{n,j}
,\ j \geq 1 \right\rbrace$ are given in (\ref{sdcco}).
\end{corollary}

\begin{proof}
Since $\widehat{C}_{n}$, with $\displaystyle \sum_{j=1}^{\infty} \widehat{C}_{n,j} = \frac{1}{n} \displaystyle \sum_{i=0}^{n-1} \displaystyle \sum_{j=1}^{\infty} X_{i,j}^{2} = \frac{1}{n} \displaystyle \sum_{i=0}^{n-1} \left\| X_i \right\|_{H}^{2}$, is in the trace class, then, under \textbf{Assumptions A1-A2},
\begin{eqnarray}
 \frac{n^{1/4}}{\left(\ln(n) \right)^{\beta}}\|\widehat{C}_{n}-C \|_{\mathcal{S}(H)} &=& \frac{n^{1/4}}{\left(\ln(n) \right)^{\beta}}\sqrt{\sum_{j\geq 1}|\widehat{C}_{n,j}-C_{j}|^{2}} = \frac{n^{1/4}}{\left(\ln(n) \right)^{\beta}}\sqrt{\sum_{j\geq 1}|C_{n}(\phi_{j})(\phi_{j})-C_{j}|^{2}}
\nonumber \\
& \leq & \frac{n^{1/4}}{\left(\ln(n) \right)^{\beta}}\sqrt{\sum_{j,l\geq 1}|C_{n}(\phi_{k})(\phi_{l})-\delta_{j,l}C_{j}|^{2}}=\frac{n^{1/4}}{\left(\ln(n) \right)^{\beta}}\left\| C_n - C \right\|_{\mathcal{S}(H)},\label{eqpqc1hh}
\end{eqnarray}

\noindent where $\delta_{j,l}$ denotes the Kronecker delta  function.
 From (\ref{eqpqc1hh}), applying  Theorem \ref{theorem2} under \textbf{Assumption A3},
\begin{eqnarray}
\frac{n^{1/4}}{\left(\ln(n) \right)^{\beta}}\sup_{j \geq 1}
|\widehat{C}_{n,j}-C_{j}| \leq \frac{n^{1/4}}{\left(\ln(n) \right)^{\beta}}\|\widehat{C}_{n}-C \|_{\mathcal{S}(H)}\leq \frac{n^{1/4}}{\left(\ln(n) \right)^{\beta}}\|C_{n}-C \|_{\mathcal{S}(H)}\to^{a.s.} 0,
\end{eqnarray}

\noindent as we wanted to prove. Equation (\ref{c11a}) is obtained in a similar way to equation (\ref{c110}), under \textbf{Assumptions A2-A3}, and keeping in mind that $\widehat{D}_{n}$ is, a.s., in the trace class, with $\left| \displaystyle \sum_{j=1}^{\infty} \widehat{D}_{n,j} \right| \leq 2 \displaystyle \sum_{j=1}^{\infty} \widehat{C}_{n,j}~a.s.$,
 \begin{eqnarray}
 \frac{n^{1/4}}{\left(\ln(n) \right)^{\beta}}\sup_{j\geq 1}
|\widehat{D}_{n,j}-D_{j}| \leq \frac{n^{1/4}}{\left(\ln(n) \right)^{\beta}}\|\widehat{D}_{n}-D \|_{\mathcal{S}(H)}\leq \frac{n^{1/4}}{\left(\ln(n) \right)^{\beta}}\|D_{n}-D \|_{\mathcal{S}(H)}\to^{a.s.} 0.\nonumber
\end{eqnarray}

\vspace{0.5cm}

From Theorem \ref{theorem2} and under \textbf{Assumption A3} for $n$ sufficiently large,  $C_{n}$  is a  Hilbert-Schmidt operator, and in particular, it is a compact operator. Thus, applying Bosq (2000, Lemma 4.2) and Theorem \ref{theorem2},  for $n\geq n_{0},$ with $n_{0}$ sufficiently large, we obtain
\begin{eqnarray}
 \frac{n^{1/4}}{\left(\ln(n) \right)^{\beta}}\sup_{k\geq 1}|C_{n,k}-C_{k}| \leq \frac{n^{1/4}}{\left(\ln(n) \right)^{\beta}}\|C_{n}-C\|_{\mathcal{L}(H)} \leq  \frac{n^{1/4}}{\left(\ln(n) \right)^{\beta}}\|C_{n}-C \|_{\mathcal{S}(H)}\to^{a.s.} 0.
\label{eqdpc1}
\end{eqnarray}

Finally, in a similar way to the derivation of (\ref{c11b}), equation (\ref{c11bhh}) is obtained, under \textbf{Assumptions A2-A3}, from Theorem \ref{theorem2} and applying Bosq (2000, Lemma 4.2).

\hfill \hfill $\blacksquare$
\end{proof}

The following lemma, which contains some assertions from Bosq (2000, Corollary 4.3), provides information on the asymptotic properties of the empirical eigenvectors.
\begin{lemma} 
\label{lemma9}
Assume that $\| X_{0}\|_{H}$ is bounded, and  if $\{k_{n}\}$ is a sequence of integers such that
$\Lambda_{k_{n}}=o\left(\left(\frac{n}{\log n}\right)^{1/2}\right),$ as $n\rightarrow \infty$, with 
\begin{equation}\Lambda_{k_n}=\sup_{1\leq j\leq k_n}(C_{j}-C_{j+1})^{-1},\quad 1 \leq j \leq k_n,\label{uee1}
\end{equation}

\noindent then, under \textbf{Assumption A1},
\begin{equation}\sup_{1\leq j\leq k_{n}}\|\phi_{n,j}^{\prime }-\phi_{n,j}\|_{H}\to^{a.s.} 0,\quad n\rightarrow \infty,\label{leq}\end{equation}
 \noindent where $\{\phi_{n,j}^{\prime },\ j\geq 1\}$ are  introduced in Section S.1.2. above.
 \end{lemma}

Let us now consider the following lemma to obtain  the strong-consistency  of $\left\{D_{n,j}, \ j\geq 1\right\}$ (see Corollary \ref{prr1} below).

\begin{lemma}
\label{lemmanew}
Assume that $\| X_{0}\|_{H}$ is bounded, and  if $\{k_{n}\}$ is a sequence of integers such that  $\Lambda_{k_{n}}=o\left(n^{1/4}(\ln(n))^{\beta -1/2} \right),$ as $n\rightarrow \infty$,  where $\Lambda_{k_n}$ is defined in equation (\ref{uee1}) under \textbf{Assumptions A1} and \textbf{A3}. The following limit then holds, for any $\beta > 1/2$,

\begin{equation}\frac{n^{1/4}}{\left(\ln(n) \right)^{\beta}}\sup_{1\leq j\leq k_{n}}\|\phi_{n,j}^{\prime }-\phi_{n,j}\|_{H}\to^{a.s.} 0,\quad n\rightarrow \infty,\label{leqkk}\end{equation}
 \noindent for any $\beta > 1/2$, where $\{\phi_{n,j}^{\prime },\ j\geq 1\}$ defined above.
\end{lemma}

\begin{proof}
From Bosq (2000, Lemma 4.3), for any $n \geq 2$ and $1 \leq j \leq k_n$,
  \begin{eqnarray}
  \left\| \phi_{n,j}^{\prime }-\phi_{n,j} \right\|_{H}  \leq a_j \left\| C_n - C \right\|_{\mathcal{L} (H)} \leq 2 \sqrt{2} \Lambda_{k_n}  \left\| C_n - C \right\|_{\mathcal{S} (H)},
  \end{eqnarray}
  which implies that
  \begin{eqnarray}
\mathcal{P} \left( \sup_{1\leq j\leq k_{n}}\| \phi_{n,j}^{\prime }-\phi_{n,j} \|_{H} \geq \eta \right) & \leq& \mathcal{P} \left( \left\| C_n - C \right\|_{\mathcal{S} (H)} \geq  \frac{\eta}{2 \sqrt{2}  \Lambda_{k_n}} \right). \label{32AA}
\end{eqnarray}

Thus, since $\left\| X_0 \right\|_{H}$ is bounded, from  Bosq (2000, Theorem 4.2), and under \textbf{Assumption A3}, for any $\eta > 0$, and $\beta  > 1/2,$
\begin{eqnarray}
\mathcal{P} \left( \frac{n^{1/4}}{\left(\ln(n) \right)^{\beta}} \sup_{1\leq j\leq k_{n}}\| \phi_{n,j}^{\prime }-\phi_{n,j} \|_{H} \geq \eta \right)
&
\leq &\mathcal{P} \left(\left\| C_n - C \right\|_{\mathcal{S} (H)} \geq  \frac{\eta}{2 \sqrt{2} \Lambda_{k_n}}  \frac{\left(\ln(n) \right)^{\beta}}{n^{1/4}}  \right) \nonumber \\
& \leq &  4 \exp \left(- \frac{n \dfrac{\eta^2}{8 \Lambda_{k_n}^{2}} \dfrac{\left( \ln(n) \right)^{2\beta}}{n^{1/2}}}{\gamma_1 + \delta_1 \dfrac{\eta}{2 \sqrt{2} \Lambda_{k_n}}  \dfrac{\left(\ln(n) \right)^{\beta}}{n^{1/4}} } \right)  \nonumber \\
&= &  \mathcal{O}\left(n^{-\frac{\eta ^{2}}{\gamma_1+\eta \delta_1\left(\frac{\ln(n)}{n}\right)^{1/2}}}\right),\quad n\rightarrow \infty.
\label{fi}
\end{eqnarray}
\noindent Thus, taking  $\eta^{2}> \gamma_1+\delta_1 \eta,$ sequence (\ref{fi}) is summable, and  applying Borel-Cantelli Lemma we arrive to the desired result.\hfill \hfill $\blacksquare$
\end{proof}

\begin{corollary}
\label{prr1}
Under the conditions of Lemma \ref{lemmanew},
considering now \textbf{Assumptions A1-A3},     for
$\beta > \frac{1}{2},$ and $n$ sufficiently large,
\begin{equation}
\frac{n^{1/4}}{\left(\ln(n) \right)^{\beta}}\sup_{j \geq 1} \left| D_{n,j}- D_j \right| \to^{a.s.} 0,\quad n\rightarrow \infty,
\label{resscc}
\end{equation}
 \noindent where  $\left\lbrace D_{n,j},\ j \geq 1 \right\rbrace$  are defined in equation (\ref{141}).
\end{corollary}

\begin{proof}
From Theorem \ref{theorem2}, under \textbf{Assumption A3}, there exists an $n_{0}$ such that for $n\geq n_{0},$ $D_{n}$  is a  Hilbert-Schmidt operator. Then, for  $n\geq n_{0},$ for every $j\geq 1,$ applying orthonormality of the empirical eigenvectors
$\{\phi_{n,j},\ j\geq 1\},$ under \textbf{Assumptions A1-A2},
\begin{eqnarray}
\frac{n^{1/4}}{\left(\ln(n) \right)^{\beta}}\left| D_{n,j} - D_j \right|&=&\frac{n^{1/4}}{\left(\ln(n) \right)^{\beta}} 
  \left|D_{n}(\phi_{n,j})(\phi_{n,j})-D_{n}(\phi_{n,j})(\phi_{j})
+D_{n}(\phi_{n,j})(\phi_{j})\right.\nonumber\\ 
&&\left.\hspace*{2cm}-D(\phi_{n,j})(\phi_{j})+D(\phi_{n,j})(\phi_{j})
-D(\phi_{j})(\phi_{j})\right|\nonumber \\
& &\leq \frac{n^{1/4}}{\left(\ln(n) \right)^{\beta}}\left[\|D_{n}(\phi_{n,j})\|_{H}\|\phi_{n,j}-\phi_{j}\|_{H}+\|(D_{n}-D)(\phi_{n,j})\|_{H}\|\phi_{j}\|_{H}\right.
\nonumber\\
& &
\left.\hspace*{1cm}+\|D(\phi_{n,j}-\phi_{j})\|_{H}\|\phi_{j}\|_{H}\right]\nonumber\\
& &\leq \frac{n^{1/4}}{\left(\ln(n) \right)^{\beta}}\left[\|D_{n}\|_{\mathcal{L}(H)}\|\phi_{n,j}-\phi_{j}\|_{H}+\|D_{n}-D\|_{\mathcal{L}(H)}\right.\nonumber\\
& &
\left.\hspace*{1cm}+\|D\|_{\mathcal{L}(H)}\|\phi_{n,j}-\phi_{j}\|_{H}\right].
\label{dd2bb}
\end{eqnarray}

From Theorem \ref{theorem2}, under \textbf{Assumption A3},
\begin{equation}
\frac{n^{1/4}}{\left(\ln(n) \right)^{\beta}}\|D_{n}-D\|_{\mathcal{L}(H)}\leq \frac{n^{1/4}}{\left(\ln(n) \right)^{\beta}}
\|D_{n}-D\|_{\mathcal{S}(H)}\to^{a.s.} 0,
\label{sc}
\end{equation}

\noindent and,  for $n$ sufficiently large, $\|D_{n}\|_{\mathcal{L}(H)} <\infty.$ Furthermore,
from Lemma \ref{lemmanew} (see equation (\ref{leqkk})), 
\begin{equation}
\frac{n^{1/4}}{\left(\ln(n) \right)^{\beta}}\sup_{1\leq j\leq k_{n}}\|\phi_{n,j}-\phi_{j}\|_{H}\to^{a.s.} 0.
\label{eiv}
\end{equation}
\noindent Hence, from equations (\ref{sc})-(\ref{eiv}), taking the supremum in $j$ at the left-hand side of equation (\ref{dd2bb}),
we obtain equation (\ref{resscc}).
\hfill \hfill $\blacksquare$

\end{proof}

\begin{remark}
\label{remark__4}
Let us now consider the sequence  $\left\lbrace a_j, \ j \geq 1 \right\rbrace$ given by
\begin{equation}
a_1 = 2 \sqrt{2} \frac{1}{C_1 - C_2}, \quad a_j = 2 \sqrt{2} \displaystyle \max \left(\frac{1}{C_{j-1} - C_{j}},\frac{1}{C_j - C_{j+1}} \right),~j \geq 2,
\label{saj}
\end{equation}

If $C_j > C_{j+1}$, when $1 \leq j \leq k_n$, hence $a_j > 0$ for any $1 \leq j \leq k_n$, and then $a_{k_n} < \displaystyle \sum_{j=1}^{k_n} a_j$, for a truncation parameter $\displaystyle \lim_{n \to \infty} k_n = \infty$, with $k_n < n$. Moreover, there exists an integer $j_0$ large enough such that, for any $j \geq j_0$, $a_j > 1$. In particular, if $k_n$ is large enough,
\begin{equation}
\frac{1}{C_{k_n}} < \frac{1}{C_{k_n} - C_{k_n + 1}} < a_{k_n} < \displaystyle \sum_{j=1}^{k_n} a_j, \quad \displaystyle \sum_{j=1}^{k_n} a_j > 1.
\end{equation}
\end{remark}

% *********************************

% *********************************

\section*{S.3. One-sided upper a.s. asymptotic estimate of the $\mathcal{S}(H)$ norm of the error associated with $\boldsymbol{\widetilde{\rho}_{k_{n}}}$}

In this section, $\rho$ does not admit a diagonal spectral decomposition in terms of the eigenvectors of $C$, being $\rho$ not positive, nor trace operator, but it is a Hilbert-Schmidt operator. In this more general framework, an asymptotically almost surely one-sided upper estimate of the $\mathcal{S}(H)$ norm of the error associated with $\widehat{\rho}_{k_{n}}$ is derived. See Ruiz-Medina \& \'Alvarez-Li\'ebana (2017b), where sufficient conditions for the strong-consistency, in the trace norm, of the autocorrelation operator of an ARH(1) process, when it is a positive trace operator which does not admit a diagonal spectral decomposition, are provided.

\begin{proposition}
\textit{Let us assume that $\rho$ is a Hilbert-Schmidt, but not positive nor trace operator. Under \textbf{Assumption A5}, and conditions imposed in Lemma \ref{lemmanew},
\begin{equation}
\left\| \widetilde{\rho}_{k_n} - \rho \right\|_{\mathcal{S}(H)}^{2} \leq \left\| \rho \right\|_{\mathcal{S}(H)}^{2} - \displaystyle \sum_{j=1}^{\infty} \left( \rho \left( \phi_j \right) \left( \phi_j \right) \right)^{2} = \displaystyle \sum_{j \neq k}^{\infty} \left( \frac{D \left( \phi_j \right)\left( \phi_k \right) }{C_j} \right)^{2} < \infty.
\end{equation}
In particular, for $n$ sufficiently large,
\begin{equation}
\left\| \widetilde{\rho}_{k_{n}}-\rho \right\|_{\mathcal{S}(H)}^{2}\leq 
\sum_{j\neq k}^{\infty}\left[\frac{D_{n}(\phi_{n,j})(\phi_{n,k})}{C_{n,j}}\right]^{2}\quad \mbox{a.s}.
\end{equation}}
\end{proposition}

\begin{proof}

Let us consider the eigenvectors $\left\lbrace \phi_{n,j}, \ j \geq 1 \right\rbrace$ of $C_n$. Applying Parseval's identity, we obtain
\begin{eqnarray}
\left\| \widetilde{\rho}_{k_{n}}-\rho \right\|_{\mathcal{S}(H)}^{2}&=&\|(\widetilde{\rho}_{k_{n}}-\rho)^{*}(\widetilde{\rho}_{k_{n}}-\rho) \|_{1} = \sum_{j=1}^{\infty}
\left\langle (\widetilde{\rho}_{k_{n}}-\rho)(\phi_{n,j}),(\widetilde{\rho}_{k_{n}}-\rho)(\phi_{n,j})\right\rangle_{H}\nonumber\\
&=& \sum_{j=1}^{\infty}\| (\widetilde{\rho}_{k_{n}}-\rho)(\phi_{n,j})\|_{H}^{2} = \sum_{j=1}^{\infty}\sum_{k=1}^{\infty}
\left[\left\langle \widetilde{\rho}_{k_{n}}(\phi_{n,j}),\phi_{n,k}\right\rangle_{H}-
\left\langle \rho(\phi_{n,j}), \phi_{n,k}\right\rangle_{H}\right]^{2}\nonumber\\
&=&\sum_{j=1}^{\infty}\sum_{k=1}^{\infty}\left[\left\langle \widetilde{\rho}_{k_{n}}(\phi_{n,j}),\phi_{n,k}\right\rangle_{H}\right]^{2}+\sum_{j=1}^{\infty}\sum_{k=1}^{\infty}\left[\left\langle \rho(\phi_{n,j}), \phi_{n,k}\right\rangle_{H}\right]^{2}\nonumber\\
&-&\sum_{j=1}^{\infty}\sum_{k=1}^{\infty}2\left\langle \widetilde{\rho}_{k_{n}}(\phi_{n,j}),\phi_{n,k}\right\rangle_{H}\left\langle \rho(\phi_{n,j}), \phi_{n,k}\right\rangle_{H} \leq \sum_{j=1}^{\infty}\sum_{k=1}^{k_{n}}\delta_{j,k}[D_{n,j}C_{n,j}^{-1}]^{2} \nonumber \\
&+&\sum_{j=1}^{\infty}\sum_{k=1}^{\infty}
[\left\langle D C^{-1}(\phi_{n,j}),\phi_{n,k}\right\rangle_{H}]^{2} - 2\sum_{j=1}^{\infty}\sum_{k=1}^{k_{n}}\delta_{j,k}D_{n,j}C_{n,j}^{-1}\left\langle C^{-1}(\phi_{n,j}),  D^{*}(\phi_{n,k})\right\rangle_{H}\nonumber\\
&=&\sum_{j=1}^{\infty}[D_{n,j}C_{n,j}^{-1}]^{2}-2D_{n,j}C_{n,j}^{-1}\left\langle 
DC^{-1}(\phi_{n,j}), \phi_{n,j}\right\rangle_{H} + \sum_{j=1}^{\infty}\left[\left\langle DC^{-1}(\phi_{n,j}), \phi_{n,j}\right\rangle_{H}\right]^{2}\nonumber\\
&+&
\sum_{j\neq k}^{\infty}
\left\langle \left[DC^{-1}(\phi_{n,j}),\phi_{n,k}\right\rangle_{H}\right]^{2}  = \sum_{j=1}^{\infty}[D_{n,j}C_{n,j}^{-1}-DC^{-1}(\phi_{n,j})(\phi_{n,j})]^{2}\nonumber\\
&+&
\sum_{j\neq k}^{\infty}
\left[\left\langle DC^{-1}(\phi_{n,j}),\phi_{n,k}\right\rangle_{H}\right]^{2},
\label{eq1s4}
\end{eqnarray}

\noindent where $\delta_{j,k}$ denotes the Kronecker delta function, and $\left\| \cdot \right\|_1$ represents the trace operator norm. From Theorem \ref{theorem2}, under \textbf{Assumption A3},
\begin{eqnarray}
\|D_{n}C_{n}^{-1}-DC^{-1}\|_{\mathcal{S}(H)} &=& 
\| D_{n}C_{n}^{-1}-DC_{n}^{-1}+D C_{n}^{-1}-DC^{-1}\|_{\mathcal{S}(H)} \leq \| D_{n}C_{n}^{-1}-DC_{n}^{-1}\|_{\mathcal{S}(H)} \nonumber\\
&+&
\|DC_{n}^{-1}-DC^{-1}\|_{\mathcal{S}(H)} =  \|(D_{n}-D)C_{n}^{-1}\|_{\mathcal{S}(H)}+
\|D(C_{n}^{-1}-C^{-1})\|_{\mathcal{S}(H)}, \nonumber \\
\label{eqproof1}
\end{eqnarray}

\noindent leading to
\begin{equation}
\lim_{n\rightarrow \infty}\sum_{j=1}^{\infty}[D_{n,j}C_{n,j}^{-1}-DC^{-1}(\phi_{n,j})(\phi_{n,j})]^{2}=0\quad a.s.
\label{eqes1}
\end{equation}

From equations (\ref{eq1s4}) and (\ref{eqes1}), and from Lemma \ref{lemmanew},
\begin{eqnarray}
\lim_{n\rightarrow \infty}\|\widetilde{\rho}_{k_{n}}-\rho\|_{\mathcal{S}(H)}^{2}&=& \lim_{n\rightarrow \infty}\|(\widetilde{\rho}_{k_{n}}-\rho)^{*}(\widetilde{\rho}_{k_{n}}-\rho) \|_{1}=\lim_{n\rightarrow \infty}\sum_{j\neq k}^{\infty}
\left[\left\langle DC^{-1}(\phi_{n,j}),\phi_{n,k}\right\rangle_{H}\right]^{2}\nonumber\\
&=&\lim_{n\rightarrow \infty}\sum_{j\neq k}^{\infty} \left[ DC^{-1}(\phi_{n,j})(\phi_{n,k})-DC^{-1}(\phi_{n,j})(\phi_{k}) + 
DC^{-1}(\phi_{n,j})(\phi_{k})\right.
\nonumber\\
&&\left.-
\sum_{j\neq k}^{\infty}DC^{-1}(\phi_{j})(\phi_{k})+DC^{-1}(\phi_{j})(\phi_{k})\right]^{2} \nonumber \\
&\leq & \lim_{n\rightarrow \infty}
\sum_{j\neq k}^{\infty}\left[\|DC^{-1}(\phi_{n,j})\|_{H}\|\phi_{n,k}-\phi_{k}\|_{H}\right.\nonumber\\
&&\left.+
\|DC^{-1}\|_{\mathcal{L}(H)}\|\phi_{n,j}-\phi_{j}\|_{H}+DC^{-1}(\phi_{j})(\phi_{k})\right]^{2} = \sum_{j\neq k}^{\infty}[DC^{-1}(\phi_{j})(\phi_{k})]^{2}\nonumber \\
& \leq & 
\|\rho\|_{\mathcal{S}(H)}^{2}\quad \mbox{a.s}
\end{eqnarray}

Therefore, when $\rho$ is not positive, nor trace operator,  but it is Hilbert-Schmidt operator, the norm of the error associated with $\widetilde{\rho}_{k_{n}},$ in the space of Hilbert-Schmidt operators, is  a.s. asymptotically upper bounded by the following quantity: \begin{equation}\|\rho\|_{\mathcal{S}(H)}^{2}-\sum_{j=1}^{\infty}[\rho(\phi_{j})(\phi_{j})]^{2}=\sum_{j\neq k}^{\infty}\left[\frac{D(\phi_{j})(\phi_{k})}{C_{j}}\right]^{2}<\infty.\label{eqbsss}
\end{equation}

Equation (\ref{eqbsss}) can be approximated by the empirical quantity:
$$\sum_{j\neq k}^{\infty}\left[\frac{D_{n}(\phi_{n,j})(\phi_{n,k})}{C_{n,j}}\right]^{2}.$$

Thus, for $n$ sufficiently large,
\begin{equation}
\left\| \widetilde{\rho}_{k_{n}}-\rho \right\|_{\mathcal{S}(H)}^{2}\leq 
\sum_{j\neq k}^{\infty}\left[\frac{D_{n}(\phi_{n,j})(\phi_{n,k})}{C_{n,j}}\right]^{2}\quad \mbox{a.s}. \label{eq__35}
\end{equation}

\hfill \hfill $\blacksquare$

\end{proof}

% *********************************

% *********************************

\section*{S.4. Simulation study: large-sample behavior of $\widetilde{\rho}_{k_n}$ when eigenvectors of $C$ are unknown}

A brief simulation study is undertaken to illustrate the theoretical results on the strong-consistency of the formulated diagonal componentwise estimator of $\rho$ in Section 8.3 of the main paper, when $\left\lbrace \phi_j, \ j\geq 1 \right\rbrace$ are unknown and a Gaussian diagonal data generation  is achieved. An almost sure rate of convergence is fitted as well.

The zero-mean Gaussian ARH(1) process $X$ generated here has covariance operator $C$ given by
\begin{equation}
 C(f) (g) = \sum_{j=1}^{\infty} C_j \langle \phi_j, f \rangle_H \langle \phi_j, g \rangle_H, \quad \phi_{j} \left( x\right)= \sqrt{\frac{2}{b-a}} \sin \left(
\frac{\pi j x}{b-a} \right), \quad f,g\in H = L^{2}((a,b)),~x\in (a,b).\label{aco}
\end{equation}

In the remaining, we fix $(a,b)=(0,4)$ and, under \textbf{Assumption A1}, $C_j = c_1 j^{-\delta_1}$, being $c_1$ a positive constant, for any $j \geq 1$. Different rates of convergence to zero of the eigenvalues of $C$ are studied, corresponding to the values of the shape parameter $\delta_{1}\in (1,2)$. In a diagonal context, autocorrelation operator and covariance operator of the error term are approximated as follows, for $M=50$:
\begin{equation}
\rho \left(X \right) \left( t \right) \simeq \displaystyle \sum_{j=1}^{M} \rho_{j,j} \langle \phi_j, X \rangle_H \phi_j (t), \quad C_\varepsilon \left(X \right) \left( t \right) \simeq \displaystyle \sum_{j=1}^{M} \sigma_{j,j}^{2} \langle \phi_j, X \rangle_H \phi_j(t),
\end{equation}

\noindent where $\rho_{j,j} = c_2j ^{-\delta_2}$ and $\sigma_{j,j}^{2} = C_j \left( 1 - \rho_{j,j} \right)$, for any $j \geq 1$, being $\delta_2 \in \left(1/2, 2 \right)$, and $c_2$ a constant which belongs to $\left[0,1\right]$. Thus, $\rho $ is a diagonal self-adjoint Hilbert-Schmidt operator, with $\left\| \rho \right\|_{\mathcal{L}(H)} = \displaystyle \sup_{j \geq 1} \left| \rho_j \right| < 1$, under \textbf{Assumption A2}. Simulations are then performed under \textbf{Assumptions A1--A3} and \textbf{A5}, and the empirical version of the upper-bound derived in (\ref{uppbound}) is considered:
\begin{equation}
UB \left(k_n, l \right) = \displaystyle \sup_{1 \leq j \leq k_n} \left| \widetilde{\rho}_{n,j}^{l} - \frac{D_{n,j}^{l}}{C_j} \right| + \displaystyle \sup_{1 \leq j \leq k_n} \left|  \frac{D_{n,j}^{l}}{C_j} - \rho_j \right| + 2\displaystyle \sum_{j=1}^{k_n} \frac{\left| D_{n,j}^{l} \right|}{C_j} \left\| \phi_{n,j}^{l} - \phi_{n,j}^{\prime,l} \right\|_{H} +\displaystyle \sup_{j > k_n} \left|  \rho_j \right|,\label{error3}
\end{equation}
\noindent for $k_n = \lceil \ln(n) \rceil$, for which conditions in Proposition \ref{proposition5} are held (see Example 8.6 in Bosq, 2000). In equation (\ref{error3}), superscript $l$ denotes the estimator computed based on the $l$th generation of the values $\widetilde{X}_{i,j}^{l} = \langle X_{i}^{l}, \phi_{n,j}^{l} \rangle$, for $l=1,\ldots,N$, $j=1,\ldots,k_n$ and $i=0,\ldots,n-1.$ Here, $N=500$ realizations have been generated, with shape parameters $\delta_{1} = 61/60,~3/2,~9/5$ and $\delta_2 = 11/10$. Discretization step $\Delta t = 0.06$ has been adopted. For sample sizes $n_t = 35000 + 40000 \left(t -1 \right),~t=1,\ldots,10$,
 \begin{eqnarray}
E \left(k_n, n_t, \beta \right)= \left(\displaystyle \sum_{l=1}^{N} \mathbf{1}_{\left(\xi_{n_t,\beta},\infty \right)} \left(UB \left(k_n, l \right)  \right) \right)/N, \quad \xi_{n_t,\beta} = \frac{\left(\ln(n_t)\right)^{\beta}}{n_{t}^{1/3}}, \label{count3} 
\end{eqnarray}

\noindent values are reflected in Table \ref{tab:Table9}, in which the curve $\xi_{n_t,95/100}$ is fitted as the almost sure rate of convergence. In equation (\ref{count3}), $\mathbf{1}_{\left(\xi_{n_t,\beta},\infty \right)}$ denotes the indicator function over the interval $\left(\xi_{n_t,\beta},\infty \right)$.

\setlength{\heavyrulewidth}{0.3em}
\setlength{\lightrulewidth}{0.15em}

\begin{table}[H]
 \caption{{\small $E \left(k_n, n_t, \beta \right)$ values defined in  (\ref{count3}), for $\beta =  95/100$  and $N = 500$ realizations,  with $\delta_2 = 11/10,$ and $\delta_{1} =61/60,~3/2,~9/5,$  considering  $n_{t}= 35000+40000(t-1),~t=1,\ldots, 10,$ and $k_{n} = \lceil \ln(n) \rceil $.}}
\centering
\vspace{-0.17cm}
\begin{tabular}{cc@{\hspace{1.8em}}c@{\hspace{1.8em}}c@{\hspace{1.8em}}c}
\toprule
  $n_t$ &  $k_n$ & $\delta_1 = 61/60$ &  $\delta_1 = 3/2$ & $\delta_1 = 9/5$ \\
 \midrule
35000 & 10  & $\frac{33}{500}$& $\frac{28}{500}$ & $\frac{20}{500}$  \\
75000 & 11  &$\frac{19}{500}$ & $\frac{15}{500}$ & $\frac{11}{500}$ \\
 115000 &  11  &$\frac{10}{500}$ & $\frac{8}{500}$ & $\frac{5}{500}$  \\
155000 & 11  & $\frac{4}{500}$ & $\frac{3}{500}$ & $\frac{2}{500}$  \\
195000 & 12  &$\frac{5}{500}$ & $\frac{4}{500}$ & $\frac{2}{500}$  \\
235000  & 12  & $\frac{3}{500}$ & $\frac{1}{500}$ & $\frac{1}{500}$  \\
 275000 & 12  & $\frac{3}{500}$ & 0 & 0  \\
  315000 & 12 &  0& $\frac{1}{500}$ & 0  \\
  355000 & 12 &$\frac{1}{500}$ & 0 & 0 \\
  395000 & 12 & 0 & 0 & 0 \\
\bottomrule
\end{tabular}
  \label{tab:Table9}
\end{table}

The convergence to zero of the empirical mean of $\left\| \widehat{\rho}_{k_n} - \rho \right\|_{\mathcal{L}(H)}$ is numerically illustrated in Figure \ref{fig:F_0} below, which displays the empirical mean of values $UB \left(k_n, l \right)$, against the curve $\xi_{n_t,95/100}$, for each $l=1,\ldots,N$ realizations, with $N=500$ and $k_n = \lceil \ln(n) \rceil$.

\begin{figure}[H]
\centering
 \includegraphics[width=15cm,height=6.5cm]{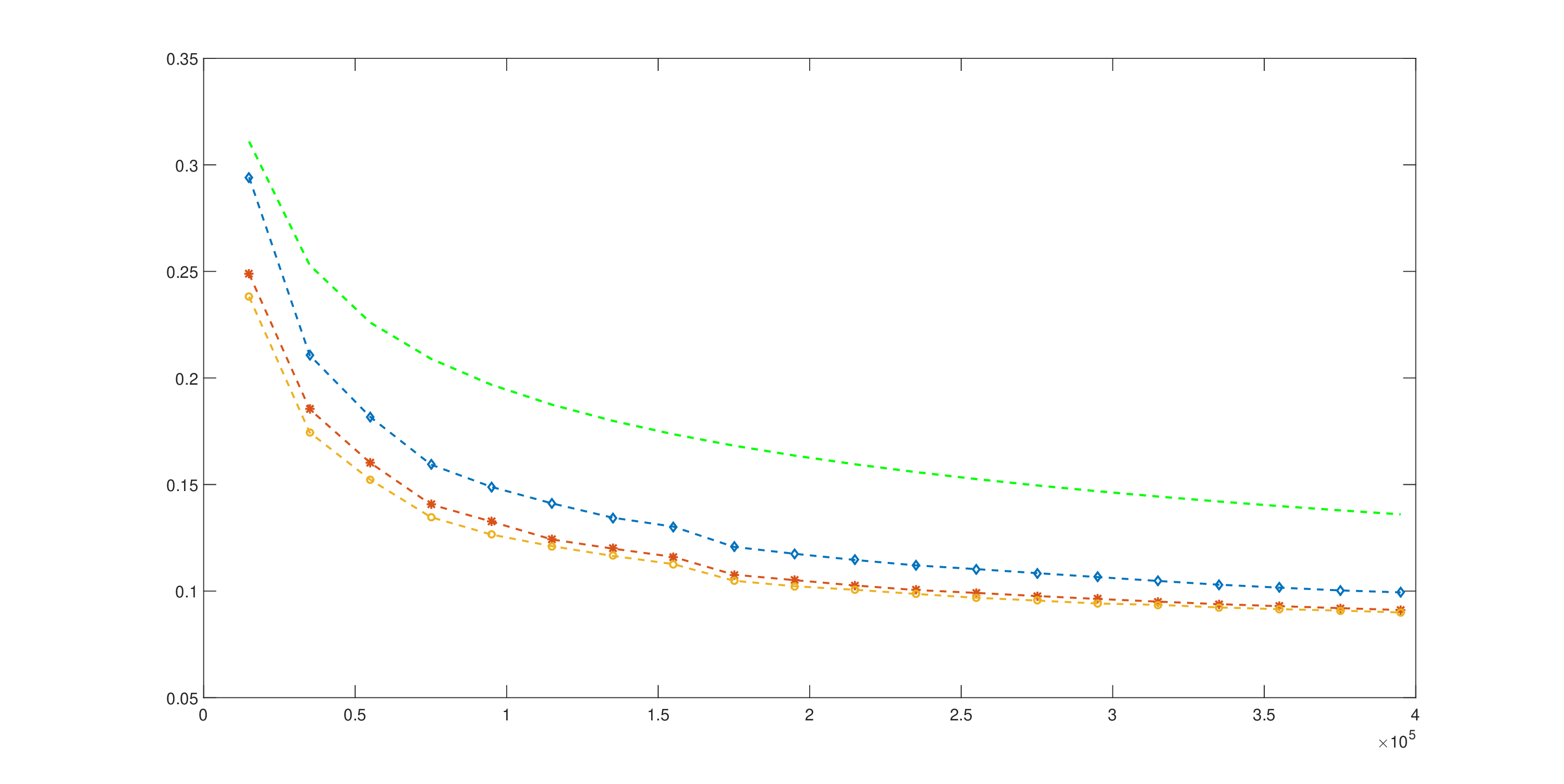}
\vspace{-0.6cm}
\caption{{\small Empirical mean of $UB \left(k_n, l \right)$, for $l=1,\ldots,N$, with $N = 500$, $\delta_2 = 11/10$, $n_{t}= 15000+20000(t-1),~t=1,\ldots, 20,$ and $k_{n} = \lceil \ln(n) \rceil $. Shape parameters $\delta_{1} =61/60,~3/2,~9/5$  are considered (blue diamond, red star and yellow circle dotted lines, respectively). The curve $\xi_{n_t,\beta} = \frac{\left(\ln(n_t) \right)^{\beta} }{n_{t}^{1/3} }$, with $\beta = 95/100$, is also drawn (green dotted line).}} \label{fig:F_0}
\end{figure}

A theoretical almost sure rate of convergence for the diagonal componentwise estimator $\widetilde{\rho}_{k_n}$ has not been derived in Proposition \ref{proposition5}. However, when a diagonal data generation is performed, under different rates of convergence to zero of the eigenvalues of $C$, the curve $\xi_{n_t,95/100} = \frac{\left( \ln(n) \right)^{95/100}}{n^{1/3}}$ is numerically fitted, when large samples sizes are considered. As expected, for the largest shape parameter value $\delta_1$, corresponding to the fastest decay velocity  of the eigenvalues of the autocovariance operator, we obtain the fastest convergence to zero of $\left\|  \widetilde{\rho}_{k_n} - \rho \right\|_{\mathcal{L}(H)}$. From results displayed in Figure \ref{fig:F_0}, the empirical mean of the upper bound in (\ref{error3}), computed from $N = 500$ realizations, is showed that can be upper bounded by the curve $\xi_{n_t,95/100} = \frac{\left( \ln(n) \right)^{95/100}}{n^{1/3}}$, for the parameters adopted.

\section*{S.5. Comparative study: numerical results}

Tables \ref{tab:TableA4}-\ref{tab:TableA9a} and \ref{tab:TableA5}-\ref{tab:TableA10} of this section reflect, in more detail, the numerical results obtained in the comparative study performed in Section 9, where a summary of such results is displayed in Figures 2-7. Details about the comparative study, and about which conditions are held under any scenario for each approach considered, can be found in the referred Section 9. As remarked at the beginning of Section 9, the ARH(1) diagonal componentwise plug-in predictor established in Section 8.3 will be only considered under diagonal scenarios. Strong-consistency results for the estimator $\widetilde{\rho}_{k_n}$, in the trace norm, when $\rho$ is a positive and trace operator, which does not admit a diagonalization in terms of the eigenvectors of $C$, have been recently provided in Ruiz-Medina \& \'Alvarez-Li\'ebana (2017b).

The large-sample behaviour of the empirical-eigenvector-based componentwise plug-in predictor formulated in Section 8, as well as those ones in Bosq (2000) and Guillas (2001), will be firstly displayed in Tables \ref{tab:TableA4}-\ref{tab:TableA9a}. For sample sizes $n_t = 35000 + 40000 \left(t -1 \right),~t=1,\ldots,10$, the following values are computed
\begin{equation}
F \left(k_n, n_t, \beta \right) = \left( \displaystyle \sum_{l=1}^{N} \mathbf{1}_{\left(\xi_{n_t,\beta},\infty \right)} \left( \left\| \left(\rho - \overline{\rho}_{k_n}^{l} \right) \left(X_{n-1}^{l} \right) \right\|_{H}^{k_n} \right) \right) /N, \label{a955}
\end{equation}

\noindent being $\xi_{n_t,\beta}$  the curve which numerically fits the almost sure rate of convergence of $\left\| \left(\rho - \overline{\rho}_{k_n}^{l} \right) \left(X_{n-1}^{l} \right) \right\|_{H}^{k_n}$. When the data is generated in the diagonal framework

\begin{equation}
\left\| \left(\rho - \overline{\rho}_{k_n}^{l} \right) \left(X_{n-1}^{l} \right) \right\|_{H}^{k_n} = \sqrt{\displaystyle \int_{a}^{b} \left(\displaystyle \sum_{j=1}^{k_n} \rho_j X_{n-1,j}^{l} \phi_j (t)  - \displaystyle \sum_{j=1}^{k_n} \overline{\rho}_{n,j}^{l} \left( X_{n-1}^{l} \right) \phi_{n,j}^{l}(t) \right)^2 dt}, \label{a755}
\end{equation}
\noindent is computed, being $\overline{\rho}_{k_n}^{l} \left( X_{n-1}^{l} \right)$ the corresponding predictors, for any $j=1,\ldots,k_n$, and based on the $l$th generation of the values  $\widetilde{X}_{i,j}^{l} = \langle X_{i}^{l}, \phi_{n,j}^{l} \rangle_H$, for $l=1, \ldots, N$, with $N=500$ realizations. The following scenarios will be considered, when the diagonal data generation is assumed:

\setlength{\heavyrulewidth}{0.3em}
\setlength{\lightrulewidth}{0.15em}

\begin{table}[H]
 \caption{{\small Diagonal scenarios considered (see Table \ref{tab:TableA4} below), with $\delta_2 = 11/10$,  $n_t = 35000 + 40000(t-1),~t=1,\ldots,10$, and $\xi_{n_t,\beta} = \frac{\left(\ln(n_t) \right)^{\beta}}{n_{t}^{1/2}}$, for $\beta = 65/100$.}}
\centering
\vspace{-0.17cm}
\begin{tabular}{ccc}
\toprule
  Scenario & $\delta_1$ &  $k_n$ \\
 \midrule
1 & $3/2$  & $ \lceil \ln(n) \rceil$  \\
2 & $24/10$  & $\lceil \ln(n) \rceil$ \\
3 &  $3/2$  & $ \lceil e^{\prime} n^{1/\left(8 \delta_1 + 2 \right)} \rceil,~e^{\prime} = 17/10$  \\
4& $24/10$  & $\lceil e^{\prime} n^{1/\left(8 \delta_1 + 2 \right)} \rceil,~e^{\prime} = 17/10$  \\
\bottomrule
\end{tabular}
  \label{tab:Table_Scen_A}
\end{table}

In the case of pseudo-diagonal and non-diagonal frameworks, the following truncated norm is then computed:
\begin{equation}
\left\| \left(\rho - \overline{\rho}_{k_n}^{l} \right) \left(X_{n-1}^{l} \right) \right\|_{H}^{k_n} = \sqrt{\displaystyle \int_{a}^{b} \left( \displaystyle \int_{a}^{b} \left(\displaystyle \sum_{j,k=1}^{k_n} \rho_{j,k} \phi_j (t) \phi_k (s) \right)ds  - \displaystyle \sum_{j=1}^{k_n} \overline{\rho}_{n,j}^{l} \left( X_{n-1}^{l} \right) \phi_{n,j}^{l}(t) \right)^2 dt}. \label{a766}
\end{equation}

Pseudo-diagonal  and non-diagonal scenarios are outlined as follows:

\begin{table}[H]
 \caption{{\small Pseudo-diagonal and non-diagonal scenarios considered (see Tables \ref{tab:TableA9}-\ref{tab:TableA9a} below), with $\delta_2 = 11/10$,  $n_t = 35000 + 40000(t-1),~t=1,\ldots,10$, and $\xi_{n_t,\beta} = \frac{\left(\ln(n_t) \right)^{\beta}}{n_{t}^{1/3}}$.}}
\centering
\vspace{-0.17cm}
\begin{tabular}{cccc||cccc}
\toprule
\multicolumn{4}{c||}{Pseudo-diagonal scenarios} & \multicolumn{4}{c}{Non-diagonal scenarios} \\ 
\hline
  Scenario & $\delta_1$ &  $k_n$ & $\beta$ & Scenario & $\delta_1$ &  $k_n$ & $\beta$ \\
 \midrule
5 & $3/2$  & $ \lceil \ln(n) \rceil$ & $3/10$ & 9 & $3/2$  & $ \lceil \ln(n) \rceil$ & $125/100$\\
6 & $24/10$  & $\lceil \ln(n) \rceil$ & $3/10$ & 10 &  $24/10$  & $ \lceil \ln(n) \rceil$ & $125/100$\\
7&  $3/2$  & $ \lceil (17/10) n^{1/\left(8 \delta_1 + 2 \right)} \rceil$  & $3/10$ & 11 & $3/2$  & $ \lceil(17/10)  n^{1/\left(8 \delta_1 + 2 \right)} \rceil$ & $125/100$\\
8& $24/10$  & $\lceil (17/10) n^{1/\left(8 \delta_1 + 2 \right)} \rceil$  & $3/10$ & 12 &  $24/10$  & $ \lceil (17/10) n^{1/\left(8 \delta_1 + 2 \right)} \rceil$ & $125/100$ \\
\bottomrule
\end{tabular}
  \label{tab:Table_Scen_B}
\end{table}

\begin{small}
\begin{table}[H]
 \caption{{\small
$F \left(k_n, n_t, \beta \right)$ values in
(\ref{a955})-(\ref{a755}), for scenarios 1-4. O.A. denotes the approach detailed in Section 8; B denotes the approach in Bosq (2000); G denotes the approach in Guillas (2001).}}
\vspace{-0.12cm}
\centering
\begin{tabular}{cc@{\hspace{1.5em}}ccc@{\hspace{1.6em}}ccc@{\hspace{1.7em}}c@{\hspace{0.16cm}}ccc@{\hspace{1.7em}}c@{\hspace{0.16cm}}ccc}
\toprule  
& & \multicolumn{3}{c}{\hspace{-0.4cm}\textit{Scenario 1}} & \multicolumn{3}{c}{\hspace{-0.5cm}\textit{Scenario 2}} & \multicolumn{4}{c}{\hspace{-0.3cm}\textit{Scenario 3}}& \multicolumn{4}{c}{\hspace{-0.2cm}\textit{Scenario 4}}\\
  \midrule[0.05cm]
   $n_t$  & $k_n$ & O.A. &  B & G & O.A. &  B & G & $k_n$ & O.A. &  B & G & $k_n$& O.A. &  B & G\\ 
  \midrule[0.05cm] 
35000 & 10 & $\frac{11}{500}$ & $\frac{68}{500}$ & $\frac{70}{500}$ &  $\frac{4}{500}$ & $\frac{24}{500}$ & $\frac{28}{500}$ & 3 & $\frac{13}{500}$ & $\frac{12}{500}$& $\frac{10}{500}$ & 2 & $\frac{7}{500}$& $\frac{7}{500}$ & $\frac{4}{500}$ \vspace{0.09cm}  \\ 
75000 & 11& $\frac{9}{500}$ & $\frac{62}{500}$& $\frac{66}{500}$ & $\frac{3}{500}$ & $\frac{18}{500}$ & $\frac{25}{500}$ & 3&   $\frac{9}{500}$& $\frac{9}{500}$& $\frac{6}{500}$ & 2 &  $\frac{4}{500}$& $\frac{3}{500}$& $\frac{2}{500}$ \vspace{0.09cm}  \\ 
115000 & 11&  $\frac{6}{500}$ & $\frac{59}{500}$ & $\frac{62}{500}$ & $\frac{3}{500}$ & $\frac{16}{500}$ & $\frac{22}{500}$ & 3& $\frac{6}{500}$& $\frac{5}{500}$& $\frac{5}{500}$ & 2&   $\frac{3}{500}$& $\frac{3}{500}$& $\frac{2}{500}$ \vspace{0.09cm}  \\ 
155000 & 11& $\frac{4}{500}$ & $\frac{57}{500}$ & $\frac{60}{500}$ &  $\frac{2}{500}$ & $\frac{12}{500}$ & $\frac{19}{500}$ & 3&  $\frac{5}{500}$&  $\frac{4}{500}$&$\frac{4}{500}$ & 2 &   $\frac{3}{500}$& $\frac{2}{500}$ & $\frac{1}{500}$ \vspace{0.09cm}  \\ 
  195000 & 12& $\frac{6}{500}$ & $\frac{60}{500}$ & $\frac{64}{500}$  & $\frac{4}{500}$ & $\frac{15}{500}$ & $\frac{21}{500}$ & 4& $\frac{6}{500}$&$\frac{4}{500}$ &$\frac{3}{500}$ & 3 & $\frac{4}{500}$& $\frac{2}{500}$ & $\frac{1}{500}$ \vspace{0.09cm}  \\  
   235000  & 12&$\frac{4}{500}$ &  $\frac{58}{500}$  & $\frac{61}{500}$  &0  & $\frac{14}{500}$ &  $\frac{17}{500}$ & 4& $\frac{4}{500}$ & $\frac{3}{500}$ & $\frac{2}{500}$ & 3 & $\frac{2}{500}$& $\frac{1}{500}$& $\frac{1}{500}$ \vspace{0.09cm}  \\   
    275000 & 12& $\frac{3}{500}$ & $\frac{51}{500}$ & $\frac{58}{500}$ & 0 & $\frac{13}{500}$ & $\frac{16}{500}$ & 4& $\frac{3}{500}$& $\frac{2}{500}$ & $\frac{1}{500}$ & 3&  $\frac{2}{500}$& $\frac{1}{500}$& 0 \vspace{0.09cm}  \\ 
  315000 & 12& $\frac{3}{500}$ & $\frac{50}{500}$& $\frac{55}{500}$  & $\frac{1}{500}$ & $\frac{12}{500}$ & $\frac{14}{500}$ & 4& $\frac{2}{500}$& $\frac{1}{500}$& $\frac{1}{500}$ & 3& $\frac{1}{500}$& 0 & 0 \vspace{0.09cm}  \\  
     355000 & 12& $\frac{2}{500}$ & $\frac{47}{500}$ & $\frac{53}{500}$ & 0 & $\frac{12}{500}$ & $\frac{13}{500}$ & 4& $\frac{2}{500}$ &  $\frac{1}{500}$& 0 & 3&  $\frac{1}{500}$ &0 & 0 \vspace{0.09cm}  \\ 
395000 & 12&$\frac{2}{500}$ & $\frac{44}{500}$ & $\frac{51}{500}$ & 0 & $\frac{11}{500}$  & $\frac{13}{500}$ & 4&  $\frac{2}{500}$ & 0& 0 & 3 & $\frac{1}{500}$ & 0& 0 \vspace{0.05cm}  \\ 
  \toprule 
\end{tabular}     
  \label{tab:TableA4}
\end{table}
\end{small}

   \begin{table}[H]
 \caption{{\small
$F \left(k_n, n_t, \beta \right)$ values in
(\ref{a955}) and (\ref{a766}), for scenarios 5-8. B denotes the approach in Bosq (2000); G denotes the approach in Guillas (2001).}}
\vspace{-0.12cm}
\centering
\begin{small}
   \begin{tabular}{c@{\hspace{1.5em}}ccc@{\hspace{1.6em}}ccc@{\hspace{1.7em}}@{\hspace{0.16cm}}ccc@{\hspace{1.7em}}@{\hspace{0.16cm}}ccc}
\toprule  
& \multicolumn{3}{c}{\hspace{-0.4cm}\textit{Scenario 5}} & \multicolumn{3}{c}{\hspace{-0.5cm}\textit{Scenario 6}} & \multicolumn{3}{c}{\hspace{-0.3cm}\textit{Scenario 7}}& \multicolumn{3}{c}{\hspace{-0.35cm}\textit{Scenario 8}}\\
  \midrule[0.05cm]
   $n_t$  & $k_n$ & B & G  & $k_n$ &  B & G & $k_n$ &  B & G & $k_n$&   B & G \\
  \midrule[0.05cm]  
35000 & 10 & $\frac{32}{500}$ & $\frac{33}{500}$ & 10 &   $\frac{25}{500}$ & $\frac{29}{500}$ & 3 & $\frac{28}{500}$& $\frac{26}{500}$ & 2  & $\frac{27}{500}$ & $\frac{24}{500}$ \vspace{0.09cm}  \\
 
75000 & 11&  $\frac{29}{500}$& $\frac{31}{500}$  & 11
 &  $\frac{21}{500}$ & $\frac{23}{500}$  & 3& $\frac{26}{500}$& $\frac{24}{500}$ & 2 &  $\frac{22}{500}$& $\frac{19}{500}$ \vspace{0.09cm}  \\ 
 
115000 & 11  & $\frac{26}{500}$ & $\frac{28}{500}$ & 11 &    $\frac{18}{500}$ & $\frac{20}{500}$ & 3 & $\frac{23}{500}$& $\frac{21}{500}$ & 2 & $\frac{18}{500}$ & $\frac{15}{500}$ \vspace{0.09cm}  \\ 

155000 & 11 & $\frac{24}{500}$ & $\frac{26}{500}$ & 11 &   $\frac{14}{500}$ & $\frac{17}{500}$ & 3 & $\frac{19}{500}$& $\frac{17}{500}$ & 2 & $\frac{16}{500}$& $\frac{12}{500}$ \vspace{0.09cm}  \\ 
  195000 & 12  & $\frac{19}{500}$ & $\frac{21}{500}$ & 12  & $\frac{10}{500}$ & $\frac{13}{500}$ & 4 & $\frac{14}{500}$&$\frac{12}{500}$ & 3 & $\frac{11}{500}$& $\frac{9}{500}$ \vspace{0.09cm}  \\  
   235000  & 12 &  $\frac{16}{500}$  & $\frac{16}{500}$  & 12 & $\frac{12}{500}$ &  $\frac{14}{500}$  & 4 & $\frac{15}{500}$& $\frac{10}{500}$ & 3 &  $\frac{13}{500}$& $\frac{10}{500}$ \vspace{0.09cm}  \\   
    275000 & 12  & $\frac{12}{500}$ & $\frac{13}{500}$ & 10  & $\frac{8}{500}$ & $\frac{10}{500}$  & 4& $\frac{9}{500}$& $\frac{7}{500}$ & 3 & $\frac{7}{500}$ & $\frac{6}{500}$ \vspace{0.09cm}  \\ 
  315000 & 12 & $\frac{9}{500}$& $\frac{15}{500}$ & 12 & $\frac{5}{500}$ & $\frac{7}{500}$ & 4& $\frac{5}{500}$& $\frac{4}{500}$ & 3 & $\frac{4}{500}$& $\frac{3}{500}$ \vspace{0.09cm}  \\  
     355000 & 12 & $\frac{8}{500}$ & $\frac{11}{500}$ & 12 & $\frac{3}{500}$ & $\frac{5}{500}$  & 4 & $\frac{3}{500}$ & $\frac{3}{500}$ & 3 & $\frac{2}{500}$& $\frac{2}{500}$ \vspace{0.09cm}  \\ 
395000 & 12  & $\frac{6}{500}$ & $\frac{9}{500}$& 12 &  $\frac{3}{500}$ & $\frac{5}{500}$  & 4  & $\frac{2}{500}$&  $\frac{1}{500}$ & 3 & $\frac{1}{500}$& 0 \vspace{0.05cm}  \\ 
  \toprule 
\end{tabular}  
   \end{small}   
  \label{tab:TableA9}
\end{table}

\begin{table}[H]
 \caption{{\small
$F \left(k_n, n_t, \beta \right)$ values in
(\ref{a955}) and (\ref{a766}), for scenarios 9-12. B denotes the approach in Bosq (2000); G denotes the approach in Guillas (2001).}}
\vspace{-0.16cm}
\centering
\begin{small}
   \begin{tabular}{c@{\hspace{1.5em}}ccc@{\hspace{1.6em}}ccc@{\hspace{1.7em}}@{\hspace{0.16cm}}ccc@{\hspace{1.7em}}@{\hspace{0.16cm}}ccc}
\toprule  
& \multicolumn{3}{c}{\hspace{-0.4cm}\textit{Scenario 9}} & \multicolumn{3}{c}{\hspace{-0.5cm}\textit{Scenario 10}} & \multicolumn{3}{c}{\hspace{-0.3cm}\textit{Scenario 11}}& \multicolumn{3}{c}{\hspace{-0.35cm}\textit{Scenario 12}}\\
  \midrule[0.05cm]
   $n_t$  & $k_n$ & B & G  & $k_n$ &  B & G & $k_n$ &  B & G & $k_n$&   B & G \\
  \midrule[0.05cm] 
35000 & 10 & $\frac{67}{500}$ & $\frac{71}{500}$  & 10 & $\frac{59}{500}$ & $\frac{62}{500}$ & 3 & $\frac{55}{500}$& $\frac{47}{500}$ & 2 &  $\frac{44}{500}$ & $\frac{40}{500}$ \vspace{0.09cm}  \\ 

75000 & 11 & $\frac{44}{500}$& $\frac{50}{500}$ & 11 & $\frac{38}{500}$ & $\frac{45}{500}$  & 3& $\frac{36}{500}$& $\frac{31}{500}$ & 2 &   $\frac{34}{500}$& $\frac{30}{500}$ \vspace{0.09cm}  \\ 

115000 & 11 & $\frac{47}{500}$ & $\frac{52}{500}$ & 11 & $\frac{32}{500}$ & $\frac{40}{500}$ & 3 & $\frac{30}{500}$& $\frac{21}{500}$ & 2 & $\frac{27}{500}$ & $\frac{20}{500}$ \vspace{0.09cm}  \\ 
155000 & 11  & $\frac{51}{500}$ & $\frac{55}{500}$ & 11 & $\frac{27}{500}$ & $\frac{34}{500}$ & 3 & $\frac{27}{500}$& $\frac{25}{500}$ & 2 & $\frac{23}{500}$& $\frac{17}{500}$ \vspace{0.09cm}  \\ 
  195000 & 12 & $\frac{39}{500}$ & $\frac{44}{500}$ & 12  & $\frac{22}{500}$ & $\frac{29}{500}$ & 4 & $\frac{21}{500}$&$\frac{14}{500}$ & 3 & $\frac{16}{500}$& $\frac{13}{500}$ \vspace{0.09cm}  \\  
   235000  & 12  &  $\frac{40}{500}$  & $\frac{42}{500}$ & 12    & $\frac{29}{500}$ &  $\frac{33}{500}$  & 4 & $\frac{18}{500}$& $\frac{16}{500}$ & 3 &  $\frac{12}{500}$& $\frac{9}{500}$ \vspace{0.09cm}  \\   
    275000 & 12 & $\frac{35}{500}$ & $\frac{37}{500}$ & 12 & $\frac{24}{500}$ & $\frac{28}{500}$  & 4 & $\frac{19}{500}$& $\frac{13}{500}$ & 3 & $\frac{9}{500}$ & $\frac{5}{500}$ \vspace{0.09cm}  \\ 
  315000 & 12  & $\frac{24}{500}$& $\frac{28}{500}$ & 12 &  $\frac{17}{500}$ & $\frac{19}{500}$ & 4 & $\frac{11}{500}$& $\frac{8}{500}$ & 3 & $\frac{6}{500}$& $\frac{3}{500}$ \vspace{0.09cm}  \\  
     355000 & 12&   $\frac{21}{500}$ & $\frac{25}{500}$& 12 & $\frac{12}{500}$ & $\frac{15}{500}$  & 4 & $\frac{7}{500}$ & $\frac{4}{500}$ & 3 & $\frac{5}{500}$& $\frac{2}{500}$ \vspace{0.09cm}  \\ 
395000 & 12 & $\frac{18}{500}$ & $\frac{21}{500}$  & 12& $\frac{9}{500}$ & $\frac{12}{500}$  & 4  & $\frac{6}{500}$&  $\frac{3}{500}$ & 3  &$\frac{4}{500}$& $\frac{2}{500}$ \vspace{0.05cm}  \\ 
  \toprule 
\end{tabular}   
   \end{small}   
  \label{tab:TableA9a}
\end{table}

The rate of convergence of the empirical mean of $\left\| \left(\rho - \overline{\rho}_{k_n}^{l} \right) \left(X_{n-1}^{l} \right) \right\|_{H}^{k_n}$, for $l=1,\ldots,N$, with $N = 500$ realizations, can be also numerically fitted (see Figures \ref{fig:F1}-\ref{fig:F3} below).

\begin{figure}[H]
  \hspace{-1.06cm} \includegraphics[width=9.3cm,height=5.9cm]{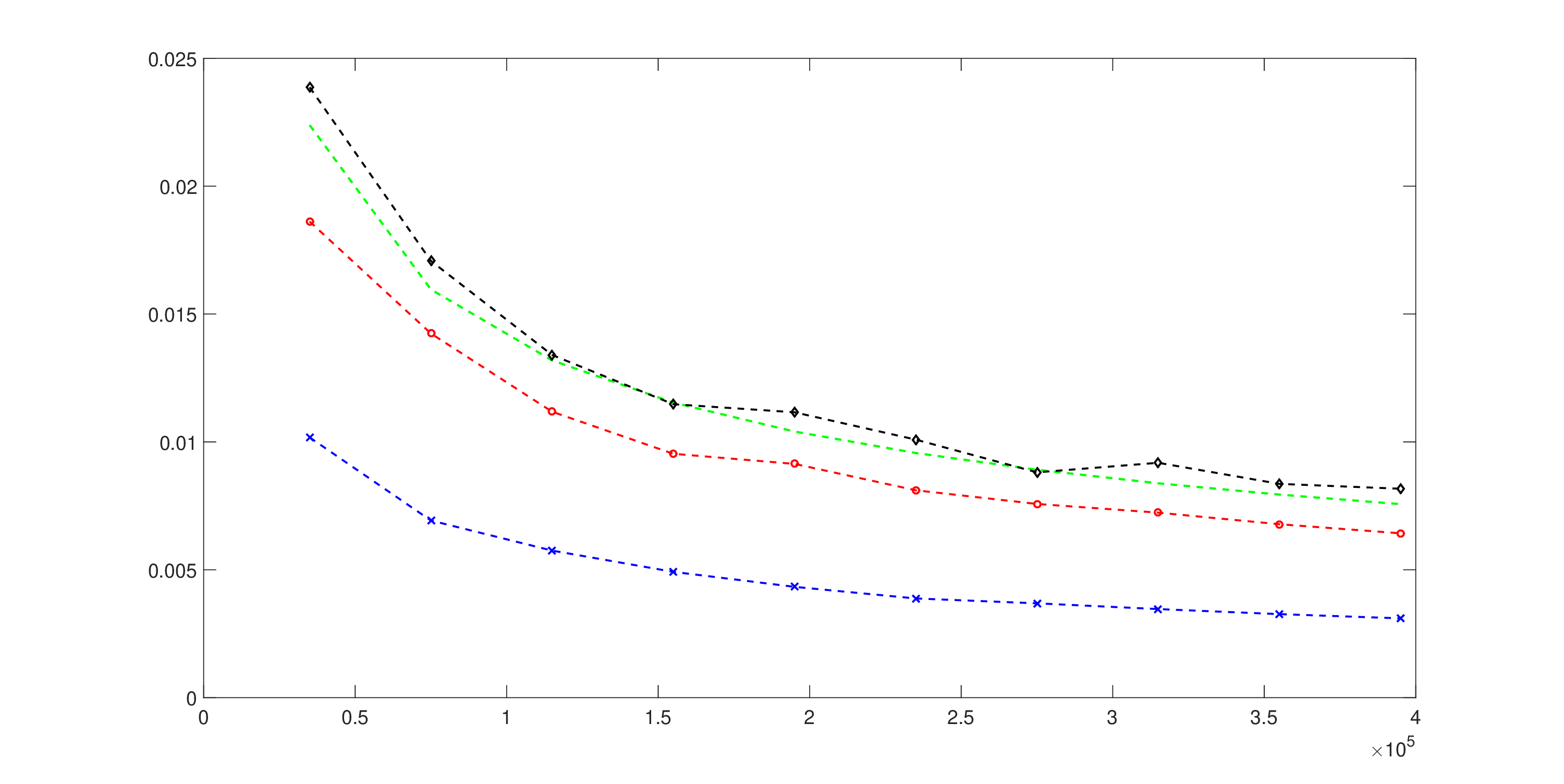}
 \hspace{-1cm}  \includegraphics[width=9.3cm,height=5.9cm]{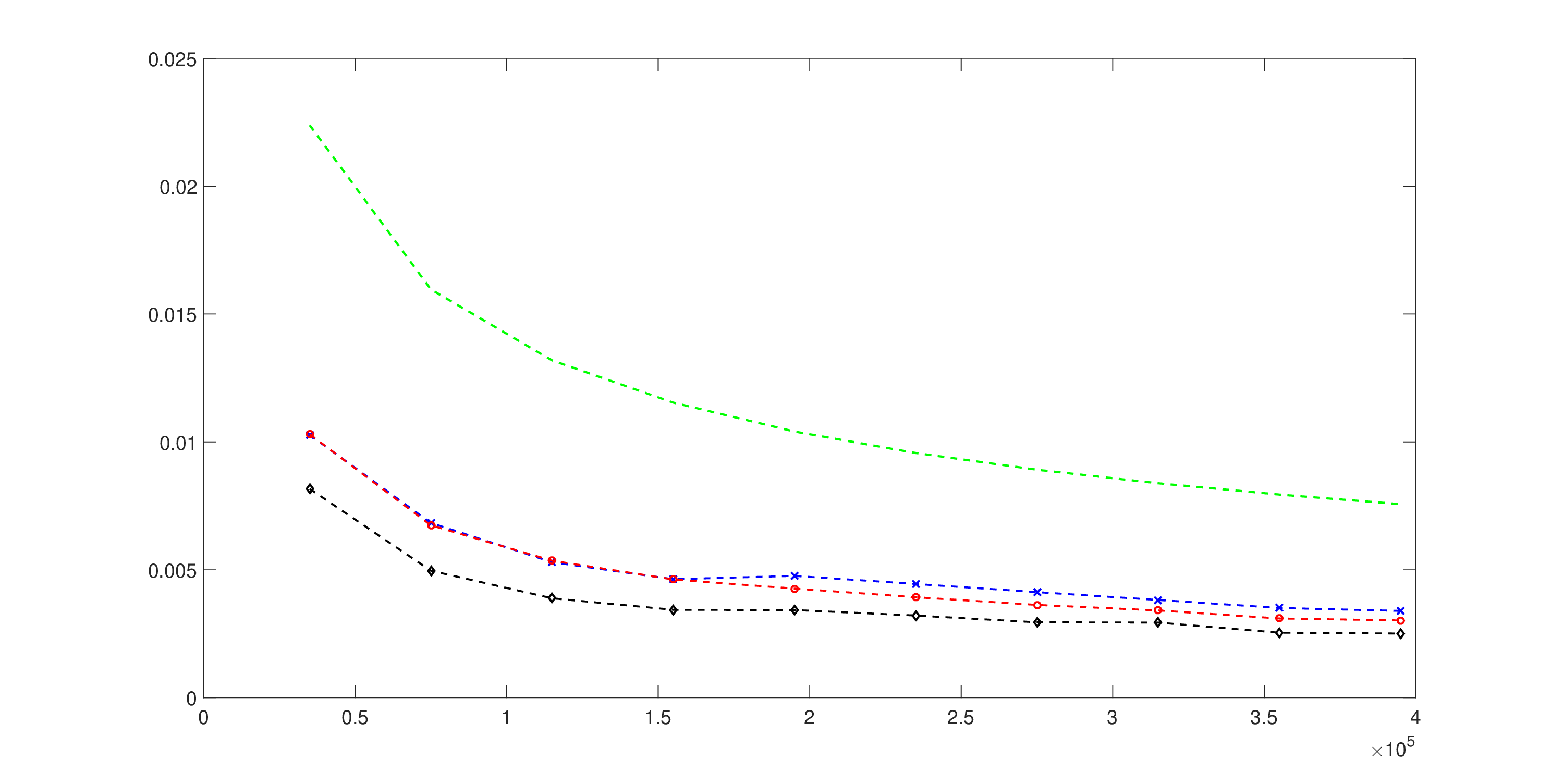}
\vspace{-1.22cm}
\caption{{\small Empirical mean of $\left\| \left(\rho - \overline{\rho}_{k_n}^{l} \right) \left(X_{n-1}^{l} \right) \right\|_{H}^{k_n}$, for scenario 1 (on left) and scenario 3 (on right), for our approach (blue star dotted line) and those one presented in Bosq (2000) (red circle dotted line) and Guillas (2001) (black diamond dotted line). The curve $\xi_{n_t,\beta} = \frac{ \left(\ln(n_t) \right)^{\beta} }{n_{t}^{1/2}}$, with $\beta = 65/100$, is drawn (green dotted line).}} \label{fig:F1}
\end{figure}

\begin{figure}[H]
 \hspace{-1.06cm} \includegraphics[width=9.3cm,height=5.9cm]{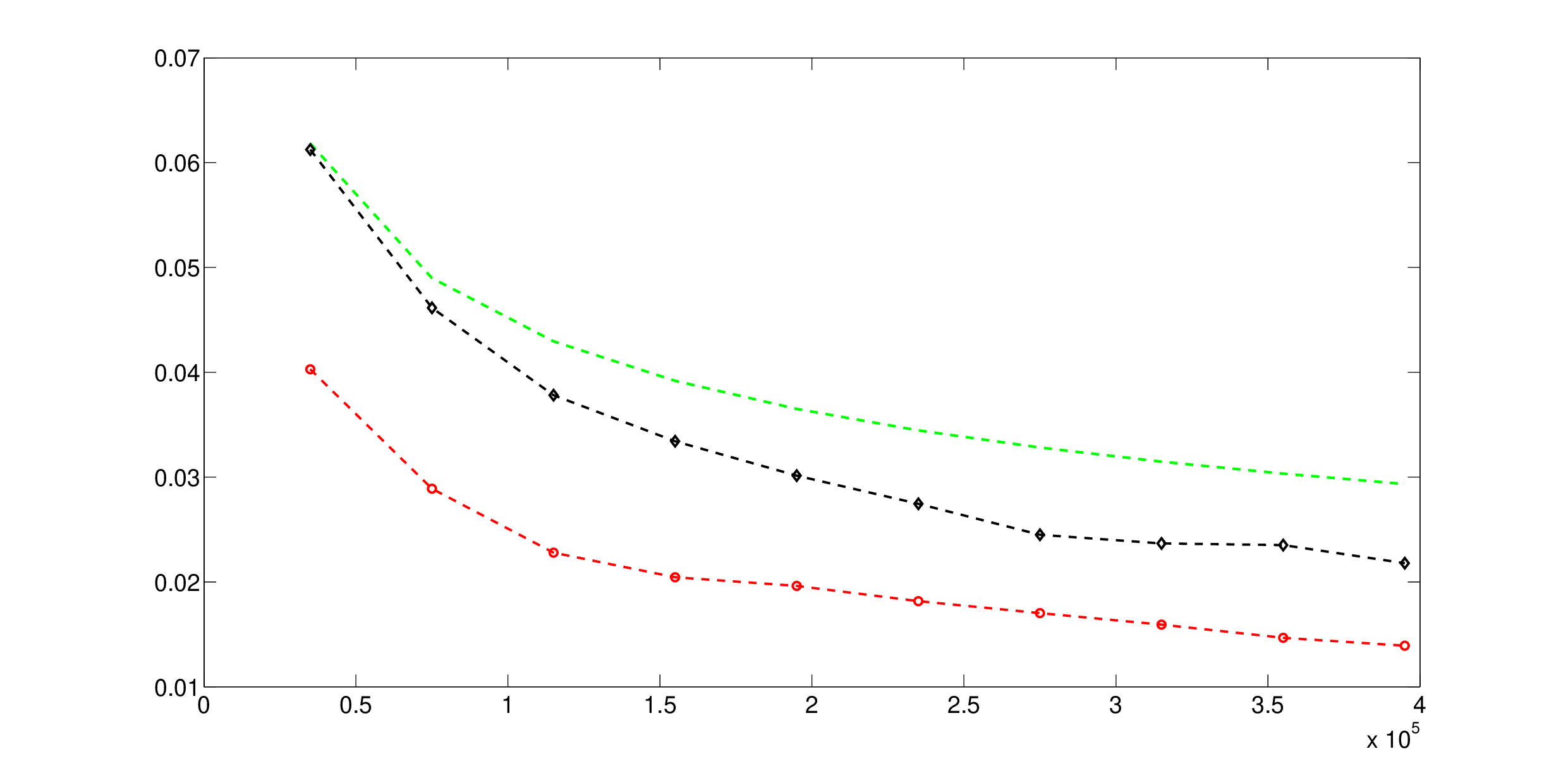}
 \hspace{-1cm}  \includegraphics[width=9.3cm,height=5.9cm]{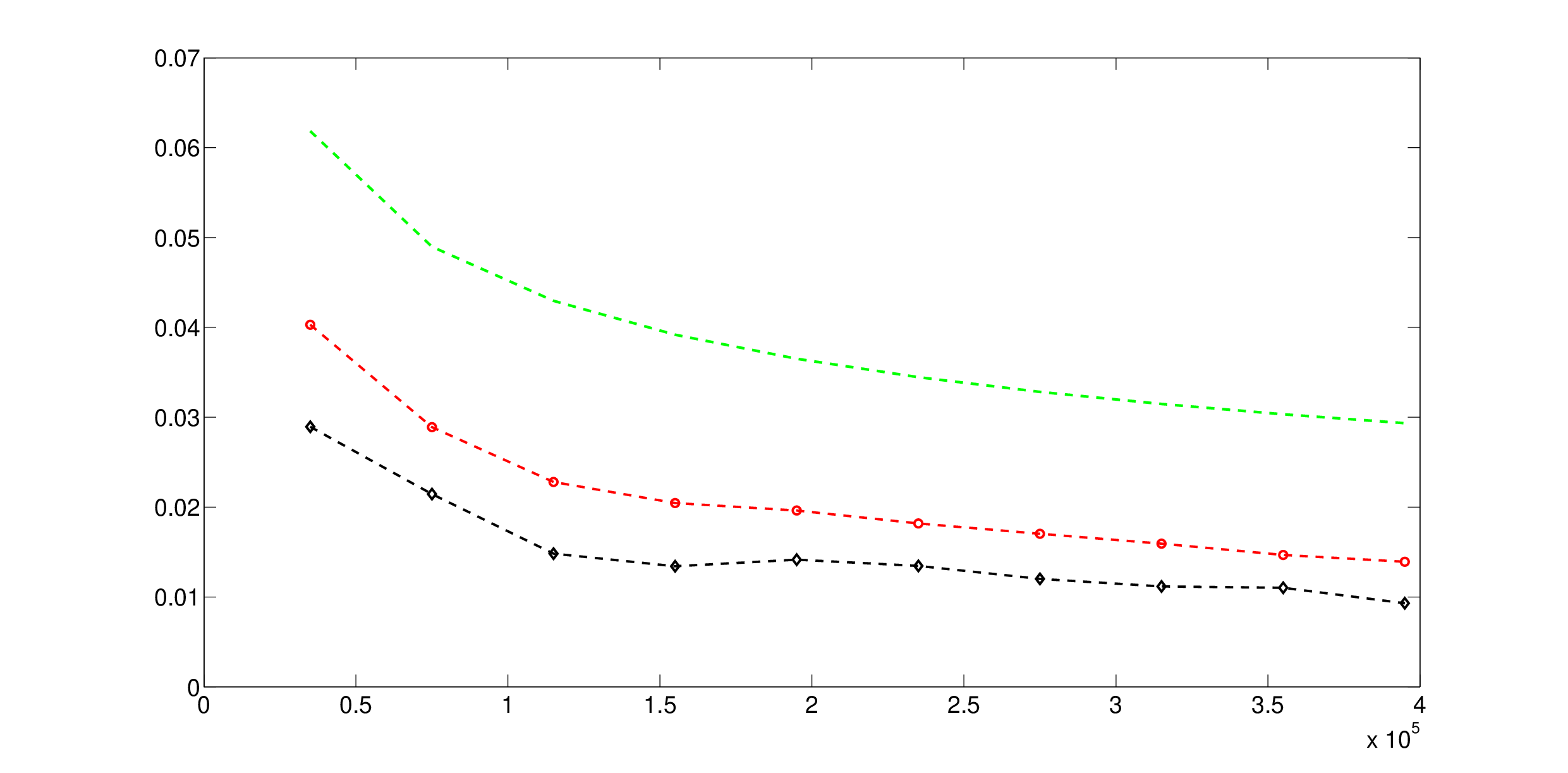}
\vspace{-1.22cm}
\caption{{\small Empirical mean of $\left\| \left(\rho - \overline{\rho}_{k_n}^{l} \right) \left(X_{n-1}^{l} \right) \right\|_{H}^{k_n}$, for scenario 5 (on left) and scenario 7 (on right), for approaches presented in Bosq (2000) (red circle dotted line) and Guillas (2001) (black diamond dotted line). The curve $\xi_{n_t,\beta} = \frac{\left(\ln(n_t) \right)^{\beta} }{n_{t}^{1/3} }$, with $\beta = 3/10$, is drawn (green dotted line).}} \label{fig:F2}
\end{figure}
\begin{figure}[H]
 \hspace{-1.06cm} \includegraphics[width=9.2cm,height=5.9cm]{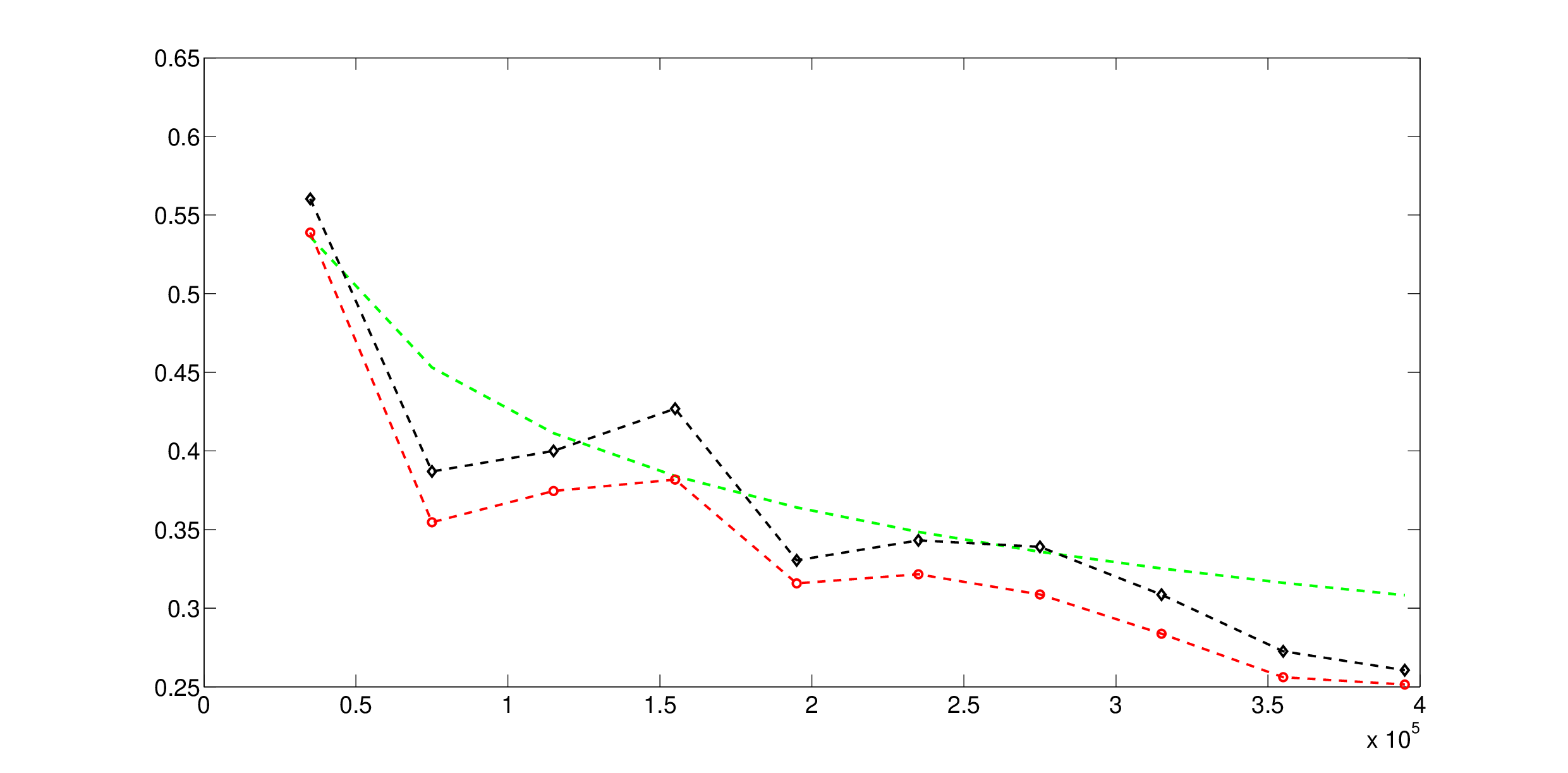}
 \hspace{-1.02cm}  \includegraphics[width=9.2cm,height=5.9cm]{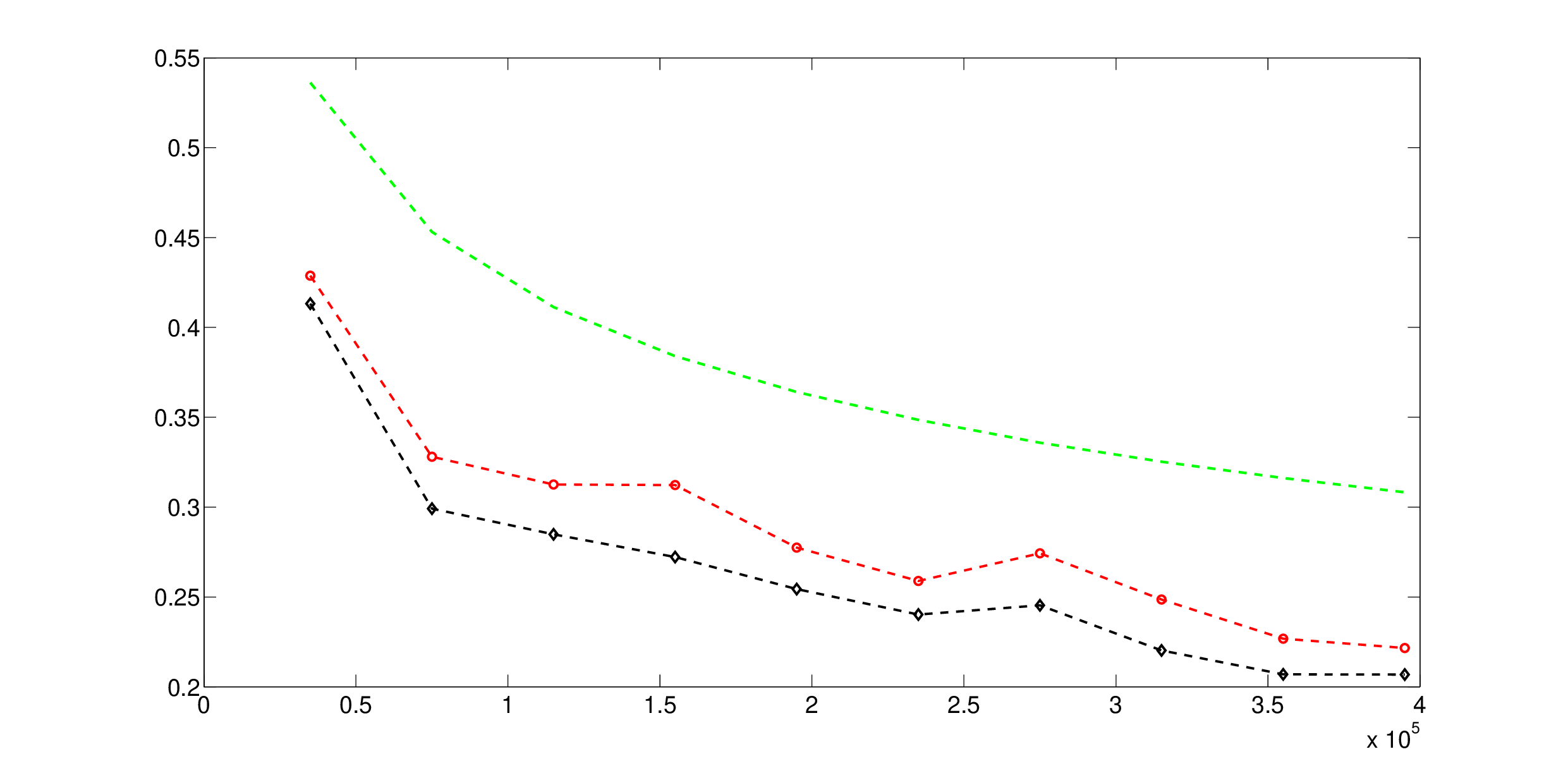}
\vspace{-1.22cm}
\caption{{\small Empirical mean of $\left\| \left(\rho - \overline{\rho}_{k_n}^{l} \right) \left(X_{n-1}^{l} \right) \right\|_{H}^{k_n}$, for scenario 9 (on left) and scenario 11 (on right), for approaches presented in Bosq (2000) (red circle dotted line) and Guillas (2001) (black diamond dotted line). The curve $\xi_{n_t,\beta} = \frac{\left(\ln(n_t) \right)^{\beta} }{n_{t}^{1/3} }$, with $\beta = 125/100$, is drawn (green dotted line).}} \label{fig:F3}
\end{figure}

\vspace{0.5cm}

When approaches formulated in Antoniadis \& Sapatinas (2003) and Besse et al. (2000) are included, smaller sample sizes must be considered due to computational limitations. Hence, a small-sample comparative study is shown in Tables \ref{tab:TableA5}-\ref{tab:TableA10}. The following alternative  norm replaces the  norm reflected in (\ref{a755})-(\ref{a766}), for values $F \left( k_n, n_t, \beta \right)$, when approaches in Besse et al. (2000) are compared:
\begin{eqnarray}
\left\| \left(\rho - \overline{\rho}_{k_n}^{l} \right) \left(X_{n-1}^{l} \right) \right\|_{H} &=& \sqrt{\displaystyle \int_{a}^{b} \left( \rho \left(X_{n-1}^{l} \right) (t) - \overline{\rho}_{k_n}^{l} \left(X_{n-1}^{l} \right) (t) \right)^2 dt }, \quad l=1,\ldots,N.\label{new2}
\end{eqnarray}

Hence, the following diagonal scenarios are regarded:

\begin{table}[H]
 \caption{{\small Diagonal scenarios considered (see Tables \ref{tab:TableA5}-\ref{tab:TableA6} below), with $M=50$, $q=10$, $\delta_2 = 11/10$,  $n_t = 750 + 500(t-1),~t=1,\ldots,13$, and $\xi_{n_t,\beta} = \frac{\left(\ln(n_t) \right)^{\beta}}{n_{t}^{1/2}},~\beta = 65/100$.}}
\centering
\vspace{-0.17cm}
\begin{small}
\begin{tabular}{cccc}
\toprule
  Scenario & $\delta_1$ &  $k_n$ & $h_n$ \\
 \midrule
13 & $3/2$  & $ \lceil \ln(n) \rceil$ & $0.15,~0.25$ \\
14 & $24/10$  & $\lceil \ln(n) \rceil$ & $0.15,~0.25$ \\
15 &  $3/2$  & $ \lceil n^{1/\alpha} \rceil,~\alpha = 6.5$  \\
16 & $24/10$  & $ \lceil n^{1/\alpha} \rceil,~\alpha = 10$ \\
\bottomrule
\end{tabular}
\end{small}
  \label{tab:Table_Scen_C}
\end{table}

Remark that, since both approaches formulated in Besse et al. (2000) not depend on the truncation parameter $k_n$ adopted, we only performed them for scenarios 13 and 14, where different decay rates are considered. Pseudo-diagonal and non-diagonal scenarios are detailed as follows:
\begin{table}[H]
 \caption{{\small Pseudo-diagonal and non-diagonal scenarios considered (see Tables \ref{tab:TableA7}-\ref{tab:TableA10} below), with $\delta_2 = 11/10$,  $n_t = 750 + 500(t-1),~t=1,\ldots,13$, and $\xi_{n_t,\beta} = \frac{\left(\ln(n_t) \right)^{\beta}}{n_{t}^{1/3}}$.}}
\centering
\vspace{-0.17cm}
\begin{small}
\begin{tabular}{ccccc||ccccc}
\toprule
\multicolumn{5}{c||}{Pseudo-diagonal scenarios} & \multicolumn{5}{c}{Non-diagonal scenarios} \\ 
\hline
  Scenario & $\delta_1$ &  $k_n$ & $\beta$ & $h_n$ & Scenario & $\delta_1$ &  $k_n$ & $\beta$ & $h_n$  \\
 \midrule
17 & $3/2$  & $ \lceil \ln(n) \rceil$ & $3/10$ & $1.2,~1.7$ & 21 & $3/2$  & $ \lceil \ln(n) \rceil$ & $125/100$ & $1.2,~1.7$\\
18 & $24/10$  & $\lceil \ln(n) \rceil$ & $3/10$ & $1.2,~1.7$ & 22 &  $24/10$  & $ \lceil \ln(n) \rceil$ & $125/100$& $1.2,~1.7$\\
19&  $3/2$  & $ \lceil n^{1/\alpha} \rceil,~\alpha = 6.5$   & $3/10$ & & 23 & $3/2$  & $ \lceil n^{1/\alpha} \rceil$ & $125/100$ & \\
20& $24/10$  & $ \lceil n^{1/\alpha} \rceil,~\alpha = 10$  & $3/10$ & & 24 &  $24/10$  & $ \lceil n^{1/\alpha} \rceil$ & $125/100$&  \\
\bottomrule
\end{tabular}
\end{small}
  \label{tab:Table_Scen_D}
\end{table}

\begin{table}[H]
 \caption{{\small
$F\left(k_n,n_t, \beta\right)$ values in
(\ref{a955})-(\ref{a755}), for scenarios 13-16. O.A., B, G and AS denote the approaches detailed in Section 8, Bosq (2000), Guillas (2001) and Antoniadis \& Sapatinas (2003), respectively.}}
\vspace{-0.16cm}
\centering
\begin{small}
\begin{tabular}{c@{\hspace{1.5em}}c@{\hspace{0.16cm}}c@{\hspace{0.14cm}}c@{\hspace{0.14cm}}c@{\hspace{0.14cm}}c@{\hspace{1.3em}}c@{\hspace{0.14cm}}c@{\hspace{0.14cm}}c@{\hspace{0.14cm}}c@{\hspace{1.3em}}c@{\hspace{0.14cm}}c@{\hspace{0.14cm}}c@{\hspace{0.14cm}}c@{\hspace{0.14cm}}c@{\hspace{1.3em}}c@{\hspace{0.14cm}}c@{\hspace{0.14cm}}c@{\hspace{0.14cm}}c@{\hspace{0.14cm}}c}
\toprule  
& & \multicolumn{4}{c}{\hspace{-0.4cm}\textit{Scenario 13}} & \multicolumn{4}{c}{\hspace{-0.25cm}\textit{Scenario 14}} & \multicolumn{5}{c}{\hspace{-0.2cm}\textit{Scenario 15}}& \multicolumn{5}{c}{\hspace{-0.2cm}\textit{Scenario 16}}\\
  \midrule[0.05cm]
   $n_t$  & $k_n$ & O.A. &  B & G & AS & O.A. &  B & G & AS & $k_n$ & O.A. &  B & G  & AS &  $k_n$& O.A. &  B & G & AS \\
  \midrule[0.05cm] 
750 & 6 & $\frac{42}{500}$ & $\frac{83}{500}$ & $\frac{90}{500}$ & $\frac{84}{500}$ & $\frac{31}{500}$  &  $\frac{76}{500}$ & $\frac{79}{500}$ & $\frac{72}{500}$ & 2 & $\frac{30}{500}$& $\frac{33}{500}$& $\frac{34}{500}$ & $\frac{21}{500}$& 1 & $\frac{19}{500}$& $\frac{24}{500}$ & $\frac{24}{500}$& $\frac{14}{500}$ \vspace{0.09cm}  \\ 

1250 & 7& $\frac{28}{500}$ & $\frac{74}{500}$ &  $\frac{88}{500}$ & $\frac{76}{500}$& $\frac{29}{500}$  
& $\frac{74}{500}$ & $\frac{78}{500}$ & $\frac{76}{500}$  & 2& $\frac{27}{500}$& $\frac{29}{500}$& $\frac{30}{500}$& $\frac{20}{500}$ & 2 &  $\frac{17}{500}$&  $\frac{21}{500}$& $\frac{22}{500}$ & $\frac{13}{500}$ \vspace{0.09cm}  \\ 
 
1750 & 7&$\frac{27}{500}$ & $\frac{70}{500}$ &   $\frac{84}{500}$ & $\frac{75}{500}$ & $\frac{28}{500}$ &    $\frac{71}{500}$ & $\frac{73}{500}$ & $\frac{70}{500}$ & 2& $\frac{25}{500}$& $\frac{25}{500}$& $\frac{27}{500}$& $\frac{17}{500}$ & 2&   $\frac{16}{500}$& $\frac{19}{500}$ & $\frac{21}{500}$  & $\frac{11}{500}$\vspace{0.09cm}  \\ 

2250  & 7& $\frac{26}{500}$ & $\frac{66}{500}$ & $\frac{81}{500}$ & $\frac{71}{500}$ & $\frac{25}{500}$ &  $\frac{68}{500}$ & $\frac{67}{500}$ & $\frac{65}{500}$ & 3&  $\frac{24}{500}$& $\frac{23}{500}$& $\frac{26}{500}$ & $\frac{15}{500}$& 2 &   $\frac{13}{500}$& $\frac{17}{500}$& $\frac{20}{500}$  & $\frac{9}{500}$\vspace{0.09cm}  \\ 

  2750 & 7& $\frac{28}{500}$ & $\frac{68}{500}$ & $\frac{82}{500}$ & $\frac{70}{500}$ & $\frac{24}{500}$  & $\frac{63}{500}$ & $\frac{62}{500}$ & $\frac{59}{500}$ & 3& $\frac{21}{500}$& $\frac{21}{500}$&$\frac{26}{500}$& $\frac{12}{500}$ & 2 & $\frac{12}{500}$& $\frac{15}{500}$& $\frac{18}{500}$ & $\frac{8}{500}$ \vspace{0.09cm}  \\  
  
 3250  & 8& $\frac{25}{500}$ & $\frac{66}{500}$ & $\frac{76}{500}$ &  $\frac{72}{500}$  & $\frac{21}{500}$   & $\frac{58}{500}$  & $\frac{59}{500}$ &  $\frac{55}{500}$  & 3& $\frac{20}{500}$& $\frac{20}{500}$& $\frac{25}{500}$ & $\frac{10}{500}$& 2 & $\frac{10}{500}$&  $\frac{14}{500}$& $\frac{17}{500}$ & $\frac{7}{500}$ \vspace{0.09cm}  \\   
  
3750 & 8& $\frac{23}{500}$ & $\frac{60}{500}$ & $\frac{72}{500}$ & $\frac{72}{500}$ & $\frac{21}{500}$ &  $\frac{53}{500}$ & $\frac{54}{500}$ & $\frac{54}{500}$  & 3& $\frac{17}{500}$ & $\frac{18}{500}$& $\frac{24}{500}$ & $\frac{9}{500}$& 2&  $\frac{9}{500}$& $\frac{11}{500}$ & $\frac{14}{500}$ &  $\frac{7}{500}$\vspace{0.09cm}  \\
 
4250 & 8&$\frac{23}{500}$ & $\frac{59}{500}$ &  $\frac{70}{500}$ & $\frac{71}{500}$& $\frac{20}{500}$  & $\frac{49}{500}$ & $\frac{51}{500}$ & $\frac{48}{500}$ & 3& $\frac{14}{500}$& $\frac{16}{500}$& $\frac{18}{500}$& $\frac{9}{500}$ & 2& $\frac{8}{500}$ & $\frac{10}{500}$& $\frac{11}{500}$ & $\frac{6}{500}$ \vspace{0.09cm}  \\  

4750 & 8& $\frac{21}{500}$ & $\frac{56}{500}$ & $\frac{67}{500}$ & $\frac{69}{500}$ & $\frac{18}{500}$ &$\frac{47}{500}$ & $\frac{49}{500}$ & $\frac{45}{500}$  & 3 & $\frac{13}{500}$& $\frac{13}{500}$ & $\frac{15}{500}$ & $\frac{8}{500}$& 2& $\frac{7}{500}$& $\frac{8}{500}$& $\frac{9}{500}$ & $\frac{5}{500}$ \vspace{0.09cm}  \\ 

5250 & 8&$\frac{18}{500}$ & $\frac{55}{500}$ & $\frac{65}{500}$ & $\frac{68}{500}$ & $\frac{15}{500}$ &  $\frac{47}{500}$ & $\frac{48}{500}$ & $\frac{44}{500}$  & 3& $\frac{12}{500}$ & $\frac{10}{500}$&  $\frac{13}{500}$ & $\frac{8}{500}$ & 2 & $\frac{7}{500}$& $\frac{7}{500}$& $\frac{7}{500}$ & $\frac{3}{500}$ \vspace{0.05cm}  \\ 

5750 & 8&$\frac{20}{500}$ & $\frac{58}{500}$ & $\frac{66}{500}$ & $\frac{68}{500}$ & $\frac{16}{500}$ &  $\frac{45}{500}$ & $\frac{50}{500}$ & $\frac{47}{500}$  & 3& $\frac{11}{500}$ & $\frac{9}{500}$&  $\frac{11}{500}$ & $\frac{7}{500}$& 2 & $\frac{6}{500}$& $\frac{7}{500}$& $\frac{6}{500}$ & $\frac{2}{500}$ \vspace{0.05cm}  \\ 

6250 & 8&$\frac{16}{500}$ & $\frac{57}{500}$ & $\frac{62}{500}$ & $\frac{67}{500}$ & $\frac{11}{500}$ &  $\frac{42}{500}$ & $\frac{47}{500}$ & $\frac{52}{500}$  & 3& $\frac{9}{500}$ & $\frac{8}{500}$&  $\frac{10}{500}$ & $\frac{7}{500}$ & 2 & $\frac{5}{500}$& $\frac{5}{500}$ & $\frac{5}{500}$& $\frac{2}{500}$ \vspace{0.05cm}  \\ 

6750 & 8&$\frac{14}{500}$ & $\frac{54}{500}$ & $\frac{59}{500}$ & $\frac{67}{500}$ & $\frac{9}{500}$ &  $\frac{41}{500}$ & $\frac{45}{500}$ & $\frac{42}{500}$  & 3& $\frac{7}{500}$ & $\frac{8}{500}$&  $\frac{8}{500}$ & $\frac{6}{500}$& 2 & $\frac{3}{500}$& $\frac{4}{500}$& $\frac{4}{500}$ & 0 \vspace{0.05cm}  \\ 
  \toprule 
\end{tabular} 
\end{small}     
  \label{tab:TableA5}
\end{table}

\begin{table}[H]
 \caption{{\small
$F \left(k_n,n_t, \beta \right)$ values in
(\ref{a955}) and (\ref{new2}), for scenarios 13-14. $B_{0.15}$ and $B_{0.25}$ denotes the kernel-based approach in Besse et al. (2000), for $h_n = 0.15,~0.25$, respectively. $B_q$ denotes its penalized prediction approach.}}
\vspace{-0.16cm}
\centering
\begin{small}
     \begin{tabular}{c@{\hspace{1.5em}}ccc@{\hspace{1.6em}}ccc}
\toprule  
& \multicolumn{3}{c}{\hspace{-0.4cm}\textit{Scenario 13}} & \multicolumn{3}{c}{\hspace{-0.25cm}\textit{Scenario 14}} \\
  \midrule[0.05cm]
   $n_t$   & $B_{0.15}$ &  $B_{0.25}$ & $B_q$ &  $B_{0.15}$ & $B_{0.25}$ & $B_q$ \\
  \midrule[0.05cm] 
750 &  $\frac{85}{500}$ & $\frac{88}{500}$ & $\frac{6}{500}$ & $\frac{76}{500}$ & $\frac{80}{500}$  &  $\frac{3}{500}$  \vspace{0.09cm}  \\ 

1250 & $\frac{80}{500}$ & $\frac{79}{500}$ & $\frac{6}{500}$ & $\frac{75}{500}$ & $\frac{73}{500}$  &  $\frac{3}{500}$  \vspace{0.09cm}  \\ 

1750 & $\frac{76}{500}$ & $\frac{71}{500}$ & $\frac{5}{500}$ & $\frac{73}{500}$ & $\frac{67}{500}$  &  $\frac{2}{500}$  \vspace{0.09cm}  \\ 

2250  & $\frac{78}{500}$ & $\frac{60}{500}$ & $\frac{4}{500}$ & $\frac{72}{500}$ & $\frac{57}{500}$  &  $\frac{3}{500}$  \vspace{0.09cm}  \\

  2750 & $\frac{73}{500}$ & $\frac{57}{500}$ & $\frac{4}{500}$ & $\frac{70}{500}$ & $\frac{53}{500}$  &  $\frac{3}{500}$  \vspace{0.09cm}  \\ 
3250  & $\frac{75}{500}$ & $\frac{53}{500}$ & $\frac{2}{500}$ & $\frac{67}{500}$ & $\frac{51}{500}$  &  $\frac{2}{500}$  \vspace{0.09cm}  \\   
3750 & $\frac{70}{500}$ & $\frac{49}{500}$ & $\frac{2}{500}$ & $\frac{67}{500}$ & $\frac{43}{500}$  &  $\frac{1}{500}$  \vspace{0.09cm}  \\ 

4250 & $\frac{72}{500}$ & $\frac{44}{500}$ & $\frac{1}{500}$ & $\frac{65}{500}$ & $\frac{41}{500}$  &  0 \vspace{0.09cm}  \\ 

4750 &$\frac{68}{500}$ & $\frac{39}{500}$ & $\frac{3}{500}$ & $\frac{63}{500}$ & $\frac{38}{500}$  &  $\frac{1}{500}$  \vspace{0.09cm}  \\ 

5250 & $\frac{65}{500}$ & $\frac{496}{500}$ & $\frac{3}{500}$ & $\frac{62}{500}$ & $\frac{33}{500}$  &  $\frac{2}{500}$  \vspace{0.09cm}  \\ 

5750 & $\frac{62}{500}$ & $\frac{34}{500}$ & $\frac{2}{500}$ & $\frac{60}{500}$ & $\frac{31}{500}$  &  $\frac{2}{500}$  \vspace{0.09cm}  \\ 

6250 & $\frac{60}{500}$ & $\frac{33}{500}$ & $\frac{3}{500}$ & $\frac{60}{500}$ & $\frac{28}{500}$  &  $\frac{1}{500}$  \vspace{0.09cm}  \\

6750 & $\frac{59}{500}$ & $\frac{33}{500}$ & $\frac{3}{500}$ & $\frac{57}{500}$ & $\frac{24}{500}$  &  $\frac{1}{500}$  \vspace{0.09cm}  \\ 
  \toprule 
\end{tabular} 
     \end{small}     
  \label{tab:TableA6}
\end{table}

\begin{table}[H]
 \caption{{\small
$F\left(k_n,n_t, \beta\right)$ values in
(\ref{a955}) and (\ref{a766}), for scenarios 17-20. B and G denote the approaches in Bosq (2000) and Guillas (2001), respectively; AS denotes the approach in Antoniadis \& Sapatinas (2003).}}
\vspace{-0.16cm}
\centering
\begin{small}
    \begin{tabular}{c@{\hspace{1.5em}}c@{\hspace{0.16cm}}c@{\hspace{0.14cm}}c@{\hspace{0.14cm}}c@{\hspace{1.3em}}c@{\hspace{0.14cm}}c@{\hspace{0.14cm}}c@{\hspace{0.14cm}}c@{\hspace{1.3em}}c@{\hspace{0.14cm}}c@{\hspace{0.14cm}}c@{\hspace{0.14cm}}c@{\hspace{1.3em}}c@{\hspace{0.14cm}}c@{\hspace{0.14cm}}c@{\hspace{0.14cm}}c}
\toprule  
& \multicolumn{4}{c}{\hspace{-0.4cm}\textit{Scenario 17}} & \multicolumn{4}{c}{\hspace{-0.25cm}\textit{Scenario 18}} & \multicolumn{4}{c}{\hspace{-0.2cm}\textit{Scenario 19}}& \multicolumn{4}{c}{\hspace{-0.2cm}\textit{Scenario 20}}\\
  \midrule[0.05cm]
   $n_t$  & $k_n$ &  B & G & AS & $k_n$ &  B & G & AS & $k_n$  &  B & G  & AS &  $k_n$ &  B & G & AS \\
  \midrule[0.05cm]  
750 & 6  & $\frac{135}{500}$ & $\frac{146}{500}$ & $\frac{176}{500}$ & 6 &   $\frac{123}{500}$ & $\frac{129}{500}$ & $\frac{180}{500}$ & 2  & $\frac{62}{500}$& $\frac{48}{500}$ & $\frac{41}{500}$& 1  & $\frac{50}{500}$ & $\frac{39}{500}$ & $\frac{38}{500}$\vspace{0.09cm}  \\ 

1250 & 7& $\frac{124}{500}$ &  $\frac{130}{500}$ & $\frac{166}{500}$ & 7
 & $\frac{117}{500}$ & $\frac{120}{500}$ & $\frac{175}{500}$  & 2 & $\frac{60}{500}$& $\frac{42}{500}$& $\frac{37}{500}$ & 2  &  $\frac{47}{500}$& $\frac{36}{500}$ & $\frac{33}{500}$ \vspace{0.09cm}  \\ 

1750 & 7&  $\frac{113}{500}$ &   $\frac{122}{500}$ & $\frac{159}{500}$ & 7 &    $\frac{104}{500}$ & $\frac{110}{500}$ & $\frac{168}{500}$ & 2 & $\frac{53}{500}$& $\frac{36}{500}$& $\frac{34}{500}$ & 2 & $\frac{41}{500}$ & $\frac{30}{500}$ & $\frac{31}{500}$ \vspace{0.09cm}  \\ 

2250  & 7  & $\frac{89}{500}$ & $\frac{115}{500}$ & $\frac{153}{500}$ & 7 &  $\frac{86}{500}$ & $\frac{91}{500}$ & $\frac{164}{500}$ & 3 & $\frac{49}{500}$& $\frac{31}{500}$ & $\frac{29}{500}$& 2  & $\frac{35}{500}$& $\frac{28}{500}$ & $\frac{30}{500}$ \vspace{0.09cm}  \\

 2750 & 7  & $\frac{80}{500}$ & $\frac{100}{500}$ & $\frac{133}{500}$ & 7  & $\frac{76}{500}$ & $\frac{83}{500}$ & $\frac{149}{500}$ & 3 & $\frac{44}{500}$&$\frac{28}{500}$& $\frac{27}{500}$ & 2  & $\frac{32}{500}$& $\frac{28}{500}$ & $\frac{28}{500}$ \vspace{0.09cm}  \\  
  
  3250  & 8  & $\frac{99}{500}$ & $\frac{104}{500}$ &  $\frac{139}{500}$ & 8   & $\frac{71}{500}$  & $\frac{78}{500}$ &  $\frac{153}{500}$  & 3 & $\frac{40}{500}$&$\frac{26}{500}$& $\frac{26}{500}$ & 2  & $\frac{27}{500}$& $\frac{27}{500}$ & $\frac{25}{500}$ \vspace{0.09cm}  \\  
    
3750 & 8  & $\frac{67}{500}$ & $\frac{78}{500}$ & $\frac{136}{500}$ & 8 &  $\frac{62}{500}$ & $\frac{67}{500}$ & $\frac{142}{500}$  & 3  & $\frac{35}{500}$& $\frac{24}{500}$ & $\frac{25}{500}$& 2 & $\frac{24}{500}$ & $\frac{26}{500}$ & $\frac{23}{500}$ \vspace{0.09cm}  \\ 

4250 & 8 & $\frac{65}{500}$ &  $\frac{74}{500}$ & $\frac{129}{500}$ & 8 & $\frac{60}{500}$ & $\frac{63}{500}$ & $\frac{133}{500}$ & 3 & $\frac{30}{500}$& $\frac{23}{500}$& $\frac{22}{500}$ & 2  & $\frac{22}{500}$& $\frac{22}{500}$ & $\frac{19}{500}$ \vspace{0.09cm}  \\  

4750 & 8  & $\frac{61}{500}$ & $\frac{63}{500}$ & $\frac{127}{500}$ & 8 &$\frac{55}{500}$ & $\frac{60}{500}$ & $\frac{126}{500}$  & 3  & $\frac{28}{500}$ & $\frac{19}{500}$ & $\frac{20}{500}$& 2 & $\frac{20}{500}$& $\frac{16}{500}$ & $\frac{13}{500}$ \vspace{0.09cm}  \\ 

5250 & 8 & $\frac{48}{500}$ & $\frac{51}{500}$ & $\frac{125}{500}$ & 8 &  $\frac{46}{500}$ & $\frac{49}{500}$ & $\frac{122}{500}$  & 3  & $\frac{25}{500}$&  $\frac{17}{500}$ & $\frac{16}{500}$ & 2  & $\frac{17}{500}$& $\frac{12}{500}$  & $\frac{10}{500}$\vspace{0.05cm}  \\ 

5750 & 8  & $\frac{4}{500}$ & $\frac{49}{500}$ & $\frac{122}{500}$ & 8 &  $\frac{39}{500}$ & $\frac{42}{500}$ & $\frac{113}{500}$  & 3 & $\frac{20}{500}$&  $\frac{14}{500}$ & $\frac{13}{500}$& 2  & $\frac{15}{500}$& $\frac{7}{500}$ & $\frac{5}{500}$ \vspace{0.05cm}  \\ 

6250 & 8  & $\frac{38}{500}$ & $\frac{45}{500}$ & $\frac{118}{500}$ & 8 &  $\frac{33}{500}$ & $\frac{35}{500}$ & $\frac{108}{500}$  & 3 & $\frac{19}{500}$&  $\frac{13}{500}$ & $\frac{10}{500}$ & 2  & $\frac{13}{500}$& $\frac{7}{500}$ & $\frac{3}{500}$ \vspace{0.05cm}  \\ 

6750 & 8  & $\frac{36}{500}$ & $\frac{40}{500}$ & $\frac{114}{500}$ & 8 &  $\frac{29}{500}$ & $\frac{31}{500}$ & $\frac{101}{500}$  & 3  & $\frac{13}{500}$&  $\frac{12}{500}$ & $\frac{9}{500}$& 2  & $\frac{10}{500}$& $\frac{8}{500}$ & $\frac{3}{500}$  \vspace{0.05cm}  \\ 
  \toprule 
\end{tabular}  
    \end{small}    
  \label{tab:TableA7}
\end{table}

\begin{table}[H]
 \caption{{\small
$F\left(k_n,n_t, \beta\right)$ values in
(\ref{a955}) and (\ref{new2}), for scenarios 17-18. $B_{1.2}$ and $B_{1.7}$ denotes the kernel-based approach in Besse et al. (2000), for $h_n = 1.2,~ 1.7$, respectively. $B_q$ denotes its penalized prediction approach.}}
\vspace{-0.16cm}
\centering
\begin{small}
     \begin{tabular}{c@{\hspace{1.5em}}ccc@{\hspace{1.6em}}ccc}
\toprule  
& \multicolumn{3}{c}{\hspace{-0.4cm}\textit{Scenario 17}} & \multicolumn{3}{c}{\hspace{-0.25cm}\textit{Scenario 18}} \\
  \midrule[0.05cm]
   $n_t$   & $B_{1.2}$ &  $B_{1.7}$ & $B_q$  & $B_{1.2}$ &  $B_{1.7}$ & $B_q$ \\  
  \midrule[0.05cm] 
750 &  $\frac{174}{500}$ & $\frac{233}{500}$ & $\frac{18}{500}$ & $\frac{167}{500}$ & $\frac{180}{500}$  &  $\frac{10}{500}$  \vspace{0.09cm}  \\ 
1250 & $\frac{158}{500}$ & $\frac{214}{500}$ & $\frac{10}{500}$ & $\frac{151}{500}$ & $\frac{169}{500}$  &  $\frac{7}{500}$  \vspace{0.09cm}  \\ 
1750 & $\frac{149}{500}$ & $\frac{199}{500}$ & $\frac{9}{500}$ & $\frac{133}{500}$ & $\frac{155}{500}$  &  $\frac{6}{500}$  \vspace{0.09cm}  \\ 
2250  & $\frac{146}{500}$ & $\frac{185}{500}$ & $\frac{7}{500}$ & $\frac{130}{500}$ & $\frac{146}{500}$  &  $\frac{4}{500}$  \vspace{0.09cm}  \\ 
  2750 & $\frac{131}{500}$ & $\frac{190}{500}$ & $\frac{6}{500}$ & $\frac{127}{500}$ & $\frac{140}{500}$  &  $\frac{3}{500}$  \vspace{0.09cm}  \\ 
 3250  & $\frac{129}{500}$ & $\frac{193}{500}$ & $\frac{5}{500}$ & $\frac{119}{500}$ & $\frac{135}{500}$  &  $\frac{3}{500}$  \vspace{0.09cm}  \\   
3750 & $\frac{125}{500}$ & $\frac{162}{500}$ & $\frac{6}{500}$ & $\frac{115}{500}$ & $\frac{130}{500}$  &  $\frac{4}{500}$  \vspace{0.09cm}  \\ 
4250 & $\frac{138}{500}$ & $\frac{160}{500}$ & $\frac{4}{500}$ & $\frac{109}{500}$ & $\frac{121}{500}$  &  $\frac{2}{500}$  \vspace{0.09cm}  \\ 
4750 &$\frac{133}{500}$ & $\frac{162}{500}$ & $\frac{2}{500}$ & $\frac{108}{500}$ & $\frac{117}{500}$  &  $\frac{2}{500}$  \vspace{0.09cm}  \\ 

5250 & $\frac{120}{500}$ & $\frac{154}{500}$ & $\frac{1}{500}$ & $\frac{107}{500}$ & $\frac{114}{500}$  &  $\frac{1}{500}$  \vspace{0.09cm}  \\ 
5750 & $\frac{118}{500}$ & $\frac{156}{500}$ & $\frac{2}{500}$ & $\frac{104}{500}$ & $\frac{111}{500}$  &  $\frac{1}{500}$  \vspace{0.09cm}  \\ 
6250 & $\frac{116}{500}$ & $\frac{144}{500}$ & $\frac{1}{500}$ & $\frac{99}{500}$ & $\frac{103}{500}$  &  0  \vspace{0.09cm}  \\ 

6750 & $\frac{111}{500}$ & $\frac{135}{500}$ & 0 & $\frac{94}{500}$ & $\frac{100}{500}$  &  0  \vspace{0.09cm}  \\ 
  \toprule 
\end{tabular} 
     \end{small}     
  \label{tab:TableA8}
\end{table}

\begin{table}[H]
 \caption{{\small
$F\left(k_n,n_t, \beta\right)$ values in
(\ref{a955}) and (\ref{a766}), for scenarios 21-24. B and G denote the approaches in  Bosq (2000) and Guillas (2001), resp.; AS denotes the approach in Antoniadis \& Sapatinas (2003).}}
\vspace{-0.16cm}
\centering
\begin{small}
    \begin{tabular}{c@{\hspace{1.5em}}c@{\hspace{0.16cm}}c@{\hspace{0.14cm}}c@{\hspace{0.14cm}}c@{\hspace{1.3em}}c@{\hspace{0.14cm}}c@{\hspace{0.14cm}}c@{\hspace{0.14cm}}c@{\hspace{1.3em}}c@{\hspace{0.14cm}}c@{\hspace{0.14cm}}c@{\hspace{0.14cm}}c@{\hspace{1.3em}}c@{\hspace{0.14cm}}c@{\hspace{0.14cm}}c@{\hspace{0.14cm}}c}
\toprule  
& \multicolumn{4}{c}{\hspace{-0.4cm}\textit{Scenario 21}} & \multicolumn{4}{c}{\hspace{-0.25cm}\textit{Scenario 22}} & \multicolumn{4}{c}{\hspace{-0.2cm}\textit{Scenario 23}}& \multicolumn{4}{c}{\hspace{-0.2cm}\textit{Scenario 24}}\\
  \midrule[0.05cm]
   $n_t$  & $k_n$ &  B & G & AS & $k_n$ &  B & G & AS & $k_n$  &  B & G  & AS &  $k_n$ &  B & G & AS \\
  \midrule[0.05cm]  
750 & 6  & $\frac{86}{500}$ & $\frac{90}{500}$ & $\frac{88}{500}$ & 6  &  $\frac{80}{500}$ & $\frac{84}{500}$ & $\frac{83}{500}$ & 2  & $\frac{73}{500}$& $\frac{66}{500}$ & $\frac{75}{500}$& 1  & $\frac{55}{500}$ & $\frac{42}{500}$  & $\frac{60}{500}$\vspace{0.09cm}  \\ 

1250 & 7  & $\frac{81}{500}$ &  $\frac{84}{500}$ & $\frac{86}{500}$ & 7
 & $\frac{78}{500}$ & $\frac{81}{500}$ & $\frac{85}{500}$  & 2 & $\frac{69}{500}$& $\frac{64}{500}$& $\frac{71}{500}$ & 2  &  $\frac{48}{500}$& $\frac{39}{500}$ & $\frac{51}{500}$ \vspace{0.09cm}  \\

1750 & 7  & $\frac{77}{500}$ &   $\frac{80}{500}$ & $\frac{85}{500}$ & 7&    $\frac{73}{500}$ & $\frac{87}{500}$ & $\frac{86}{500}$ & 2 & $\frac{64}{500}$& $\frac{60}{500}$& $\frac{70}{500}$ & 2 & $\frac{46}{500}$ & $\frac{32}{500}$ & $\frac{50}{500}$ \vspace{0.09cm}  \\ 

2250  & 7  & $\frac{73}{500}$ & $\frac{77}{500}$ & $\frac{86}{500}$ & 7 &  $\frac{68}{500}$ & $\frac{72}{500}$ & $\frac{84}{500}$ & 3 & $\frac{59}{500}$& $\frac{56}{500}$ & $\frac{63}{500}$& 2  & $\frac{41}{500}$& $\frac{31}{500}$ & $\frac{46}{500}$\vspace{0.09cm}  \\ 

  2750 & 7 & $\frac{70}{500}$ & $\frac{73}{500}$ & $\frac{83}{500}$ & 7  & $\frac{55}{500}$ & $\frac{70}{500}$ & $\frac{80}{500}$ & 3 & $\frac{50}{500}$&$\frac{54}{500}$& $\frac{55}{500}$ & 2  & $\frac{37}{500}$& $\frac{27}{500}$& $\frac{45}{500}$  \vspace{0.09cm}  \\ 
   
 3250  & 8  & $\frac{65}{500}$ & $\frac{68}{500}$ &  $\frac{82}{500}$  & 8  & $\frac{47}{500}$  & $\frac{60}{500}$ &  $\frac{78}{500}$  & 3& $\frac{47}{500}$& $\frac{50}{500}$ & $\frac{51}{500}$& 2 &  $\frac{35}{500}$& $\frac{25}{500}$ & $\frac{41}{500}$ \vspace{0.09cm}  \\ 
    
3750 & 8  & $\frac{54}{500}$ & $\frac{59}{500}$ & $\frac{80}{500}$ & 8 &  $\frac{43}{500}$ & $\frac{53}{500}$ & $\frac{75}{500}$  & 3 & $\frac{45}{500}$& $\frac{43}{500}$ & $\frac{48}{500}$& 2 & $\frac{31}{500}$ & $\frac{24}{500}$ & $\frac{37}{500}$ \vspace{0.09cm}  \\ 

4250 & 8 & $\frac{51}{500}$ &  $\frac{57}{500}$ & $\frac{77}{500}$ & 8 & $\frac{39}{500}$ & $\frac{46}{500}$ & $\frac{72}{500}$ & 3 & $\frac{42}{500}$& $\frac{38}{500}$& $\frac{40}{500}$ & 2  & $\frac{27}{500}$& $\frac{21}{500}$& $\frac{35}{500}$ \vspace{0.09cm}  \\  

4750 & 8  & $\frac{45}{500}$ & $\frac{51}{500}$ & $\frac{79}{500}$ & 8 &$\frac{37}{500}$ & $\frac{41}{500}$ & $\frac{73}{500}$  & 3  & $\frac{35}{500}$ & $\frac{33}{500}$ & $\frac{38}{500}$& 2 & $\frac{23}{500}$& $\frac{17}{500}$ & $\frac{32}{500}$ \vspace{0.09cm}  \\ 

5250 & 8  & $\frac{40}{500}$ & $\frac{49}{500}$ & $\frac{73}{500}$ & 8 &  $\frac{33}{500}$ & $\frac{36}{500}$ & $\frac{72}{500}$  & 3  & $\frac{37}{500}$&  $\frac{35}{500}$ & $\frac{41}{500}$ & 2  & $\frac{24}{500}$& $\frac{19}{500}$ & $\frac{34}{500}$ \vspace{0.05cm}  \\ 

5750 & 8  & $\frac{38}{500}$ & $\frac{43}{500}$ & $\frac{74}{500}$ & 8 &  $\frac{32}{500}$ & $\frac{34}{500}$ & $\frac{59}{500}$  & 3 & $\frac{33}{500}$&  $\frac{32}{500}$ & $\frac{37}{500}$& 2  & $\frac{19}{500}$& $\frac{13}{500}$ & $\frac{29}{500}$  \vspace{0.05cm}  \\ 

6250 & 8  & $\frac{34}{500}$ & $\frac{37}{500}$ & $\frac{70}{500}$ & 8 &  $\frac{27}{500}$ & $\frac{30}{500}$ & $\frac{69}{500}$  & 3  & $\frac{30}{500}$&  $\frac{30}{500}$ & $\frac{36}{500}$ & 2  & $\frac{16}{500}$& $\frac{10}{500}$ & $\frac{25}{500}$  \vspace{0.05cm}  \\ 

6750 & 8  & $\frac{30}{500}$ & $\frac{33}{500}$ & $\frac{68}{500}$ & 8 &  $\frac{25}{500}$ & $\frac{29}{500}$ & $\frac{66}{500}$  & 3  & $\frac{29}{500}$&  $\frac{25}{500}$ & $\frac{35}{500}$& 2 & $\frac{12}{500}$& $\frac{9}{500}$ & $\frac{21}{500}$\vspace{0.05cm}  \\ 
  \toprule 
\end{tabular} 
     \end{small}     
  \label{tab:TableA99}
\end{table}

\begin{table}[H]
 \caption{{\small
$F\left(k_n,n_t, \beta\right)$ values in
(\ref{a955}) and (\ref{new2}), for scenarios 21-22. $B_{1.2}$ and $B_{1.7}$ denotes the kernel-based approach in Besse et al. (2000), for $h_n = 1.2,~ 1.7$, respectively. $B_q$ denotes its penalized prediction approach.}}
\vspace{-0.16cm}
\centering
\begin{small}
    \begin{tabular}{c@{\hspace{1.5em}}ccc@{\hspace{1.6em}}ccc}
\toprule  
& \multicolumn{3}{c}{\hspace{-0.4cm}\textit{Scenario 21}} & \multicolumn{3}{c}{\hspace{-0.25cm}\textit{Scenario 22}} \\
  \midrule[0.05cm]
   $n_t$   & $B_{1.2}$ &  $B_{1.7}$ & $B_q$ & $B_{1.2}$ &  $B_{1.7}$ & $B_q$ \\
  \midrule[0.05cm] 
750 &  $\frac{449}{500}$ & $\frac{281}{500}$ & $\frac{7}{500}$ & $\frac{377}{500}$ & $\frac{222}{500}$  &  $\frac{5}{500}$  \vspace{0.09cm}  \\ 
1250 & $\frac{434}{500}$ & $\frac{225}{500}$ & $\frac{5}{500}$ & $\frac{355}{500}$ & $\frac{209}{500}$  &  $\frac{4}{500}$  \vspace{0.09cm}  \\ 
1750 & $\frac{436}{500}$ & $\frac{164}{500}$ & $\frac{5}{500}$ & $\frac{330}{500}$ & $\frac{196}{500}$  &  $\frac{4}{500}$  \vspace{0.09cm}  \\ 
2250  & $\frac{426}{500}$ & $\frac{142}{500}$ & $\frac{4}{500}$ & $\frac{309}{500}$ & $\frac{162}{500}$  &  $\frac{3}{500}$  \vspace{0.09cm}  \\ 
  2750 & $\frac{422}{500}$ & $\frac{123}{500}$ & $\frac{3}{500}$ & $\frac{292}{500}$ & $\frac{130}{500}$  &  $\frac{3}{500}$  \vspace{0.09cm}  \\ 
  3250  & $\frac{417}{500}$ & $\frac{105}{500}$ & $\frac{3}{500}$ & $\frac{281}{500}$ & $\frac{107}{500}$  &  $\frac{2}{500}$  \vspace{0.09cm}  \\   
3750 & $\frac{376}{500}$ & $\frac{97}{500}$ & $\frac{3}{500}$ & $\frac{269}{500}$ & $\frac{83}{500}$  &  $\frac{2}{500}$  \vspace{0.09cm}  \\ 

4250 & $\frac{358}{500}$ & $\frac{80}{500}$ & $\frac{2}{500}$ & $\frac{252}{500}$ & $\frac{72}{500}$  &  $\frac{1}{500}$  \vspace{0.09cm}  \\ 

4750 &$\frac{345}{500}$ & $\frac{71}{500}$ & $\frac{1}{500}$ & $\frac{241}{500}$ & $\frac{69}{500}$  &  0 \vspace{0.09cm}  \\ 

5250 & $\frac{313}{500}$ & $\frac{61}{500}$ & 0 & $\frac{230}{500}$ & $\frac{56}{500}$  &  $\frac{1}{500}$  \vspace{0.09cm}  \\ 

5750 & $\frac{262}{500}$ & $\frac{55}{500}$ & $\frac{1}{500}$ & $\frac{215}{500}$ & $\frac{45}{500}$  &  $\frac{1}{500}$  \vspace{0.09cm}  \\ 

6250 & $\frac{240}{500}$ & $\frac{52}{500}$ & $\frac{1}{500}$ & $\frac{203}{500}$ & $\frac{37}{500}$  &  0  \vspace{0.09cm}  \\ 

6750 & $\frac{230}{500}$ & $\frac{46}{500}$ & 0 & $\frac{195}{500}$ & $\frac{32}{500}$  &  0  \vspace{0.09cm}  \\ 
  \toprule 
\end{tabular}  
    \end{small}    
  \label{tab:TableA10}
  \vspace{-0.4cm}
\end{table}

As above, the empirical mean of $\left\| \left(\rho - \overline{\rho}_{k_n}^{l} \right) \left(X_{n-1}^{l} \right) \right\|_{H}^{k_n}$, for $l=1,\ldots,N$, with $N = 500$ realizations, for small-sample scenarios considered, will be illustrated in Figures \ref{fig:F5}-\ref{fig:F7} below.

\begin{figure}[H]
 \hspace{-1.06cm} \includegraphics[width=9.2cm,height=5.9cm]{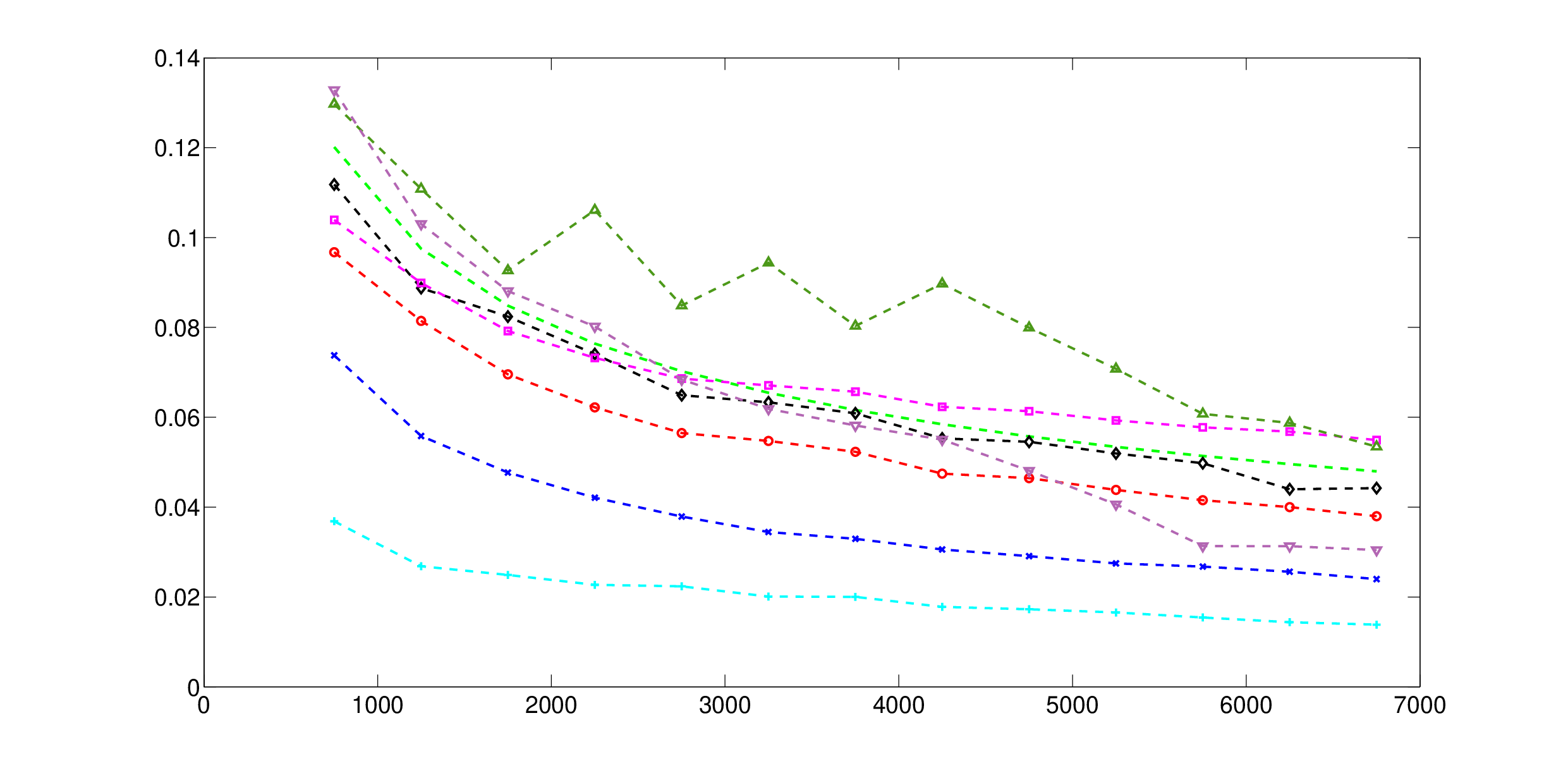}
 \hspace{-1.02cm}  \includegraphics[width=9.2cm,height=5.9cm]{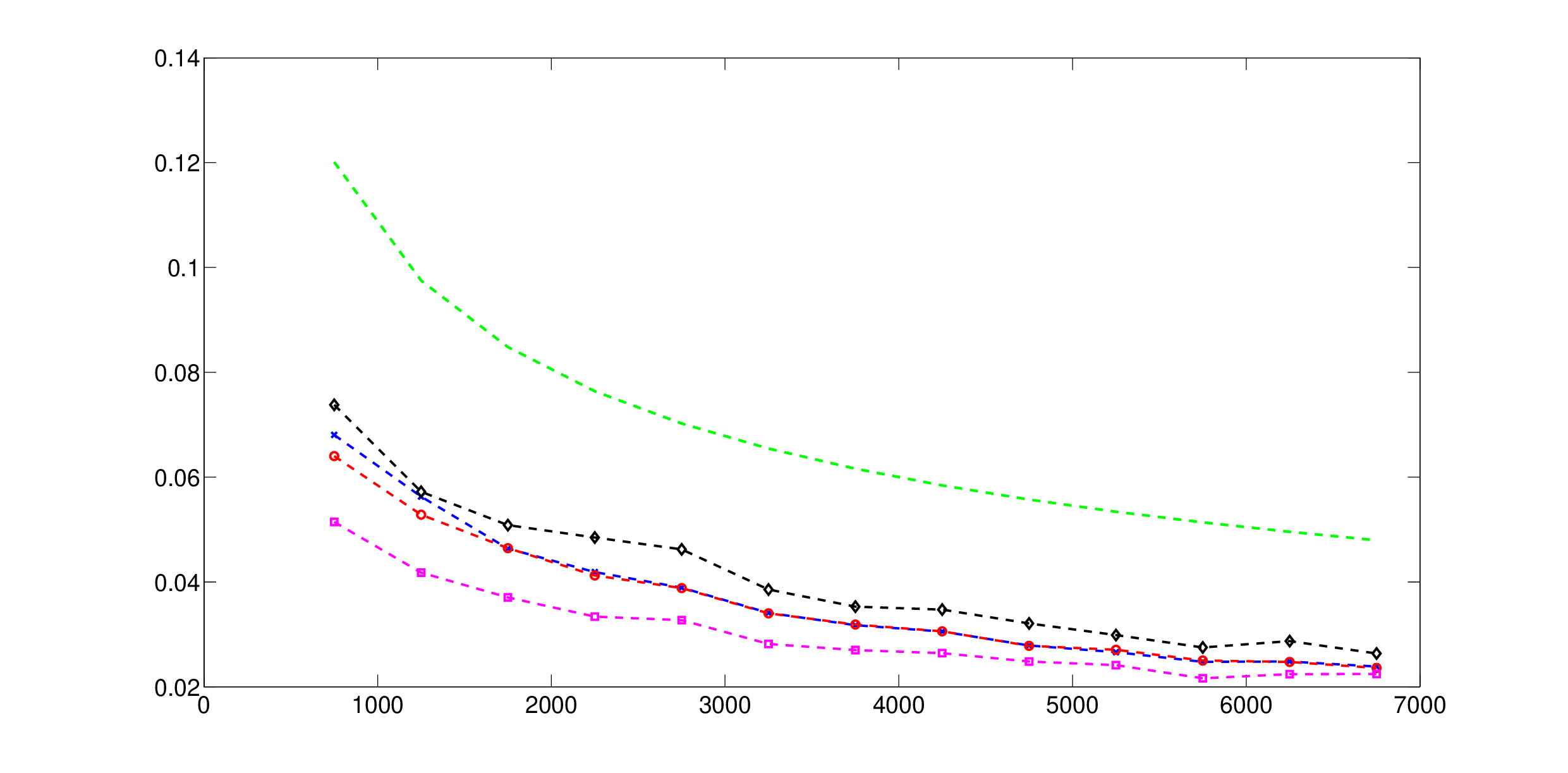}
\vspace{-1.22cm}
\caption{{\small Empirical mean of $\left\| \left(\rho - \overline{\rho}_{k_n}^{l} \right) \left(X_{n-1}^{l} \right) \right\|_{H}$ and $\left\| \left(\rho - \overline{\rho}_{k_n}^{l} \right) \left(X_{n-1}^{l} \right) \right\|_{H}^{k_n}$, for scenario 13 (on left) and scenario 15 (on right), for our approach (blue star dotted line) and those one presented in Antoniadis and Sapatinas (2001) (pink square dotted line), Besse et al. (2000) (cyan blue plus dotted line for penalized prediction; dark green upward-pointing triangle  and purple downward-pointing triangle dotted lines, for kernel-based prediction, for $h_n = 0.15$ and $h_n = 0.25$, respectively), Bosq (2000) (red circle dotted line) and Guillas (2001) (black diamond dotted line). The curve $\xi_{n_t,\beta} = \frac{\left(\ln(n_t) \right)^{\beta} }{n_{t}^{1/2} }$, with $\beta = 65/100$, is drawn (light green dotted line).}} \label{fig:F5}
\end{figure}

\begin{figure}[H]
 \hspace{-1.06cm} \includegraphics[width=9.2cm,height=5.9cm]{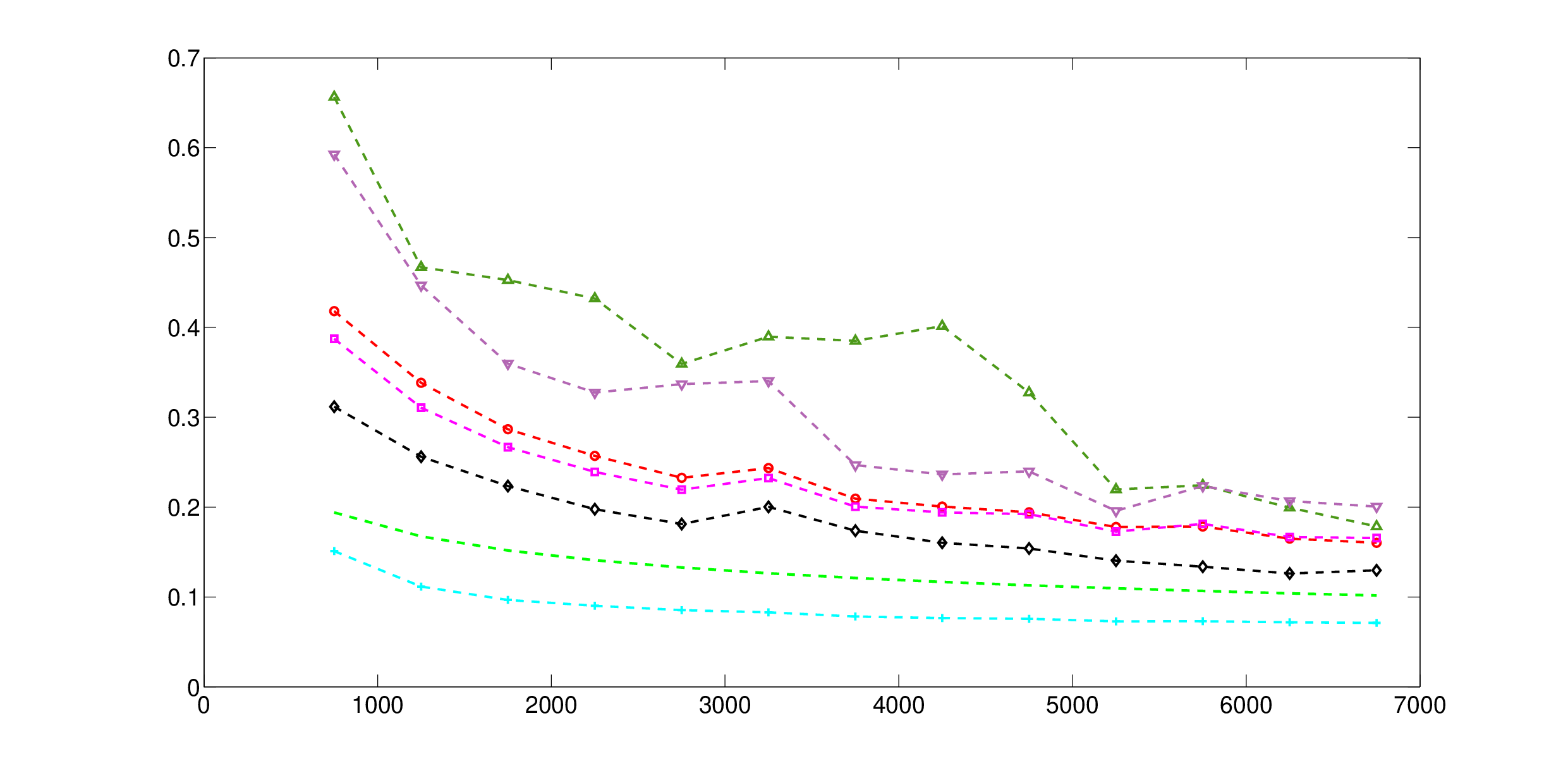}
 \hspace{-1.02cm}  \includegraphics[width=9.2cm,height=5.9cm]{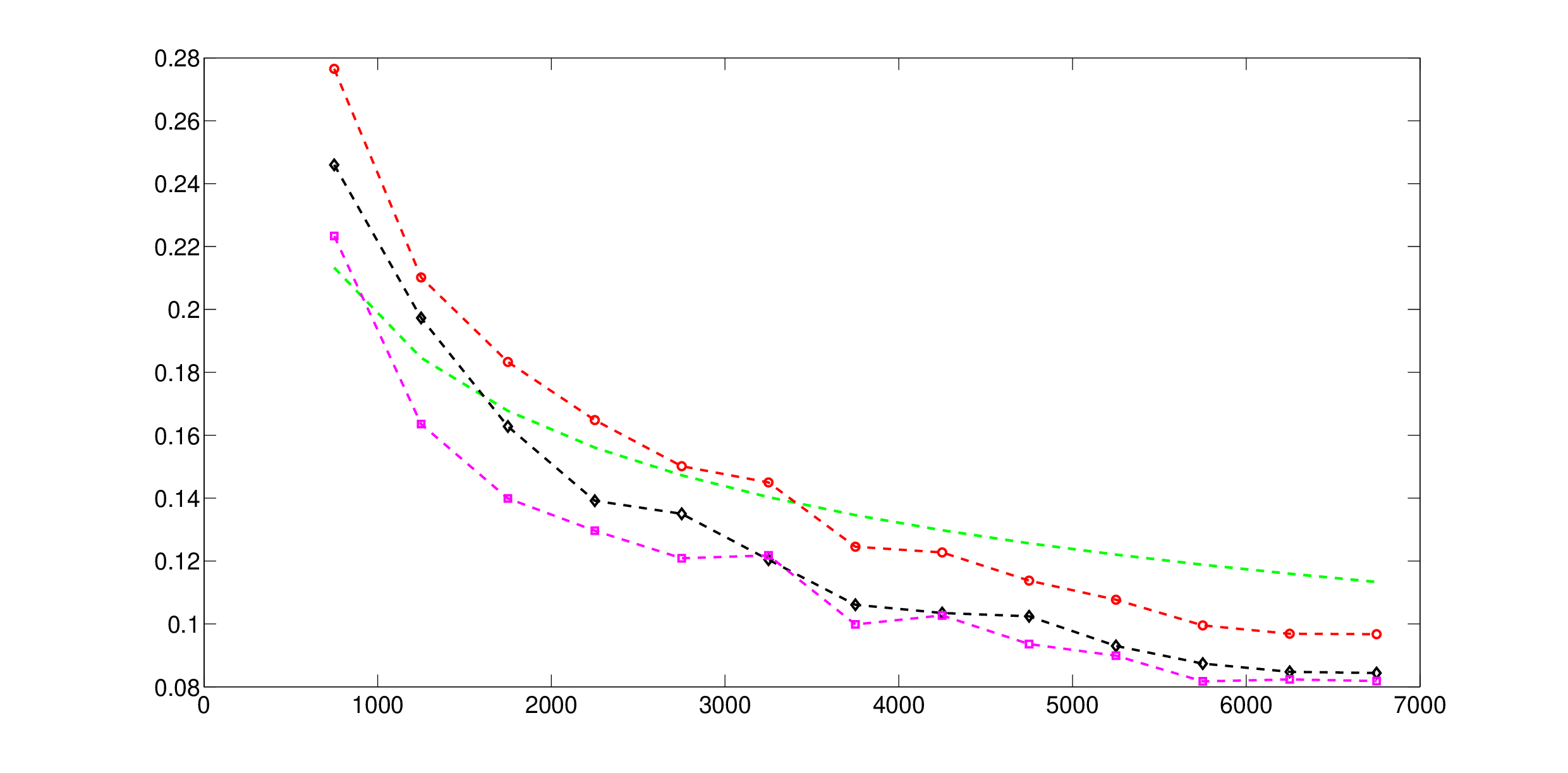}
\vspace{-1.22cm}
\caption{{\small Empirical mean of $\left\| \left(\rho - \overline{\rho}_{k_n}^{l} \right) \left(X_{n-1}^{l} \right) \right\|_{H}$ and $\left\| \left(\rho - \overline{\rho}_{k_n}^{l} \right) \left(X_{n-1}^{l} \right) \right\|_{H}^{k_n}$, for scenario 17 (on left) and scenario 19 (on right), for approaches presented in Antoniadis \& Sapatinas (2003) (pink square dotted line), Besse et al. (2000) (cyan blue plus dotted line for penalized prediction; dark green upward-pointing triangle  and purple downward-pointing triangle dotted lines, for kernel-based prediction, for $h_n = 1.2$ and $h_n = 1.7$, respectively), Bosq (2000) (red circle dotted line) and Guillas (2001) (black diamond dotted line). The curve $\xi_{n_t,\beta} = \frac{\left(\ln(n_t) \right)^{\beta} }{n_{t}^{1/3} }$, with $\beta = 3/10$, is drawn (light green dotted line).}} \label{fig:F6}
\end{figure}

\begin{figure}[H]
 \hspace{-1.06cm} \includegraphics[width=9.2cm,height=5.9cm]{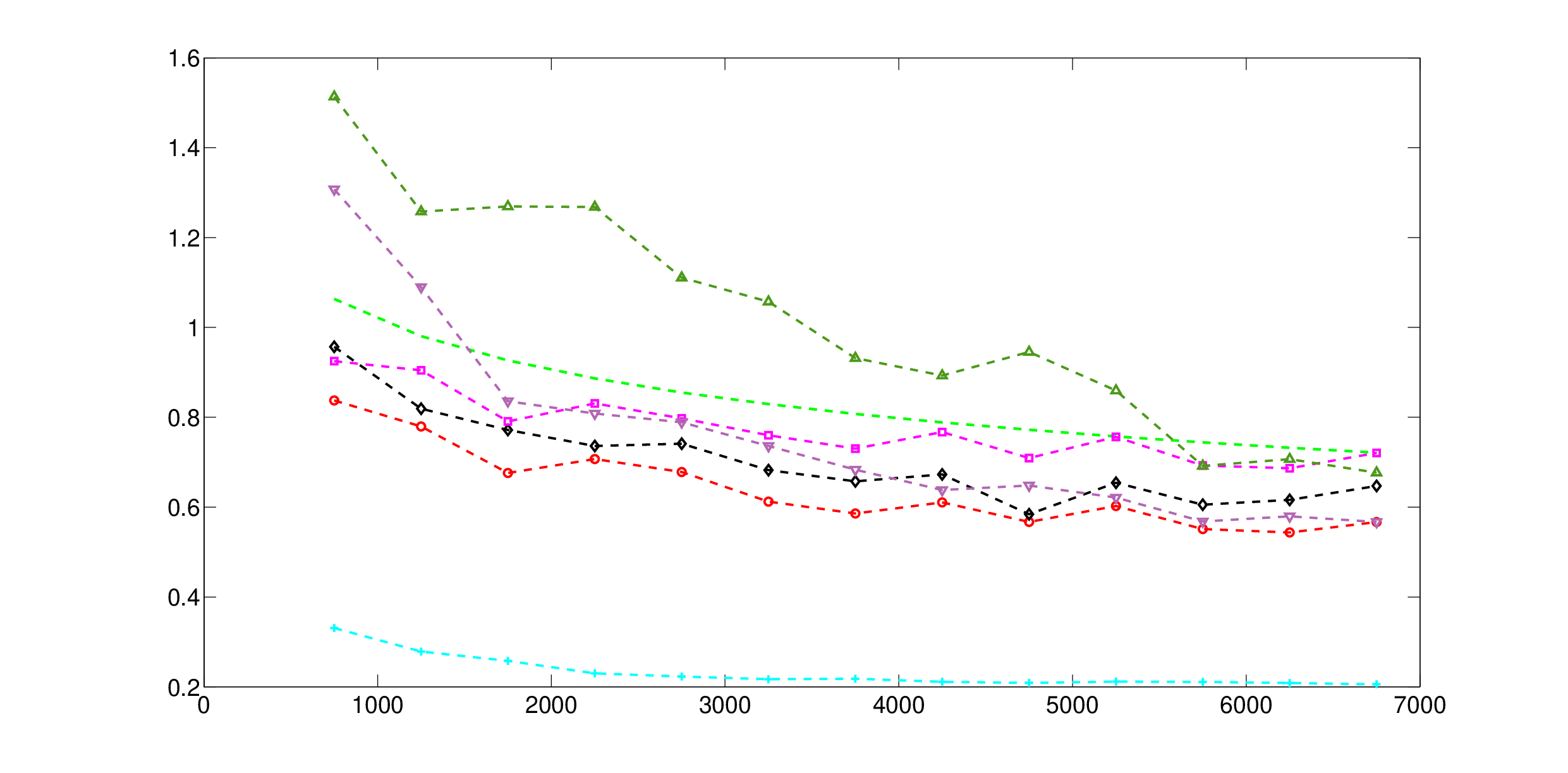}
 \hspace{-1.02cm}  \includegraphics[width=9.2cm,height=5.9cm]{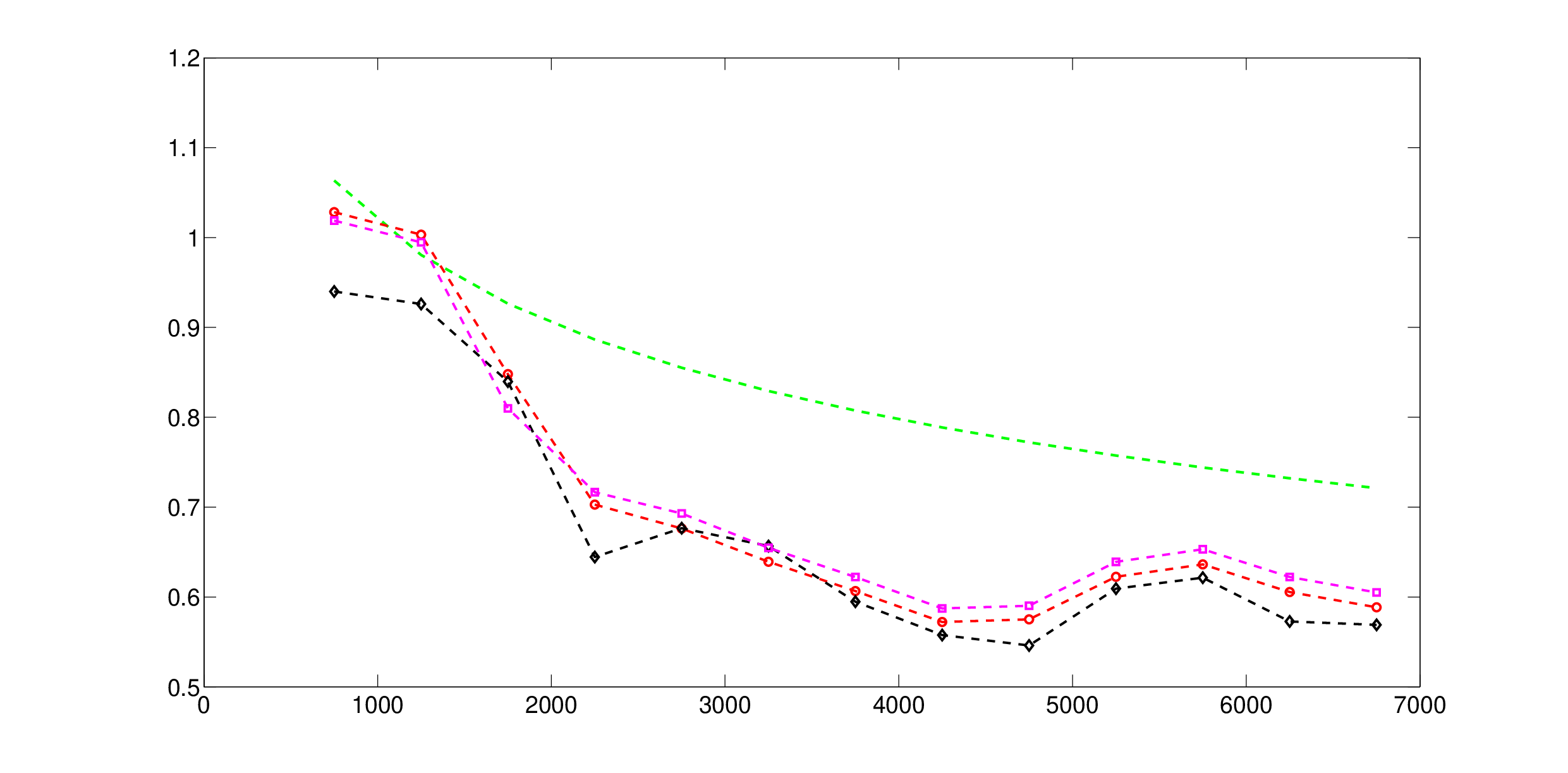}
\vspace{-1.22cm}
\caption{{\small Empirical mean of $\left\| \left(\rho - \overline{\rho}_{k_n}^{l} \right) \left(X_{n-1}^{l} \right) \right\|_{H}$ and $\left\| \left(\rho - \overline{\rho}_{k_n}^{l} \right) \left(X_{n-1}^{l} \right) \right\|_{H}^{k_n}$, for scenario 21 (on left) and scenario 23 (on right), for approaches presented in Antoniadis \& Sapatinas (2003) (pink square dotted line), Besse et al. (2000) (cyan blue plus dotted line for penalized prediction; dark green upward-pointing triangle  and purple downward-pointing triangle dotted lines, for kernel-based prediction, for $h_n = 1.2$ and $h_n = 1.7$, respectively), Bosq (2000) (red circle dotted line) and Guillas (2001) (black diamond dotted line). The curve $\xi_{n_t,\beta} = \frac{\left(\ln(n_t) \right)^{\beta} }{n_{t}^{1/3} }$, with $\beta = 125/100$, is drawn (light green dotted line).}} \label{fig:F7}
\end{figure}

Results displayed in Tables 4-6 and 9-14, and Figures 2-7, are discussed in Section 9 of the main paper.

% *********************************

% *********************************

\vspace*{0.5cm}

 \noindent {\bfseries {\large Acknowledgments}}

This work has been supported in part by project MTM2015--71839--P (co-funded by Feder funds), of the DGI, MINECO, Spain.